\documentclass[11pt]{article}

\usepackage[letterpaper,margin=1in]{geometry}
\usepackage{amsmath,amssymb,amsthm,mathtools}
\usepackage{microtype}
\usepackage{xcolor}
\usepackage[colorlinks=true,hyperfootnotes=false,linkcolor=blue!60!black,citecolor=blue,urlcolor=blue!60!black]{hyperref}
\hypersetup{
 pdftitle={The Berry--Esseen Constant Conjecture is Eventually True}
}

\newtheorem{theorem}{Theorem}[section]
\newtheorem{lemma}[theorem]{Lemma}
\newtheorem{proposition}[theorem]{Proposition}
\newtheorem{corollary}[theorem]{Corollary}
\theoremstyle{definition}

\theoremstyle{remark}

\usepackage{hyperref}
\newcommand{\PP}{\mathbb P}
\newcommand{\EE}{\mathbb E}
\newcommand{\R}{\mathbb R}

\newcommand{\Pthree}{\mathcal P_3}
\newcommand{\ce}{c_{\mathrm E}}
\newcommand{\pe}{p_{\mathrm E}}
\newcommand{\qe}{q_{\mathrm E}}
\newcommand{\se}{\sigma_{\mathrm E}}
\newcommand{\be}{\beta_{\mathrm E}}
\newcommand{\ke}{\kappa_{\mathrm E}}
\newcommand{\Pe}{P_{\mathrm E}}
\newcommand{\aE}{a_{\mathrm E}}
\newcommand{\bE}{b_{\mathrm E}}
\newcommand{\ind}{\mathbf 1}
\DeclareMathOperator{\Var}{Var}
\DeclareMathOperator{\supp}{supp}
\DeclareMathOperator{\dist}{dist}
\DeclareMathOperator{\sgn}{sgn}
\DeclareMathOperator{\sinc}{sinc}

\allowdisplaybreaks

\title{The Berry-Esseen Constant Conjecture is Eventually True}
\author{Hengzhi He\thanks{Department of Statistics and Data Science, University of California, Los Angeles. Contact: guangcheng@ucla.edu; hengzhihe@ucla.edu}
\qquad Guang Cheng\footnotemark[1]}
\date{September 2026}

\begin{document}
\maketitle

\begin{abstract}
We prove that the Esseen's constant $0.4097...$ is valid in the iid Berry--Esseen inequality as long as the sample sizes $n\ge N$. The threshold $N$ is universal over all real summand distributions with positive variance and a finite third absolute moment. In particular, one may take $N=2\left\lceil\exp(10^{17})\right\rceil.
$ 
\end{abstract}

\begingroup
\renewcommand{\thefootnote}{}
\footnotetext{\textit{2020 Mathematics Subject Classification.}
Primary 60F05; Secondary 60E15.}
\endgroup

\noindent\textbf{Keywords:} Berry--Esseen inequality; Edgeworth expansion; Esseen constant; Extremal distributions; Lattice distributions; Normal approximation.

% ===== introduction =====
\section{Introduction}\label{sec:introduction}

Let \(X_1,X_2,\ldots\) be independent identically distributed real random variables with mean zero, variance one and finite third absolute moment \(\beta=\EE|X_1|^3\). The Berry--Esseen theorem asserts that
\begin{equation}\label{eq:berry-esseen}
 \sup_{x\in\R}\left|
 \PP\left\{\frac{X_1+\cdots+X_n}{\sqrt n}\le x\right\}-\Phi(x)
 \right|\le C\frac{\beta}{\sqrt n}
\end{equation}
for an absolute constant \(C\), where \(\Phi\) is the standard normal distribution function~\cite{Berry1941,Esseen1942}. Determining the smallest possible constant $C$ in \eqref{eq:berry-esseen} is a classical problem in normal approximation. Esseen proved that any such constant is at least
\begin{equation}\label{eq:ce}
 \ce=\frac{3+\sqrt{10}}{6\sqrt{2\pi}}=0.4097321837\ldots,
\end{equation}
and identified the two-point distributions that attain this value asymptotically~\cite{Esseen1956}. 

The conjecture is that \eqref{eq:berry-esseen} holds with \(C=\ce\)
for every summand distribution and every positive integer \(n\) \cite{Zolotarev1966,ZolotukhinNagaevChebotarev2018}.
We establish this bound for every sample size above a universal threshold.

Earlier work on refining Berry--Esseen bounds and their remainder
terms includes Prawitz~\cite{Prawitz1975} and Bentkus~\cite{Bentkus1994}.
Esseen's fixed-distribution asymptotics were subsequently
strengthened to uniform asymptotic results by
Chistyakov~\cite{Chistyakov1999,Chistyakov2002a,
Chistyakov2002b,Chistyakov2003}; see also the account in
Mattner and Shevtsova~\cite{MattnerShevtsova2019}. In the identically distributed case, a consequence of Shevtsova's estimates is
\[
 \Delta_n\le \ce\frac{\beta}{\sqrt n}+3\frac{\beta^2}{n},
\]
where \(\Delta_n\) denotes the left side of \eqref{eq:berry-esseen}~\cite[Corollary~4.18, p.~303]{Shevtsova2012}. For two-point distributions, Schulz proved the strict bound \(\Delta_n<\ce\beta/\sqrt n\) for every \(n\)~\cite[Theorem~1, p.~1]{Schulz2016}.
Complementary analytic and computational estimates in the binomial
case were developed in \cite{ZolotukhinNagaevChebotarev2018}.

For \(n=1\), Bentkus and Kir\v{s}a~\cite{BentkusKirsa1989}
proved the stronger bound \(\Delta_1\le 0.3704\,\beta<\ce\beta\)
for all distributions satisfying the above moment assumptions.

The Theorem below summarizes our main result of this paper. 

\begin{theorem}\label{thm:main}
There exists an integer \(N\ge1\) such that, for every probability distribution \(P\) with mean zero, variance one and finite third absolute moment \(\beta\), every integer \(n\ge N\), and independent random variables \(X_1,\ldots,X_n\) with common law \(P\),
\[
 \sup_{x\in\R}\left|
 \PP\left\{\frac{X_1+\cdots+X_n}{\sqrt n}\le x\right\}-\Phi(x)
 \right|
 \le \ce\frac{\beta}{\sqrt n}.
\]
The integer \(N\) is independent of the summand distribution. In particular, one may take
\begin{equation}\label{eq:explicit-threshold}
 N=2\left\lceil\exp(10^{17})\right\rceil.
\end{equation}
\end{theorem}

The theorem allows the summand distribution to vary with \(n\), and imposes no common bound on the third absolute moment. The coefficient \(\ce\) cannot be decreased, even if \(N\) is increased, because a fixed Esseen two-point distribution attains it in the limit. The explicit threshold \eqref{eq:explicit-threshold} is a sufficient bound, with no attempt at optimization. Establishing the validity of the Berry--Esseen bound with constant \(\ce\) for every \(n\geq 1\) remains an open problem.

The distinction between uniform asymptotic bounds, such as that of
Shevtsova~\cite[Corollary~4.18, p.~303]{Shevtsova2012},
and Theorem~\ref{thm:main} is that the former do not rule out
arbitrarily small excesses above \(\ce\).

The proof combines a variational analysis of extremal
distributions with three classical ingredients:
Esseen's smoothing method, short Edgeworth expansions, and
moment inequality~\cite{Esseen1945,Esseen1956};
Shevtsova's explicit normal approximation
bounds~\cite{Shevtsova2012,Shevtsova2013};
and Schulz's exact bound for two-point
distributions~\cite{Schulz2016}.
We use the first two ingredients to control and identify
limits of extremal distributions.
Schulz's result supplies the finite-sample bound for the
Bernoulli laws used in the local comparison.

A Lean formalization of the main theorem, conditional on explicitly stated previous results, is available at
\href{https://github.com/UCLA-Trustworthy-AI-Lab/BerryEsseenLean}{https://github.com/UCLA-Trustworthy-AI-Lab/BerryEsseenLean}.

Much of modern statistical inference relies on asymptotic approximations, yet achieving adequate accuracy may require substantial sample sizes. Despite important advances in non-asymptotic inference \cite{YangShangCheng2020,ZhengCheng2021}, the practical usefulness of existing guarantees can be limited by insufficiently sharp bounds on Gaussian approximation errors. Sharpening the Berry--Esseen bound~(\ref{eq:berry-esseen})---a cornerstone of Gaussian approximation theory---could therefore enable tighter finite-sample guarantees and more accurately calibrated statistical procedures.

\section{Proof Sketch}

The main idea of our proof is to study the worst counterexample
and use its optimality to obtain additional information about
its distribution. This information will eventually allow us
to rule out such counterexamples at large sample sizes.

To carry out this idea, we argue by contradiction.
Suppose that the Berry--Esseen bound with constant \(\ce\)
fails for arbitrarily large sample sizes.
For each such sample size, we choose a distribution that
maximizes the normalized Berry--Esseen error; the existence
of a maximizer is justified below.
Since no admissible perturbation can increase this error,
the maximizing distribution satisfies first-order optimality
conditions. These conditions yield exact identities at every
point in its support, which are the main tool in the argument.

By choosing the sample sizes appropriately and using these
identities, we obtain a sequence of worst counterexamples
whose supports all lie in one fixed bounded interval.
We can then pass to a subsequence that converges to a limiting
distribution.
A smoothing argument, together with Esseen's moment inequality
and its equality characterization, identifies the limit,
up to reflection, as the standardized two-point distribution
appearing in the equality case of Esseen's theorem
\cite[p.~161, Theorem]{Esseen1956}.
We refer to this distribution as the \emph{Esseen two-point law}
throughout the paper.

This convergence alone does not rule out small violations at
finite sample sizes. In particular, additional support points
may remain in regions whose probabilities tend to zero.
We therefore return to the optimality conditions and prove
the stronger statement that the entire support must lie in
shrinking neighborhoods of the two atoms of the limiting
distribution.
The probabilities of the two clusters also converge to the
corresponding limiting masses.
A local comparison then gives the exact bound with constant
\(\ce\) for every sufficiently late member of the selected
sequence.
This contradicts their being counterexamples and proves the
existence of a universal finite threshold.

In the following subsections, we present the argument step by step.
Each subsection begins with a brief summary in bold,
explaining its goal and its role in the proof.

\subsection{Extremal distributions and first-order conditions}
\label{subsec:outline-extremizers}

\textbf{We first make the choice of a worst counterexample precise.
We then perturb its distribution and derive the identities
imposed by its optimality.}

For a mean-zero, variance-one law \(P\), write
\begin{equation}\label{eq:outline-definitions}
 \beta(P)=\int |x|^3P(dx),\qquad
 R_n(P)=\frac{\sqrt n\,\Delta_n(P)}{\beta(P)},\qquad
 C_n=\sup_P R_n(P),
\end{equation}
where \(\Delta_n(P)\) denotes the error in \eqref{eq:berry-esseen}
for summands with law \(P\), and the supremum is over all such
laws with finite third absolute moment.
We use the convergence and uniform bound established in
\eqref{eq:cn-limit}:
\begin{equation}\label{eq:outline-cn-bounds}
 C_n\longrightarrow\ce,
 \qquad C_n\le C_{\mathrm{BE}}:=0.4690\quad(n\ge1).
\end{equation}
Suppose that \(C_n>\ce\) for arbitrarily large \(n\).
Fix any such sample size \(n\ge2\).
Lemma~\ref{lem:attainment} supplies a standardized law \(P\) and a
threshold \(t\) such that
\begin{equation}\label{eq:outline-maximizing-pair}
 \frac{\sqrt n}{\beta(P)}
 \left\{\Pr(X_1+\cdots+X_n\le t)-\Phi(t/\sqrt n)\right\}
 =R_n(P)=C_n,
 \qquad 1\le\beta(P)<B_*,
\end{equation}
where \(X_1,\ldots,X_n\) are iid with law \(P\).
The expression on the left is the normalized signed error.
Here \(B_*<\infty\) is the universal moment cutoff in
\eqref{eq:moment-cutoff}: every standardized law satisfying
\(R_n(P)>\ce\) has third absolute moment less than \(B_*\).

To exploit maximality, consider the perturbation
\begin{equation}\label{eq:outline-contamination}
 P_\varepsilon=(1-\varepsilon)P+\varepsilon\delta_y,
\end{equation}
where \(\delta_y\) is a point mass at \(y\).
Write
\begin{equation}\label{eq:outline-perturbed-moments}
\begin{aligned}
 m_\varepsilon
 &=\int x\,P_\varepsilon(dx),\qquad
 \sigma_\varepsilon^2
 =\int (x-m_\varepsilon)^2P_\varepsilon(dx),\\
 \beta_\varepsilon
 &=\sigma_\varepsilon^{-3}
   \int |x-m_\varepsilon|^3P_\varepsilon(dx).
\end{aligned}
\end{equation}
Let \(F_{n,\varepsilon}\) denote the distribution function of the
sum of \(n\) independent random variables with law \(P_\varepsilon\), namely 
\[
 F_{n,\varepsilon}(u)
 =\Pr\left(\sum_{i=1}^n X_{i,\varepsilon}\le u\right),
 \qquad
 X_{1,\varepsilon},\ldots,X_{n,\varepsilon}
 \overset{\mathrm{iid}}{\sim}P_\varepsilon.
\]

For the unperturbed law, we write \(F_{n,P}=F_{n,0}\). At the fixed threshold \(t\), define the normalized signed error by
\begin{equation}\label{eq:outline-signed-error}
 r_y(\varepsilon)
 =\frac{\sqrt n}{\beta_\varepsilon}
 \left\{
 F_{n,\varepsilon}(t)
 -\Phi\left(\frac{t-nm_\varepsilon}
                  {\sigma_\varepsilon\sqrt n}\right)
 \right\}.
\end{equation}
Let \(I(y)=r_y'(0+)\) be its right derivative with respect to
\(\varepsilon\) at zero. (We call \(I(y)\) the \emph{contamination derivative} at \(y\).
It measures the first-order change in the normalized signed
error under the perturbation
\(P_\varepsilon=(1-\varepsilon)P+\varepsilon\delta_y\).)
By standardizing the perturbed law, the definition of \(C_n\) in
\eqref{eq:outline-definitions} and the maximizing identity
\eqref{eq:outline-maximizing-pair} imply
\(r_y(\varepsilon)\le C_n=r_y(0)\) for all sufficiently small
\(\varepsilon\ge0\).
Thus \(I(y)\le0\).
Lemma~\ref{lem:influence} also gives equality at every support point:
\begin{equation}\label{eq:outline-contact}
 I(y)\le0\quad(y\in\mathbb R),
 \qquad I(x)=0\quad(x\in\operatorname{supp}P).
\end{equation}
The equality in \eqref{eq:outline-contact} supplies an exact
identity at every support point, even when very little probability
lies near that point.
These identities provide the additional information needed to
constrain the support of a worst counterexample.

For the maximizing pair in \eqref{eq:outline-maximizing-pair}, write
\[
 \beta=\beta(P),\qquad
 z=\frac{t}{\sqrt n},\qquad
 M=\int x|x|\,P(dx),\qquad
 G(u)=\Pr(X_1+\cdots+X_{n-1}\le u),
\]
where the \(X_i\) are iid with law \(P\).
For the contamination \eqref{eq:outline-contamination},
differentiation with respect to the mixing weight gives
\begin{equation}\label{eq:outline-cdf-derivative}
 \left.\frac{d}{d\varepsilon}F_{n,\varepsilon}(t)
 \right|_{\varepsilon=0+}
 =n\{G(t-y)-F_{n,0}(t)\}.
\end{equation}
The factor \(n\) counts the possible positions of the contaminated
summand: to first order, one summand is replaced by \(y\), while
the remaining \(n-1\) summands retain law \(P\).
Differentiating \eqref{eq:outline-signed-error}, using
\eqref{eq:outline-cdf-derivative} and differentiating the moments
in \eqref{eq:outline-perturbed-moments}, gives the exact influence formula
\eqref{eq:influence} in Lemma~\ref{lem:influence}, namely 
\begin{equation}\label{eq:outline-influence}
\begin{aligned}
 \beta I(y)
 ={}& n^{3/2}\{G(t-y)-F_{n,P}(t)\}
      +ny\phi(z)\\
 &+\frac{\sqrt n}{2}z\phi(z)(y^2-1)
      +\frac32 C_n\beta(y^2-1)\\
 &-C_n\{|y|^3-\beta-3yM\}.
\end{aligned}
\end{equation}

\subsection{Obtaining a uniform support bound}
\label{subsec:outline-bounded-support}

\textbf{We now use the optimality conditions to select sample
sizes \(n_j\to\infty\) and corresponding maximizing pairs
\((P_j,t_j)\) such that the supports of \(P_j\) all lie in
one fixed bounded interval.
For each pair, write \(I_j\) for the contamination derivative
introduced in the preceding subsection.
We will show that, for all sufficiently large \(j\),
\(I_j(y)<0\) whenever \(|y|>L\), where \(L\) is independent
of \(j\).
Such points cannot belong to the support of \(P_j\), because
the contact condition \eqref{eq:outline-contact} requires
\(I_j(y)=0\) at every support point.
Thus, after discarding finitely many initial terms, all the
supports lie in \([-L,L]\).}

\textbf{The upper bound for \(I_j\) will contain a term involving
\(n_j(C_{n_j-1}-C_{n_j})_+\).
We therefore choose the sample sizes so that this quantity
tends to zero while the corresponding distributions remain
counterexamples. More precisely, the convergence
\(C_n\to\ce\), together with our assumption that
\(C_n>\ce\) for arbitrarily large \(n\), allows us to apply
Lemma~\ref{lem:selection} and choose
\[
 C_{n_j}>\ce,
 \qquad
 n_j(C_{n_j-1}-C_{n_j})_+\longrightarrow0.
\]
Here \(a_+=\max\{a,0\}\).
For these sample sizes, the negative cubic term in the
influence bound will dominate the remaining terms when
\(|y|\) is large, yielding the desired uniform support bound.}

Since the supremum defining \(C_{n-1}\) in
\eqref{eq:outline-definitions} is taken over all standardized
laws, we have \(R_{n-1}(P)\le C_{n-1}\).
Applying this bound at the threshold \(t-y\) gives
\begin{equation}\label{eq:outline-previous-error}
 G(t-y)-\Phi\left(\frac{t-y}{\sqrt{n-1}}\right)
 \le \frac{\beta C_{n-1}}{\sqrt{n-1}}.
\end{equation}
Since \(F_{n,0}(t)=\Pr(X_1+\cdots+X_n\le t)\),
\eqref{eq:outline-maximizing-pair} gives
\(F_{n,0}(t)-\Phi(z)=\beta C_n/\sqrt n\).
We collect the normal-distribution terms in
\begin{equation}\label{eq:outline-gaussian-H}
\begin{aligned}
 H_n(z,y):={}&
 n^{3/2}\left\{
 \Phi\left(\frac{z-y/\sqrt n}{\sqrt{1-1/n}}\right)
 -\Phi(z)\right\}\\
 &+ny\phi(z)
 +\frac{\sqrt n}{2}z\phi(z)(y^2-1),
\end{aligned}
\end{equation}
which is the quantity defined in \eqref{eq:gaussian-H}.
Substituting the maximizing identity
\eqref{eq:outline-maximizing-pair} and the bound
\eqref{eq:outline-previous-error} into
\eqref{eq:outline-influence} then yields
\begin{equation}\label{eq:outline-influence-step}
\begin{split}
 \beta I(y)\le{}&
 H_n(z,y)
 +\beta\left\{
 \frac{n^{3/2}}{\sqrt{n-1}}C_{n-1}-nC_n
 \right\}\\
 &+\frac32 C_n\beta(y^2-1)
 -C_n\{|y|^3-\beta-3yM\}.
\end{split}
\end{equation}
Here \(\phi=\Phi'\) is the standard normal density.

To bound the expression in braces in
\eqref{eq:outline-influence-step}, use the uniform bound
\(C_m\le C_{\mathrm{BE}}\) from \eqref{eq:outline-cn-bounds}.
Since
\[
 0\le n\left(\sqrt{\frac n{n-1}}-1\right)\le1
 \qquad(n\ge2),
\]
we have
\begin{equation}\label{eq:outline-recursive-bound}
\begin{split}
 \frac{n^{3/2}}{\sqrt{n-1}}C_{n-1}-nC_n
 &=
 n(C_{n-1}-C_n)
 +n\left(\sqrt{\frac n{n-1}}-1\right)C_{n-1}\\
 &\le n(C_{n-1}-C_n)_++C_{\mathrm{BE}},
\end{split}
\end{equation}
where \(a_+=\max\{a,0\}\).
Moreover, \(|M|\le\int x^2P(dx)=1\), and
\(1\le\beta<B_*\) by \eqref{eq:outline-maximizing-pair}.
Lemma~\ref{lem:gaussian-expansion} gives the Gaussian bound
\eqref{eq:gaussian-global}: for all integers \(n\ge2\) and
all \(y,z\in\mathbb R\),
\begin{equation}\label{eq:outline-gaussian-global}
 |H_n(z,y)|
 \le \frac{\phi(0)}6|y|^3+A_0(1+|y|),
\end{equation}
where \(A_0<\infty\) is an absolute constant.
Substituting \eqref{eq:outline-recursive-bound} and
\eqref{eq:outline-gaussian-global} into
\eqref{eq:outline-influence-step}, and collecting terms, gives
\[
\begin{aligned}
 \beta I(y)\le{}&
 -\left(C_n-\frac{\phi(0)}6\right)|y|^3
 +\beta n(C_{n-1}-C_n)_+\\
 &+\frac32 C_n\beta y^2+3C_nMy
 +\beta\left(C_{\mathrm{BE}}-\frac{C_n}{2}\right)
 +A_0(1+|y|).
\end{aligned}
\]
Using \(1\le\beta<B_*\) from
\eqref{eq:outline-maximizing-pair}, \(|M|\le1\), and
\(C_n\le C_{\mathrm{BE}}\) from
\eqref{eq:outline-cn-bounds}, together with
\(|y|\le(1+y^2)/2\), we bound the last four terms by
\(A(1+y^2)\). Thus
\begin{equation}\label{eq:outline-cubic-bound}
 \beta I(y)
 \le
 -\left(C_n-\frac{\phi(0)}6\right)|y|^3
 +A(1+y^2)
 +\beta n(C_{n-1}-C_n)_+,
\end{equation}
where \(A<\infty\) is independent of \(n,P,t,y\).
Since \(C_n>\ce>\phi(0)/6\), the coefficient of \(|y|^3\)
is bounded above by the fixed negative constant
\(-(\ce-\phi(0)/6)\).

We now choose the sample sizes so that the last term in
\eqref{eq:outline-cubic-bound} tends to zero.
By \eqref{eq:outline-cn-bounds}, the sequence
\(a_n=C_n-\ce\) tends to zero; by our contradiction assumption,
it is positive at arbitrarily large indices.
Lemma~\ref{lem:selection} therefore provides integers
\(n_j\to\infty\) such that
\begin{equation}\label{eq:outline-selected-sizes}
 C_{n_j}>\ce,\qquad
 d_j:=n_j(C_{n_j-1}-C_{n_j})_+\longrightarrow0.
\end{equation}
At each selected sample size, choose a maximizing pair
\((P_j,t_j)\) satisfying \eqref{eq:outline-maximizing-pair}, and write
\(\beta_j=\beta(P_j)\) and \(I_j\) for their contamination
derivatives.
For each selected maximizing pair \((P_j,t_j)\),
\eqref{eq:outline-selected-sizes} gives \(C_{n_j}>\ce\),
while \eqref{eq:outline-maximizing-pair} gives
\(1\le\beta_j<B_*\).
Since \(d_j=n_j(C_{n_j-1}-C_{n_j})_+\ge0\), we have
\[
 C_{n_j}-\frac{\phi(0)}6
 \ge \ce-\frac{\phi(0)}6,
 \qquad
 \beta_jd_j\le B_*d_j.
\]
Applying \eqref{eq:outline-cubic-bound} with
\((n,P,t)=(n_j,P_j,t_j)\) therefore yields
\begin{equation}\label{eq:outline-selected-influence}
 \beta_j I_j(y)
 \le -\gamma|y|^3+A(1+y^2)+B_*d_j,
 \qquad
 \gamma:=\ce-\frac{\phi(0)}6>0.
\end{equation}
Choose \(L<\infty\) such that
\[
 \gamma r^3>A(1+r^2)+1
 \qquad\text{whenever }r>L.
\]
Since \(d_j\to0\) by \eqref{eq:outline-selected-sizes}
and \(B_*\) is a fixed finite constant, we have
\(B_*d_j\le1\) for all sufficiently large \(j\).
Equation~\eqref{eq:outline-selected-influence} then implies
\(I_j(y)<0\) whenever \(|y|>L\).
The support contact condition \eqref{eq:outline-contact},
namely \(I_j(x)=0\) for every \(x\in\operatorname{supp}P_j\),
therefore forces
\begin{equation}\label{eq:outline-common-support}
 \operatorname{supp}P_j\subseteq[-L,L].
\end{equation}
Thus there exists \(j_0\) such that
\eqref{eq:outline-common-support} holds for every \(j\ge j_0\).
The tail sequence \((n_j,P_j,t_j)_{j\ge j_0}\) retains
the properties in \eqref{eq:outline-selected-sizes}
and has the common support bound required in
Proposition~\ref{prop:bounded-extremizers}.
For the remainder of this qualitative argument, we restrict
attention to \(j\ge j_0\).
The explicit sample-size threshold is established separately
in Appendix~\ref{app:effective}.

\subsection{Bounding the error by smoothing}
\label{subsec:outline-smoothing}

\textbf{We now use the uniform support bound to pass to a limiting
distribution \(P\). To identify this limit, we first need to
relate the approximation errors along the selected sequence
to properties of \(P\).
For this purpose, we derive a first-order expansion that is
uniform in the threshold, adding independent uniform noise
when \(P\) is lattice.
We then compare the original and smoothed distribution
functions to obtain the bound}
\[
 \limsup_{j\to\infty}\sqrt{n_j}\,\Delta_{n_j}(P_j)
 \le \frac{\phi(0)}6\{3h(P)+|\kappa(P)|\},
\]
\textbf{where \(h(P)\) is the maximal lattice span, taken to be
zero for a nonlattice law, and \(\kappa(P)=\int x^3P(dx)\).
In the next subsection, we will combine this bound with the
maximizing property of \(P_j\) to force equality in Esseen's
moment inequality and thereby identify \(P\).}

The common compact support in \eqref{eq:outline-common-support}
allows us to pass to a subsequence for which \(P_j\) converges
weakly to a mean-zero, variance-one law \(P\), and
\(\beta_j\to\beta(P)\).
Consequently, \(W_3(P_j,P)\to0\).\footnote{
For probability laws \(\mu,\nu\) on \(\mathbb R\) with finite
third absolute moments, the Wasserstein distance of order three is
\(W_3(\mu,\nu)
=\bigl(\inf_{\pi\in\Pi(\mu,\nu)}
\int_{\mathbb R^2}|x-y|^3\,\pi(dx,dy)\bigr)^{1/3}\),
where \(\Pi(\mu,\nu)\) consists of all probability measures
on \(\mathbb R^2\) with marginals \(\mu\) and \(\nu\).
See \cite[Section~5.1, p.~179]{Santambrogio2015}.
Convergence in this distance is equivalent to weak convergence
together with convergence of third absolute moments;
see \cite[Theorem~5.11, pp.~185--187]{Santambrogio2015}.}

We use this convergence to obtain a first-order error bound
through the smoothing expansion in Lemma~\ref{lem:jitter}.
The smoothing variable is determined by the limiting law \(P\). We call \(P\) lattice if
\(\operatorname{supp}P\subseteq a+h\mathbb Z\) for some
\(a\in\mathbb R\) and \(h>0\), where
\(a+h\mathbb Z=\{a+kh:k\in\mathbb Z\}\).
Its maximal lattice span \(h(P)\) is the largest such \(h\);
we set \(h(P)=0\) when \(P\) is nonlattice.
Let \(U\) be independent of the summands and uniformly distributed
on \([-h(P)/2,h(P)/2]\) when \(h(P)>0\), and set \(U=0\)
when \(h(P)=0\).

For each \(j\), write
\[
 S_j=\frac{1}{\sqrt{n_j}}\sum_{r=1}^{n_j}X_{j,r},
 \qquad
 X_{j,1},\ldots,X_{j,n_j}\overset{\mathrm{iid}}{\sim}P_j.
\]
Put \(\mu_j=\int x^3P_j(dx)\).
Since \(n_j\to\infty\) and \(W_3(P_j,P)\to0\),
Lemma~\ref{lem:jitter} gives the first-order expansion
\eqref{eq:jitter-expansion}, namely,
\begin{equation}\label{eq:outline-jitter-expansion}
 \sup_{x\in\mathbb R}
 \left|
 \PP\{S_j+U/\sqrt{n_j}\le x\}
 -\Phi(x)
 -\frac{\mu_j}{6\sqrt{n_j}}(1-x^2)\phi(x)
 \right|
 =o(n_j^{-1/2}).
\end{equation}
Put \(a_j=h(P)/(2\sqrt{n_j})\).
Since \(|U|/\sqrt{n_j}\le a_j\), we have
\[
 \PP\{S_j+U/\sqrt{n_j}\le x-a_j\}
 \le \PP\{S_j\le x\}
 \le \PP\{S_j+U/\sqrt{n_j}\le x+a_j\}.
\]
Applying \eqref{eq:outline-jitter-expansion} at \(x+a_j\)
to the right-hand side and using a uniform Taylor expansion
gives the following upper bound, as in
Corollary~\ref{cor:jitter-envelopes}:
\eqref{eq:jitter-upper}:
\begin{equation}\label{eq:outline-jitter-upper}
 \PP\{S_j\le x\}-\Phi(x)
 \le
 \frac{1}{\sqrt{n_j}}
 \left[
   \frac{h(P)}2+\frac{\mu_j}{6}(1-x^2)
 \right]\phi(x)
 +o(n_j^{-1/2}),
\end{equation}
uniformly in \(x\in\mathbb R\).
The term \(h(P)/2\) arises because the upper bound evaluates
the smoothed distribution function at \(x+h(P)/(2\sqrt{n_j})\);
the corresponding normal term satisfies
\[
 \Phi\left(x+\frac{h(P)}{2\sqrt{n_j}}\right)-\Phi(x)
 =
 \frac{h(P)}{2\sqrt{n_j}}\phi(x)+O(n_j^{-1}),
\]
uniformly in \(x\).

To bound the absolute approximation error
\(\Delta_{n_j}(P_j)\), we also need a lower bound
on \(\PP\{S_j\le x\}-\Phi(x)\).
For this, we use the other side of the preceding
distribution-function comparison:
\[
 \PP\{S_j\le x\}
 \ge \PP\{S_j+U/\sqrt{n_j}\le x-a_j\}.
\]
Applying \eqref{eq:outline-jitter-expansion} at \(x-a_j\)
to the right-hand side and using a uniform Taylor
expansion gives
\[
 \PP\{S_j\le x\}-\Phi(x)
 \ge
 \frac{1}{\sqrt{n_j}}
 \left[
   -\frac{h(P)}2+\frac{\mu_j}{6}(1-x^2)
 \right]\phi(x)
 +o(n_j^{-1/2}),
\]
uniformly in \(x\in\mathbb R\).

We continue with the subsequence \((P_j)\) converging
to \(P\) in \(W_3\) that was used in
\eqref{eq:outline-jitter-expansion}.
Together with \eqref{eq:outline-selected-sizes} and
\eqref{eq:outline-common-support}, this gives
\begin{equation}\label{eq:outline-array-convergence}
 n_j\longrightarrow\infty,\qquad
 P_j\longrightarrow P\text{ in }W_3,\qquad
 \operatorname{supp}P_j\subseteq[-L,L].
\end{equation}
Here the sample sizes and the support bound are those in
\eqref{eq:outline-selected-sizes} and
\eqref{eq:outline-common-support}, respectively.
For subsequent reference, put \(a_j=h(P)/(2\sqrt{n_j})\)
and recall the distribution-function comparison
\begin{equation}\label{eq:outline-cdf-comparison}
 \PP\{S_j+U/\sqrt{n_j}\le x-a_j\}
 \le\PP\{S_j\le x\}
 \le\PP\{S_j+U/\sqrt{n_j}\le x+a_j\}.
\end{equation}
Applying the smoothing expansion
\eqref{eq:outline-jitter-expansion} to the right and left
sides of \eqref{eq:outline-cdf-comparison}, respectively,
gives the upper bound \eqref{eq:outline-jitter-upper} and
the corresponding lower bound. We record both estimates here:
\begin{align}
 \PP\{S_j\le x\}-\Phi(x)
 &\le \frac1{\sqrt{n_j}}
 \left[\frac{h(P)}2+\frac{\mu_j}{6}(1-x^2)\right]\phi(x)
 +o(n_j^{-1/2}),\label{eq:outline-cdf-upper}\\
 \PP\{S_j\le x\}-\Phi(x)
 &\ge \frac1{\sqrt{n_j}}
 \left[-\frac{h(P)}2+\frac{\mu_j}{6}(1-x^2)\right]\phi(x)
 +o(n_j^{-1/2}),\label{eq:outline-cdf-lower}
\end{align}
both uniformly in \(x\in\mathbb R\).
We then combine \eqref{eq:outline-cdf-upper} and
\eqref{eq:outline-cdf-lower} into a bound for the absolute error.
This bound will allow us to identify the limiting distribution
in the next subsection.

By the uniform smoothing expansion
\eqref{eq:outline-jitter-expansion} and the uniform Taylor
estimates used to derive \eqref{eq:outline-cdf-upper}
and \eqref{eq:outline-cdf-lower}, the remainder terms
in both bounds are \(o(n_j^{-1/2})\), uniformly in \(x\), after multiplication by \(\sqrt{n_j}\)
they can be bounded in absolute value by a common sequence
\(\varepsilon_j\ge0\), with \(\varepsilon_j\to0\).
Thus, for every \(x\in\mathbb R\),
\begin{equation}\label{eq:outline-centered-two-sided-bound}
\begin{aligned}
 -\frac{h(P)}2\phi(x)-\varepsilon_j
 &\le
 \sqrt{n_j}\bigl[\PP\{S_j\le x\}-\Phi(x)\bigr]
 -\frac{\mu_j}{6}(1-x^2)\phi(x)\\
 &\le \frac{h(P)}2\phi(x)+\varepsilon_j.
\end{aligned}
\end{equation}
The two inequalities in \eqref{eq:outline-centered-two-sided-bound}
imply
\begin{equation}\label{eq:outline-corrected-error-bound}
 \left|
 \sqrt{n_j}\bigl[\PP\{S_j\le x\}-\Phi(x)\bigr]
 -\frac{\mu_j}{6}(1-x^2)\phi(x)
 \right|
 \le \frac{h(P)}2\phi(x)+\varepsilon_j.
\end{equation}
We write the original error as the sum of the corrected error
in \eqref{eq:outline-corrected-error-bound} and the third-moment
correction term. The triangle inequality and
\eqref{eq:outline-corrected-error-bound} therefore give
\begin{equation}\label{eq:outline-pointwise-absolute-bound}
\begin{aligned}
 &\sqrt{n_j}\left|\PP\{S_j\le x\}-\Phi(x)\right|\\
 &\quad=
 \left|
 \left(
 \sqrt{n_j}\bigl[\PP\{S_j\le x\}-\Phi(x)\bigr]
 -\frac{\mu_j}{6}(1-x^2)\phi(x)
 \right)
 +\frac{\mu_j}{6}(1-x^2)\phi(x)
 \right|\\
 &\quad\le
 \left|
 \sqrt{n_j}\bigl[\PP\{S_j\le x\}-\Phi(x)\bigr]
 -\frac{\mu_j}{6}(1-x^2)\phi(x)
 \right|
 +\frac{|\mu_j|}{6}|1-x^2|\phi(x)\\
 &\quad\le
 \frac{h(P)}2\phi(x)
 +\frac{|\mu_j|}{6}|1-x^2|\phi(x)
 +\varepsilon_j.
\end{aligned}
\end{equation}
By definition,
\begin{equation}\label{eq:outline-kolmogorov-error}
 \Delta_{n_j}(P_j)
 =\sup_{x\in\mathbb R}
   \left|\PP\{S_j\le x\}-\Phi(x)\right|.
\end{equation}
We also have the elementary identities
\begin{equation}\label{eq:outline-normal-suprema}
 \sup_{x\in\mathbb R}\phi(x)
 =\sup_{x\in\mathbb R}|1-x^2|\phi(x)
 =\phi(0).
\end{equation}
For the second identity in \eqref{eq:outline-normal-suprema},
on \(|x|\le1\) the expression is at
most \(\phi(0)\), with equality at zero. On \(|x|\ge1\),
differentiating \((x^2-1)\phi(x)\) shows that its maximum is
\(2\phi(\sqrt3)=2e^{-3/2}\phi(0)<\phi(0)\).
Taking the supremum over \(x\) in
\eqref{eq:outline-pointwise-absolute-bound}, and using the
definition \eqref{eq:outline-kolmogorov-error} and the identities
\eqref{eq:outline-normal-suprema}, yields
\begin{equation}\label{eq:outline-supremum-error-bound}
\begin{aligned}
 \sqrt{n_j}\,\Delta_{n_j}(P_j)
 &\le
 \frac{h(P)}2\sup_{x\in\mathbb R}\phi(x)
 +\frac{|\mu_j|}{6}\sup_{x\in\mathbb R}|1-x^2|\phi(x)
 +\varepsilon_j\\
 &=\frac{\phi(0)}6\{3h(P)+|\mu_j|\}+\varepsilon_j.
\end{aligned}
\end{equation}
Since \(\varepsilon_j\to0\),
\eqref{eq:outline-supremum-error-bound} proves
\begin{equation}\label{eq:outline-absolute-error-bound}
 \sqrt{n_j}\,\Delta_{n_j}(P_j)
 \le\frac{\phi(0)}6\{3h(P)+|\mu_j|\}+o(1).
\end{equation}

Write \(\kappa(P)=\int x^3P(dx)\).
The common compact support and convergence in
\eqref{eq:outline-array-convergence} give the moment limits
\begin{equation}\label{eq:outline-limiting-moments}
 \mu_j=\int x^3P_j(dx)\longrightarrow\kappa(P),
 \qquad \beta_j=\int|x|^3P_j(dx)\longrightarrow\beta(P).
\end{equation}
Indeed, both \(x^3\) and \(|x|^3\) are continuous and bounded
on \([-L,L]\).
Taking the upper limit in \eqref{eq:outline-absolute-error-bound}
and using \eqref{eq:outline-limiting-moments} gives
\begin{equation}\label{eq:outline-triangular-upper}
 \limsup_{j\to\infty}\sqrt{n_j}\,\Delta_{n_j}(P_j)
 \le\frac{\phi(0)}6\{3h(P)+|\kappa(P)|\}.
\end{equation}
The bound \eqref{eq:outline-triangular-upper} is recorded later
as \eqref{eq:triangular-upper}
in Corollary~\ref{cor:jitter-envelopes}.

\subsection{Identifying the limiting distribution}
\label{subsec:outline-limit-identification}

\textbf{We now combine the upper bound
\eqref{eq:outline-triangular-upper} with the maximizing identity
\eqref{eq:outline-maximizing-pair} and Esseen's moment inequality.
These relations force equality in the limit and restrict
\(P\) to the Esseen two-point law or its reflection.
Recall that, without loss of generality, the maximizing pairs
were chosen to have positive signed error, as justified by
Lemma~\ref{lem:attainment} and recorded in
\eqref{eq:outline-maximizing-pair}.
We use this choice first to exclude the reflected law and identify
the limit as the Esseen two-point law, and then to show that
\(t_j/\sqrt{n_j}\to0\).}

The maximizing identity \eqref{eq:outline-maximizing-pair},
applied at the selected sample sizes, also gives
\begin{equation}\label{eq:outline-maximal-error-limit}
 \sqrt{n_j}\,\Delta_{n_j}(P_j)
 =C_{n_j}\beta_j\longrightarrow\ce\beta(P).
\end{equation}
Indeed, \(C_{n_j}\to\ce\) by \eqref{eq:outline-cn-bounds},
and \(\beta_j\to\beta(P)\) by
\eqref{eq:outline-limiting-moments}.
For the limiting mean-zero, variance-one law \(P\),
Esseen's moment inequality, stated in
\eqref{eq:esseen-moment} below, gives
\begin{equation}\label{eq:outline-esseen-moment}
 3h(P)+|\kappa(P)|\le(3+\sqrt{10})\beta(P).
\end{equation}
Since \(\ce=(3+\sqrt{10})\phi(0)/6\), combining
\eqref{eq:outline-triangular-upper},
\eqref{eq:outline-maximal-error-limit}, and
\eqref{eq:outline-esseen-moment} yields
\begin{equation}\label{eq:outline-limit-equality}
 \ce\beta(P)
 =\lim_{j\to\infty}\sqrt{n_j}\,\Delta_{n_j}(P_j)
 \le\frac{\phi(0)}6\{3h(P)+|\kappa(P)|\}
 \le\ce\beta(P).
\end{equation}
The endpoints of \eqref{eq:outline-limit-equality} coincide,
so both inequalities in that display are equalities.
In particular, equality holds in
\eqref{eq:outline-esseen-moment}. The equality characterization in Esseen's moment inequality
\cite[Lemma~1, pp.~162--163]{Esseen1956}
therefore implies that \(P\) is the Esseen two-point law
or its reflection.

To specify these two possibilities more explicitly, we state here
the parameters and law given in
\eqref{eq:esseen-parameters}--\eqref{eq:esseen-law} below:
\begin{equation}\label{eq:outline-esseen-parameters}
 \pe=\frac{4-\sqrt{10}}2,\qquad
 \qe=1-\pe,\qquad \se=\sqrt{\pe\qe},\qquad
 \aE=\frac{\pe}{\se},\quad \bE=\frac{\qe}{\se},
\end{equation}
and
\begin{equation}\label{eq:outline-esseen-law}
 \Pe=\qe\delta_{-\aE}+\pe\delta_{\bE}.
\end{equation}
Here \(0<\pe<1/2<\qe<1\).
The moment identities in \eqref{eq:esseen-identities} take
the following form:
\begin{equation}\label{eq:outline-esseen-identities}
\begin{aligned}
 h_{\mathrm E}:=h(\Pe)&=\aE+\bE=\frac1{\se},\\
 \beta_{\mathrm E}:=\beta(\Pe)&=\frac{\pe^2+\qe^2}{\se},
 &\kappa_{\mathrm E}:=\kappa(\Pe)&=\frac{\qe-\pe}{\se}>0,\\
 M_{\mathrm E}:=\int x|x|\,\Pe(dx)&=\qe-\pe=\sqrt{10}-3,
 &\ce\beta_{\mathrm E}&=
 \phi(0)\left(\frac{h_{\mathrm E}}2+
                   \frac{\kappa_{\mathrm E}}6\right).
\end{aligned}
\end{equation}
Reflection preserves the span and third absolute moment and
reverses the third moment. Thus the two candidates for \(P\)
have span \(h(P)=h_{\mathrm E}>0\), third absolute moment
\(\beta(P)=\beta_{\mathrm E}\), and third moments
\(\kappa_{\mathrm E}\) and \(-\kappa_{\mathrm E}\), respectively.

To distinguish the two possibilities, we use the fact that
the maximizing pairs \((P_j,t_j)\) were chosen with positive
signed error, as specified in
\eqref{eq:outline-maximizing-pair}. Put \(z_j=t_j/\sqrt{n_j}\) and define
\begin{equation}\label{eq:outline-positive-envelope}
 E_P(z):=\phi(z)\left[
           \frac{h(P)}2+\frac{\kappa(P)}6(1-z^2)\right].
\end{equation}
Multiplying the upper bound \eqref{eq:outline-cdf-upper}
by \(\sqrt{n_j}\) and using the moment convergence in
\eqref{eq:outline-limiting-moments} to replace \(\mu_j\)
by \(\kappa(P)\) gives
\begin{equation}\label{eq:outline-uniform-positive-envelope}
 \sqrt{n_j}\bigl[\PP\{S_j\le z\}-\Phi(z)\bigr]
 \le E_P(z)+o(1),
 \qquad z\in\mathbb R,
\end{equation}
uniformly in \(z\). To justify this replacement, the additional
term is bounded, using \eqref{eq:outline-normal-suprema} and
\eqref{eq:outline-limiting-moments}, by
\begin{equation}\label{eq:outline-third-moment-replacement}
 \frac{|\mu_j-\kappa(P)|}{6}
 \sup_{z\in\mathbb R}|1-z^2|\phi(z)
 =\frac{\phi(0)}6|\mu_j-\kappa(P)|\longrightarrow0.
\end{equation}
Equation~\eqref{eq:outline-uniform-positive-envelope} is the
upper envelope stated later in \eqref{eq:limit-upper-envelope}.
The uniformity of \eqref{eq:outline-uniform-positive-envelope}
allows us to evaluate it at the moving
threshold \(z_j\), whether or not \(z_j\) is bounded at this
stage. The positive signed maximizing identity
\eqref{eq:outline-maximizing-pair} then yields
\begin{equation}\label{eq:outline-limit-upper-envelope}
 C_{n_j}\beta_j
 =\sqrt{n_j}\bigl[\PP\{S_j\le z_j\}-\Phi(z_j)\bigr]
 \le E_P(z_j)+o(1).
\end{equation}

Put \(A=h_{\mathrm E}/2\) and \(B=\kappa_{\mathrm E}/6\).
The identities in \eqref{eq:outline-esseen-identities} give
\(h_{\mathrm E}>\kappa_{\mathrm E}>0\), and hence \(A>3B>0\).
According to whether \(P\) is \(\Pe\) or its reflection,
the upper envelope \(E_P\) is
\begin{equation}\label{eq:outline-candidate-envelopes}
 E_+(z)=\phi(z)(A+B-Bz^2)
 \quad\text{or}\quad
 E_-(z)=\phi(z)(A-B+Bz^2),
\end{equation}
respectively.

We first exclude the reflected law. Direct differentiation gives
\[
 E_-'(z)=z\phi(z)(3B-A-Bz^2).
\]
Since \(A>3B\), the function \(E_-\) is even and strictly
decreasing on \((0,\infty)\). Therefore
\begin{equation}\label{eq:outline-reflected-envelope-gap}
 \sup_{z\in\mathbb R}E_-(z)
 =\phi(0)(A-B)
 <\phi(0)(A+B)
 =\ce\beta_{\mathrm E}.
\end{equation}
If \(P\) were the reflection of \(\Pe\), then
\eqref{eq:outline-limit-upper-envelope} would imply
\[
 C_{n_j}\beta_j
 \le\sup_{z\in\mathbb R}E_-(z)+o(1).
\]
However, \eqref{eq:outline-maximal-error-limit} and
\(\beta(P)=\beta_{\mathrm E}\) give
\(C_{n_j}\beta_j\to\ce\beta_{\mathrm E}\), contradicting
\eqref{eq:outline-reflected-envelope-gap}.
Thus \(P=\Pe\).

It remains to show that \(z_j\to0\).
Now \(E_P=E_+\), and \(E_+\) has its unique global maximum
at zero. Indeed, for \(z\ne0\), if \(A+B-Bz^2\le0\), then
\(E_+(z)\le0<E_+(0)\). Otherwise, both the positive factor
\(A+B-Bz^2\) and \(\phi(z)\) are strictly smaller than their
values at zero, so again \(E_+(z)<E_+(0)\).
Moreover, \(E_+\) is continuous and tends to zero as
\(|z|\to\infty\). Consequently, for every \(\delta>0\),
\begin{equation}\label{eq:outline-positive-envelope-gap}
 \sup_{|z|\ge\delta}E_+(z)
 <E_+(0)=\ce\beta_{\mathrm E}.
\end{equation}
If \(|z_j|\ge\delta\) along an infinite subsequence, then
\eqref{eq:outline-limit-upper-envelope} would give, along
that subsequence,
\[
 C_{n_j}\beta_j
 \le\sup_{|z|\ge\delta}E_+(z)+o(1),
\]
contradicting \eqref{eq:outline-positive-envelope-gap}
and \(C_{n_j}\beta_j\to\ce\beta_{\mathrm E}\).
Hence \(z_j\to0\).

Together with \(P_j\to P\) in \(W_3\), this proves the
conclusion of Lemma~\ref{lem:limit-extremizer}:
\begin{equation}\label{eq:outline-limit-law}
 P_j\longrightarrow\Pe\quad\text{in }W_3,
 \qquad z_j=\frac{t_j}{\sqrt{n_j}}\longrightarrow0.
\end{equation}
In particular, \eqref{eq:outline-limit-law} and
\eqref{eq:outline-esseen-identities} give
\begin{equation}\label{eq:outline-esseen-moment-limits}
 \beta_j\longrightarrow\beta_{\mathrm E},\qquad
 \mu_j\longrightarrow\kappa_{\mathrm E},\qquad
 M_j:=\int x|x|\,P_j(dx)\longrightarrow M_{\mathrm E}.
\end{equation}
The last convergence in \eqref{eq:outline-esseen-moment-limits}
follows from weak convergence and continuity of \(x|x|\) on
the common compact support \eqref{eq:outline-common-support}.
The convergences in \eqref{eq:outline-limit-law} and
\eqref{eq:outline-esseen-moment-limits} persist along every
further subsequence.

\subsection{Locating the possible support limits}
\label{subsec:outline-support-interval}
\textbf{Although \(P_j\) converges to \(\Pe\), this convergence
does not exclude additional support points away from the two
Esseen atoms: neighborhoods of such points may carry
probabilities tending to zero.
We now return to the pointwise optimality condition
\eqref{eq:outline-contact}, which holds at every support point,
to restrict where such points can occur.
Our next goal is to prove that every convergent sequence
of support points has its limit in the following interval:}
\[
 y_j\in\operatorname{supp}P_j,\qquad y_j\to y
 \quad\Longrightarrow\quad
 y\in\left[-\frac{\sqrt{10}}2\aE,\,
                \frac{\sqrt{10}}2\bE\right].
\]
\textbf{To obtain this bound, we combine \(I_j(y_j)=0\)
with the influence-function upper bound
\eqref{eq:outline-influence-step}.
After controlling the remaining terms, we pass to the limit
and obtain a necessary inequality for \(y\), which gives
the endpoints of this interval.
This is the argument of Lemma~\ref{lem:support-interval}.
The interval still contains points other than the two Esseen
atoms. In the next subsection, we will derive a restriction
on the distances between distinct support limits.
We will then combine these two restrictions in
Subsection~\ref{subsec:outline-confinement} to exclude
all additional support limits.}

Fix any convergent sequence
\(y_j\in\operatorname{supp}P_j\), and write \(y_j\to y\).
Recall from \eqref{eq:outline-common-support},
\eqref{eq:outline-limit-law}, and
\eqref{eq:outline-selected-sizes} that the selected sequence satisfies
\begin{equation}\label{eq:outline-support-limit-assumptions}
\begin{gathered}
 \operatorname{supp}P_j\subseteq[-L,L],\qquad
 P_j\longrightarrow\Pe\text{ in }W_3,\qquad z_j\longrightarrow0,
 \\
 d_j:=n_j(C_{n_j-1}-C_{n_j})_+\longrightarrow0.
\end{gathered}
\end{equation}
In particular, \(\beta_j\to\beta_{\mathrm E}\) and
\(M_j\to M_{\mathrm E}\), as recorded in
\eqref{eq:outline-esseen-moment-limits}.

For the function \(H_n\) defined in
\eqref{eq:outline-gaussian-H}, the compact expansion in
Lemma~\ref{lem:gaussian-expansion} states that
\begin{equation}\label{eq:outline-compact-gaussian-expansion}
 H_n(z,u)
 =\phi''(z)\left(\frac u2-\frac{u^3}{6}\right)
   +O_L(n^{-1/2}),
 \qquad |u|\le L,\quad z\in\mathbb R,
\end{equation}
where the remainder is uniform over these \(u,z\).
Since \(\phi''(0)=-\phi(0)\), evaluating
\eqref{eq:outline-compact-gaussian-expansion} at
\((n,z,u)=(n_j,z_j,y_j)\) and using
\eqref{eq:outline-support-limit-assumptions} gives
\begin{equation}\label{eq:outline-sharp-gaussian-limit}
 H_{n_j}(z_j,y_j)
 \longrightarrow
 \frac{\phi(0)}6(y^3-3y).
\end{equation}

The support contact condition \eqref{eq:outline-contact} gives
\(I_j(y_j)=0\).
Thus \eqref{eq:outline-influence-step}, after expanding its last
two terms and collecting the constant terms, implies
\begin{equation}\label{eq:outline-support-contact-upper}
\begin{aligned}
 0=\beta_j I_j(y_j)\le{}&
 H_{n_j}(z_j,y_j)-C_{n_j}|y_j|^3
 +\frac32 C_{n_j}\beta_j y_j^2+3C_{n_j}M_jy_j\\
 &+\beta_j\left\{
 \frac{n_j^{3/2}}{\sqrt{n_j-1}}C_{n_j-1}
 -n_jC_{n_j}-\frac12C_{n_j}\right\}.
\end{aligned}
\end{equation}
To retain the constant-order cancellation in the second line of
\eqref{eq:outline-support-contact-upper}, use
\begin{equation}\label{eq:outline-square-root-expansion}
 n\left(\sqrt{\frac n{n-1}}-1\right)
 =\frac12+O(n^{-1}).
\end{equation}
Combining \eqref{eq:outline-square-root-expansion} with the uniform
bounds on \(\beta_j\) and \(C_n\) in
\eqref{eq:outline-maximizing-pair} and
\eqref{eq:outline-cn-bounds}, respectively, gives
\begin{equation}\label{eq:outline-recursive-cancellation}
\begin{aligned}
 &\beta_j\left\{
 \frac{n_j^{3/2}}{\sqrt{n_j-1}}C_{n_j-1}
 -n_jC_{n_j}-\frac12C_{n_j}\right\}\\
 &\quad=
 \beta_jn_j(C_{n_j-1}-C_{n_j})
 +\frac12\beta_j(C_{n_j-1}-C_{n_j})+O(n_j^{-1})\\
 &\quad\le \beta_jd_j+o(1)=o(1).
\end{aligned}
\end{equation}
Here \eqref{eq:outline-cn-bounds} gives
\(C_{n_j-1}\to\ce\) and \(C_{n_j}\to\ce\), so the
second term in the middle line of
\eqref{eq:outline-recursive-cancellation} tends to zero.
The signed first term need not tend to zero; its upper bound
\(\beta_jd_j\) tends to zero by
\eqref{eq:outline-maximizing-pair} and
\eqref{eq:outline-selected-sizes}.
In particular, the contribution
\(\beta_jC_{n_j-1}/2\) from the recursive term cancels, in the
limit, the constant \(-\beta_jC_{n_j}/2\) in
\eqref{eq:outline-support-contact-upper}.

Substituting the upper bound in
\eqref{eq:outline-recursive-cancellation} into
\eqref{eq:outline-support-contact-upper} gives
\begin{equation}\label{eq:outline-support-contact-prelimit}
\begin{aligned}
 0\le{}&
 H_{n_j}(z_j,y_j)-C_{n_j}|y_j|^3
 +\frac32 C_{n_j}\beta_j y_j^2
 +3C_{n_j}M_jy_j\\
 &+\beta_jd_j+o(1).
\end{aligned}
\end{equation}
The moment bound in \eqref{eq:outline-maximizing-pair}
and the sample-size selection in
\eqref{eq:outline-selected-sizes} imply
\(0\le\beta_jd_j\le B_*d_j\to0\).
Moreover, \(y_j\to y\) by our choice of the support points,
\(C_{n_j}\to\ce\) by \eqref{eq:outline-cn-bounds}, and
\(\beta_j\to\beta_{\mathrm E}\) and \(M_j\to M_{\mathrm E}\)
by \eqref{eq:outline-esseen-moment-limits}.
Together with the Gaussian limit
\eqref{eq:outline-sharp-gaussian-limit}, these convergences
allow us to pass to the limit in
\eqref{eq:outline-support-contact-prelimit}, yielding
\begin{equation}\label{eq:outline-limiting-contact-inequality}
 0\le
 \frac{\phi(0)}6(y^3-3y)-\ce|y|^3
 +\frac32\ce\beta_{\mathrm E}y^2+3\ce M_{\mathrm E}y.
\end{equation}
Put \(c_*:=6\ce/\phi(0)=3+\sqrt{10}\).
The Esseen parameters and identities in
\eqref{eq:outline-esseen-parameters}--\eqref{eq:outline-esseen-identities}
give
\begin{equation}\label{eq:outline-contact-constant-identities}
 c_*M_{\mathrm E}
 =(3+\sqrt{10})(\sqrt{10}-3)=1,
 \qquad
 c_*\beta_{\mathrm E}=\frac{\sqrt{10}}{\se}.
\end{equation}
Multiplying \eqref{eq:outline-limiting-contact-inequality} by
\(6/\phi(0)\) and using
\eqref{eq:outline-contact-constant-identities} cancels its linear terms:
\begin{equation}\label{eq:outline-contact-polynomial-expansion}
\begin{aligned}
 0
 &\le y^3-3y-c_*|y|^3
       +\frac32c_*\beta_{\mathrm E}y^2+3c_*M_{\mathrm E}y\\
 &=y^3-c_*|y|^3+\frac32c_*\beta_{\mathrm E}y^2.
\end{aligned}
\end{equation}
Rearranging \eqref{eq:outline-contact-polynomial-expansion} gives
\begin{equation}\label{eq:outline-contact-polynomial}
 c_*|y|^3-y^3-\frac32c_*\beta_{\mathrm E}y^2\le0.
\end{equation}
For \(y>0\), the left side of
\eqref{eq:outline-contact-polynomial} is
\(y^2\{(c_*-1)y-3c_*\beta_{\mathrm E}/2\}\), so
\begin{equation}\label{eq:outline-positive-support-root}
 y\le\frac{3c_*\beta_{\mathrm E}}{2(c_*-1)}
   =\frac{\sqrt{10}}2\bE.
\end{equation}
For \(y<0\), the left side of
\eqref{eq:outline-contact-polynomial} is
\(y^2\{(c_*+1)(-y)-3c_*\beta_{\mathrm E}/2\}\), so
\begin{equation}\label{eq:outline-negative-support-root}
 y\ge-\frac{3c_*\beta_{\mathrm E}}{2(c_*+1)}
   =-\frac{\sqrt{10}}2\aE.
\end{equation}
The equalities in \eqref{eq:outline-positive-support-root} and
\eqref{eq:outline-negative-support-root} use
\eqref{eq:outline-contact-constant-identities} together with
the following consequences of
\eqref{eq:outline-esseen-parameters} and
\eqref{eq:outline-esseen-law}:
\begin{equation}\label{eq:outline-support-root-identities}
 (c_*-1)\qe=3,\qquad (c_*+1)\pe=3,\qquad
 \aE=\pe/\se,\qquad \bE=\qe/\se.
\end{equation}
The case \(y=0\) also lies between the endpoints in
\eqref{eq:outline-positive-support-root} and
\eqref{eq:outline-negative-support-root}.
Together, these two bounds prove the following local form of
the limiting-support bound:
\begin{equation}\label{eq:outline-limiting-support-interval}
 y_j\in\operatorname{supp}P_j,\quad y_j\to y
 \quad\Longrightarrow\quad
 y\in\left[-\frac{\sqrt{10}}2\aE,
                \frac{\sqrt{10}}2\bE\right].
\end{equation}
The hypotheses in \eqref{eq:outline-support-limit-assumptions}
remain valid after taking a subsequence, so
\eqref{eq:outline-limiting-support-interval} applies to limits
along any subsequence as well.

\subsection{Separating limits of support points}
\label{subsec:outline-support-separation}

\textbf{The location bound
\eqref{eq:outline-limiting-support-interval} confines all
possible support limits to a fixed interval, but this interval
still contains points other than the two Esseen atoms.
We now derive an additional restriction on the distances
between distinct support limits.
For this purpose, we again use the pointwise optimality
condition \eqref{eq:outline-contact}.
Substituting \(I_j(x)=0\) into the influence formula
\eqref{eq:outline-influence} gives an exact relation between
the distribution function of the \((n_j-1)\)-term sum at
\(t_j-x\) and that of the \(n_j\)-term sum at \(t_j\).
These are the contact equations used below.
They hold at every support point, even when the probability
near that point is arbitrarily small.}

\textbf{We combine these equations with the smoothing
expansion to prove that, if
\(x_j,y_j\in\operatorname{supp}P_j\) converge to distinct
limits \(x\) and \(y\), then \(|x-y|\ge h_{\mathrm E}\).
The proof compares the contact equations at \(x_j\) and \(y_j\)
and shows that a limiting distance strictly between zero and
\(h_{\mathrm E}\) would contradict the smoothing estimate.
In the next subsection, we will combine this separation
property with the location bound
\eqref{eq:outline-limiting-support-interval} to show that
only the two Esseen atoms can occur as support limits.}

Throughout this argument, we work along the subsequence in
\eqref{eq:outline-limit-law}. Every convergence statement used below
also holds along any further subsequence.
Put \(h=h_{\mathrm E}\), and recall the identities in
\eqref{eq:outline-esseen-identities}.
For each \(j\), define
\begin{equation}\label{eq:outline-contact-notation}
 m_j=n_j-1,\qquad
 G_j(u)=F_{m_j,P_j}(u)
 =\PP\left\{\sum_{r=1}^{m_j}X_{j,r}\le u\right\}.
\end{equation}
Thus \(G_j\) is the distribution function of the unnormalized sum
of \(m_j\) iid variables with law \(P_j\). Since these variables
have mean zero and variance one, the corresponding normal
distribution function at the raw threshold \(u\) is
\(\Phi(u/\sqrt{m_j})\).

For every \(v\in\operatorname{supp}P_j\), the contact condition
\eqref{eq:outline-contact} gives \(I_j(v)=0\). The exact contact
equation \eqref{eq:contact-equation} therefore reads
\begin{equation}\label{eq:outline-contact-equation}
\begin{aligned}
 &n_j\{G_j(t_j-v)-F_{n_j,P_j}(t_j)\}\\
 &\quad=
 -\sqrt{n_j}\,\phi(z_j)v
 -\frac{z_j\phi(z_j)}2(v^2-1)\\
 &\qquad\quad+
 \frac{C_{n_j}}{\sqrt{n_j}}
 \left\{
 |v|^3-\beta_j-3M_jv-\frac32\beta_j(v^2-1)
 \right\}.
\end{aligned}
\end{equation}
Here \(z_j=t_j/\sqrt{n_j}\). The common support bound
\eqref{eq:outline-common-support}, the moment bound in
\eqref{eq:outline-maximizing-pair}, and the constant bound in
\eqref{eq:outline-cn-bounds} give
\begin{equation}\label{eq:outline-contact-uniform-bounds}
 |v|\le L,\qquad
 |M_j|\le\int x^2P_j(dx)=1,\qquad
 \beta_j<B_*,\qquad C_{n_j}\le C_{\mathrm{BE}}.
\end{equation}
After division by \(n_j\), write the last two terms of
\eqref{eq:outline-contact-equation} as
\begin{equation}\label{eq:outline-contact-remainder}
 \mathcal R_j(v):=
 -\frac{z_j\phi(z_j)}{2n_j}(v^2-1)
 +\frac{C_{n_j}}{n_j^{3/2}}
 \left\{
 |v|^3-\beta_j-3M_jv-\frac32\beta_j(v^2-1)
 \right\}.
\end{equation}
Applying \eqref{eq:outline-contact-uniform-bounds} to the remainder
defined in \eqref{eq:outline-contact-remainder} gives
\begin{equation}\label{eq:outline-contact-remainder-bound}
\begin{aligned}
 \sup_{v\in\operatorname{supp}P_j}|\mathcal R_j(v)|
 &\le
 \frac{L^2+1}{2n_j}\sup_{z\in\mathbb R}|z\phi(z)|\\
 &\quad+
 \frac{C_{\mathrm{BE}}}{n_j^{3/2}}
 \left\{L^3+B_*+3L+\frac32B_*(L^2+1)\right\}
 =O(n_j^{-1}).
\end{aligned}
\end{equation}
In particular, choose any sequence
\(x_j\in\operatorname{supp}P_j\), and put
\begin{equation}\label{eq:outline-contact-thresholds}
 u_j=t_j-x_j,\qquad w_j=\frac{u_j}{\sqrt{m_j}}.
\end{equation}
Substituting \eqref{eq:outline-contact-thresholds} into
\eqref{eq:outline-contact-equation} and using the definition
\eqref{eq:outline-contact-remainder} and bound
\eqref{eq:outline-contact-remainder-bound} gives
\begin{equation}\label{eq:outline-contact-cdf}
\begin{aligned}
 G_j(u_j)
 &=F_{n_j,P_j}(t_j)
   -\frac{\phi(z_j)x_j}{\sqrt{n_j}}
   +\mathcal R_j(x_j)\\
 &=F_{n_j,P_j}(t_j)
   -\frac{\phi(z_j)x_j}{\sqrt{n_j}}
   +O(n_j^{-1}).
\end{aligned}
\end{equation}

The definitions \eqref{eq:outline-contact-notation} and
\eqref{eq:outline-contact-thresholds} give the normalization
\begin{equation}\label{eq:outline-contact-standardization}
\begin{aligned}
 w_j
 &=\frac{z_j-x_j/\sqrt{n_j}}{\sqrt{1-1/n_j}}\\
 &=z_j-\frac{x_j}{\sqrt{n_j}}+O(n_j^{-1}).
\end{aligned}
\end{equation}
Indeed, \((1-1/n_j)^{-1/2}=1+O(n_j^{-1})\),
\(|x_j|\le L\) by \eqref{eq:outline-contact-uniform-bounds}, and
\(z_j\to0\) by \eqref{eq:outline-limit-law}.
In particular, \eqref{eq:outline-contact-standardization} gives
\(w_j-z_j=O(n_j^{-1/2})\).
Using \eqref{eq:outline-contact-standardization} in Taylor's
theorem, together with boundedness of \(\phi\) and \(\phi'\), gives
\begin{equation}\label{eq:outline-contact-normal}
\begin{aligned}
 \Phi(w_j)
 &=\Phi(z_j)+\phi(z_j)(w_j-z_j)
   +O\bigl((w_j-z_j)^2\bigr)\\
 &=\Phi(z_j)-\frac{\phi(z_j)x_j}{\sqrt{n_j}}
   +O(n_j^{-1}).
\end{aligned}
\end{equation}
The remainders in \eqref{eq:outline-contact-cdf},
\eqref{eq:outline-contact-standardization}, and
\eqref{eq:outline-contact-normal} are uniform over
\(x_j\in\operatorname{supp}P_j\).

Subtracting \eqref{eq:outline-contact-normal} from
\eqref{eq:outline-contact-cdf} cancels the common term
\(-\phi(z_j)x_j/\sqrt{n_j}\). The difference of the two
\(O(n_j^{-1})\) remainders is still \(O(n_j^{-1})\), so
\begin{equation}\label{eq:outline-contact-error-difference}
 G_j(u_j)-\Phi(w_j)
 =
 F_{n_j,P_j}(t_j)-\Phi(z_j)+O(n_j^{-1}).
\end{equation}
Thus the normal approximation error for the \(m_j\)-term sum
at \(u_j\) differs from that for the \(n_j\)-term sum at
\(t_j\) by only \(O(n_j^{-1})\).

The positive signed maximizing identity
\eqref{eq:outline-maximizing-pair}, together with
\(z_j=t_j/\sqrt{n_j}\), gives the exact value
\begin{equation}\label{eq:outline-contact-maximizing-error}
 F_{n_j,P_j}(t_j)-\Phi(z_j)
 =
 \frac{C_{n_j}\beta_j}{\sqrt{n_j}}.
\end{equation}
Multiplying \eqref{eq:outline-contact-error-difference}
by \(\sqrt{m_j}\) and then substituting
\eqref{eq:outline-contact-maximizing-error}, we obtain
\begin{equation}\label{eq:outline-contact-saturation}
\begin{aligned}
 \sqrt{m_j}\{G_j(u_j)-\Phi(w_j)\}
 &=
 \sqrt{m_j}
 \{F_{n_j,P_j}(t_j)-\Phi(z_j)\}
 +O\!\left(\frac{\sqrt{m_j}}{n_j}\right)\\
 &=
 \sqrt{m_j}\,
 \frac{C_{n_j}\beta_j}{\sqrt{n_j}}
 +O\!\left(\frac{\sqrt{m_j}}{n_j}\right)\\
 &=
 \sqrt{\frac{m_j}{n_j}}\,C_{n_j}\beta_j
 +O(n_j^{-1/2})\\
 &\longrightarrow \ce\beta_{\mathrm E}.
\end{aligned}
\end{equation}
Here the last equality uses \(m_j=n_j-1\), which implies
\(\sqrt{m_j}/n_j=O(n_j^{-1/2})\).
The limit in \eqref{eq:outline-contact-saturation} uses
\(m_j/n_j=1-1/n_j\to1\) from
\eqref{eq:outline-contact-notation},
\(C_{n_j}\to\ce\) from \eqref{eq:outline-cn-bounds}, and
\(\beta_j\to\beta_{\mathrm E}\) from the \(W_3\) convergence in
\eqref{eq:outline-limit-law}.
Thus \eqref{eq:outline-contact-saturation} shows that
the signed error at each such threshold, multiplied by
\(\sqrt{m_j}\), converges to \(\ce\beta_{\mathrm E}\).

We now compare the error in \eqref{eq:outline-contact-saturation}
with the error for the smoothed sum.
Let \(U\) be independent of the summands and uniform on
\([-h/2,h/2]\), and define its raw distribution function by
\begin{equation}\label{eq:outline-contact-jitter-cdf}
 J_j(u)=
 \PP\left\{\sum_{r=1}^{m_j}X_{j,r}+U\le u\right\}.
\end{equation}
By \eqref{eq:outline-contact-notation}, \(m_j\to\infty\), and
\eqref{eq:outline-limit-law} gives \(P_j\to\Pe\) in \(W_3\).
Lemma~\ref{lem:jitter}, applied with sample size \(m_j\), gives
\begin{equation}\label{eq:outline-previous-jitter-expansion}
 \sup_{w\in\mathbb R}
 \left|
 J_j(\sqrt{m_j}w)-\Phi(w)
 -\frac{\mu_j}{6\sqrt{m_j}}(1-w^2)\phi(w)
 \right|
 =o(m_j^{-1/2}).
\end{equation}
Here the third moment is still
\(\mu_j=\int x^3P_j(dx)\), since only the number of summands
has changed. The \(W_3\) convergence in
\eqref{eq:outline-limit-law} gives
\(\mu_j\to\kappa_{\mathrm E}\), and
\eqref{eq:outline-contact-standardization}, together with
\(z_j\to0\), gives \(w_j\to0\).

Put \(\delta_j=h/(2\sqrt{m_j})\).
By \eqref{eq:outline-contact-thresholds}, the raw threshold
\(u_j+h/2\) corresponds to the standardized threshold
\(w_j+\delta_j\). Since \(\phi'\) and the derivative
of \((1-w^2)\phi(w)\) are bounded, we have
\begin{equation}\label{eq:outline-contact-normal-shift}
 \Phi(w_j+\delta_j)-\Phi(w_j)
 =\delta_j\phi(w_j)+O(\delta_j^2)
 =\frac{h}{2\sqrt{m_j}}\phi(w_j)+O(m_j^{-1}),
\end{equation}
For the third-moment correction term, put
\(g(w)=(1-w^2)\phi(w)\).
Its derivative \(g'(w)=(w^3-3w)\phi(w)\) is bounded
on \(\mathbb R\). Hence the mean value theorem gives
\[
 |g(w_j+\delta_j)-g(w_j)|
 \le \|g'\|_\infty\,\delta_j.
\]
Moreover, the moment convergence in
\eqref{eq:outline-esseen-moment-limits} gives
\(\mu_j\to\kappa_{\mathrm E}\), and therefore
\(\sup_j|\mu_j|<\infty\).
Multiplying the preceding estimate by
\(|\mu_j|/(6\sqrt{m_j})\) and using
\(\delta_j=h/(2\sqrt{m_j})\), we obtain
\begin{equation}\label{eq:outline-contact-correction-shift}
\begin{aligned}
 &\frac{\mu_j}{6\sqrt{m_j}}
   \{1-(w_j+\delta_j)^2\}\phi(w_j+\delta_j)\\
 &\qquad=
 \frac{\mu_j}{6\sqrt{m_j}}(1-w_j^2)\phi(w_j)
 +O\!\left(\frac{\delta_j}{\sqrt{m_j}}\right)\\
 &\qquad=
 \frac{\mu_j}{6\sqrt{m_j}}(1-w_j^2)\phi(w_j)
 +O(m_j^{-1}).
\end{aligned}
\end{equation}
Thus shifting the threshold changes the third-moment
correction by only \(O(m_j^{-1})\), which becomes \(o(1)\)
after multiplication by \(\sqrt{m_j}\).
Evaluating \eqref{eq:outline-previous-jitter-expansion} at
\(w=w_j+\delta_j\) and substituting
\eqref{eq:outline-contact-normal-shift} and
\eqref{eq:outline-contact-correction-shift} yields
\begin{equation}\label{eq:outline-smoothed-contact}
\begin{aligned}
 &\sqrt{m_j}\{J_j(u_j+h/2)-\Phi(w_j)\}\\
 &\qquad=
 \left[\frac h2+\frac{\mu_j}{6}(1-w_j^2)\right]\phi(w_j)
 +o(1)\\
 &\qquad\longrightarrow
 \phi(0)\left(\frac h2+\frac{\kappa_{\mathrm E}}6\right)
 =\ce\beta_{\mathrm E}.
\end{aligned}
\end{equation}
The last equality in \eqref{eq:outline-smoothed-contact} is the
Esseen identity recorded in
\eqref{eq:outline-esseen-identities}.
Subtracting \eqref{eq:outline-contact-saturation} from
\eqref{eq:outline-smoothed-contact} gives
\begin{equation}\label{eq:outline-jitter-saturation}
 \sqrt{m_j}\{J_j(u_j+h/2)-G_j(u_j)\}\longrightarrow0.
\end{equation}

To use \eqref{eq:outline-jitter-saturation}, we express
the smoothed distribution function as an average of the
original distribution function.
Write \(T_j=\sum_{r=1}^{m_j}X_{j,r}\). By
\eqref{eq:outline-contact-notation} and
\eqref{eq:outline-contact-jitter-cdf},
\(G_j(a)=\PP\{T_j\le a\}\) and
\(J_j(a)=\PP\{T_j+U\le a\}\).

For a given value \(U=v\), the event \(T_j+U\le a\)
becomes \(T_j\le a-v\). Since \(T_j\) and \(U\) are
independent, the conditional distribution of \(T_j\)
is unchanged, so the corresponding conditional probability
is \(G_j(a-v)\).
To obtain the unconditional probability \(J_j(a)\),
we then average these conditional probabilities over
the distribution of \(U\).
The law of total probability and the density \(1/h\)
of \(U\) on \([-h/2,h/2]\) therefore give
\begin{equation}\label{eq:outline-jitter-convolution}
\begin{aligned}
 J_j(u_j+h/2)
 &=
 \mathbb E\!\left[
   \PP\{T_j+U\le u_j+h/2\mid U\}
 \right]\\
 &=
 \mathbb E\!\left[G_j(u_j+h/2-U)\right]\\
 &=
 \frac1h\int_{-h/2}^{h/2}
 G_j(u_j+h/2-v)\,dv\\
 &=
 \frac1h\int_0^h G_j(u_j+s)\,ds.
\end{aligned}
\end{equation}
The last equality in \eqref{eq:outline-jitter-convolution}
uses the change of variables \(s=h/2-v\). Subtracting
\(G_j(u_j)=h^{-1}\int_0^h G_j(u_j)\,ds\)
from \eqref{eq:outline-jitter-convolution} gives
\begin{equation}\label{eq:outline-jitter-contact-integral}
 J_j(u_j+h/2)-G_j(u_j)
 =
 \frac1h\int_0^h
 \{G_j(u_j+s)-G_j(u_j)\}\,ds.
\end{equation}

The integrand in \eqref{eq:outline-jitter-contact-integral} is
nonnegative because \(s\ge0\) and \(G_j\) is nondecreasing.
For any sequence
\(\ell_j\to\ell\in(0,h)\), we have \(0<\ell_j<h\) for all
sufficiently large \(j\). Restricting the integral in
\eqref{eq:outline-jitter-contact-integral} to \([\ell_j,h]\)
and applying monotonicity gives
\begin{equation}\label{eq:outline-monotone-integral-bound}
\begin{aligned}
 J_j(u_j+h/2)-G_j(u_j)
 &\ge\frac1h\int_{\ell_j}^{h}
       \{G_j(u_j+s)-G_j(u_j)\}\,ds\\
 &\ge\frac1h\int_{\ell_j}^{h}
       \{G_j(u_j+\ell_j)-G_j(u_j)\}\,ds\\
 &=\frac{h-\ell_j}{h}
       \{G_j(u_j+\ell_j)-G_j(u_j)\}.
\end{aligned}
\end{equation}
Rearranging \eqref{eq:outline-monotone-integral-bound} yields
\begin{equation}\label{eq:outline-short-interval-bound}
 0\le G_j(u_j+\ell_j)-G_j(u_j)
 \le\frac{h}{h-\ell_j}
       \{J_j(u_j+h/2)-G_j(u_j)\}.
\end{equation}
The factor \(h/(h-\ell_j)\) in
\eqref{eq:outline-short-interval-bound} tends to the finite number
\(h/(h-\ell)\). Multiplying
\eqref{eq:outline-short-interval-bound} by \(\sqrt{m_j}\) and applying
\eqref{eq:outline-jitter-saturation} therefore gives
\begin{equation}\label{eq:outline-flat-contact-interval}
\begin{aligned}
 &\sqrt{m_j}\,
 \PP\left\{u_j<\sum_{r=1}^{m_j}X_{j,r}\le u_j+\ell_j\right\}\\
 &\qquad=
 \sqrt{m_j}\{G_j(u_j+\ell_j)-G_j(u_j)\}
 \longrightarrow0.
\end{aligned}
\end{equation}
The conclusion \eqref{eq:outline-flat-contact-interval} holds
for every choice of support points
\(x_j\), with \(u_j=t_j-x_j\), and every such sequence
\(\ell_j\). No convergence of \(x_j\) was needed.

To prove separation, suppose that
\(x_j,y_j\in\operatorname{supp}P_j\),
\(x_j\to x\), \(y_j\to y\), and \(0<y-x<h\).
Set \(\ell_j=y_j-x_j\), so that \(\ell_j\to y-x\in(0,h)\).
Apply \eqref{eq:outline-flat-contact-interval} with the support
point \(y_j\), whose contact threshold is \(t_j-y_j\).
Since \((t_j-y_j)+\ell_j=t_j-x_j\), this gives
\begin{equation}\label{eq:outline-flat-increment}
 \sqrt{m_j}\{G_j(t_j-x_j)-G_j(t_j-y_j)\}
 \longrightarrow0.
\end{equation}

On the other hand, subtracting
\eqref{eq:outline-contact-equation} at \(v=y_j\) from the same
equation at \(v=x_j\) cancels the common term
\(F_{n_j,P_j}(t_j)\) and gives
\begin{equation}\label{eq:outline-contact-subtraction}
\begin{aligned}
 &n_j\{G_j(t_j-x_j)-G_j(t_j-y_j)\}\\
 &\quad=
 -\sqrt{n_j}\,\phi(z_j)(x_j-y_j)
 -\frac{z_j\phi(z_j)}2(x_j^2-y_j^2)\\
 &\qquad\quad+
 \frac{C_{n_j}}{\sqrt{n_j}}
 \left\{
 |x_j|^3-|y_j|^3-3M_j(x_j-y_j)
 -\frac32\beta_j(x_j^2-y_j^2)
 \right\}.
\end{aligned}
\end{equation}
The second term on the right of
\eqref{eq:outline-contact-subtraction} is \(O(1)\), by
\(|x_j|,|y_j|\le L\) from
\eqref{eq:outline-contact-uniform-bounds} and boundedness of
\(z\phi(z)\). The expression in braces in
\eqref{eq:outline-contact-subtraction} is uniformly bounded by
\eqref{eq:outline-contact-uniform-bounds}, so its last term is
\(O(n_j^{-1/2})\). Multiplying
\eqref{eq:outline-contact-subtraction} by \(\sqrt{m_j}/n_j\)
and using \(m_j/n_j\to1\) from
\eqref{eq:outline-contact-notation} and \(z_j\to0\) from
\eqref{eq:outline-limit-law}, we obtain
\begin{equation}\label{eq:outline-positive-increment-limit}
\begin{aligned}
 &\sqrt{m_j}\{G_j(t_j-x_j)-G_j(t_j-y_j)\}\\
 &\qquad=
 \sqrt{\frac{m_j}{n_j}}\,\phi(z_j)(y_j-x_j)
 +O(n_j^{-1/2})\\
 &\qquad\longrightarrow\phi(0)(y-x)>0.
\end{aligned}
\end{equation}
The positive limit in \eqref{eq:outline-positive-increment-limit}
contradicts \eqref{eq:outline-flat-increment}.
If \(0<x-y<h\), the same argument applies after interchanging
the two sequences. This proves the separation property in
Lemma~\ref{lem:support-separation}:
\begin{equation}\label{eq:outline-support-separation}
 x_j,y_j\in\operatorname{supp}P_j,\quad x_j\to x,\quad y_j\to y
 \quad\Longrightarrow\quad
 x=y\ \text{or}\ |x-y|\ge h_{\mathrm E}.
\end{equation}

\subsection{Confining the support near the two atoms}
\label{subsec:outline-confinement}
\textbf{We now combine the location bound
\eqref{eq:outline-limiting-support-interval} with the separation
property \eqref{eq:outline-support-separation} to show that
\(-\aE\) and \(\bE\) are the only possible limits of support
points. We will then use the uniform support bound
\eqref{eq:outline-common-support} to deduce
that the entire support approaches these two atoms, namely,}
\[
 \sup_{x\in\operatorname{supp}P_j}
 \operatorname{dist}\bigl(x,\{-\aE,\bE\}\bigr)
 \longrightarrow0.
\]

\textbf{To see how the two bounds work together, we change
coordinates by \(v=\pe+\se x\).
The two Esseen atoms become \(0\) and \(1\).
Since \(\se h_{\mathrm E}=1\), the separation property
\eqref{eq:outline-support-separation} requires distinct
support limits in the new coordinates to be at least
one unit apart.
Meanwhile, the location bound
\eqref{eq:outline-limiting-support-interval} places all
possible support limits in these coordinates inside
\((-1,2)\).
Weak convergence supplies support sequences approaching
both \(0\) and \(1\).
Every other point in this interval has distance strictly
between zero and one from at least one of these two points,
so the separation property excludes it.
This leaves only the two Esseen atoms as possible support
limits.}

We continue with the sequence of maximizing laws \(P_j\)
in \eqref{eq:outline-limit-law}.
For each \(j\), let \(X_j\sim P_j\), and define \(Q_j\) to be
the distribution of
\[
 Y_j:=\pe+\se X_j.
\]
Thus \(Q_j\) is obtained from \(P_j\) by applying the same affine
transformation \(v=\pe+\se x\) for every \(j\);
the parameters \(\pe\) and \(\se>0\) do not depend on \(j\).
By \eqref{eq:outline-esseen-law} and
\eqref{eq:outline-esseen-identities}, we have
\(\aE=\pe/\se\), \(\bE=\qe/\se\), and
\(h_{\mathrm E}=1/\se\). Consequently,
\begin{equation}\label{eq:outline-affine-atom-identities}
 \pe-\se\aE=\pe-\pe=0,\qquad
 \pe+\se\bE=\pe+\qe=1,\qquad
 \se h_{\mathrm E}=1.
\end{equation}

Thus the affine map \(x\mapsto\pe+\se x\) sends the two
support points \(-\aE\) and \(\bE\) of \(\Pe\) to \(0\)
and \(1\), respectively. In particular, if \(X\sim\Pe\),
then \(\pe+\se X\) has distribution
\(\qe\delta_0+\pe\delta_1\).
Since \(\se>0\), the supports satisfy
\begin{equation}\label{eq:outline-affine-support}
 \operatorname{supp}Q_j=\pe+\se\operatorname{supp}P_j.
\end{equation}
Combining \eqref{eq:outline-affine-support} and
\eqref{eq:outline-limiting-support-interval} places every
limit of transformed support points in the explicit interval
\begin{equation}\label{eq:outline-transformed-support-interval}
 \left[\pe-\frac{\sqrt{10}}2\pe,
       \pe+\frac{\sqrt{10}}2\qe\right]
 =\left[-\pe\qe,\frac92-\sqrt{10}\right]
 \subset(-1,2).
\end{equation}
By \eqref{eq:outline-affine-atom-identities} and
\eqref{eq:outline-affine-support}, the separation property
\eqref{eq:outline-support-separation} transforms into
\begin{equation}\label{eq:outline-transformed-separation}
\begin{gathered}
 v_j,\widetilde v_j\in\operatorname{supp}Q_j,\qquad
 v_j\to v,\quad \widetilde v_j\to\widetilde v\\
 \Longrightarrow\qquad v=\widetilde v\quad\text{or}\quad
 |v-\widetilde v|\ge1.
\end{gathered}
\end{equation}
Indeed, define the corresponding support points of \(P_j\) by
\[
 x_j:=\frac{v_j-\pe}{\se},
 \qquad
 \widetilde x_j:=\frac{\widetilde v_j-\pe}{\se}.
\]
By \eqref{eq:outline-affine-support},
\(x_j,\widetilde x_j\in\operatorname{supp}P_j\), and
\[
 x_j\longrightarrow x:=\frac{v-\pe}{\se},
 \qquad
 \widetilde x_j\longrightarrow
 \widetilde x:=\frac{\widetilde v-\pe}{\se}.
\]
If \(v\ne\widetilde v\), then \(x\ne\widetilde x\).
The separation property \eqref{eq:outline-support-separation}
and the identity \(\se h_{\mathrm E}=1\) in
\eqref{eq:outline-affine-atom-identities} therefore give
\[
 |v-\widetilde v|
 =\se|x-\widetilde x|
 \ge \se h_{\mathrm E}
 =1.
\]

The weak convergence in \eqref{eq:outline-limit-law}, together
with the affine identities \eqref{eq:outline-affine-atom-identities},
implies
\begin{equation}\label{eq:outline-transformed-weak-limit}
 Q_j\Longrightarrow Q_{\mathrm E}:=\qe\delta_0+\pe\delta_1.
\end{equation}
By the Portmanteau theorem
\cite[Theorem~2.1(iv)]{Billingsley1999},
the weak convergence in
\eqref{eq:outline-transformed-weak-limit} implies
\begin{equation}\label{eq:outline-portmanteau-open}
 \liminf_{j\to\infty}Q_j(G)\ge Q_{\mathrm E}(G)
 \qquad\text{for every open set }G\subseteq\mathbb R.
\end{equation}
For each fixed \(\varepsilon>0\), apply
\eqref{eq:outline-portmanteau-open} to the open intervals
\((-\varepsilon,\varepsilon)\) and
\((1-\varepsilon,1+\varepsilon)\).
Since \(Q_{\mathrm E}=\qe\delta_0+\pe\delta_1\), we obtain
\begin{equation}\label{eq:outline-principal-neighborhood-mass}
\begin{aligned}
 \liminf_{j\to\infty}Q_j((-\varepsilon,\varepsilon))
 &\ge Q_{\mathrm E}((-\varepsilon,\varepsilon))
 \ge\qe>0,\\
 \liminf_{j\to\infty}Q_j((1-\varepsilon,1+\varepsilon))
 &\ge Q_{\mathrm E}((1-\varepsilon,1+\varepsilon))
 \ge\pe>0.
\end{aligned}
\end{equation}
By \eqref{eq:outline-principal-neighborhood-mass}, both intervals
meet \(\operatorname{supp}Q_j\) for all sufficiently large \(j\).
Since \eqref{eq:outline-principal-neighborhood-mass} holds for every
\(\varepsilon>0\), the distances from \(0\) and \(1\) to
\(\operatorname{supp}Q_j\) tend to zero.  The support \(\operatorname{supp}Q_j\) is compact by
\eqref{eq:outline-common-support} and
\eqref{eq:outline-affine-support}.
Hence the distances from \(0\) and \(1\) to
\(\operatorname{supp}Q_j\) are each attained at a point
of \(\operatorname{supp}Q_j\).
We may consequently
choose support points \(v_j^{(0)},v_j^{(1)}\) such that
\begin{equation}\label{eq:outline-principal-support-sequences}
 v_j^{(0)},v_j^{(1)}\in\operatorname{supp}Q_j,
 \qquad v_j^{(0)}\to0,\quad v_j^{(1)}\to1.
\end{equation}
We now show that every convergent sequence of support points
of \(Q_j\), taken along any subsequence, has its limit in
\(\{0,1\}\).
Let \(j_1<j_2<\cdots\) be arbitrary indices, and suppose that
\[
 v_{j_k}\in\operatorname{supp}Q_{j_k},
 \qquad v_{j_k}\longrightarrow v.
\]
By \eqref{eq:outline-principal-support-sequences}, the points
\(v_{j_k}^{(0)}\) and \(v_{j_k}^{(1)}\) also belong to
\(\operatorname{supp}Q_{j_k}\), and they converge to \(0\)
and \(1\), respectively.
Thus the separation property
\eqref{eq:outline-transformed-separation} can be applied
to \(v_{j_k}\) together with either of these two sequences,
using the same indices \(j_k\).

The location bound
\eqref{eq:outline-transformed-support-interval} gives
\(v\in(-1,2)\).
Suppose, for a contradiction, that \(v\notin\{0,1\}\).
There are then three cases:
\begin{equation}\label{eq:outline-transformed-gap-cases}
\begin{array}{lll}
 -1<v<0 &: &0<|v-0|=-v<1,\\
 0<v<1 &: &0<|v-0|=v<1,\\
 1<v<2 &: &0<|v-1|=v-1<1.
\end{array}
\end{equation}
In the first two cases of
\eqref{eq:outline-transformed-gap-cases}, apply
\eqref{eq:outline-transformed-separation} to
\(v_{j_k}\) and \(v_{j_k}^{(0)}\).
Their limits are \(v\) and \(0\), which are distinct,
so this separation property requires \(|v-0|\ge1\).
This contradicts \eqref{eq:outline-transformed-gap-cases}.
In the third case, apply
\eqref{eq:outline-transformed-separation} to
\(v_{j_k}\) and \(v_{j_k}^{(1)}\).
Their distinct limits \(v\) and \(1\) must satisfy
\(|v-1|\ge1\), again contradicting
\eqref{eq:outline-transformed-gap-cases}.
Consequently, \(v\in\{0,1\}\).

We next express this conclusion in terms of the original
distributions \(P_j\).
If \(x_k\in\operatorname{supp}P_{j_k}\) and \(x_k\to x\),
then \eqref{eq:outline-affine-support} gives
\[
 \pe+\se x_k\in\operatorname{supp}Q_{j_k},
 \qquad
 \pe+\se x_k\longrightarrow\pe+\se x.
\]
The conclusion just proved therefore implies
\(\pe+\se x\in\{0,1\}\).
By \eqref{eq:outline-affine-atom-identities}, this is
equivalent to \(x\in\{-\aE,\bE\}\).
Since the indices \(j_k\) and the convergent support point
sequence were arbitrary, we have proved
\begin{equation}\label{eq:outline-only-support-limits}
 x_k\in\operatorname{supp}P_{j_k},\quad x_k\to x
 \quad\Longrightarrow\quad x\in\{-\aE,\bE\}
 \qquad\text{for every }j_k\to\infty.
\end{equation}

Combining the common compact support
\eqref{eq:outline-common-support} with
\eqref{eq:outline-only-support-limits} gives uniform confinement:
\begin{equation}\label{eq:outline-support-confinement}
 \sup_{x\in\operatorname{supp}P_j}
 \operatorname{dist}\bigl(x,\{-\aE,\bE\}\bigr)
 \longrightarrow0.
\end{equation}
Indeed, if \eqref{eq:outline-support-confinement} failed,
there would be \(\varepsilon>0\),
indices \(j_k\to\infty\), and
\(x_k\in\operatorname{supp}P_{j_k}\) such that
\begin{equation}\label{eq:outline-confinement-countersequence}
 \operatorname{dist}\bigl(x_k,\{-\aE,\bE\}\bigr)
 \ge\varepsilon\qquad\text{for every }k.
\end{equation}
By \eqref{eq:outline-common-support}, one has \(x_k\in[-L,L]\),
so a further subsequence would converge to some \(x\).
Continuity of the distance function and
\eqref{eq:outline-confinement-countersequence} would give
\(\operatorname{dist}(x,\{-\aE,\bE\})\ge\varepsilon\),
contradicting \eqref{eq:outline-only-support-limits}.
By \eqref{eq:outline-affine-atom-identities} and
\eqref{eq:outline-affine-support}, the confinement
\eqref{eq:outline-support-confinement} is equivalently expressed
for the affine images as
\begin{equation}\label{eq:outline-transformed-confinement}
\begin{aligned}
 &\sup_{v\in\operatorname{supp}Q_j}
       \operatorname{dist}(v,\{0,1\})\\
 &\qquad=
 \se\sup_{x\in\operatorname{supp}P_j}
       \operatorname{dist}\bigl(x,\{-\aE,\bE\}\bigr)
 \longrightarrow0.
\end{aligned}
\end{equation}

\subsection{Local comparison and the final contradiction}
\label{subsec:outline-local-comparison}

\textbf{We now complete the contradiction by applying the local
two-cluster comparison in Proposition~\ref{prop:clusters}.
By \eqref{eq:outline-transformed-confinement}, the support of
\(Q_j\) eventually lies in the required small intervals around
\(0\) and \(1\).
The comparison also requires the probability of the cluster
near \(1\) to be close to \(\pe\).
We will verify this using the weak convergence in
\eqref{eq:outline-transformed-weak-limit}.}

\textbf{To transfer the resulting bound from \(Q_j\) back to
the original maximizing distribution \(P_j\), we also verify
that the affine transformation preserves the normalized
Berry--Esseen error.
Here each distribution is normalized using its own mean
and variance.
Thus, for all sufficiently large \(j\), the local comparison
will give}
\[
 C_{n_j}=R_{n_j}(P_j)=R_{n_j}(Q_j)\le\ce,
\]
\textbf{contradicting \(C_{n_j}>\ce\) in
\eqref{eq:outline-selected-sizes}.
This will complete the proof of the existence of a universal
finite threshold.}

More precisely, Proposition~\ref{prop:clusters} states that
there exist fixed constants \(\eta\in(0,1/4)\) and
an integer \(N_{\mathrm{loc}}<\infty\) such that every law
\(Q\) satisfying
\begin{equation}\label{eq:outline-local-cluster-hypotheses}
 \operatorname{supp}Q
 \subseteq[-\eta,\eta]\cup[1-\eta,1+\eta],
 \qquad
 \bigl|Q([1-\eta,1+\eta])-\pe\bigr|<\eta
\end{equation}
obeys
\begin{equation}\label{eq:outline-local-cluster-conclusion}
 R_n(Q)\le\ce
 \qquad\text{for every integer }n\ge N_{\mathrm{loc}}.
\end{equation}

To relate the bound \eqref{eq:outline-local-cluster-conclusion}
for \(Q_j\) to the maximizing identity
\eqref{eq:outline-maximizing-pair} for \(P_j\), we verify that
\(R_n(Q_j)=R_n(P_j)\).

For a distribution \(Q\) that is not necessarily standardized,
the quantities in \eqref{eq:outline-definitions} are defined
using its own mean and variance.
Specifically, if \(Q\) has mean \(m_Q\), variance \(v_Q>0\),
and finite third absolute central moment, the normalization
in \eqref{eq:raw-normalization} is
\begin{equation}\label{eq:outline-affine-normalization}
\begin{aligned}
 \beta(Q)
 &=\frac{\int|x-m_Q|^3\,Q(dx)}{v_Q^{3/2}},\\
 \Delta_n(Q)
 &=\sup_{t\in\mathbb R}
 \left|Q^{*n}((-\infty,t])
 -\Phi\left(\frac{t-nm_Q}{\sqrt{nv_Q}}\right)\right|,\\
 R_n(Q)&=\frac{\sqrt n\,\Delta_n(Q)}{\beta(Q)}.
\end{aligned}
\end{equation}

Since \(Q_j\) is the distribution of \(\pe+\se X\) for
\(X\sim P_j\), and \(P_j\) has mean zero and variance one,
we have \(m_{Q_j}=\pe\) and \(v_{Q_j}=\se^2\).
Consequently,
\begin{equation}\label{eq:outline-affine-moment-identities}
 \beta(Q_j)
 =\frac{\int|\se x|^3\,P_j(dx)}{\se^3}
 =\beta_j.
\end{equation}
For every integer \(n\ge1\), the sum of \(n\) independent
variables with distribution \(Q_j\) has the same distribution
as \(n\pe+\se\sum_{r=1}^n X_r\), where
\(X_1,\ldots,X_n\) are independent with distribution \(P_j\).
Thus the substitution \(t=n\pe+\se u\) gives
\[
 Q_j^{*n}((-\infty,t])=P_j^{*n}((-\infty,u]).
\]
Since \(\se>0\), this substitution is a bijection of
\(\mathbb R\).
Using it in \eqref{eq:outline-affine-normalization},
together with \eqref{eq:outline-affine-moment-identities},
we obtain
\begin{equation}\label{eq:outline-affine-invariance}
\begin{aligned}
 \Delta_n(Q_j)
 &=\sup_{u\in\mathbb R}
 \left|P_j^{*n}((-\infty,u])-\Phi(u/\sqrt n)\right|
 =\Delta_n(P_j),\\
 R_n(Q_j)
 &=\frac{\sqrt n\,\Delta_n(P_j)}{\beta_j}
 =R_n(P_j).
\end{aligned}
\end{equation}

Fix the constants \(\eta,N_{\mathrm{loc}}\) in
\eqref{eq:outline-local-cluster-hypotheses}--\eqref{eq:outline-local-cluster-conclusion}.
By \eqref{eq:outline-transformed-confinement}, the whole
support of \(Q_j\) lies in the prescribed two intervals for
all sufficiently large \(j\). The boundary of \([1-\eta,1+\eta]\) consists of the two
endpoints \(1-\eta\) and \(1+\eta\).
Since \(Q_{\mathrm E}\) is concentrated on \(\{0,1\}\)
and neither endpoint belongs to \(\{0,1\}\), we have
\[
 Q_{\mathrm E}\bigl(\partial[1-\eta,1+\eta]\bigr)=0.
\]
Thus \([1-\eta,1+\eta]\) is a continuity set for
\(Q_{\mathrm E}\). The weak convergence in
\eqref{eq:outline-transformed-weak-limit} then implies
\begin{equation}\label{eq:outline-upper-cluster-mass-limit}
 Q_j([1-\eta,1+\eta])
 \longrightarrow Q_{\mathrm E}([1-\eta,1+\eta])=\pe.
\end{equation}
Combining \eqref{eq:outline-upper-cluster-mass-limit} with
\eqref{eq:outline-transformed-confinement} and the convergence
\(n_j\to\infty\) in \eqref{eq:outline-selected-sizes}, we obtain,
for all sufficiently large \(j\),
\begin{equation}\label{eq:outline-cluster-conditions}
 \begin{gathered}
 \operatorname{supp}Q_j
 \subseteq[-\eta,\eta]\cup[1-\eta,1+\eta],\\
 \bigl|Q_j([1-\eta,1+\eta])-\pe\bigr|<\eta,
 \qquad n_j\ge N_{\mathrm{loc}}.
 \end{gathered}
\end{equation}
The conditions in \eqref{eq:outline-cluster-conditions}
are precisely the hypotheses needed for
\eqref{eq:outline-local-cluster-conclusion}.
Applying that conclusion and using
\eqref{eq:outline-affine-invariance} and the maximizing identity
\eqref{eq:outline-maximizing-pair} gives
\begin{equation}\label{eq:outline-final-contradiction}
 C_{n_j}=R_{n_j}(P_j)=R_{n_j}(Q_j)\le\ce,
\end{equation}
contrary to \(C_{n_j}>\ce\) in
\eqref{eq:outline-selected-sizes}.
The contradiction between \eqref{eq:outline-final-contradiction}
and \eqref{eq:outline-selected-sizes} shows that
\(C_n>\ce\) for at most finitely many positive
integers \(n\).  Taking an integer \(N\) larger than every
such index, or \(N=1\) if there are none, proves
\begin{equation}\label{eq:outline-eventual-bound}
 C_n\le\ce\qquad\text{for every integer }n\ge N.
\end{equation}
By the definition of \(C_n\) in \eqref{eq:outline-definitions},
\eqref{eq:outline-eventual-bound} establishes the existence of
the universal threshold in
Theorem~\ref{thm:main}.  

The argument above does not specify its
numerical value; the explicit threshold in
\eqref{eq:explicit-threshold} is obtained separately in
Appendix~\ref{app:effective}.

Section~\ref{sec:smoothing} establishes the uniform smoothing
expansion and its binomial consequences.
Section~\ref{sec:clusters} proves the local two-cluster
comparison.  In Section~\ref{sec:extremizers}, we derive the
variational identities and extract a sequence of maximizing
laws with uniformly bounded support.
Section~\ref{sec:completion} proves the limiting-support bound,
the separation property, and the existence assertion in
Theorem~\ref{thm:main}.
Appendix~\ref{app:effective} gives quantitative versions of the
required estimates and proves the explicit threshold
\eqref{eq:explicit-threshold}.

\section{Notation and related works}\label{subsec:notation}

We collect the notation used in the proof sketch and the classical
results needed below. The main inputs are Esseen's smoothing method
and moment inequality~\cite{Esseen1945,Esseen1956},
Shevtsova's explicit Berry--Esseen
estimates~\cite{Shevtsova2012,Shevtsova2013,Shevtsova2017},
and Schulz's finite-sample bound for two-point
laws~\cite{Schulz2016}.
We first recall the notation and convergence facts used below.
We then derive a common third-moment bound for standardized
counterexamples and establish \(C_n\to\ce\).
Finally, we describe the other classical results and methods
used in the subsequent proofs.

\subsection{Notation and convergence}
\label{subsec:notation-convergence}

We write \(F(t-)=\lim_{s\uparrow t}F(s)\) for the left limit
of a distribution function \(F\).
We use \(x_+=\max(x,0)\), and write
\(\|f\|_\infty=\sup_x|f(x)|\).

We continue to write \(\phi=\Phi'\), and put \(\phi_0=\phi(0)\).
Let \(\Pthree\) be the set of probability laws \(P\) on \(\R\)
satisfying
\[
 \int x\,P(dx)=0,\qquad \int x^2\,P(dx)=1,\qquad
 \beta(P)=\int |x|^3\,P(dx)<\infty.
\]
For a law \(Q\) with mean \(m\), variance \(v>0\), and finite
third absolute central moment \(\rho\), recall the normalization
\begin{equation}\label{eq:raw-normalization}
\begin{split}
 \Delta_n(Q)&=\sup_{t\in\R}
 \left|Q^{*n}(({-\infty},t])-
 \Phi\left(\frac{t-nm}{\sqrt{nv}}\right)\right|,\\
 \beta(Q)&=\rho/v^{3/2},\qquad
 R_n(Q)=\frac{\sqrt n\,\Delta_n(Q)}{\beta(Q)}.
\end{split}
\end{equation}
These quantities are invariant under nonzero affine changes
of the summands. For \(P\in\Pthree\) and an integer \(k\ge1\), let
\[
 F_{k,P}(t):=P^{*k}((-\infty,t])
 =\PP\{X_1+\cdots+X_k\le t\},\qquad t\in\R,
\]
where \(X_1,\ldots,X_k\) are independent with common law \(P\).
Thus \(F_{k,P}\) is the distribution function of the
unnormalized sum.
The extremal constant \(C_n\) can be written as
\begin{equation}\label{eq:cn}
 C_n=\sup_{P\in\Pthree}R_n(P)
 =\sup_{P\in\Pthree,\,t\in\R}
 \frac{\sqrt n\,[F_{n,P}(t)-\Phi(t/\sqrt n)]}{\beta(P)}.
\end{equation}

We may restrict attention to positive discrepancies because
any negative discrepancy can be matched arbitrarily closely
in magnitude by a positive discrepancy after reflecting the
summand distribution.
Indeed, let \(Q\) be the law of \(-X\) when \(X\) has law \(P\).
Then
\[
 F_{n,Q}((-t)-)=1-F_{n,P}(t),
\]
and the continuity and symmetry of \(\Phi\) imply
\[
 \lim_{u\uparrow -t}
 \bigl[F_{n,Q}(u)-\Phi(u/\sqrt n)\bigr]
 =\Phi(t/\sqrt n)-F_{n,P}(t).
\]
Thus, whenever the discrepancy for \(P\) at \(t\) is negative,
the positive discrepancies for \(Q\) at thresholds approaching
\(-t\) from below approach its absolute value.
Since \(Q\in\Pthree\) and \(\beta(Q)=\beta(P)\), allowing
negative discrepancies cannot increase the supremum.
This proves the second equality in \eqref{eq:cn}.

We use the third Wasserstein distance \(W_3\) introduced in
Section~\ref{subsec:outline-smoothing}.
For laws \(P_j,P\) with finite third absolute moments, convergence
in \(W_3\) is equivalent to weak convergence together with
\(\int |x|^3\,P_j(dx)\to\int |x|^3\,P(dx)\);
see \cite[Theorem~5.11, pp.~185--187]{Santambrogio2015}.
In particular, weak convergence of laws supported in one common
compact interval implies convergence in \(W_3\).
The compactness argument in Lemma~\ref{lem:attainment} also uses
tightness, the Portmanteau theorem, and convergence of expectations
under uniform integrability; see
\cite[Theorems~5.1, 2.1, and 3.5]{Billingsley1999}.

Constants denoted by \(C\) or \(c\) may change from one estimate
to the next; their dependence on parameters is indicated where
needed.

For a fixed sequence of laws and sample sizes satisfying the
stated convergence assumptions, an \(o(1)\) term tends to zero
along that sequence. Uniformity in any additional variables
is stated explicitly. No common convergence rate over all
such sequences is asserted.

When an estimate is stated uniformly over a class of laws,
its constants may depend on the stated fixed parameters,
but not on the particular law in that class.

\subsection{Classical bounds and the extremal constant}
\label{subsec:classical-bounds}

We recall Shevtsova's explicit Berry--Esseen
bounds~\cite{Shevtsova2012,Shevtsova2013} and Esseen's
moment inequality and fixed-law asymptotic
formula~\cite{Esseen1956}.
From these classical results, we derive a uniform bound on
the third absolute moments of standardized counterexamples
and prove that \(C_n\to\ce\).
In Section~\ref{sec:extremizers}, the moment bound will allow
us to prove that a maximizing distribution exists whenever
\(C_n>\ce\), while the convergence of \(C_n\) will allow us
to select suitable sample sizes for the contradiction argument.

We begin with two estimates of Shevtsova~\cite{Shevtsova2013,Shevtsova2012}.
For \(P\in\Pthree\) and \(n\ge1\), these estimates are
\begin{equation}\label{eq:structural}
 \Delta_n(P)\le\frac1{\sqrt n}
 \min\{0.4690\beta(P),\;0.3031(\beta(P)+0.646)\},
\end{equation}
and
\begin{equation}\label{eq:asymptotic}
 \Delta_n(P)\le \ce\frac{\beta(P)}{\sqrt n}
                  +3\frac{\beta(P)^2}{n}.
\end{equation}
The first follows from~\cite[pp.~124--125]{Shevtsova2013}: for
standardized iid summands, the quantities in that paper are
\(\ell_n=\beta(P)/\sqrt n\) and \(\tau_n=1/\sqrt n\).
The coefficient \(3\) in the second is a convenient enlargement of
\(2.5786\) in~\cite[Corollary~4.18, p.~303]{Shevtsova2012}.

For \(P\in\Pthree\), recall that
\(\kappa(P)=\int x^3\,P(dx)\), and that \(h(P)\) denotes its
maximal lattice span, with \(h(P)=0\) for a nonlattice law.
Esseen's moment inequality and his asymptotic formula for a
fixed summand law \(P\) are
\begin{equation}\label{eq:esseen-moment}
 |\kappa(P)|+3h(P)\le(3+\sqrt{10})\beta(P),
\end{equation}
and
\begin{equation}\label{eq:esseen-limit}
 \lim_{n\to\infty}\sqrt n\,\Delta_n(P)
 =\frac{\phi_0}{6}\bigl(|\kappa(P)|+3h(P)\bigr).
\end{equation}
The inequality \eqref{eq:esseen-moment} and its equality
characterization are due to Esseen
\cite[Lemma~1, pp.~162--163]{Esseen1956}.
The fixed-law formula \eqref{eq:esseen-limit} is also due to
Esseen~\cite[p.~162, following equations~(6)--(7)]{Esseen1956}; see
\cite[Section~1.1, equations~(1.5)--(1.9), p.~491]{MattnerShevtsova2019}
for both results.
Equality in \eqref{eq:esseen-moment} holds exactly for the
standardizations of Bernoulli laws with parameters \(\pe\)
and \(1-\pe\).
We recall the parameters
\begin{equation}\label{eq:esseen-parameters}
 \pe=\frac{4-\sqrt{10}}2,\qquad \qe=1-\pe,
 \qquad \se=\sqrt{\pe\qe}.
\end{equation}
The positive-skew representative used in the proof sketch is
\begin{equation}\label{eq:esseen-law}
 \Pe=\qe\delta_{-\aE}+\pe\delta_{\bE},
 \qquad \aE=\pe/\se,\quad \bE=\qe/\se.
\end{equation}
Since \(\Pe\) has only two support points, its maximal lattice
span is the distance between them:
\[
 h_{\mathrm E}=\aE+\bE=1/\se.
\]
We also have
\begin{equation}\label{eq:esseen-identities}
\begin{aligned}
 \be=\beta(\Pe)&=\frac{\pe^2+\qe^2}{\se},\qquad
 \ke=\kappa(\Pe)=\frac{\qe-\pe}{\se},\\
 M_{\mathrm E}=\int x|x|\,\Pe(dx)
 &=\qe-\pe=\sqrt{10}-3.
\end{aligned}
\end{equation}

The equality case of \eqref{eq:esseen-moment} will identify
the possible limiting extremizers in Lemma~\ref{lem:limit-extremizer}.
The fixed-law formula \eqref{eq:esseen-limit}, applied to
\(\Pe\), supplies the lower asymptotic bound for \(C_n\).

The proof of Lemma~\ref{lem:effective-lattice-stability}
in Appendix~\ref{app:effective} expresses the deficit in
Esseen's moment inequality using a product coupling of
the normalized first-moment measures of the positive
and negative parts.
For related representations of zero-mean laws through
couplings of these measures, see
Pinelis~\cite[equation~(3.12), Proposition~3.23,
and the following remark]{Pinelis2009}.

We now derive the moment cutoff and convergence of \(C_n\)
used in the proof sketch.
By \eqref{eq:structural}, every \(P\in\Pthree\) satisfying
\(R_n(P)>\ce\) must satisfy
\begin{equation}\label{eq:moment-cutoff}
 \beta(P)<B_*:=\frac{0.3031\cdot0.646}{\ce-0.3031}
 =1.83624299\ldots.
\end{equation}

The same cutoff was also used by Schulz~\cite[p.~20]{Schulz2016}.

For \(\beta(P)\ge B_*\), the same bound gives \(R_n(P)\le\ce\).
For \(\beta(P)<B_*\), equation~\eqref{eq:asymptotic} gives
\(R_n(P)\le\ce+3B_*/\sqrt n\).
Thus \(C_n\le\ce+3B_*/\sqrt n\).
On the other hand, \eqref{eq:esseen-limit} applied to \(\Pe\)
gives \(R_n(\Pe)\to\ce\), and \(C_n\ge R_n(\Pe)\).
Together with \eqref{eq:structural}, these bounds imply
\begin{equation}\label{eq:cn-limit}
 C_n\longrightarrow\ce,\qquad C_n\le0.4690.
\end{equation}

\subsection{Further background results}
\label{subsec:classical-tools}

The smoothing argument in Section~\ref{sec:smoothing} follows
Esseen's classical approach to short Edgeworth expansions.
Esseen's smoothing inequality
\cite[Chapter~II, Theorem~2a, p.~32]{Esseen1945}
is stated for a comparison function \(G\) of bounded variation
whose derivative exists everywhere and satisfies
\(|G'(x)|\le M\) for all \(x\in\R\).
His proof uses the derivative bound only through the
Lipschitz estimate
\[
 |G(x)-G(y)|\le M|x-y|,
\]
which is the regularity assumption used in
Lemma~\ref{lem:signed-smoothing}.
This Lipschitz form also follows directly from
\cite[Chapter~V, Section~1, Theorem~1, p.~104]{Petrov1975}.
We include a proof of Lemma~\ref{lem:signed-smoothing}
for completeness.

In Lemma~\ref{lem:jitter}, we use uniform smoothing to obtain
a first-order Edgeworth expansion when the summand distribution
varies with the sample size.
For standardized laws \(P_j\in\Pthree\) satisfying
\(P_j\to P\) in \(W_3\), along sample sizes \(n_j\to\infty\),
we add to each unnormalized row sum an independent uniform
variable on \([-h(P)/2,h(P)/2]\).
The smoothing width is determined by the limiting law \(P\);
when \(h(P)=0\), the added variable is identically zero.
Corollary~\ref{cor:jitter-envelopes} then uses this expansion
to obtain upper and lower bounds for the distribution functions
of the original sums.

For a fixed lattice summand law, this smoothing removes
the first-order lattice correction.
A precise one-dimensional statement is recorded in
\cite[Theorem~A, p.~29]{BabuSingh1989},
where it is attributed to Feller.
The underlying nonlattice and lattice expansions are
classical results of Esseen
\cite[Chapter~IV, Theorem~2, p.~49, and Theorem~3, p.~56]{Esseen1945}.

Booth, Hall and Wood
\cite[Theorem~2.2, pp.~124--125]{BoothHallWood1994}
established threshold-uniform Edgeworth expansions for discrete
resampling laws with a growing number of atoms, almost surely
in iid atom locations drawn from a parent density, under
assumptions on that density and the weights.
In contrast, Lemma~\ref{lem:jitter} applies to every sequence of
standardized laws \(P_j\to P\) in \(W_3\), including lattice limits,
after adding independent uniform noise on
\([-h(P)/2,h(P)/2]\) to each unnormalized row sum.
Here \(h(P)\) is the maximal lattice span of \(P\), with zero
noise when \(P\) is nonlattice.
G\"otze and Hipp~\cite[Theorem~3.6, p.~71]{GoetzeHipp1978}
obtained expansions for fixed smooth test functions under moment
and Lindeberg-type conditions.
Their theorem does not directly yield the threshold-uniform
\(o(n_j^{-1/2})\) remainder in Lemma~\ref{lem:jitter},
since the relevant test functions here are indicators or
\(n_j\)-dependent piecewise-linear ramps.

To complete the comparison with a Bernoulli law in
Lemma~\ref{lem:small-cluster-variance}, we use Schulz's
finite-sample theorem
\cite[Theorem~1, p.~1]{Schulz2016}.
Together with affine invariance, this theorem gives
\begin{equation}\label{eq:schulz-two-point}
 R_n(Q)<\ce,\qquad n\ge1,
\end{equation}
for every nondegenerate two-point law \(Q\).

The same comparison also uses the binomial local limit
and the expansions at the two sides of each binomial jump
established in Lemma~\ref{lem:binomial-estimates}.
There we derive these formulas from the smoothing expansion,
uniformly for \(p\) in compact subsets of \((0,1)\).
The local limit theorem is classical; see, for example,
\cite[Chapter~VII, Section~1, Theorem~1, p.~187]{Petrov1975}.
An explicit bound on the error in approximating binomial
point probabilities by the normal density is given in
\cite[Lemma~5, p.~393]{ZolotukhinNagaevChebotarev2018}.

To prove the lower bound of order \(n^{-1/2}\) for interval
probabilities in Lemma~\ref{lem:two-cluster-local-mass},
we approximate the noise sum conditional on the cluster
counts by a normal law.
This conditional noise sum can be represented as a sum of
independent centered variables whose distributions may differ.
We therefore use the Berry--Esseen inequality for independent,
not necessarily identically distributed summands.
If \(Y_1,\ldots,Y_k\) are independent and centered, have finite
third absolute moments, and satisfy
\(s^2=\sum_{i=1}^k\EE Y_i^2>0\), then
\begin{equation}\label{eq:independent-berry-esseen}
 \sup_{x\in\R}
 \left|\PP\left\{\frac{\sum_{i=1}^kY_i}{s}\le x\right\}
       -\Phi(x)\right|
 \le C_{\mathrm{ind}}\,
       \frac{\sum_{i=1}^k\EE|Y_i|^3}{s^3}.
\end{equation}
We take \(C_{\mathrm{ind}}=1\), a convenient enlargement of
the constant \(0.5583\) in
\cite[pp.~124--125]{Shevtsova2013}.

The proof of the explicit threshold in Appendix~\ref{app:effective}
also requires control of the approximation error at thresholds
far from the mean.
For this purpose, we use the nonuniform Berry--Esseen bound
of Nagaev--Bikelis type~\cite{Nagaev1965,Bikelis1966}
for iid summands, with constant \(17.36\) proved in
\cite[Corollary~3.4.5, p.~90, and Table~3.3, p.~55]{Shevtsova2017}.
In our normalization, this bound is stated in
\eqref{eq:effective-nonuniform}, whose right-hand side decreases
as the standardized threshold moves away from zero.

The proof of Lemma~\ref{lem:effective-global-jitter} uses an integer
representation of the lattice span to control the characteristic
function away from its resonance frequencies.
For related quantitative bounds, see Buhler, Gamst, Graham and
Hales~\cite{BuhlerGamstGrahamHales2018}.

For the variational argument in Section~\ref{sec:extremizers},
we perturb a law \(P\) by replacing it with
\((1-\varepsilon)P+\varepsilon\delta_y\).
This is the classical point-mass contamination used to define
influence functions; see
\cite[Section~2.1, p.~383]{Hampel1974}.
Related variational contact conditions appear in
Mattner~\cite[Theorem~1 and Example~1]{Mattner1992}
and, for tail probabilities of iid sums, in
Meester~\cite[Section~2]{Meester2008}.
Lemma~\ref{lem:influence} computes the derivative of the
normalized signed discrepancy at a fixed threshold and derives
the support contact identities used in the proof.

% ===== smoothing =====
\section{Uniform smoothing along convergent arrays}\label{sec:smoothing}

We now prove the smoothing expansion for convergent arrays
described in Section~\ref{subsec:classical-tools}.
We begin with Esseen's smoothing inequality and then establish
the first-order expansion in Lemma~\ref{lem:jitter}, for
standardized laws \(P_j\to P\) in \(W_3\) and sample sizes
\(n_j\to\infty\).
From this expansion, we obtain upper and lower bounds for
the distribution functions of the original sums, as well as
the binomial estimates needed in Section~\ref{sec:clusters}.

For related explicit finite-sample bounds on the distance
to the first-order Edgeworth expansion under fourth-moment
assumptions, see Derumigny, Girard and
Guyonvarch~\cite{DerumignyGirardGuyonvarch2024}.

Throughout this section, we use the Fourier convention
\(\widehat\nu(t)=\int e^{itx}\,d\nu(x)\).

\begin{lemma}[Esseen's smoothing inequality, Lipschitz form]
\label{lem:signed-smoothing}
Let \(F\) be the distribution function of a probability measure with finite
first absolute moment.  Let \(\nu\) be a finite signed measure of total mass
one such that \(\int |x|\,d|\nu|(x)<\infty\), and suppose that
\(G(x)=\nu(( -\infty,x])\) is globally Lipschitz with constant \(M\).
There are numerical constants \(C_1,C_2\) such that, for every \(L>0\),
\begin{equation}\label{eq:signed-smoothing}
 \|F-G\|_\infty
 \le C_1\int_{-L}^{L}
       \frac{|\widehat F(t)-\widehat\nu(t)|}{|t|}\,dt
       +\frac{C_2M}{L}.
\end{equation}
Here \(\widehat F\) denotes the characteristic function of the probability
measure associated with \(F\).
\end{lemma}

\begin{proof}
We first bound the difference \(F-G\) in terms of a smoothed
version of this difference. We then estimate the smoothed
difference by Fourier inversion. 

Choose the even probability density
\[
 K(x)=c\left(\frac{\sin(x/4)}{x/4}\right)^4,
 \qquad
 c^{-1}=\int_{\R}\left(\frac{\sin(x/4)}{x/4}\right)^4dx.
\]
The ratio at \(x=0\) is defined by continuity.
Since \(K(x)=O(|x|^{-4})\) as \(|x|\to\infty\), this density
has a finite first absolute moment.
Its characteristic function is supported in \([-1,1]\). Indeed,
\[
 \frac{\sin(x/4)}{x/4}=2\int_{-1/4}^{1/4}e^{iux}\,du,
\]
so its fourth power is the Fourier transform of the convolution of four
uniform densities on \([-1/4,1/4]\); Fourier inversion gives the asserted
support.  Also \(|\widehat K|\le1\).  For \(Z\) with density \(K\), fix a
numerical \(a>0\) such that
\[
 \varepsilon=\PP(Z>a)=\PP(Z<-a)<\tfrac14,
 \qquad C_a=\EE(a-Z)_+<\infty.
\]
By symmetry, \(\EE(a+Z)_+=C_a\).

Set \(K_L(y)=LK(Ly)\), \(d=F-G\), and
\[
 D_+=\sup_xd(x),\qquad D_-=\sup_x[-d(x)],\qquad
 D=\max(D_+,D_-),\qquad R=\|d*K_L\|_\infty.
\]
The quantities \(D_+,D_-\) are nonnegative because \(d\) tends to zero at
both infinities.  For \(y\ge x\), monotonicity of \(F\) and the Lipschitz
bound for \(G\) give \(d(y)\ge d(x)-M(y-x)\).  Splitting the following
expectation according to \(Z\le a\) consequently gives
\[
 \begin{split}
 (d*K_L)(x+a/L)
 &=\EE\,d\bigl(x+(a-Z)/L\bigr)\\
 &\ge(1-\varepsilon)d(x)-\varepsilon D_- -MC_a/L.
 \end{split}
\]
Taking a sequence of points at which \(d\) tends to \(D_+\) yields
\[
 (1-\varepsilon)D_+-\varepsilon D_-\le R+MC_a/L.
\]
For the other side, \(y\le x\) implies
\(d(y)\le d(x)+M(x-y)\). Write
\[
 (d*K_L)(x-a/L)
 =\EE\,d\bigl(x-(a+Z)/L\bigr).
\]
On the event \(Z\ge-a\), the argument of \(d\) is at most \(x\),
so the preceding inequality applies. On the complementary
event, we use \(d\le D_+\). By symmetry,
\(\EE(a+Z)_+=C_a\), and hence
\[
 (d*K_L)(x-a/L)
 \le(1-\varepsilon)d(x)+\varepsilon D_++MC_a/L.
\]
Since the left-hand side is at least \(-R\), letting
\(d(x)\) approach \(-D_-\) gives
\[
 (1-\varepsilon)D_--\varepsilon D_+\le R+MC_a/L.
\]

If \(D=D_+\), then
\[
 (1-\varepsilon)D_+-\varepsilon D_-
 \ge(1-2\varepsilon)D.
\]
The same conclusion follows from the other inequality when
\(D=D_-\). Thus, in either case,
\[
 (1-2\varepsilon)D\le R+MC_a/L.
\]

We now estimate \(R\) by Fourier inversion.
Let \(\mu=dF-\nu\). Since both \(dF\) and \(\nu\) have total
mass one, \(\mu(\R)=0\). Moreover,
\[
 \int_{\R}|u|\,d|\mu|(u)
 \le \int_{\R}|u|\,dF(u)
    +\int_{\R}|u|\,d|\nu|(u)<\infty.
\]
Using \(\mu(\R)=0\), we may write
\[
 d(x)=\int_{\R}
 \bigl(\mathbf 1_{\{u\le x\}}-\mathbf 1_{\{0\le x\}}\bigr)
 \,\mu(du).
\]
For each \(u\in\R\),
\[
 \int_{\R}
 \left|\mathbf 1_{\{u\le x\}}-\mathbf 1_{\{0\le x\}}\right|
 \,dx=|u|.
\]
Tonelli's theorem therefore gives
\[
 \int_{\R}|d(x)|\,dx
 \le\int_{\R}|u|\,d|\mu|(u)<\infty.
\]

The same integrability justifies Fubini's theorem in the
following calculation. For \(t\ne0\),
\[
 \int_{\R}e^{itx}
 \bigl(\mathbf 1_{\{u\le x\}}-\mathbf 1_{\{0\le x\}}\bigr)
 \,dx
 =\frac{1-e^{itu}}{it}.
\]
Consequently,
\[
 \widehat d(t)
 =\int_{\R}\frac{1-e^{itu}}{it}\,\mu(du)
 =\frac{\widehat\nu(t)-\widehat F(t)}{it}.
\]

Since \(d\) and \(K_L\) are integrable,
\[
 \widehat{d*K_L}(t)=\widehat d(t)\widehat K(t/L).
\]
The function \(\widehat d\) is bounded, and
\(\widehat K(t/L)\) vanishes outside \([-L,L]\).
Thus \(\widehat{d*K_L}\) is integrable.
Also, \(d*K_L\) is continuous because \(d\in L^1(\R)\)
and \(K_L\) is bounded and continuous.
Fourier inversion therefore holds at every \(x\in\R\):
\[
 (d*K_L)(x)
 =\frac1{2\pi}\int_{-L}^{L}
 e^{-itx}\widehat d(t)\widehat K(t/L)\,dt.
\]
Using \(|\widehat K|\le1\), we obtain
\[
 R\le\frac1{2\pi}\int_{-L}^{L}
 \frac{|\widehat F(t)-\widehat\nu(t)|}{|t|}\,dt.
\]

Combining this estimate with
\((1-2\varepsilon)D\le R+MC_a/L\) proves the result with
\[
 C_1=\frac1{2\pi(1-2\varepsilon)},
 \qquad
 C_2=\frac{C_a}{1-2\varepsilon}.
\]
Since \(K\) and \(a\) were fixed independently of \(F\), \(\nu\),
and \(L\), these are numerical constants.
\end{proof}

We now establish the first-order expansion for the smoothed
sums described in Section~\ref{subsec:classical-tools},
when the summand laws satisfy \(P_j\to P\) in \(W_3\)
and the sample sizes satisfy \(n_j\to\infty\).

\begin{lemma}\label{lem:jitter}
Let \(P_j\) be distributions of mean zero and variance one, and suppose
that \(P_j\to P\) in \(W_3\).  Let \(n_j\to\infty\) be positive integers,
and let \(X_{j,1},\ldots,X_{j,n_j}\) be independent with common law
\(P_j\).  Write
\[
 S_j=\frac1{\sqrt{n_j}}\sum_{r=1}^{n_j}X_{j,r},\qquad
 \mu_j=\int x^3\,dP_j(x).
\]
Let \(h=h(P)\), with \(h(P)=0\) when \(P\) is nonlattice,
and let \(U\), independent of the row, be uniform on \([-h/2,h/2]\);
when \(h=0\), set \(U=0\).  Then
\begin{equation}\label{eq:jitter-expansion}
 \sup_{x\in\R}\left|
 \PP\{S_j+U/\sqrt{n_j}\le x\}
 -\Phi(x)-\frac{\mu_j}{6\sqrt{n_j}}(1-x^2)\phi(x)
 \right|=o(n_j^{-1/2}).
\end{equation}
\end{lemma}

\begin{proof}
We first describe the argument and the estimates to be proved.
Write \(n=n_j\), and let
\[
 \mathcal J_j(x)=\PP\{S_j+U/\sqrt n\le x\},
 \qquad
 G_j(x)=\Phi(x)+\frac{\mu_j}{6\sqrt n}(1-x^2)\phi(x).
\]
Our goal is
\[
 \|\mathcal J_j-G_j\|_\infty=o(n^{-1/2}).
\]
Define
\[
 f_j(u)=\int_{\R}e^{iux}\,P_j(dx),
 \qquad
 f(u)=\int_{\R}e^{iux}\,P(dx),
 \qquad
 H(u)=\EE e^{iuU}=\sinc(hu/2),
\]
where
\[
 \sinc(v)=
 \begin{cases}
  \sin(v)/v,&v\ne0,\\
  1,&v=0.
 \end{cases}
\]
In particular, \(H\equiv1\) when \(h=0\).
As verified below, the Fourier transforms of the law associated
with \(\mathcal J_j\) and the signed measure
\(dG_j(x)=G_j'(x)\,dx\) are, respectively,
\[
 \begin{aligned}
 a_j(t)&=H(t/\sqrt n)f_j(t/\sqrt n)^n,\\
 b_j(t)&=e^{-t^2/2}
          \left(1+\frac{\mu_j(it)^3}{6\sqrt n}\right).
 \end{aligned}
\]

Following Esseen's Fourier approach to first-order expansions
\cite[Chapter~IV]{Esseen1945}, we estimate
\[
 \mathcal I_j(T)
 :=\int_{-T\sqrt n}^{T\sqrt n}
       \frac{|a_j(t)-b_j(t)|}{|t|}\,dt.
\]
We will verify the hypotheses of
Lemma~\ref{lem:signed-smoothing}, including the uniform bound
\[
 M_g:=\sup_j\|G_j'\|_\infty<\infty.
\]
Applying that lemma with \(L=T\sqrt n\) will then give
\[
 \sqrt n\,\|\mathcal J_j-G_j\|_\infty
 \le C_1\sqrt n\,\mathcal I_j(T)+\frac{C_2M_g}{T},
\]
where \(C_1,C_2\) and \(M_g\) are independent of \(j\) and \(T\).
Thus it suffices to prove
\(\mathcal I_j(T)=o(n^{-1/2})\) for arbitrarily large fixed
cutoffs \(T\).

We first estimate the integral near zero.
Using \(W_3\) convergence, we prove the uniform expansion
\[
 \log f_j(u)
 =-\frac{u^2}{2}+\frac{\mu_j(iu)^3}{6}
   +\widetilde r_j(u),
 \qquad
 \sup_j|\widetilde r_j(u)|=o(|u|^3)
 \quad(u\to0).
\]
Here the logarithms are taken on the principal branch in
one common neighborhood of zero.
Substituting \(u=t/\sqrt n\) and multiplying by \(n\) gives
\[
 n\log f_j(t/\sqrt n)
 =-\frac{t^2}{2}+\frac{\mu_j(it)^3}{6\sqrt n}
   +n\widetilde r_j(t/\sqrt n).
\]
This expansion, together with the bounds
\[
 |f_j(u)|\le e^{-u^2/4}
 \quad (|u|\le\delta),
 \qquad
 |H(u)-1|\le\frac{h^2u^2}{24},
\]
will yield
\[
 \int_{|t|\le\delta\sqrt n}
       \frac{|a_j(t)-b_j(t)|}{|t|}\,dt
 =o(n^{-1/2})
\]
for a sufficiently small fixed \(\delta>0\).

For the remaining frequencies, we use the variable
\(u=t/\sqrt n\).
We first identify the set of limiting resonances:
\[
 \mathcal R:=\{u:|f(u)|=1\}
 =
 \begin{cases}
  \{0\},&h=0,\\
  \{2\pi k/h:k\in\mathbb Z\},&h>0.
 \end{cases}
\]
For every nonempty compact set \(K\) disjoint from
\(\mathcal R\), we show that there exist
\(\rho_K\in(0,1)\) and an index \(j_K\) such that
\[
 \sup_{u\in K}|f_j(u)|\le\rho_K,
 \qquad j\ge j_K.
\]
Consequently,
\[
 \sup_{u\in K}|f_j(u)|^n\le\rho_K^n,
 \qquad j\ge j_K,
\]
so these frequencies contribute an exponentially small error.

In the lattice case, the nonzero resonances
\(u_k=2\pi k/h\), \(k\ne0\), require a separate estimate.
After decreasing \(\delta\), if necessary, we have
\(\mathcal R\cap[-\delta,\delta]=\{0\}\).
Fix \(T>\delta\) with \(|f(\pm T)|<1\).
There are only finitely many nonzero resonances in
\((-T,T)\).
Around each of them, we choose a sufficiently small fixed
closed interval
\[
 I_k\subset\{u:\delta<|u|<T\},
\]
with these intervals pairwise disjoint.
We show that, for all sufficiently large \(j\),
\(q_j:=|f_j|^2\) has a unique maximizer \(u_{j,k}\) on
\(I_k\), and that
\[
 u_{j,k}\longrightarrow u_k.
\]
The curvature of \(q_j\) near this maximum gives a Gaussian
bound for \(|f_j|^n\).
The smoothing factor \(H\) vanishes at each nonzero limiting resonance:
\[
 H(u_k)=\sinc(\pi k)=0.
\]
Combining these facts, we obtain
\[
 \int_{I_k}\frac{|H(u)|\,|f_j(u)|^n}{|u|}\,du
 \le C\left(
       \frac1n+\frac{|u_{j,k}-u_k|}{\sqrt n}
       \right)
 =o(n^{-1/2}),
\]
where \(C\) is independent of \(j\).
Thus convergence of the peak locations is sufficient;
no rate of convergence is needed.

Combining the estimate near zero, the estimates near the
nonzero resonances, the exponential bound on the remaining
compact set, and the Gaussian tail estimate for \(b_j\),
we prove
\[
 \mathcal I_j(T)=o(n^{-1/2})
\]
for every fixed \(T>\delta\) with \(|f(\pm T)|<1\).
The smoothing inequality then gives
\[
 \limsup_{j\to\infty}
 \sqrt{n_j}\,\|\mathcal J_j-G_j\|_\infty
 \le\frac{C_2M_g}{T}.
\]
Since \(C_2M_g\) is independent of \(T\), and such cutoffs
can be arbitrarily large, letting \(T\to\infty\) proves
\eqref{eq:jitter-expansion}.

\textbf{Step~1. Preparing the comparison measure.}
We first record the consequences of \(W_3\) convergence
and compute the Fourier transforms of the two measures
to be compared.
For the signed density \(g_j=G_j'\), we will verify
\[
 \int_{\R}g_j(x)\,dx=1,\qquad
 \sup_j\int_{\R}(1+|x|)|g_j(x)|\,dx<\infty,
 \qquad
 M_g=\sup_j\|g_j\|_\infty<\infty.
\]
These facts provide the mass, moment, and Lipschitz
conditions on \(G_j\) required by
Lemma~\ref{lem:signed-smoothing}.

Write \(n=n_j\) within the proof, and put
\begin{equation}\label{eq:jitter-proof-characteristic-definitions}
 \begin{aligned}
 f_j(u)&=\int_{\R}e^{iux}\,P_j(dx),&
 f(u)&=\int_{\R}e^{iux}\,P(dx),\\
 q_j(u)&=|f_j(u)|^2,& q(u)&=|f(u)|^2.
 \end{aligned}
\end{equation}
All limits indexed by \(j\) below are along the sequence in the
statement, so \(n=n_j\to\infty\).

We first explain why the first three moments pass to the limit.
By the definition of \(W_3\), for each \(j\) we can couple
\(\xi_j\sim P_j\) and \(\eta_j\sim P\) so that
\[
 \|\xi_j-\eta_j\|_3
 \le W_3(P_j,P)+j^{-1}\longrightarrow0,
 \qquad
 \|\eta_j\|_3=\left(\int_{\R}|x|^3\,P(dx)\right)^{1/3}<\infty.
\]
Here \(\|Y\|_p=(\EE|Y|^p)^{1/p}\), and only the coupling for
each fixed \(j\) is needed. The triangle inequality for the
\(L^3\) norm gives
\(\|\xi_j\|_3\le\|\eta_j\|_3+\|\xi_j-\eta_j\|_3\), so
\(\sup_j\|\xi_j\|_3<\infty\).
For \(k=1,2,3\), the factorization of \(x^k-y^k\) gives
\[
 |x^k-y^k|\le k|x-y|(|x|+|y|)^{k-1}.
\]
Applying the triangle inequality when \(k=1\), and H\"older's
inequality with exponents \(k\) and \(k/(k-1)\) when \(k=2,3\),
we obtain
\[
 \begin{aligned}
 \left|\int_{\R}x^k\,P_j(dx)-\int_{\R}x^k\,P(dx)\right|
 &\le k\|\xi_j-\eta_j\|_k
       \bigl(\|\xi_j\|_k+\|\eta_j\|_k\bigr)^{k-1}\\
 &\le k\|\xi_j-\eta_j\|_3
       \bigl(\|\xi_j\|_3+\|\eta_j\|_3\bigr)^{k-1}
 \longrightarrow0.
 \end{aligned}
\]
The second inequality uses \(\|Y\|_k\le\|Y\|_3\) on a probability
space for \(k\le3\). Since every \(P_j\) has mean zero and second
moment one, this proves
\begin{equation}\label{eq:jitter-proof-moments}
 \begin{gathered}
 \int_{\R}x\,P_j(dx)=\int_{\R}x\,P(dx)=0,\\
 \int_{\R}x^2\,P_j(dx)=\int_{\R}x^2\,P(dx)=1,\\
 \mu_j\longrightarrow\mu:=\int_{\R}x^3\,P(dx),\qquad
 M_\mu:=\sup_j|\mu_j|<\infty.
 \end{gathered}
\end{equation}
In particular, \(P\) has variance one and is nondegenerate.

Recall that
\begin{equation}\label{eq:jitter-proof-jitter-factor}
 H(u):=\EE e^{iuU}=\sinc(hu/2),
\end{equation}
where \(\sinc(v)=\sin(v)/v\) for \(v\ne0\) and
\(\sinc(0)=1\). Indeed, for \(h>0\) and \(u\ne0\),
\[
 \EE e^{iuU}
 =\frac1h\int_{-h/2}^{h/2}e^{iuv}\,dv
 =\frac{e^{iuh/2}-e^{-iuh/2}}{ihu}
 =\frac{\sin(hu/2)}{hu/2}.
\]
At \(u=0\) both sides equal one. When \(h=0\), the definition
\(U=0\) gives \(H\equiv1\), also in agreement with
\eqref{eq:jitter-proof-jitter-factor}. As a characteristic function,
\(H\) satisfies \(|H(u)|\le\EE|e^{iuU}|=1\).
By the independence of \(U,X_{j,1},\ldots,X_{j,n}\), the
characteristic function of \(S_j+U/\sqrt n\) is
\begin{equation}\label{eq:jitter-proof-sum-transform}
 \begin{aligned}
 a_j(t)&:=\EE e^{it(S_j+U/\sqrt n)}\\
 &=\EE e^{itU/\sqrt n}
       \prod_{r=1}^{n}\EE e^{itX_{j,r}/\sqrt n}\\
 &=H(t/\sqrt n)f_j(t/\sqrt n)^n.
 \end{aligned}
\end{equation}
The last equality uses the common law \(P_j\) of the summands
and definitions \eqref{eq:jitter-proof-characteristic-definitions}
and \eqref{eq:jitter-proof-jitter-factor}.

For the application of Lemma~\ref{lem:signed-smoothing}, recall the comparison function
\begin{equation}\label{eq:jitter-proof-comparison-function}
 G_j(x):=\Phi(x)+\frac{\mu_j}{6\sqrt n}(1-x^2)\phi(x).
\end{equation}
Our goal in \eqref{eq:jitter-expansion} is to show that \(G_j\)
approximates the distribution function of \(S_j+U/\sqrt n\)
with uniform error \(o(n^{-1/2})\).
Using \(\Phi'=\phi\) and \(\phi'(x)=-x\phi(x)\), we calculate
\[
 \frac{d}{dx}\bigl[(1-x^2)\phi(x)\bigr]
 =-2x\phi(x)-x(1-x^2)\phi(x)
 =(x^3-3x)\phi(x).
\]
Thus differentiating \eqref{eq:jitter-proof-comparison-function}
gives the signed density
\begin{equation}\label{eq:jitter-proof-comparison-density}
 g_j(x):=G_j'(x)
 =\phi(x)\left[1+\frac{\mu_j}{6\sqrt n}(x^3-3x)\right].
\end{equation}
The identities
\(\phi''(x)=(x^2-1)\phi(x)\) and
\(\phi'''(x)=-(x^3-3x)\phi(x)\) follow by differentiating
\(\phi'(x)=-x\phi(x)\) twice.
For clarity, the Gaussian Fourier transform can also be obtained
from this derivative identity. Let
\(\Gamma(t)=\int_{\R}e^{itx}\phi(x)\,dx\).
Since \(\int_{\R}|x|\phi(x)\,dx<\infty\), differentiation
under the integral and one integration by parts give
\[
 \Gamma'(t)
 =i\int_{\R}xe^{itx}\phi(x)\,dx
 =-i\int_{\R}e^{itx}\phi'(x)\,dx
 =-t\Gamma(t).
\]
Together with \(\Gamma(0)=1\), this differential equation gives
\(\Gamma(t)=e^{-t^2/2}\).
Three integrations by parts therefore give
\[
 \begin{aligned}
 \int_{\R}e^{itx}(x^3-3x)\phi(x)\,dx
 &=-\int_{\R}e^{itx}\phi'''(x)\,dx\\
 &=it\int_{\R}e^{itx}\phi''(x)\,dx\\
 &=-(it)^2\int_{\R}e^{itx}\phi'(x)\,dx\\
 &=(it)^3\int_{\R}e^{itx}\phi(x)\,dx
 =(it)^3e^{-t^2/2}.
 \end{aligned}
\]
All boundary terms vanish because \(\phi,\phi',\phi''\) tend
to zero at both infinities, and their products with \(e^{itx}\)
have the same absolute values. The last equality uses the
Fourier transform of the standard normal density.
Combining this calculation with
\eqref{eq:jitter-proof-comparison-density}, we find
\begin{equation}\label{eq:jitter-proof-comparison-transform}
 b_j(t):=\int_{\R}e^{itx}g_j(x)\,dx
 =e^{-t^2/2}\left[1+\frac{\mu_j(it)^3}{6\sqrt n}\right].
\end{equation}

The correction term \((x^3-3x)\phi(x)\) in
\eqref{eq:jitter-proof-comparison-density} is odd and integrable.
Its integral over \(\R\) is consequently zero, and hence
\begin{equation}\label{eq:jitter-proof-comparison-mass}
 \int_{\R}g_j(x)\,dx
 =\int_{\R}\phi(x)\,dx
  +\frac{\mu_j}{6\sqrt n}\int_{\R}(x^3-3x)\phi(x)\,dx=1.
\end{equation}
The bound on \(M_\mu\) in \eqref{eq:jitter-proof-moments},
together with \(n=n_j\ge1\), gives
\[
 \frac{|\mu_j|}{6\sqrt n}\le\frac{|\mu_j|}{6}
 \le\frac{M_\mu}{6}.
\]
Using this coefficient bound and the triangle inequality in
\eqref{eq:jitter-proof-comparison-density}, we obtain
\[
 \begin{aligned}
 |g_j(x)|
 &\le\phi(x)\left[1+\frac{|\mu_j|}{6\sqrt n}|x^3-3x|\right]\\
 &\le\phi(x)\left[1+\frac{M_\mu}{6}(|x|^3+3|x|)\right].
 \end{aligned}
\]
For \(|x|\le1\), \(|x|\le1\), whereas for \(|x|>1\),
\(|x|\le|x|^3\). Thus \(|x|\le1+|x|^3\) for all \(x\), and
\[
 |x|^3+3|x|\le3+4|x|^3\le4(1+|x|^3).
\]
Substituting this inequality in the preceding bound yields
\begin{equation}\label{eq:jitter-proof-density-bound}
 |g_j(x)|\le C(1+|x|^3)\phi(x),\qquad
 C:=1+\frac23M_\mu<\infty.
\end{equation}
This choice of \(C\) is independent of \(j\).
The function \((1+|x|^3)\phi(x)\) is continuous and tends to
zero as \(|x|\to\infty\), so it is bounded. Taking suprema in
\eqref{eq:jitter-proof-density-bound} gives
\begin{equation}\label{eq:jitter-proof-density-sup}
 M_g:=\sup_j\|g_j\|_\infty
 \le C\sup_{x\in\R}(1+|x|^3)\phi(x)<\infty.
\end{equation}
To verify the required absolute moments as well, recall that
\(\phi'(x)=-x\phi(x)\) gives
\[
 \begin{aligned}
 \int_{\R}|x|\phi(x)\,dx
 &=-2\int_0^\infty\phi'(x)\,dx
   =2\phi(0)=\sqrt{2/\pi},\\
 \int_{\R}|x|^3\phi(x)\,dx
 &=-2\int_0^\infty x^2\phi'(x)\,dx
   =4\int_0^\infty x\phi(x)\,dx=2\sqrt{2/\pi},\\
 \int_{\R}x^4\phi(x)\,dx
 &=-\int_{\R}x^3\phi'(x)\,dx
   =3\int_{\R}x^2\phi(x)\,dx=3.
 \end{aligned}
\]
The integrations by parts have zero boundary terms because
each polynomial times \(\phi\) tends to zero at infinity;
at zero the factors \(x^2\) and \(x^3\) vanish.
Multiplying \eqref{eq:jitter-proof-density-bound} by \(1+|x|\)
and using these Gaussian integrals gives
\begin{equation}\label{eq:jitter-proof-density-moments}
 \begin{aligned}
 \int_{\R}(1+|x|)|g_j(x)|\,dx
 &\le C\int_{\R}(1+|x|+|x|^3+|x|^4)\phi(x)\,dx\\
 &=C\left(4+3\sqrt{2/\pi}\right)<\infty.
 \end{aligned}
\end{equation}
The two terms on the left of
\eqref{eq:jitter-proof-density-moments} are the total variation
and the first absolute moment of the signed measure
\(g_j(x)\,dx\), respectively. This signed measure therefore is
finite and has finite first absolute moment, in addition to
having total mass one by \eqref{eq:jitter-proof-comparison-mass}.

\textbf{Step~2. A uniform expansion near zero.}
We next use \(W_3\) convergence to control the third-order
Taylor remainder uniformly in \(j\).
Our aim is to obtain the logarithmic expansion
\[
 \log f_j(u)
 =-\frac{u^2}{2}+\frac{\mu_j(iu)^3}{6}
   +\widetilde r_j(u),
 \qquad
 |\widetilde r_j(u)|
 \le |u|^3\widetilde\omega(|u|),
 \qquad
 \widetilde\omega(r)\longrightarrow0
 \quad(r\downarrow0).
\]
Here the remainder bound uses a single function
\(\widetilde\omega\) for all \(j\), and the logarithms
will be defined on the principal branch in one common
neighborhood of zero.
We will also choose a fixed \(\delta>0\) such that
\[
 |f_j(u)|\le e^{-u^2/4},
 \qquad |u|\le\delta,\quad\text{for every }j.
\]
The expansion controls the approximation error, while
the Gaussian bound will allow us to integrate the
resulting estimates.
These two conclusions are established in
\eqref{eq:jitter-proof-log-expansion} and
\eqref{eq:jitter-proof-gaussian-decay}. 

By the characterization of \(W_3\) convergence in
Section~\ref{subsec:notation-convergence},
\[
 \int |x|^3\,P_j(dx)\longrightarrow\int |x|^3\,P(dx).
\]
In particular,
\begin{equation}\label{eq:jitter-proof-third-moment-bound}
 B_3:=\sup_j\int |x|^3\,P_j(dx)<\infty.
\end{equation}
We also need uniform integrability of the random variables
\(|X_{j,1}|^3\).
For \(A>0\), the function \(x\mapsto\min\{|x|^3,A^3\}\) is
bounded and continuous. Weak convergence and convergence of the
third absolute moments therefore give
\[
 \int (|x|^3-A^3)_+\,P_j(dx)
 \longrightarrow
 \int (|x|^3-A^3)_+\,P(dx).
\]
The expression on the right tends to zero as \(A\to\infty\)
by dominated convergence, since its integrand tends to zero
and is bounded by the \(P\)-integrable function \(|x|^3\).
Given any positive tolerance, first choose \(A\) so that this
expression is small. The same bound, with an arbitrarily small
enlargement, then holds for all sufficiently large \(j\).
Increasing \(A\) further handles the finitely many remaining laws,
since each has a finite third absolute moment. Thus
\[
 \lim_{A\to\infty}\sup_j
 \int (|x|^3-A^3)_+\,P_j(dx)=0.
\]
For \(|x|>2A\), we have
\(|x|^3\le2(|x|^3-A^3)_+\). Consequently,
\begin{equation}\label{eq:jitter-proof-third-moment-tails}
 \lim_{A\to\infty}\sup_j
 \int_{|x|>A}|x|^3\,P_j(dx)=0.
\end{equation}

We use the classical Taylor remainder bound for the complex
exponential; see, for example, \cite{Shevtsova2014Taylor}.
Writing
\[
 T_m(y):=e^{iy}-\sum_{k=0}^{m-1}\frac{(iy)^k}{k!},
\]
the bound is
\begin{equation}\label{eq:jitter-proof-taylor-remainder}
 |T_m(y)|\le\frac{|y|^m}{m!},
 \qquad y\in\R,\quad m=1,2,\ldots.
\end{equation}
Taking \(m=4\) gives
\[
 |T_4(y)|\le\frac{|y|^4}{24}.
\]
Also, since \(T_4(y)=T_3(y)-(iy)^3/6\), the case \(m=3\)
and the triangle inequality give
\[
 |T_4(y)|
 \le |T_3(y)|+\frac{|y|^3}{6}
 \le\frac{|y|^3}{3}.
\]
Combining these two bounds, we obtain
\begin{equation}\label{eq:jitter-proof-taylor-truncated}
 \left|e^{iy}-1-iy-\frac{(iy)^2}{2}-\frac{(iy)^3}{6}\right|
 \le
 \begin{cases}
  |y|^4/24,& |y|\le1,\\
  |y|^3/3,& |y|>1,
 \end{cases}
 \le\frac13|y|^3\min\{|y|,1\}.
\end{equation}
Set \(y=uX_{j,1}\) and take expectations. The mean and variance
identities in \eqref{eq:jitter-proof-moments} give
\begin{equation}\label{eq:jitter-proof-characteristic-expansion}
 \begin{split}
 f_j(u)&=1-u^2/2+\mu_j(iu)^3/6+r_j(u),\\
 r_j(u)&:=\EE T_4(uX_{j,1}).
 \end{split}
\end{equation}
By \eqref{eq:jitter-proof-taylor-truncated},
\[
 |r_j(u)|\le\frac{|u|^3}{3}
 \EE\bigl[|X_{j,1}|^3\min\{|u|\,|X_{j,1}|,1\}\bigr].
\]
We may therefore use the common remainder function
\begin{equation}\label{eq:jitter-proof-remainder-modulus}
 \omega(r):=\frac13\sup_j
 \EE\bigl[|X_{j,1}|^3\min\{r|X_{j,1}|,1\}\bigr],
 \qquad r\ge0,
\end{equation}
so that
\begin{equation}\label{eq:jitter-proof-characteristic-remainder}
 |r_j(u)|\le |u|^3\omega(|u|).
\end{equation}
Definition~\eqref{eq:jitter-proof-remainder-modulus} and
\eqref{eq:jitter-proof-third-moment-bound} show that
\(\omega\) is nondecreasing and finite, with
\(0\le\omega(r)\le B_3/3\). To check that \(\omega(r)\to0\)
as \(r\downarrow0\), split the expectation at \(|X_{j,1}|=A\).
On \(|X_{j,1}|\le A\), use
\(\min\{r|X_{j,1}|,1\}\le rA\); on the complement, use the
bound by one. This gives
\[
 3\omega(r)
 \le rAB_3+
 \sup_j\EE\bigl[|X_{j,1}|^3
                  \mathbf 1_{\{|X_{j,1}|>A\}}\bigr].
\]
First choose \(A\) large using
\eqref{eq:jitter-proof-third-moment-tails}, and then let \(r\)
tend to zero. This proves the claimed limit.

We now pass from \(f_j\) to its logarithm on one fixed interval
around zero. By \eqref{eq:jitter-proof-taylor-remainder} with
\(m=2\), together with \(\EE X_{j,1}=0\) and
\(\EE X_{j,1}^2=1\),
\[
 |f_j(u)-1|
 =\left|\EE\bigl(e^{iuX_{j,1}}-1-iuX_{j,1}\bigr)\right|
 \le u^2/2.
\]
Thus, for \(|u|\le1\), every \(f_j(u)\) lies in the disk
\(\{w:|w-1|\le1/2\}\), where the principal logarithm is
defined. For \(|w|\le1/2\), its power series gives
\[
 |\log(1+w)-w|
 \le\sum_{k=2}^{\infty}|w|^k
 =\frac{|w|^2}{1-|w|}
 \le2|w|^2.
\]
Taking \(w=f_j(u)-1\) therefore shows that
\[
 |\log f_j(u)-(f_j(u)-1)|\le u^4/2,
 \qquad |u|\le1.
\]
Combining this estimate with
\eqref{eq:jitter-proof-characteristic-expansion} and
\eqref{eq:jitter-proof-characteristic-remainder}, we obtain
\begin{equation}\label{eq:jitter-proof-log-expansion}
 \begin{split}
 \log f_j(u)
 &=-u^2/2+\mu_j(iu)^3/6+\widetilde r_j(u),\\
 |\widetilde r_j(u)|
 &\le |u|^3\widetilde\omega(|u|),\qquad |u|\le1,
 \end{split}
\end{equation}
where we may take
\[
 \widetilde\omega(r):=\omega(r)+r/2.
\]
This function is nondecreasing and tends to zero as
\(r\downarrow0\); the term \(r/2\) accounts for the uniform
\(O(u^4)\) error from taking the logarithm.

Fix \(\delta\in(0,1]\), independent of \(j\), sufficiently
small that
\begin{equation}\label{eq:jitter-proof-delta}
 \delta\widetilde\omega(\delta)\le\tfrac18.
\end{equation}
Since \(\mu_j(iu)^3/6\) is purely imaginary for real \(u\),
\eqref{eq:jitter-proof-log-expansion} implies, for
\(|u|\le\delta\),
\[
 \Re\log f_j(u)
 \le-u^2/2+|u|^3\widetilde\omega(|u|)
 \le-u^2/2+u^2\delta\widetilde\omega(\delta)
 \le-3u^2/8.
\]
Since \(|f_j(u)|=\exp(\Re\log f_j(u))\), this implies
\begin{equation}\label{eq:jitter-proof-gaussian-decay}
 |f_j(u)|\le e^{-u^2/4},\qquad |u|\le\delta.
\end{equation}

\textbf{Step~3. Estimating the Fourier integral near zero.}
We now apply the preceding expansion with
\(u=t/\sqrt n\).
The Fourier error separates into the error for the
unsmoothed sum and the contribution of the added
uniform variable:
\[
 \begin{split}
 a_j(t)-b_j(t)
 &=
 \bigl[f_j(t/\sqrt n)^n-b_j(t)\bigr]\\
 &\quad+
 \bigl[H(t/\sqrt n)-1\bigr]f_j(t/\sqrt n)^n.
 \end{split}
\]
We will show that these two terms contribute
\(o(n^{-1/2})\) and \(O(n^{-1})\), respectively, after
integration against \(dt/|t|\) over
\(|t|\le\delta\sqrt n\).
Together they give
\[
 \int_{|t|\le\delta\sqrt n}
       \frac{|a_j(t)-b_j(t)|}{|t|}\,dt
 =o(n^{-1/2}),
\]
which is recorded in
\eqref{eq:jitter-proof-near-zero}.

For \(|t|\le\delta\sqrt n\), set
\[
 z=\frac{\mu_j(it)^3}{6\sqrt n},
 \qquad r=n\widetilde r_j(t/\sqrt n).
\]
The number \(z\) is purely imaginary, and
\eqref{eq:jitter-proof-log-expansion}--\eqref{eq:jitter-proof-delta}
give
\begin{equation}\label{eq:jitter-proof-scaled-log-remainder}
 |r|\le\frac{|t|^3}{\sqrt n}
          \widetilde\omega(|t|/\sqrt n)
 \le t^2\delta\widetilde\omega(\delta)
 \le t^2/8.
\end{equation}
Exponentiating \eqref{eq:jitter-proof-log-expansion}, and using
\eqref{eq:jitter-proof-comparison-transform}, yields
\[
 \begin{split}
 f_j(t/\sqrt n)^n-b_j(t)
 &=e^{-t^2/2}\bigl(e^{z+r}-1-z\bigr)\\
 &=e^{-t^2/2}
   \bigl[e^z(e^r-1)+(e^z-1-z)\bigr].
 \end{split}
\]
For any complex \(r\), the identity
\(e^r-1=r\int_0^1e^{sr}\,ds\) gives
\(|e^r-1|\le|r|e^{|r|}\).
Because \(z\) is purely imaginary, \(|e^z|=1\), and
\eqref{eq:jitter-proof-taylor-remainder} with \(m=2\) gives
\(|e^z-1-z|\le|z|^2/2\).
Using these bounds, \eqref{eq:jitter-proof-scaled-log-remainder},
and \(|\mu_j|\le M_\mu\), we obtain
\begin{equation}\label{eq:jitter-proof-low-frequency-integrand}
 \begin{split}
 |f_j(t/\sqrt n)^n-b_j(t)|
 &\le e^{-t^2/2}\bigl(|r|e^{|r|}+|z|^2/2\bigr)\\
 &\le e^{-3t^2/8}\frac{|t|^3}{\sqrt n}
          \widetilde\omega(|t|/\sqrt n)
       +\frac{M_\mu^2}{72}e^{-t^2/2}\frac{|t|^6}{n}.
 \end{split}
\end{equation}
To integrate \eqref{eq:jitter-proof-low-frequency-integrand},
multiply by \(\sqrt n/|t|\)
on \(0<|t|\le\delta\sqrt n\), and set the resulting
integrand equal to zero outside this interval and at \(t=0\).
For each fixed \(t\ne0\), its upper bound tends to zero,
since \(\widetilde\omega(|t|/\sqrt n)\to0\) and
\(n^{-1/2}\to0\). By monotonicity of \(\widetilde\omega\)
and \(n\ge1\), it is bounded on the whole real line by
\[
 e^{-3t^2/8}t^2\widetilde\omega(\delta)
 +\frac{M_\mu^2}{72}e^{-t^2/2}|t|^5.
\]
This function is integrable and independent of \(j\).
Dominated convergence therefore gives
\begin{equation}\label{eq:jitter-proof-low-frequency-unsmoothed}
 \int_{|t|\le\delta\sqrt n}
 \frac{|f_j(t/\sqrt n)^n-b_j(t)|}{|t|}\,dt
 =o(n^{-1/2}).
\end{equation}

It remains to include the factor \(H\) in
\eqref{eq:jitter-proof-sum-transform}.
For \(h>0\), symmetry of the uniform density gives \(\EE U=0\),
and direct integration gives
\[
 \EE U^2=\frac1h\int_{-h/2}^{h/2}v^2\,dv=\frac{h^2}{12}.
\]
These identities also hold when \(h=0\), since then \(U=0\).
Using \eqref{eq:jitter-proof-taylor-remainder} with \(m=2\),
we obtain
\begin{equation}\label{eq:jitter-proof-jitter-remainder}
 |H(u)-1|
 =\left|\EE\bigl(e^{iuU}-1-iuU\bigr)\right|
 \le\frac{u^2}{2}\EE U^2
 =\frac{h^2u^2}{24}.
\end{equation}
Combining \eqref{eq:jitter-proof-jitter-remainder} with
\eqref{eq:jitter-proof-gaussian-decay}, we obtain
\[
 \begin{split}
 &\int_{|t|\le\delta\sqrt n}
 \frac{|H(t/\sqrt n)-1|\,|f_j(t/\sqrt n)|^n}{|t|}\,dt\\
 &\qquad\le\frac{h^2}{24n}
      \int_{\R}|t|e^{-t^2/4}\,dt
 =\frac{h^2}{6n}=O_h(n^{-1}).
 \end{split}
\]
Here
\(\int_{\R}|t|e^{-t^2/4}\,dt
=2\int_0^\infty t e^{-t^2/4}\,dt=4\),
by the substitution \(s=t^2/4\).
This contribution is zero when \(h=0\).
Finally,
\[
 a_j(t)-b_j(t)
 =f_j(t/\sqrt n)^n-b_j(t)
  +\bigl(H(t/\sqrt n)-1\bigr)f_j(t/\sqrt n)^n.
\]
The triangle inequality and
\eqref{eq:jitter-proof-low-frequency-unsmoothed} consequently give
\begin{equation}\label{eq:jitter-proof-near-zero}
 \int_{|t|\le\delta\sqrt n}
 \frac{|a_j(t)-b_j(t)|}{|t|}\,dt=o(n^{-1/2}).
\end{equation}

\textbf{Step~4. Compact convergence and limiting resonances.}
To estimate the remaining frequencies, we first prove
that \(q_j=|f_j|^2\) and its second derivative converge
uniformly on every fixed compact interval.
Writing \(q=|f|^2\), the required conclusions are
\[
 \sup_{|u|\le T}|q_j(u)-q(u)|\longrightarrow0,
 \qquad
 \sup_{|u|\le T}|q_j''(u)-q''(u)|\longrightarrow0
 \quad\text{for every fixed }T>0.
\]
We then identify the points where \(q=1\):
\[
 \mathcal R=\{u:q(u)=1\}
 =
 \begin{cases}
  \{0\},&h=0,\\
  \{u_k=2\pi k/h:k\in\mathbb Z\},&h>0.
 \end{cases}
\]
The characterization of lattice laws by the existence
of a nonzero point with \(|f(u)|=1\) is classical;
see \cite[Chapter~I, Theorem~5, p.~26]{Esseen1945}.
We give the argument below and, in the lattice case,
also establish
\[
 q''(u_k)=-2.
\]
The convergence of \(q_j\) will control frequencies
away from \(\mathcal R\), while the convergence of
\(q_j''\) and this curvature identity will control
the peaks near the nonzero resonances.

For each \(j\), choose a coupling of \(X_j\sim P_j\) and \(X\sim P\)
such that
\[
 \EE|X_j-X|^3\le W_3(P_j,P)^3+\frac1j.
\]
The coupling may depend on \(j\), but the marginal law of \(X\)
is always \(P\). The definition of the Wasserstein distance
permits this choice, and the right side tends to zero.
H\"older's inequality gives
\[
 \EE|X_j-X|^2
 \le \bigl(\EE|X_j-X|^3\bigr)^{2/3},
\]
so the coupled difference also tends to zero in \(L^2\).
Take an independent copy \((X'_j,X')\) of this pair and put
\begin{equation}\label{eq:jitter-proof-coupled-differences}
 D_j=X_j-X'_j,\qquad D=X-X'.
\end{equation}
The definitions in \eqref{eq:jitter-proof-coupled-differences},
independence of each variable and its copy, and the first two
moments recorded in \eqref{eq:jitter-proof-moments} give
\[
 \EE D_j^2=\EE D^2=2.
\]
Moreover,
\begin{equation}\label{eq:jitter-proof-difference-l2}
 \begin{aligned}
 \EE|D_j-D|^2
 &\le 2\EE|X_j-X|^2+2\EE|X'_j-X'|^2\\
 &=4\EE|X_j-X|^2\longrightarrow0.
 \end{aligned}
\end{equation}
The Cauchy--Schwarz inequality now gives
\begin{equation}\label{eq:jitter-proof-squares-l1}
 \begin{aligned}
 \EE|D_j^2-D^2|
 &=\EE\bigl[|D_j-D|\,|D_j+D|\bigr]\\
 &\le\bigl(\EE|D_j-D|^2\bigr)^{1/2}
       \bigl(\EE|D_j+D|^2\bigr)^{1/2}\\
 &\le\sqrt8\,\bigl(\EE|D_j-D|^2\bigr)^{1/2}
 \longrightarrow0,
 \end{aligned}
\end{equation}
where the last bound uses
\(\EE|D_j+D|^2\le2\EE D_j^2+2\EE D^2=8\).

Independence of the copied pair gives
\[
 q_j(u)=f_j(u)\overline{f_j(u)}
       =\EE e^{iuD_j}=\EE\cos(uD_j),
 \qquad
 q(u)=\EE\cos(uD).
\]
The cosine representations hold because \(D_j\) and \(D\) have
symmetric distributions: interchanging each variable with its
independent copy changes the difference's sign without changing
its law. The fundamental theorem of calculus gives
\( |\cos r-\cos s|=|\int_s^r\sin v\,dv|\le|r-s| \).
For every fixed \(T<\infty\), this bound implies
\begin{equation}\label{eq:jitter-proof-q-compact}
 \sup_{|u|\le T}|q_j(u)-q(u)|
 \le T\EE|D_j-D|
 \le T\bigl(\EE|D_j-D|^2\bigr)^{1/2}
 \longrightarrow0.
\end{equation}
The finite second moments justify differentiating the cosine
representations twice under the expectation. They give
\begin{equation}\label{eq:jitter-proof-q-second-derivative}
 q_j''(u)=-\EE[D_j^2\cos(uD_j)],
 \qquad
 q''(u)=-\EE[D^2\cos(uD)].
\end{equation}
For example, the first and second derivatives of the integrand
are bounded in absolute value by \(|D_j|\) and \(D_j^2\),
respectively, and both bounds are integrable. Dominated
convergence also shows that these second derivatives are continuous.
Subtracting the two expressions in
\eqref{eq:jitter-proof-q-second-derivative} gives
\[
 \begin{split}
 q_j''(u)-q''(u)
 &=-\EE[(D_j^2-D^2)\cos(uD_j)]\\
 &\quad-\EE[D^2(\cos(uD_j)-\cos(uD))].
 \end{split}
\]
Using \(|\cos(uD_j)|\le1\) and
\(|\cos(uD_j)-\cos(uD)|\le\min\{2,|u|\,|D_j-D|\}\)
therefore yields
\begin{equation}\label{eq:jitter-proof-q-second-compact-bound}
 \begin{aligned}
 \sup_{|u|\le T}|q_j''(u)-q''(u)|
 &\le \EE|D_j^2-D^2|\\
 &\quad+\EE\!\left[
       D^2\min\{2,T|D_j-D|\}\right].
 \end{aligned}
\end{equation}
To control the second term without a fourth-moment assumption,
fix \(A>0\) and split according to \(|D|>A\). Then
\begin{equation}\label{eq:jitter-proof-q-second-truncation}
 \begin{aligned}
 \EE\!\left[D^2\min\{2,T|D_j-D|\}\right]
 &\le 2\EE[D^2\ind_{\{|D|>A\}}]
       +A^2T\EE|D_j-D|.
 \end{aligned}
\end{equation}
For fixed \(A,T\), the last term tends to zero by
\eqref{eq:jitter-proof-difference-l2}. The first term can be
made arbitrarily small by increasing \(A\), because \(\EE D^2=2\);
the marginal law of \(D\) does not depend on \(j\).
Combining \eqref{eq:jitter-proof-q-second-compact-bound},
\eqref{eq:jitter-proof-q-second-truncation}, and
\eqref{eq:jitter-proof-squares-l1} therefore proves
\begin{equation}\label{eq:jitter-proof-q-second-compact}
 \sup_{|u|\le T}|q_j''(u)-q''(u)|\longrightarrow0
 \qquad(T<\infty).
\end{equation}
Thus both compact convergence statements use only second moments
and the \(L^2\) approximation supplied by \(W_3\) convergence.

We now justify the characterization of the limiting resonance set stated above. Recall that
\begin{equation}\label{eq:jitter-proof-resonance-set}
 \mathcal R=\{u\in\R:q(u)=1\}.
\end{equation}
Since \(q(u)=\EE\cos(uD)\) and \(1-\cos(uD)\ge0\),
\[
 u\in\mathcal R
 \quad\Longleftrightarrow\quad
 \cos(uD)=1\ \hbox{almost surely}
 \quad\Longleftrightarrow\quad
 e^{iuD}=1\ \hbox{almost surely}.
\]
This characterization shows that \(0\in\mathcal R\), that
\(-u\in\mathcal R\) whenever \(u\in\mathcal R\), and that
\(u+v\in\mathcal R\) whenever \(u,v\in\mathcal R\).
For the last assertion, intersect the two probability-one events
and multiply \(e^{iuD}\) and \(e^{ivD}\).
Continuity of \(q\) makes \(\mathcal R\) closed.
Also, \(q\) is even and differentiable, so \(q'(0)=0\).
Together with \(q(0)=1\) and
\(q''(0)=-\EE D^2=-2\), Taylor's theorem gives
\[
 q(u)=1-u^2+o(u^2)\qquad(u\to0).
\]
In particular, \(q(u)<1\) for every sufficiently small nonzero \(u\);
zero is isolated in the set \(\mathcal R\) defined in
\eqref{eq:jitter-proof-resonance-set}.

These properties imply that either \(\mathcal R=\{0\}\) or
\(\mathcal R=\alpha\mathbb Z\) for some \(\alpha>0\).
Here is the classification in this setting. If the subgroup contains
a nonzero point, it contains a positive point. Isolation of zero
makes
\[
 \alpha:=\inf\bigl(\mathcal R\cap(0,\infty)\bigr)>0.
\]
A sequence of positive points tending to this infimum and
closedness show that \(\alpha\in\mathcal R\).
For any \(u\in\mathcal R\) with \(u\ge0\), subtracting
\(\lfloor u/\alpha\rfloor\alpha\) leaves a point of \(\mathcal R\)
in \([0,\alpha)\). By the definition of \(\alpha\), this
remainder must be zero. Negative points are treated by changing
sign, proving \(\mathcal R=\alpha\mathbb Z\).

If \(\mathcal R=\alpha\mathbb Z\), then
\(\alpha\in\mathcal R\), so
\[
 |\EE e^{i\alpha X}|^2=q(\alpha)=1.
\]
Consequently,
\[
 \EE\left|e^{i\alpha X}-\EE e^{i\alpha X}\right|^2
 =1-|\EE e^{i\alpha X}|^2=0,
\]
so \(e^{i\alpha X}\) is almost surely one fixed complex number
of modulus one. Choose \(\theta\in\R\) such that this number is
\(e^{i\theta}\). Then
\[
 X\in\frac{\theta}{\alpha}+\frac{2\pi}{\alpha}\mathbb Z
 \quad\hbox{almost surely}.
\]
This lattice is closed, so it contains \(\supp P\).
Conversely, if \(\supp P\) is contained in a translate of
\(d\mathbb Z\), where \(d>0\), then \(D\in d\mathbb Z\)
almost surely and \(2\pi/d\in\mathcal R\).
In the case \(\mathcal R=\alpha\mathbb Z\), this means
\(2\pi/d=m\alpha\) for a positive integer \(m\), and hence
\(d\le2\pi/\alpha\).
Therefore \(h=2\pi/\alpha\) is the maximal lattice span.
If \(\mathcal R=\{0\}\), no positive lattice span is possible.
We have thus proved
\begin{equation}\label{eq:jitter-proof-resonances}
 \mathcal R=
 \begin{cases}
  \{0\},&h=0,\\
  \{u_k=2\pi k/h:k\in\mathbb Z\},&h>0.
 \end{cases}
\end{equation}
In the lattice case, at every resonance \(u_k\),
\(\cos(u_kD)=1\) almost surely, so
\eqref{eq:jitter-proof-q-second-derivative} also gives
\begin{equation}\label{eq:jitter-proof-resonance-curvature}
 q''(u_k)=-\EE D^2=-2.
\end{equation}

\textbf{Step~5. Estimating the Fourier integral away from zero.}
We now use the preceding convergence and resonance
properties to prove, for each fixed cutoff \(T\)
chosen below,
\[
 \int_{\delta\le|u|\le T}
       \frac{|H(u)|\,|f_j(u)|^n}{|u|}\,du
 =o(n^{-1/2}).
\]
In the lattice case, we choose small pairwise disjoint
closed intervals \(I_k\) around the nonzero resonances
\(u_k\in(-T,T)\), with
\[
 I_k\subset\{u:\delta<|u|<T\}.
\]
The curvature estimates will give \(q_j''\le-1\) on
these intervals for all sufficiently large \(j\).
We will then locate the unique maximizer \(u_{j,k}\)
of \(q_j\) on each interval \(I_k\) and prove
\[
 u_{j,k}\longrightarrow u_k,
 \qquad
 |f_j(u)|^n
 \le \exp\!\left(-\frac n4(u-u_{j,k})^2\right),
 \quad u\in I_k.
\]
We will use both the identity \(H(u_k)=0\) and the
Lipschitz estimate
\eqref{eq:jitter-proof-jitter-lipschitz}, proved below.
Together they give
\[
 |H(u)|
 =|H(u)-H(u_k)|
 \le\frac h4|u-u_k|.
\]
Combining these estimates yields
\[
 \int_{I_k}\frac{|H(u)|\,|f_j(u)|^n}{|u|}\,du
 \le
 \frac h{4\delta}
 \left(
   \frac4n+
   \frac{2\sqrt\pi}{\sqrt n}|u_{j,k}-u_k|
 \right)
 =o(n^{-1/2}).
\]
In particular, convergence of \(u_{j,k}\) to \(u_k\)
is sufficient for this estimate; no rate is needed.
On the remaining compact set, a bound
\(|f_j|\le\rho<1\), valid for all sufficiently large
\(j\), gives exponential decay.
In the nonlattice case, this latter argument covers
the whole region \(\delta\le|u|\le T\).
The resulting estimate is
\eqref{eq:jitter-proof-away-fourier-factor}.

To prove this estimate, fix \(T>\delta\) such that neither
endpoint of \([-T,T]\) is a resonance. In the lattice case, decrease \(\delta\), if
necessary, so that \([-\delta,\delta]\) contains no nonzero
resonance. This decrease preserves \(T>\delta\), all the
previous bounds on \(|u|\le\delta\), and the estimate
\eqref{eq:jitter-proof-near-zero}; for the last assertion,
the integral over the smaller interval is bounded by the
one already estimated.
There are only finitely many nonzero resonances in \((-T,T)\).
Around each of them, choose a closed interval
\[
 I_k=[u_k-r_k,u_k+r_k],\qquad 0<r_k<\tfrac12,
\]
so that the intervals are pairwise disjoint, lie in \((-T,T)\),
are disjoint from \([-\delta,\delta]\), and contain no other
resonance. By continuity of \(q''\) and
\eqref{eq:jitter-proof-resonance-curvature}, the radii can also
be chosen to satisfy
\begin{equation}\label{eq:jitter-proof-resonance-concavity}
 q''(u)\le-\tfrac32\qquad(u\in I_k).
\end{equation}
All intervals in this construction are fixed once \(P,\delta,T\)
are fixed. By \eqref{eq:jitter-proof-q-second-compact},
\(\sup_{|u|\le T}|q_j''(u)-q''(u)|\le1/2\) for all sufficiently
large \(j\). Together with
\eqref{eq:jitter-proof-resonance-concavity}, this gives
\begin{equation}\label{eq:jitter-proof-array-concavity}
 q_j''(u)\le-1\qquad(u\in I_k)
\end{equation}
on every one of them.

At either endpoint of \(I_k\), the value of \(q\) is strictly
less than \(q(u_k)=1\), because that endpoint is not a resonance.
The uniform convergence in \eqref{eq:jitter-proof-q-compact}
therefore implies that, for all sufficiently large \(j\), the
value \(q_j(u_k)\) exceeds both endpoint values of \(q_j\).
The continuous function \(q_j\) attains a maximum on \(I_k\);
the preceding strict inequalities place that maximum in the
interior. By \eqref{eq:jitter-proof-array-concavity}, \(q_j\)
is strictly concave there, so the maximizer is unique.
Denote it by \(u_{j,k}\). Since it is interior,
\(q_j'(u_{j,k})=0\).

We also have
\begin{equation}\label{eq:jitter-proof-peak-location}
 u_{j,k}\longrightarrow u_k.
\end{equation}
To verify this, fix \(\eta>0\). If
\(\{u\in I_k:|u-u_k|\ge\eta\}\) is nonempty, it is compact
and contains no resonance. Its maximum value of \(q\) is
therefore strictly below \(1\). Uniform convergence then
ensures that \(q_j(u_k)\) eventually exceeds every value
of \(q_j\) on this set. The maximizer must consequently
satisfy \(|u_{j,k}-u_k|<\eta\). If the set is empty,
this inequality holds automatically. This proves
\eqref{eq:jitter-proof-peak-location} for every fixed \(k\).

For \(u\ne u_{j,k}\) in \(I_k\), Taylor's theorem gives a point
\(v\) between \(u\) and \(u_{j,k}\) such that
\[
 q_j(u)=q_j(u_{j,k})+q_j'(u_{j,k})(u-u_{j,k})
            +\tfrac12q_j''(v)(u-u_{j,k})^2.
\]
Since \(v\in I_k\), \eqref{eq:jitter-proof-array-concavity}
and \(q_j'(u_{j,k})=0\) give
\[
 q_j(u)\le q_j(u_{j,k})-\tfrac12(u-u_{j,k})^2
          \le1-\tfrac12(u-u_{j,k})^2
 \qquad(u\in I_k).
\]
Here \(q_j\le1\) because it is the squared modulus of a
characteristic function. Our choice \(r_k<1/2\) also gives
\(|u-u_{j,k}|<1\), so the last upper bound is positive.
Using \(1-v\le e^{-v}\) for \(v\ge0\), we obtain
\begin{equation}\label{eq:jitter-proof-gaussian-peak}
 |f_j(u)|^n=q_j(u)^{n/2}
 \le \exp\!\left[-\frac n4(u-u_{j,k})^2\right]
 \qquad(u\in I_k).
\end{equation}

The jitter factor cancels these peaks at their limiting locations.
For \(k\ne0\), \eqref{eq:jitter-proof-jitter-factor} gives
\[
 H(u_k)=\sinc(hu_k/2)=\sinc(\pi k)=0.
\]
For \(h>0\),
\(\EE|U|=(2/h)\int_0^{h/2}v\,dv=h/4\), and this identity
also holds when \(h=0\).
Equation~\eqref{eq:jitter-proof-taylor-remainder} with \(m=1\)
gives \(|e^{ir}-e^{is}|=|e^{i(r-s)}-1|\le|r-s|\).
It follows that
\begin{equation}\label{eq:jitter-proof-jitter-lipschitz}
 |H(u)-H(v)|
 \le \EE|e^{iuU}-e^{ivU}|
 \le |u-v|\EE|U|=\frac h4|u-v|.
\end{equation}
Combining \eqref{eq:jitter-proof-jitter-lipschitz} with
\(H(u_k)=0\) yields
\[
 |H(u)|\le\frac h4|u-u_k|
 \le\frac h4\bigl(|u-u_{j,k}|+|u_{j,k}-u_k|\bigr).
\]
Since \(|u|\ge\delta\) on \(I_k\),
\eqref{eq:jitter-proof-gaussian-peak} implies
\begin{equation}\label{eq:jitter-proof-resonance-integral}
 \begin{aligned}
 \int_{I_k}\frac{|H(u)|\,|f_j(u)|^n}{|u|}\,du
 &\le\frac h{4\delta}\int_{\R}
       \bigl(|u-u_{j,k}|+|u_{j,k}-u_k|\bigr)
       e^{-n(u-u_{j,k})^2/4}\,du\\
 &=\frac h{4\delta}
       \left(\frac4n+
       \frac{2\sqrt\pi}{\sqrt n}|u_{j,k}-u_k|\right)
 =o(n^{-1/2}).
 \end{aligned}
\end{equation}
The displayed integrals are evaluated by the substitution
\(v=u-u_{j,k}\):
\(\int_{\R}|v|e^{-nv^2/4}\,dv=4/n\) and
\(\int_{\R}e^{-nv^2/4}\,dv=2\sqrt{\pi/n}\).
After multiplication by \(\sqrt n\), the two terms in
\eqref{eq:jitter-proof-resonance-integral} tend to zero by
\(n=n_j\to\infty\) and \eqref{eq:jitter-proof-peak-location}.
Thus convergence of the peak locations is enough; no rate of
convergence for \(u_{j,k}-u_k\) is required.

It remains to bound the part away from these intervals. Let
\[
 E_T=\{u:\delta\le|u|\le T\}
       \setminus\bigcup_k\operatorname{int}(I_k).
\]
This is a compact set with no resonances. It is nonempty because
it contains \(-T\) and \(T\). Hence continuity of \(q\) gives
\[
 \kappa_T:=1-\max_{u\in E_T}q(u)>0,
 \qquad q(u)\le1-\kappa_T\quad(u\in E_T).
\]
Equation~\eqref{eq:jitter-proof-q-compact} gives
\(q_j(u)\le1-\kappa_T/2\) there for all sufficiently large \(j\).
Since \(|H|\le1\), we consequently have
\begin{equation}\label{eq:jitter-proof-nonresonant-integral}
 \begin{aligned}
 \int_{E_T}\frac{|H(u)|\,|f_j(u)|^n}{|u|}\,du
 &\le (1-\kappa_T/2)^{n/2}
       \int_{\delta\le|u|\le T}\frac{du}{|u|}\\
 &\le 2\log(T/\delta)\,e^{-\kappa_Tn/4}
 =o(n^{-1/2}).
 \end{aligned}
\end{equation}
In the nonlattice case there are no intervals \(I_k\):
by \eqref{eq:jitter-proof-resonances}, the entire compact set
\(\{u:\delta\le|u|\le T\}\) contains no resonance.
The same argument therefore applies to that entire set.
In the lattice case there are only finitely many \(I_k\) for
the fixed \(T\). Summing
\eqref{eq:jitter-proof-resonance-integral} and
\eqref{eq:jitter-proof-nonresonant-integral} proves in either case
\begin{equation}\label{eq:jitter-proof-away-fourier-factor}
 \int_{\delta\le|u|\le T}
       \frac{|H(u)|\,|f_j(u)|^n}{|u|}\,du=o(n^{-1/2}).
\end{equation}

\textbf{Step~6. Combining the Fourier estimates.}
The estimate for \(a_j\) away from zero must be combined
with the corresponding tail estimate for \(b_j\).
We first prove
\[
 \int_{\delta\sqrt n\le|t|\le T\sqrt n}
       \frac{|b_j(t)|}{|t|}\,dt
 =o(n^{-1/2}).
\]
We then use \(t=u\sqrt n\), for which
\(dt/|t|=du/|u|\), to express
\eqref{eq:jitter-proof-away-fourier-factor}
in the \(t\)-variable.
Combining it with this tail bound and
\eqref{eq:jitter-proof-near-zero} gives
\[
 \mathcal I_j(T)
 =
 \int_{-T\sqrt n}^{T\sqrt n}
       \frac{|a_j(t)-b_j(t)|}{|t|}\,dt
 =o(n^{-1/2}),
\]
as stated in \eqref{eq:compact-fourier-jitter}.

The formula for \(b_j\) in
\eqref{eq:jitter-proof-comparison-transform}
and the bound \(|\mu_j|\le M_\mu\) give, for \(t\ne0\),
\[
 \frac{|b_j(t)|}{|t|}
 \le e^{-t^2/2}
       \left(\frac1{|t|}+\frac{M_\mu|t|^2}{6\sqrt n}\right).
\]
On \(|t|\ge\delta\sqrt n\), we have
\(1/|t|\le1/(\delta\sqrt n)\) and
\[
 e^{-t^2/2}
 =e^{-t^2/4}e^{-t^2/4}
 \le e^{-\delta^2n/4}e^{-t^2/4}.
\]
Extending the remaining integral to all of \(\R\) therefore gives
\begin{equation}\label{eq:jitter-proof-comparison-tail}
 \begin{aligned}
 \int_{\delta\sqrt n\le|t|\le T\sqrt n}
       \frac{|b_j(t)|}{|t|}\,dt
 &\le\frac{e^{-\delta^2n/4}}{\sqrt n}
       \int_{\R}e^{-t^2/4}
       \left(\frac1\delta+\frac{M_\mu t^2}{6}\right)\,dt\\
 &=\frac{2\sqrt\pi}{\sqrt n}
       \left(\frac1\delta+\frac{M_\mu}{3}\right)
       e^{-\delta^2n/4}
 =o(n^{-1/2}).
 \end{aligned}
\end{equation}

For the other Fourier term, the substitution \(t=u\sqrt n\)
has \(dt=\sqrt n\,du\) and \(|t|=\sqrt n\,|u|\).
Using \eqref{eq:jitter-proof-sum-transform}, it gives the exact identity
\[
 \int_{\delta\sqrt n\le|t|\le T\sqrt n}
       \frac{|a_j(t)|}{|t|}\,dt
 =\int_{\delta\le|u|\le T}
       \frac{|H(u)|\,|f_j(u)|^n}{|u|}\,du.
\]
The triangle inequality, \eqref{eq:jitter-proof-away-fourier-factor},
and \eqref{eq:jitter-proof-comparison-tail} thus control the
integral of \(|a_j-b_j|/|t|\) on
\(\delta\sqrt n\le|t|\le T\sqrt n\).
Adding the estimate \eqref{eq:jitter-proof-near-zero} on the
remaining interval proves that, for every fixed \(T\) chosen above,
\begin{equation}\label{eq:compact-fourier-jitter}
 \int_{-T\sqrt n}^{T\sqrt n}
       \frac{|a_j(t)-b_j(t)|}{|t|}\,dt=o(n^{-1/2}).
\end{equation}

\textbf{Step~7. Applying the smoothing inequality.}
We finally verify the remaining hypotheses of
Lemma~\ref{lem:signed-smoothing} and apply it with
\(L=T\sqrt n\).
This converts the Fourier estimate into
\[
 \sqrt n\,\|\mathcal J_j-G_j\|_\infty
 \le
 C_1\sqrt n\,\mathcal I_j(T)+\frac{C_2M_g}{T}.
\]
For each fixed \(T>\delta\) with \(|f(\pm T)|<1\),
\eqref{eq:compact-fourier-jitter} makes the first
term tend to zero as \(j\to\infty\).
We therefore obtain
\[
 \limsup_{j\to\infty}
 \sqrt{n_j}\,\|\mathcal J_j-G_j\|_\infty
 \le\frac{C_2M_g}{T}.
\]
Since \(C_2\) and \(M_g\) do not depend on \(T\),
letting \(T\to\infty\) through cutoffs satisfying
\(|f(\pm T)|<1\) will prove
\eqref{eq:jitter-expansion}.

Recall that
\[
 \mathcal J_j(x)=\PP\{S_j+U/\sqrt n\le x\}.
\]
The independence of the summands and their first two moments in
\eqref{eq:jitter-proof-moments} give
\[
 \begin{aligned}
 \EE S_j^2
 &=\frac1n\left[\sum_{r=1}^n\EE X_{j,r}^2
       +2\sum_{1\le r<s\le n}\EE(X_{j,r}X_{j,s})\right]\\
 &=\frac1n\left[n+2\sum_{1\le r<s\le n}
                  (\EE X_{j,r})(\EE X_{j,s})\right]=1.
 \end{aligned}
\]
The Cauchy--Schwarz inequality therefore gives
\[
 \EE|S_j+U/\sqrt n|
 \le\EE|S_j|+\frac{\EE|U|}{\sqrt n}
 \le1+\frac h{4\sqrt n}<\infty.
\]
Thus \(\mathcal J_j\) is the distribution function of a probability
measure with finite first absolute moment.
The density \(g_j\) has integral one by
\eqref{eq:jitter-proof-comparison-mass}.
The bound in \eqref{eq:jitter-proof-density-moments} shows that
\(\nu_j(dx)=g_j(x)\,dx\) is a finite signed measure with finite
first absolute moment.
Its primitive is \(G_j\): we have \(G_j'=g_j\) by
\eqref{eq:jitter-proof-comparison-density}, and
\eqref{eq:jitter-proof-comparison-function} gives
\(G_j(x)\to0\) as \(x\to-\infty\), because the Gaussian
correction tends to zero there. Consequently,
\[
 G_j(x)=\int_{-\infty}^x g_j(s)\,ds.
\]
Furthermore, \eqref{eq:jitter-proof-density-sup} implies
\[
 |G_j(x)-G_j(y)|
 =\left|\int_y^x g_j(s)\,ds\right|
 \le M_g|x-y|.
\]
Here \(M_g=\sup_j\|g_j\|_\infty<\infty\) is independent of \(j\)
and of the cutoff \(T\).
The Fourier transforms of the two measures are \(a_j\) and
\(b_j\), by \eqref{eq:jitter-proof-sum-transform} and
\eqref{eq:jitter-proof-comparison-transform}.

We can therefore apply Lemma~\ref{lem:signed-smoothing}
with \(L=T\sqrt n\), \(F=\mathcal J_j\), and \(\nu=\nu_j\).
Its bound \eqref{eq:signed-smoothing} gives
\begin{equation}\label{eq:jitter-proof-final-smoothing}
 \sqrt n\,\|\mathcal J_j-G_j\|_\infty
 \le C_1\sqrt n\int_{-T\sqrt n}^{T\sqrt n}
       \frac{|a_j(t)-b_j(t)|}{|t|}\,dt
       +\frac{C_2M_g}{T}.
\end{equation}
For this fixed \(T\), the integral term on the right of
\eqref{eq:jitter-proof-final-smoothing} tends to zero by
\eqref{eq:compact-fourier-jitter}. Hence
\begin{equation}\label{eq:jitter-proof-fixed-cutoff-limit}
 \limsup_{j\to\infty}
       \sqrt{n_j}\,\|\mathcal J_j-G_j\|_\infty
 \le\frac{C_2M_g}{T}.
\end{equation}
The expression on the left does not depend on \(T\).
The resonance set is either \(\{0\}\) or a discrete lattice,
so there are arbitrarily large choices of \(T\) with
nonresonant endpoints.
The constants \(C_2\) and \(M_g\) are independent of all these
choices. Letting \(T\) tend to infinity through them makes the
right side of \eqref{eq:jitter-proof-fixed-cutoff-limit} tend
to zero. Since the quantities on the left are nonnegative,
\[
 \sqrt{n_j}\,\|\mathcal J_j-G_j\|_\infty\longrightarrow0.
\]
This is \eqref{eq:jitter-expansion}.
The neighborhoods \(I_k\), the gaps \(\kappa_T\), and the
indices after which the estimates hold were all chosen for
each fixed \(T\). We first took the upper limit in \(j\) with
that cutoff fixed, and only then sent \(T\) to infinity;
no bound uniform in growing \(T\), and no interchange of
these two limiting operations, was used.
\end{proof}

We next extend Lemma~\ref{lem:jitter} to varying jitter widths
\(h_j\to h(P)\), and state the corresponding expansion for
summand laws that are not already standardized.
These extensions will be used in
Lemmas~\ref{lem:binomial-estimates}
and~\ref{lem:accumulated-cluster-variance}.
In both applications, we add an independent uniform variable
on \([-1/2,1/2]\) to the original sum.
Although this noise has fixed width before standardization,
its width after division by the summand standard deviation
varies with the summand law.

\begin{proposition}[Varying jitter widths and affine normalization]
\label{prop:jitter-width}
\textup{(i)} Under the assumptions and notation of
Lemma~\ref{lem:jitter}, let \(h_j\ge0\) satisfy \(h_j\to h\).
Let \(U_j\), independent of the row, be uniform on
\([-h_j/2,h_j/2]\), with \(U_j=0\) when \(h_j=0\). Then
\begin{equation}\label{eq:jitter-variable-width}
 \sup_{x\in\R}\left|
 \PP\{S_j+U_j/\sqrt{n_j}\le x\}
 -\Phi(x)-\frac{\mu_j}{6\sqrt{n_j}}(1-x^2)\phi(x)
 \right|=o(n_j^{-1/2}).
\end{equation}

\textup{(ii)} Let \(Q_j\to Q\) in \(W_3\), and suppose that
\(Q\) has mean \(m\) and variance \(\sigma^2>0\).
Write \(m_j\) and \(\sigma_j\) for the mean and standard
deviation of \(Q_j\). Let \(n_j\to\infty\) be positive integers,
and let \(Y_{j,1},\ldots,Y_{j,n_j}\) be independent with common
law \(Q_j\). Let \(d\) be the maximal lattice span of \(Q\),
with \(d=0\) if \(Q\) is nonlattice, and suppose that
\(d_j\ge0\) and \(d_j\to d\).
Let \(H_j\), independent of the row, be uniform on
\([-d_j/2,d_j/2]\), with \(H_j=0\) when \(d_j=0\).
Then \(\sigma_j>0\) for all sufficiently large \(j\).
For those \(j\), put
\[
 T_j=\frac{\sum_{r=1}^{n_j}Y_{j,r}-n_jm_j}
              {\sigma_j\sqrt{n_j}},
 \qquad
 \gamma_j=\frac{\int (y-m_j)^3\,Q_j(dy)}{\sigma_j^3}.
\]
Then
\begin{equation}\label{eq:jitter-raw-expansion}
 \sup_{x\in\R}\left|
 \PP\left\{T_j+\frac{H_j}{\sigma_j\sqrt{n_j}}\le x\right\}
 -\Phi(x)-\frac{\gamma_j}{6\sqrt{n_j}}(1-x^2)\phi(x)
 \right|=o(n_j^{-1/2}).
\end{equation}
\end{proposition}

\begin{proof}
\textbf{Changing the jitter width.}
For each \(j\), take \(V\) uniform on \([-1/2,1/2]\),
independently of the row, and couple the two jitters as
\(U=hV\) and \(U_j=h_jV\). Define
\[
 A_j=S_j+\frac{hV}{\sqrt{n_j}},\qquad
 B_j=S_j+\frac{h_jV}{\sqrt{n_j}},\qquad
 \delta_j=\frac{|h_j-h|}{2\sqrt{n_j}}.
\]
Since \(|V|\le1/2\), we have \(|A_j-B_j|\le\delta_j\)
almost surely. Consequently,
\[
 \{A_j\le x-\delta_j\}\subseteq\{B_j\le x\}
 \subseteq\{A_j\le x+\delta_j\}.
\]
Writing \(J_j(x)=\PP(A_j\le x)\) and
\(\widetilde J_j(x)=\PP(B_j\le x)\), we obtain
\begin{equation}\label{eq:jitter-width-sandwich}
 J_j(x-\delta_j)\le\widetilde J_j(x)
 \le J_j(x+\delta_j),\qquad x\in\R.
\end{equation}

Recall the comparison function from
\eqref{eq:jitter-proof-comparison-function}, with \(n=n_j\):
\[
 G_j(x)=\Phi(x)+\frac{\mu_j}{6\sqrt{n_j}}(1-x^2)\phi(x).
\]
By \eqref{eq:jitter-expansion},
\(r_j:=\|J_j-G_j\|_\infty=o(n_j^{-1/2})\).
Moreover, \eqref{eq:jitter-proof-comparison-density} and
\eqref{eq:jitter-proof-density-sup} give
\(M_g:=\sup_j\|G_j'\|_\infty<\infty\).
The mean value theorem therefore implies
\[
 |G_j(x\pm\delta_j)-G_j(x)|\le M_g\delta_j.
\]
Using these bounds in \eqref{eq:jitter-width-sandwich} yields
\[
 \begin{aligned}
 \widetilde J_j(x)
 &\le G_j(x+\delta_j)+r_j
 \le G_j(x)+M_g\delta_j+r_j,\\
 \widetilde J_j(x)
 &\ge G_j(x-\delta_j)-r_j
 \ge G_j(x)-M_g\delta_j-r_j.
 \end{aligned}
\]
Thus
\begin{equation}\label{eq:jitter-width-error}
 \|\widetilde J_j-G_j\|_\infty
 \le r_j+\frac{M_g|h_j-h|}{2\sqrt{n_j}}
 =o(n_j^{-1/2}),
\end{equation}
where the last equality uses \(h_j\to h\).
This proves \eqref{eq:jitter-variable-width}, including
\(h=0\).

\textbf{Centering and scaling the summand laws.}
We first verify that the standardized laws converge in
\(W_3\). For each \(j\), the definition of \(W_3\) allows
us to couple \(Z_j\sim Q_j\) and \(Z_j^*\sim Q\) so that
\[
 a_j:=\|Z_j-Z_j^*\|_3
 \le W_3(Q_j,Q)+j^{-1}\longrightarrow0,
\]
where \(\|Z\|_r=(\EE|Z|^r)^{1/r}\).
Since these variables have means \(m_j\) and \(m\),
\[
 |m_j-m|\le\EE|Z_j-Z_j^*|\le a_j.
\]
The reverse triangle inequality for the \(L^2\) norm gives
\[
 \begin{aligned}
 |\sigma_j-\sigma|
 &=\bigl|\|Z_j-m_j\|_2-\|Z_j^*-m\|_2\bigr|\\
 &\le\|(Z_j-m_j)-(Z_j^*-m)\|_2\\
 &\le\|Z_j-Z_j^*\|_2+|m_j-m|
 \le2a_j\longrightarrow0.
 \end{aligned}
\]
Here we used \(\|Z\|_2\le\|Z\|_3\) on a probability
space. In particular, \(\sigma_j\to\sigma>0\), so we may
restrict attention to indices for which \(\sigma_j>0\).
For those indices, the identity
\[
 \frac{Z_j-m_j}{\sigma_j}-\frac{Z_j^*-m}{\sigma}
 =\frac{(Z_j-Z_j^*)-(m_j-m)}{\sigma_j}
  +\left(\frac1{\sigma_j}-\frac1\sigma\right)(Z_j^*-m)
\]
and the triangle inequality for the \(L^3\) norm imply
\begin{equation}\label{eq:jitter-width-standardization}
 \left\|\frac{Z_j-m_j}{\sigma_j}
             -\frac{Z_j^*-m}{\sigma}\right\|_3
 \le\frac{2a_j}{\sigma_j}
 +\left|\frac1{\sigma_j}-\frac1\sigma\right|
       \left(\int|y-m|^3\,Q(dy)\right)^{1/3}
 \longrightarrow0.
\end{equation}
The integral is finite because \(Q\) has a finite third
absolute moment. Let \(P_j\) and \(P\) denote the laws of
\((Z_j-m_j)/\sigma_j\) and \((Z_j^*-m)/\sigma\), respectively.
They have mean zero and variance one, and
\eqref{eq:jitter-width-standardization} implies
\(W_3(P_j,P)\to0\).

Centering does not change a lattice span, and division by
\(\sigma>0\) divides every lattice span by \(\sigma\).
Indeed, support on \(a+d\mathbb Z\) is transformed into
support on \((a-m)/\sigma+(d/\sigma)\mathbb Z\), and the
inverse affine transformation gives the converse.
Thus maximality of the span is preserved under this scaling.
The same correspondence preserves the nonlattice case, so
\[
 h(P)=\frac d\sigma,\qquad
 h_j:=\frac{d_j}{\sigma_j}\longrightarrow\frac d\sigma.
\]
The variables \(X_{j,r}=(Y_{j,r}-m_j)/\sigma_j\) have law
\(P_j\), and their third moments equal \(\gamma_j\).
Also, \(H_j/\sigma_j\) is independent of the row and is
uniform on \([-h_j/2,h_j/2]\), with value zero when
\(h_j=0\). Finally,
\[
 T_j+\frac{H_j}{\sigma_j\sqrt{n_j}}
 =\frac1{\sqrt{n_j}}\sum_{r=1}^{n_j}X_{j,r}
  +\frac{H_j/\sigma_j}{\sqrt{n_j}}.
\]
Applying \eqref{eq:jitter-variable-width} to these standardized
laws and widths proves \eqref{eq:jitter-raw-expansion}.
The argument also applies when \(d=0\), since then
\(h(P)=0\) and \(h_j\to0\).
\end{proof}

We now use Lemma~\ref{lem:jitter} to bound the distribution
function of the sum without the added uniform variable.
The resulting upper and lower bounds retain the sign of the
third-moment correction and give an asymptotic bound on the
absolute approximation error.
The absolute-error bound will be used in the proof of
Proposition~\ref{prop:clusters}.
The signed upper bound will be used in
Lemma~\ref{lem:limit-extremizer} to identify the limiting
maximizing law and show that the standardized maximizing
thresholds converge to zero.

\begin{corollary}\label{cor:jitter-envelopes}
Under the assumptions of Lemma~\ref{lem:jitter}, let
\(F_j(x)=\PP(S_j\le x)\).  Uniformly in \(x\in\R\),
\begin{align}
 F_j(x)-\Phi(x)
 &\le\frac1{\sqrt{n_j}}
       \left[\frac h2+\frac{\mu_j}{6}(1-x^2)\right]\phi(x)
       +o(n_j^{-1/2}),\label{eq:jitter-upper}\\
 F_j(x)-\Phi(x)
 &\ge\frac1{\sqrt{n_j}}
       \left[-\frac h2+\frac{\mu_j}{6}(1-x^2)\right]\phi(x)
       +o(n_j^{-1/2}).\label{eq:jitter-lower}
\end{align}
Consequently, with \(\mu=\int x^3dP(x)\),
\begin{equation}\label{eq:triangular-upper}
 \limsup_j\sqrt{n_j}\,\|F_j-\Phi\|_\infty
 \le\frac{|\mu|+3h}{6\sqrt{2\pi}}.
\end{equation}
\end{corollary}

\begin{proof}
\textbf{Removing the added uniform variable.}
Let \(U\) be the independent uniform variable from
Lemma~\ref{lem:jitter}, and define
\[
 J_j(x)=\PP\{S_j+U/\sqrt{n_j}\le x\},
 \qquad a_j=\frac{h}{2\sqrt{n_j}}.
\]
Since \(|U|/\sqrt{n_j}\le a_j\) almost surely, we have
\[
 \{S_j+U/\sqrt{n_j}\le x-a_j\}
 \subseteq\{S_j\le x\}
 \subseteq\{S_j+U/\sqrt{n_j}\le x+a_j\}.
\]
Taking probabilities gives
\begin{equation}\label{eq:jitter-envelope-sandwich}
 J_j(x-a_j)\le F_j(x)\le J_j(x+a_j),
 \qquad x\in\R.
\end{equation}

\textbf{Expanding the two shifted bounds.}
Put \(\psi(x)=(1-x^2)\phi(x)\), and recall from
\eqref{eq:jitter-proof-comparison-function} that
\[
 G_j(x)=\Phi(x)+\frac{\mu_j}{6\sqrt{n_j}}\psi(x).
\]
Equation~\eqref{eq:jitter-expansion} gives
\begin{equation}\label{eq:jitter-envelope-smoothed-error}
 r_j:=\|J_j-G_j\|_\infty=o(n_j^{-1/2}).
\end{equation}
By \eqref{eq:jitter-proof-moments},
\(M_\mu:=\sup_j|\mu_j|<\infty\).
Moreover,
\[
 \phi'(x)=-x\phi(x),\qquad
 \psi'(x)=(x^3-3x)\phi(x),
\]
so both derivatives are bounded. Write
\(M_1=\|\phi'\|_\infty\) and
\(M_2=\|\psi'\|_\infty\).
Since \(\Phi''=\phi'\), Taylor's theorem and the mean
value theorem give, for \(\theta\in\{-1,1\}\),
\[
 \begin{aligned}
 |\Phi(x+\theta a_j)-\Phi(x)-\theta a_j\phi(x)|
 &\le\frac{M_1a_j^2}{2},\\
 |\psi(x+\theta a_j)-\psi(x)|
 &\le M_2a_j.
 \end{aligned}
\]
Consequently,
\begin{equation}\label{eq:jitter-envelope-shifted-expansion}
 G_j(x+\theta a_j)-\Phi(x)
 =\frac1{\sqrt{n_j}}
   \left[\frac{\theta h}{2}
         +\frac{\mu_j}{6}(1-x^2)\right]\phi(x)
   +E_{j,\theta}(x),
\end{equation}
where, uniformly in \(x\in\R\) and \(\theta\in\{-1,1\}\),
\begin{equation}\label{eq:jitter-envelope-taylor-error}
 |E_{j,\theta}(x)|
 \le\frac{M_1a_j^2}{2}
       +\frac{M_\mu M_2a_j}{6\sqrt{n_j}}
 =\frac{K}{n_j},
 \qquad
 K:=\frac{M_1h^2}{8}+\frac{M_\mu M_2h}{12}.
\end{equation}
The constant \(K\) is independent of \(j\) and \(x\).
Combining \eqref{eq:jitter-envelope-smoothed-error} with
\eqref{eq:jitter-envelope-taylor-error}, put
\[
 e_j=r_j+\frac{K}{n_j}=o(n_j^{-1/2}),
\]
where \(n_j^{-1}=o(n_j^{-1/2})\) because \(n_j\to\infty\).
In \eqref{eq:jitter-envelope-sandwich}, use
\eqref{eq:jitter-envelope-shifted-expansion} with \(\theta=1\)
for the upper bound and with \(\theta=-1\) for the lower bound.
Together with \eqref{eq:jitter-envelope-smoothed-error}, this gives
\[
 \begin{aligned}
 F_j(x)-\Phi(x)
 &\le\frac1{\sqrt{n_j}}
       \left[\frac h2+\frac{\mu_j}{6}(1-x^2)\right]\phi(x)
       +e_j,\\
 F_j(x)-\Phi(x)
 &\ge\frac1{\sqrt{n_j}}
       \left[-\frac h2+\frac{\mu_j}{6}(1-x^2)\right]\phi(x)
       -e_j.
 \end{aligned}
\]
These are \eqref{eq:jitter-upper} and \eqref{eq:jitter-lower},
with remainders uniform in \(x\).
When \(h=0\), both \(a_j\) and \(K\) are zero, so the
same argument applies.

\textbf{Bounding the absolute approximation error.}
Subtracting the common correction term from the preceding
two inequalities yields
\[
 \left|F_j(x)-\Phi(x)
       -\frac{\mu_j}{6\sqrt{n_j}}\psi(x)\right|
 \le\frac{h}{2\sqrt{n_j}}\phi(x)+e_j.
\]
The triangle inequality therefore gives
\begin{equation}\label{eq:jitter-envelope-absolute-bound}
 |F_j(x)-\Phi(x)|
 \le\frac1{\sqrt{n_j}}
       \left[\frac h2\phi(x)+\frac{|\mu_j|}{6}|\psi(x)|\right]
       +e_j.
\end{equation}
We now verify the two Gaussian suprema used in this bound.
Since \(\phi(x)=\phi(0)e^{-x^2/2}\), its maximum is
\(\phi(0)=(2\pi)^{-1/2}\).
For \(|x|\le1\),
\[
 |\psi(x)|=(1-x^2)\phi(x)\le\phi(0),
\]
with equality at \(x=0\).
For \(|x|\ge1\), evenness reduces the calculation to
\(x\ge1\), where
\[
 \frac{d}{dx}\bigl[(x^2-1)\phi(x)\bigr]
 =x(3-x^2)\phi(x).
\]
This derivative is positive for \(1\le x<\sqrt3\) and
negative for \(x>\sqrt3\). Thus the maximum on this region is
\[
 2\phi(\sqrt3)=2e^{-3/2}\phi(0)<\phi(0).
\]
We have proved
\begin{equation}\label{eq:jitter-envelope-gaussian-suprema}
 \sup_{x\in\R}\phi(x)
 =\sup_{x\in\R}|\psi(x)|
 =\phi(0)=\frac1{\sqrt{2\pi}}.
\end{equation}
Taking the supremum in \eqref{eq:jitter-envelope-absolute-bound}
and using \eqref{eq:jitter-envelope-gaussian-suprema}, we obtain
\[
 \sqrt{n_j}\,\|F_j-\Phi\|_\infty
 \le\frac{3h+|\mu_j|}{6\sqrt{2\pi}}+\sqrt{n_j}\,e_j.
\]
Finally, \(\sqrt{n_j}\,e_j\to0\), and
\(\mu_j\to\mu\) by \eqref{eq:jitter-proof-moments}.
Taking the upper limit proves \eqref{eq:triangular-upper}.
\end{proof}

We now specialize the smoothing expansion to binomial sums.
We will need estimates for point probabilities and for the
signed discrepancies on both sides of each binomial jump,
uniformly when the success probability lies in a fixed compact
subinterval of \((0,1)\).
These estimates will be used in
Lemma~\ref{lem:small-cluster-variance} to compare a perturbed
Bernoulli sum with the corresponding binomial sum, treating
central and distant thresholds separately.
The estimate for point probabilities in a central interval
will also be used in Lemma~\ref{lem:two-cluster-local-mass}
to obtain a lower bound for interval probabilities.

For a Bernoulli random variable \(B_p\) with success
probability \(p\), put
\[
 q=1-p,\qquad v=pq,\qquad
 \tau=p^2+q^2,\qquad d=1-2p.
\]
Since \(\EE B_p=p\), the centered variable \(B_p-p\) equals
\(q\) with probability \(p\) and \(-p\) with probability \(q\).
Consequently,
\[
 \begin{aligned}
 \EE(B_p-p)^2&=pq^2+qp^2=v,\\
 \EE|B_p-p|^3&=pq^3+qp^3=v\tau,\\
 \EE(B_p-p)^3&=pq^3-qp^3=vd.
 \end{aligned}
\]
Dividing the last two moments by \(v^{3/2}\) gives
\begin{equation}\label{eq:binomial-standardized-moments}
 \beta(B_p)=\frac{\tau}{\sqrt v},\qquad
 \frac{\EE(B_p-p)^3}{v^{3/2}}=\frac d{\sqrt v}.
\end{equation}

\begin{lemma}\label{lem:binomial-estimates}
Fix a compact interval \(I\subset(0,1)\).  Let \(K\sim\mathrm{Bin}(n,p)\)
with \(p\in I\), and define
\[
 b_{n,k}(p)=\PP(K=k),\qquad F_{n,p}(k)=\PP(K\le k),\qquad
 z_k=\frac{k-np}{\sqrt{nv}},\qquad k\in\mathbb Z.
\]
There is a constant \(C_I\) such that
\begin{equation}\label{eq:binomial-mass-bound}
 \sup_{p\in I}\sup_{k\in\mathbb Z}b_{n,k}(p)
 \le C_I n^{-1/2},\qquad n\ge1.
\end{equation}
As \(n\to\infty\), uniformly over \(p\in I\) and \(k\in\mathbb Z\),
\begin{equation}\label{eq:binomial-local}
 b_{n,k}(p)=\frac{\phi(z_k)}{\sqrt{nv}}+o(n^{-1/2}).
\end{equation}
In particular, for every fixed \(0<M<\infty\),
\[
 \sup_{\substack{p\in I,\ k\in\mathbb Z\\ |z_k|\le M}}
 \left|
 \frac{\sqrt{np(1-p)}\,b_{n,k}(p)}{\phi(z_k)}-1
 \right|\longrightarrow0.
\]

Define the normalized discrepancies at the two sides of a jump by
\[
 \begin{split}
 R^+_{n,k}(B_p)
 &=\frac{\sqrt n}{\beta(B_p)}
        [F_{n,p}(k)-\Phi(z_k)],\\
 R^-_{n,k}(B_p)
 &=\frac{\sqrt n}{\beta(B_p)}
        [\Phi(z_k)-F_{n,p}(k-1)].
 \end{split}
\]
Then, uniformly over \(p\in I\) and \(k\in\mathbb Z\),
\begin{equation}\label{eq:binomial-envelopes}
 R^\pm_{n,k}(B_p)=A_p^\pm(z_k)+o(1),
\end{equation}
where
\[
 A_p^+(z)=\frac{\phi(z)}{6\tau}(3+d-dz^2),\qquad
 A_p^-(z)=\frac{\phi(z)}{6\tau}(3-d+dz^2).
\]
These functions tend uniformly to zero as \(|z|\to\infty\), for
\(p\in I\).  If \(k\) is an integer nearest to \(np\), then
\begin{equation}\label{eq:binomial-central-branch}
 R^+_{n,k}(B_p)=g(p)+o(1),\qquad
 g(p)=\frac{\phi(0)(2-p)}{3[p^2+(1-p)^2]},
\end{equation}
uniformly over \(p\in I\).
\end{lemma}

\begin{proof}
\textbf{Bounding the point probabilities.}
Compactness of \(I\subset(0,1)\) gives
\begin{equation}\label{eq:binomial-parameter-bounds}
 v_I:=\min_{p\in I}p(1-p)>0,\qquad
 v_I\le v\le\frac14,\qquad
 \frac12\le\tau=1-2v\le1,\qquad |d|\le1.
\end{equation}
The binomial theorem gives
\[
 (q+pe^{it})^n=\sum_{\ell=0}^n b_{n,\ell}(p)e^{i\ell t}.
\]
For integers \(\ell,k\),
\[
 \frac1{2\pi}\int_{-\pi}^{\pi}e^{i(\ell-k)t}\,dt
 =\begin{cases}1,&\ell=k,\\0,&\ell\ne k.\end{cases}
\]
Multiplying the finite sum by \(e^{-ikt}\) and integrating
term by term therefore yields, for every \(k\in\mathbb Z\),
\begin{equation}\label{eq:binomial-fourier-inversion}
 b_{n,k}(p)
 =\frac1{2\pi}\int_{-\pi}^{\pi}
       (q+pe^{it})^n e^{-ikt}\,dt
 \le\frac1{2\pi}\int_{-\pi}^{\pi}|q+pe^{it}|^n\,dt.
\end{equation}
Here \(b_{n,k}(p)=0\) when \(k\notin\{0,\ldots,n\}\).
Direct calculation gives
\[
 |q+pe^{it}|^2
 =q^2+p^2+2pq\cos t
 =1-4v\sin^2(t/2).
\]
For \(0\le u\le1\), the inequality \(1-u\le e^{-u}\)
implies \((1-u)^{n/2}\le e^{-nu/2}\).
Also, concavity of sine on \([0,\pi/2]\) gives
\(\sin s\ge2s/\pi\) on this interval.
Applying these inequalities with
\(u=4v\sin^2(t/2)\) and \(s=|t|/2\), we obtain
\[
 |q+pe^{it}|^n
 \le e^{-2nv\sin^2(t/2)}
 \le e^{-2nvt^2/\pi^2},\qquad |t|\le\pi.
\]
Thus \eqref{eq:binomial-fourier-inversion} implies
\[
 \begin{aligned}
 b_{n,k}(p)
 &\le\frac1{2\pi}\int_{\R}e^{-2nvt^2/\pi^2}\,dt\\
 &=\frac{\sqrt\pi}{2\sqrt{2nv}}
 \le\frac{\sqrt\pi}{2\sqrt{2v_I}}\,n^{-1/2}.
 \end{aligned}
\]
This proves \eqref{eq:binomial-mass-bound}, with
\(C_I=\sqrt\pi/(2\sqrt{2v_I})\).

\textbf{Obtaining a smoothing expansion uniformly in \(p\).}
Let \(H_0\), independent of \(K\), be uniform on
\([-1/2,1/2]\), and define
\[
 \mathcal J_{n,p}(x)
 =\PP\left\{\frac{K-np+H_0}{\sqrt{nv}}\le x\right\},
 \qquad \psi(x)=(1-x^2)\phi(x).
\]
We claim that
\begin{equation}\label{eq:binomial-jitter}
 \mathcal J_{n,p}(x)
 =\Phi(x)+\frac{d}{6\sqrt{nv}}\psi(x)+o(n^{-1/2}),
\end{equation}
uniformly over \(p\in I\) and \(x\in\R\).
To prove this uniformity, set
\[
 \mathcal R_{n,p}(x)
 =\mathcal J_{n,p}(x)-\Phi(x)
  -\frac{1-2p}{6\sqrt{np(1-p)}}\psi(x),
 \qquad
 r_n=\sup_{p\in I}\|\mathcal R_{n,p}\|_\infty.
\]
If \(\sqrt n\,r_n\) did not tend to zero, there would be
\(\varepsilon>0\) and integers \(n_j\to\infty\) such that
\(\sqrt{n_j}\,r_{n_j}>2\varepsilon\).
By the definition of the supremum, we could choose
\(p_j\in I\) satisfying
\begin{equation}\label{eq:binomial-jitter-contradiction}
 \sqrt{n_j}\,\|\mathcal R_{n_j,p_j}\|_\infty>\varepsilon.
\end{equation}
After passing to a subsequence, compactness gives
\(p_j\to p_*\in I\).
For \(V\) uniform on \([0,1]\), the variables
\(\ind_{\{V\le p_j\}}\) and \(\ind_{\{V\le p_*\}}\)
have Bernoulli parameters \(p_j\) and \(p_*\), and
\[
 \EE\left|\ind_{\{V\le p_j\}}-\ind_{\{V\le p_*\}}\right|^3
 =|p_j-p_*|.
\]
Hence the definition of \(W_3\) gives
\[
 W_3\bigl(\mathrm{Bernoulli}(p_j),
          \mathrm{Bernoulli}(p_*)\bigr)
 \le |p_j-p_*|^{1/3}\longrightarrow0.
\]
The limiting Bernoulli law has positive variance
\(p_*(1-p_*)\) and maximal lattice span one.
Apply Proposition~\ref{prop:jitter-width}\textup{(ii)},
specifically \eqref{eq:jitter-raw-expansion}, with these
summand laws, sample sizes \(n_j\), and raw jitter widths
equal to one. The moment identities
\eqref{eq:binomial-standardized-moments} show that the
resulting expansion has exactly the correction term in
\(\mathcal R_{n_j,p_j}\).
It follows that
\(\sqrt{n_j}\|\mathcal R_{n_j,p_j}\|_\infty\to0\),
contradicting \eqref{eq:binomial-jitter-contradiction}.
We have therefore proved \eqref{eq:binomial-jitter} and
\begin{equation}\label{eq:binomial-uniform-jitter-error}
 r_n=o(n^{-1/2}).
\end{equation}

\textbf{Evaluating the expansion at the two sides of each jump.}
For any integer \(\ell\), the event \(K\le\ell\) implies
\(K+H_0\le\ell+1/2\).
If \(K\ge\ell+1\), the latter inequality can hold only
when \(K=\ell+1\) and \(H_0=-1/2\), an event of
probability zero. Therefore,
\[
 \PP(K+H_0\le\ell+1/2)=\PP(K\le\ell).
\]
Put \(a_{n,p}=1/(2\sqrt{nv})\).
Taking \(\ell=k\) and \(\ell=k-1\) in this identity gives
\begin{equation}\label{eq:binomial-jump-identities}
 \begin{aligned}
 F_{n,p}(k)&=\mathcal J_{n,p}(z_k+a_{n,p}),\\
 F_{n,p}(k-1)&=\mathcal J_{n,p}(z_k-a_{n,p}).
 \end{aligned}
\end{equation}

To expand these two values, write
\(M_1=\|\phi'\|_\infty\) and
\(M_2=\|\psi'\|_\infty\).
These constants are finite because
\(\phi'(z)=-z\phi(z)\) and
\(\psi'(z)=(z^3-3z)\phi(z)\).
Taylor's theorem for \(\Phi\) and the mean value theorem
for \(\psi\) give, for \(\theta\in\{-1,1\}\),
\[
 \begin{aligned}
 |\Phi(z+\theta a_{n,p})-\Phi(z)
             -\theta a_{n,p}\phi(z)|
 &\le\frac{M_1a_{n,p}^2}{2},\\
 |\psi(z+\theta a_{n,p})-\psi(z)|
 &\le M_2a_{n,p}.
 \end{aligned}
\]
The total Taylor error in the comparison function is thus
at most
\begin{equation}\label{eq:binomial-taylor-error}
 \frac{M_1a_{n,p}^2}{2}
 +\frac{|d|M_2a_{n,p}}{6\sqrt{nv}}
 \le\frac{M_1/8+M_2/12}{nv_I}
 =\frac{C_T}{n},
 \qquad C_T:=\frac{M_1/8+M_2/12}{v_I},
\end{equation}
where we used \eqref{eq:binomial-parameter-bounds}.
Combining this estimate with
\eqref{eq:binomial-uniform-jitter-error} gives
\begin{equation}\label{eq:binomial-shifted-expansion}
 \mathcal J_{n,p}(z+\theta a_{n,p})-\Phi(z)
 =\frac{\phi(z)}{\sqrt{nv}}
       \left[\frac\theta2+\frac d6(1-z^2)\right]
       +\mathcal E_{n,p,\theta}(z),
\end{equation}
with
\[
 \sup_{\substack{p\in I,\ z\in\R\\\theta\in\{-1,1\}}}
       |\mathcal E_{n,p,\theta}(z)|
 \le r_n+\frac{C_T}{n}=:e_n=o(n^{-1/2}).
\]
Using \eqref{eq:binomial-jump-identities} in
\eqref{eq:binomial-shifted-expansion}, we obtain
\begin{equation}\label{eq:binomial-jump-expansions}
 \begin{aligned}
 F_{n,p}(k)-\Phi(z_k)
 &=\frac{\phi(z_k)}{\sqrt{nv}}
       \left[\frac12+\frac d6(1-z_k^2)\right]
       +\varepsilon^+_{n,k}(p),\\
 \Phi(z_k)-F_{n,p}(k-1)
 &=\frac{\phi(z_k)}{\sqrt{nv}}
       \left[\frac12-\frac d6(1-z_k^2)\right]
       +\varepsilon^-_{n,k}(p).
 \end{aligned}
\end{equation}
Here
\(\varepsilon^+_{n,k}(p)=\mathcal E_{n,p,1}(z_k)\) and
\(\varepsilon^-_{n,k}(p)=-\mathcal E_{n,p,-1}(z_k)\), so
\[
 \sup_{p\in I,\,k\in\mathbb Z}
 \max\{|\varepsilon^+_{n,k}(p)|,
       |\varepsilon^-_{n,k}(p)|\}\le e_n.
\]

\textbf{Deriving the local limit and normalized discrepancies.}
Adding the two equations in \eqref{eq:binomial-jump-expansions}
cancels their third-moment correction terms and gives
\[
 \begin{aligned}
 b_{n,k}(p)
 &=F_{n,p}(k)-F_{n,p}(k-1)\\
 &=\frac{\phi(z_k)}{\sqrt{nv}}
   +\varepsilon^+_{n,k}(p)+\varepsilon^-_{n,k}(p).
 \end{aligned}
\]
In particular,
\begin{equation}\label{eq:binomial-local-error}
 \sup_{p\in I,\,k\in\mathbb Z}
 \left|b_{n,k}(p)-\frac{\phi(z_k)}{\sqrt{nv}}\right|
 \le2e_n=o(n^{-1/2}),
\end{equation}
which proves \eqref{eq:binomial-local}.
For fixed \(0<M<\infty\) and \(|z_k|\le M\),
\(\phi(z_k)\ge\phi(M)>0\).
Using \eqref{eq:binomial-local-error} and \(\sqrt v\le1/2\),
we obtain
\[
 \left|\frac{\sqrt{nv}\,b_{n,k}(p)}{\phi(z_k)}-1\right|
 \le\frac{2\sqrt{nv}\,e_n}{\phi(M)}
 \le\frac{\sqrt n\,e_n}{\phi(M)}\longrightarrow0,
\]
uniformly over the stated central interval.
This proves the relative local limit estimate.

By \eqref{eq:binomial-standardized-moments},
\(\sqrt n/\beta(B_p)=\sqrt{nv}/\tau\).
Multiplying \eqref{eq:binomial-jump-expansions} by this factor
therefore gives
\[
 R^\pm_{n,k}(B_p)
 =\frac{\phi(z_k)}{6\tau}\bigl[3\pm d(1-z_k^2)\bigr]
  +\frac{\sqrt{nv}}\tau\varepsilon^\pm_{n,k}(p).
\]
The first term is \(A_p^\pm(z_k)\).
By \eqref{eq:binomial-parameter-bounds},
\(\sqrt v/\tau\le1\), and hence the second term is bounded
in absolute value by \(\sqrt n\,e_n\to0\), uniformly in
\(p\) and \(k\). This proves \eqref{eq:binomial-envelopes}.

\textbf{Controlling distant and central thresholds.}
The bounds \(|d|\le1\) and \(\tau\ge1/2\) imply
\[
 |A_p^\pm(z)|
 \le\frac{\phi(z)}{6\tau}\bigl(3+|d|(1+z^2)\bigr)
 \le\frac{4+z^2}{3}\phi(z).
\]
The last expression tends to zero as \(|z|\to\infty\)
and is independent of \(p\), proving the asserted uniform
decay.

Finally, if \(k\) is an integer nearest to \(np\), then
\[
 |k-np|\le\frac12,\qquad
 |z_k|\le\frac1{2\sqrt{nv_I}}.
\]
Using \(A_p^+(z)=[3\phi(z)+d\psi(z)]/(6\tau)\), we have
\[
 \left|\frac{d}{dz}A_p^+(z)\right|
 \le\frac{3M_1+|d|M_2}{6\tau}
 \le M_1+\frac{M_2}{3}.
\]
The mean value theorem gives
\[
 |A_p^+(z_k)-A_p^+(0)|
 \le\frac{M_1+M_2/3}{2\sqrt{nv_I}}\longrightarrow0
\]
uniformly over \(p\in I\) and either choice of nearest
integer in the case of a tie. Since
\[
 A_p^+(0)=\frac{\phi(0)(3+d)}{6\tau}
 =\frac{\phi(0)(2-p)}{3[p^2+(1-p)^2]}=g(p),
\]
combining this estimate with \eqref{eq:binomial-envelopes}
proves \eqref{eq:binomial-central-branch}.
\end{proof}

% ===== clusters =====
\section{Stability near the Esseen two-point law}
\label{sec:clusters}

We work in affine coordinates in which the two atoms of the
Esseen law are \(0\) and \(1\).
We prove an exact Berry--Esseen bound for laws supported in
two sufficiently small intervals around these points.
The conditional distributions inside the intervals may be
arbitrary and may differ between the two intervals.

Related lattice-smoothing effects for symmetric Bernoulli
laws mixed with a small absolutely continuous component are
studied by Johnston~\cite[Section~1.3 and Lemma~5.1]{Johnston2023}.

\begin{proposition}\label{prop:clusters}
There exist $\eta\in(0,1/4)$ and $N<\infty$ such that, whenever a law $Q$
satisfies
\[
 \supp Q\subset[-\eta,\eta]\cup[1-\eta,1+\eta],\qquad
 \bigl|Q([1-\eta,1+\eta])-\pe\bigr|<\eta,
\]
one has $R_n(Q)\le\ce$ for every integer $n\ge N$.
The same assertion holds with $\pe$ replaced by $1-\pe$.
\end{proposition}

The proof separates small and macroscopic accumulated conditional
variance. In the first regime we compare directly with the Bernoulli law;
in the second, smoothing gives a strict first-order deficit.

To handle the small-variance regime, we first quantify the
loss caused by a centered perturbation near a binomial jump.
In Lemma~\ref{lem:small-cluster-variance}, \(W\) will be the
noise sum conditional on the number of summands in the
upper cluster, and \(f(u)\) will describe the change in the
normal distribution function, divided by the corresponding
binomial point probability.
The two inequalities below will provide the loss terms in
\eqref{eq:cluster-central-comparison}: a contribution
proportional to \(\EE|W|\), up to an error controlled by
\(\EE W^4\).

\begin{lemma}\label{lem:one-sided-loss}
There exist absolute constants $c_0,C,\kappa>0$ such that, for every centered
real random variable $W$ with $\EE W^4\le\kappa$ and every continuously
differentiable function $f$ on $[-3/5,3/5]$ satisfying
$f(0)=0$ and $3/4\le f'\le5/4$, one has, for $|u|\le1/2$,
\begin{align}
 \PP(W>u)+f(u)&\ge c_0\EE|W|-C\EE W^4,
 \label{eq:loss-upper}\\
 \PP(W<u)-f(u)&\ge c_0\EE|W|-C\EE W^4.
 \label{eq:loss-lower}
\end{align}
\end{lemma}

\begin{proof}
\textbf{First consider a centered variable with bounded support.}
Let \(a=1/10\), and suppose that \(\EE V=0\) and
\(|V|\le a\) almost surely. We first prove
\begin{equation}\label{eq:loss-bounded}
 \PP(V>u)+f(u)\ge\frac38\EE|V|,
 \qquad |u|\le\frac{11}{20}.
\end{equation}
The enlarged threshold interval will accommodate the mean
shift introduced by clipping below.
Put \(V_-:=(-V)_+\). Since
\(V=V_+-V_-\) and \(|V|=V_++V_-\), centering gives
\begin{equation}\label{eq:loss-centered-moment}
 \EE V_+=\EE V_-=\frac12\EE|V|.
\end{equation}
The assumptions on \(f\), together with \(f(0)=0\), imply
for \(0\le r\le11/20\) that
\begin{equation}\label{eq:loss-function-bounds}
 \begin{aligned}
 f(r)&=\int_0^r f'(s)\,ds\ge\frac34r,\\
 f(-r)&=-\int_{-r}^0 f'(s)\,ds\ge-\frac54r.
 \end{aligned}
\end{equation}
All integration intervals lie in the domain \([-3/5,3/5]\)
of \(f\).

For \(0\le u\le11/20\), let \(t=\PP(V>u)\).
The pointwise bound
\[
 V_+\le u+a\ind_{\{V>u\}}
\]
follows because \(V_+\le u\) on \(\{V\le u\}\) and
\(V_+\le a\) everywhere. Hence \(\EE V_+\le u+at\).
Using \eqref{eq:loss-function-bounds} and \(3a/4\le1\),
we obtain
\[
 \begin{aligned}
 \PP(V>u)+f(u)
 &\ge t+\frac34u
 \ge\frac34(u+at)
 \ge\frac34\EE V_+
 =\frac38\EE|V|,
 \end{aligned}
\]
where the last equality uses \eqref{eq:loss-centered-moment}.

Next let \(u=-r\), where \(0<r\le a\), and put
\(t=\PP(V>-r)\).
Splitting the expectation over \(\{V>-r\}\) and its
complement \(\{V\le-r\}\) gives
\[
 0=\EE V
 =\EE[V\ind_{\{V>-r\}}]+\EE[V\ind_{\{V\le-r\}}]
 \le at-r(1-t).
\]
It follows that
\begin{equation}\label{eq:loss-negative-threshold}
 r\le(a+r)t\le2at.
\end{equation}
Also, \(V_+\le a\ind_{\{V>0\}}\) and
\(\{V>0\}\subseteq\{V>-r\}\), so
\[
 \EE V_+\le a\PP(V>0)\le at\le t.
\]
Using \eqref{eq:loss-function-bounds},
\eqref{eq:loss-negative-threshold}, and \(a=1/10\), we find
\[
 \begin{aligned}
 \PP(V>-r)+f(-r)
 &\ge t-\frac54r
 \ge\left(1-\frac52a\right)t\\
 &=\frac34t
 \ge\frac34\EE V_+
 =\frac38\EE|V|.
 \end{aligned}
\]
The decomposition above includes any atom at \(-r\) in
\(\{V\le-r\}\), as required by the strict event in
\eqref{eq:loss-bounded}.

Finally, if \(u=-r\) with \(a<r\le11/20\), then
\(V\ge-a>-r\) almost surely, so \(\PP(V>-r)=1\).
By \eqref{eq:loss-function-bounds} and \(\EE|V|\le a\),
\[
 \begin{aligned}
 \PP(V>-r)+f(-r)
 &\ge1-\frac54r
 \ge1-\frac54\frac{11}{20}
 =\frac5{16}\\
 &\ge\frac38a
 \ge\frac38\EE|V|.
 \end{aligned}
\]
These three cases prove \eqref{eq:loss-bounded}.

\textbf{Clip and recenter the general variable.}
Let \(a_0=1/20\), and define
\[
 T=\max(-a_0,\min(W,a_0)),\qquad
 m=\EE T,\qquad V=T-m,\qquad M_4=\EE W^4.
\]
The clipping error satisfies the pointwise bound
\[
 |T-W|=(|W|-a_0)_+
 \le |W|\ind_{\{|W|>a_0\}}
 \le a_0^{-3}|W|^4.
\]
Since \(\EE W=0\), and
\(\ind_{\{|W|>a_0\}}\le a_0^{-4}|W|^4\), we obtain
\begin{equation}\label{eq:loss-clipping-errors}
 \begin{aligned}
 |m|&=|\EE(T-W)|
 \le\EE|T-W|\le a_0^{-3}M_4,\\
 \PP(T\ne W)&=\PP(|W|>a_0)\le a_0^{-4}M_4.
 \end{aligned}
\end{equation}
The triangle inequality also gives
\begin{equation}\label{eq:loss-clipped-first-moment}
 \begin{aligned}
 \EE|V|
 &=\EE|T-m|\ge\EE|T|-|m|\\
 &\ge\EE|W|-\EE|W-T|-|m|\\
 &\ge\EE|W|-2a_0^{-3}M_4.
 \end{aligned}
\end{equation}
Choose \(\kappa=a_0^4=1/160000\).
If \(M_4\le\kappa\), then \eqref{eq:loss-clipping-errors}
implies \(|m|\le a_0\). Thus \(V\) is centered and
\[
 |V|\le|T|+|m|\le2a_0=\frac1{10},\qquad
 |u-m|\le\frac12+a_0=\frac{11}{20}
 \quad\text{when }|u|\le\frac12.
\]
Therefore \eqref{eq:loss-bounded} applies to \(V\) at the
threshold \(u-m\).
Both \(u\) and \(u-m\), and the interval between them, lie
in \([-3/5,3/5]\), so the mean value theorem gives
\begin{equation}\label{eq:loss-function-shift}
 |f(u)-f(u-m)|\le\frac54|m|.
\end{equation}
Moreover,
\[
 \{V>u-m\}=\{T>u\}
 \subseteq\{W>u\}\cup\{T\ne W\}.
\]
Taking probabilities yields
\begin{equation}\label{eq:loss-probability-shift}
 \PP(W>u)\ge\PP(V>u-m)-\PP(T\ne W).
\end{equation}
Combining \eqref{eq:loss-probability-shift} and
\eqref{eq:loss-function-shift}, and then applying
\eqref{eq:loss-bounded}, gives
\[
 \begin{aligned}
 \PP(W>u)+f(u)
 &\ge\PP(V>u-m)+f(u-m)
       -\PP(T\ne W)-\frac54|m|\\
 &\ge\frac38\EE|V|-\PP(T\ne W)-\frac54|m|.
 \end{aligned}
\]
Using \eqref{eq:loss-clipping-errors} and
\eqref{eq:loss-clipped-first-moment}, we conclude that
\[
 \begin{aligned}
 \PP(W>u)+f(u)
 &\ge\frac38\EE|W|
 -\left(\frac34a_0^{-3}+a_0^{-4}
                    +\frac54a_0^{-3}\right)M_4\\
 &=\frac38\EE|W|-(a_0^{-4}+2a_0^{-3})M_4.
 \end{aligned}
\]
This proves \eqref{eq:loss-upper}, with the absolute constants
\[
 c_0=\frac38,\qquad
 C=a_0^{-4}+2a_0^{-3}=176000,\qquad
 \kappa=a_0^4=\frac1{160000}.
\]

\textbf{Obtain the second inequality by reflection.}
Define \(g(x)=-f(-x)\) on \([-3/5,3/5]\).
Then \(g\) is continuously differentiable and satisfies
\[
 g(0)=0,\qquad g'(x)=f'(-x),\qquad
 \frac34\le g'(x)\le\frac54.
\]
The variable \(-W\) is centered, and
\(\EE(-W)^4=\EE W^4\le\kappa\).
Since \(|-u|\le1/2\), the already proved inequality
\eqref{eq:loss-upper}, applied to \(-W\), \(g\), and \(-u\),
gives
\[
 \PP(-W>-u)+g(-u)
 \ge c_0\EE|-W|-C\EE(-W)^4.
\]
Finally,
\[
 \PP(-W>-u)=\PP(W<u),\qquad
 g(-u)=-f(u),\qquad \EE|-W|=\EE|W|.
\]
Substitution proves \eqref{eq:loss-lower}, with the same
constants and the stated strict inequality inside the probability.
\end{proof}

We now use Lemma~\ref{lem:one-sided-loss} to compare a
Bernoulli law with a perturbation \(Y=B+U\) satisfying
\(\EE(U\mid B)=0\).
The next lemma treats the case in which both \(|U|\) and
the accumulated conditional variance \(n\EE U^2\) are small.
It gives the exact comparison \(R_n(Y)\le R_n(B)\), with
strict inequality whenever \(\EE U^2>0\).
This comparison will complete the proof of
Proposition~\ref{prop:clusters}, after
Lemma~\ref{lem:accumulated-cluster-variance} has reduced
the argument to this regime.

\begin{lemma}\label{lem:small-cluster-variance}
There exist a compact interval $I\subset(0,1/2)$ containing $\pe$ in its
interior, constants $\epsilon_0,\lambda_0>0$, and $N_0<\infty$ with
the following property. Let $B\sim\mathrm{Bernoulli}(p)$, $p\in I$,
and suppose
\[
 Y=B+U,\qquad \EE(U\mid B)=0,\qquad |U|\le\epsilon_0.
\]
If $n\ge N_0$ and $n\EE U^2\le\lambda_0$, then
\[
 R_n(Y)\le R_n(B)<\ce.
\]
The first inequality is strict when $\EE U^2>0$.
\end{lemma}

\begin{proof}
\textbf{Compute the moments and fix the uniform parameter bounds.}
Fix a compact interval \(I_*\subset(0,1/2)\) containing
\(\pe\) in its interior, and put
\[
 \alpha=\min_{p\in I_*}\min(p,1-p)>0,\qquad
 v_*=\min_{p\in I_*}p(1-p)>0.
\]
We first work with \(p\in I_*\), and choose a smaller interval
\(I\) at the end. Initially require \(0<\epsilon_0<\alpha\)
and \(0<\lambda_0\le1\); further restrictions will be imposed below.
Constants denoted by \(c,C\) depend only on \(I_*\), and
those denoted by \(c_M,C_M\) may also depend on a fixed
central-window size \(M>0\).
They are independent of \(n,p\), the conditional laws of \(U\),
and the subsequently chosen \(\epsilon_0,\lambda_0\).
Write
\[
 q=1-p,\quad v=pq,\quad \tau=p^2+q^2,\quad
 s_b=\EE(U^2\mid B=b),\quad s=qs_0+ps_1=\EE U^2,\quad
 \lambda=ns.
\]
Let \((B_i,U_i)\), \(1\le i\le n\), be independent copies
of \((B,U)\), and set \(Y_i=B_i+U_i\).
Conditional centering gives
\[
 \EE U=\EE[\EE(U\mid B)]=0,\qquad
 \EE[(B-p)U]=\EE[(B-p)\EE(U\mid B)]=0.
\]
Consequently, \(\EE Y=p\) and \(\Var Y=v+s\ge v_*>0\).

Since \(|U|\le\epsilon_0<\min(p,q)\), the centered variable
\(Y-p\) is negative when \(B=0\) and positive when \(B=1\).
Thus
\[
 \rho:=\EE|Y-p|^3
 =q\,\EE[(p-U)^3\mid B=0]
  +p\,\EE[(q+U)^3\mid B=1].
\]
Expanding the cubes and using \(\EE(U\mid B)=0\) gives
\begin{equation}\label{eq:cluster-third-moment}
 \rho
 =v\tau+q\{3ps_0-\EE(U^3\mid B=0)\}
       +p\{3qs_1+\EE(U^3\mid B=1)\}.
\end{equation}
Here \(qp^3+pq^3=v\tau\). Moreover,
\[
 |\EE(U^3\mid B=b)|\le\epsilon_0s_b,\qquad b=0,1.
\]
Because \(3p-\epsilon_0>0\) and \(3q-\epsilon_0>0\),
\eqref{eq:cluster-third-moment} implies
\begin{equation}\label{eq:cluster-moment-bounds}
 0\le\rho-v\tau\le(3+\epsilon_0)s\le4s,\qquad
 \rho\ge v\tau\ge v_*/2.
\end{equation}
We used \(\tau\ge1/2\) and \(\epsilon_0<1\).
Also \(s\le\epsilon_0^2\le1\), so \(\rho\) and \(v+s\)
are bounded above uniformly.

\textbf{Control the noise conditional on the binomial count.}
Let \(K=\sum_{i=1}^nB_i\), \(b_k=\PP(K=k)\), and let
\(F_B\) denote the CDF of \(K\).
Conditional on \(K=k\), the sum \(\sum_iY_i\) has the law
of \(k+W_k\), where \(W_k\) is a sum of \(n-k\) independent
variables with law \(U\mid B=0\) and \(k\) independent
variables with law \(U\mid B=1\).
Indeed, conditional on the entire vector \((B_1,\ldots,B_n)\),
the \(U_i\) are independent with these conditional laws.
Every label vector containing \(k\) ones gives the same
convolution law for their sum, so conditioning only on \(K=k\)
gives that law as well.
In particular, \(\EE W_k=0\); independence of \(B\) and \(U\)
is not required.

Since \(p,q\ge\alpha\), one has
\(s_0+s_1\le s/\alpha\).
If \(\xi_1,\ldots,\xi_n\) are the independent centered
variables in the representation of \(W_k\), then
\[
 \begin{aligned}
 \EE W_k^4
 &=\sum_i\EE\xi_i^4
     +6\sum_{i<j}\EE\xi_i^2\,\EE\xi_j^2\\
 &\le\epsilon_0^2\sum_i\EE\xi_i^2
        +3\left(\sum_i\EE\xi_i^2\right)^2.
 \end{aligned}
\]
The equality follows by expanding the fourth power:
independence and centering make every term with an index
appearing only once vanish. The inequality uses
\(|\xi_i|\le\epsilon_0\).
It follows, uniformly for \(0\le k\le n\), that
\begin{equation}\label{eq:cluster-fourth-moment}
 \begin{aligned}
 t_k:=\EE W_k^2
 &=(n-k)s_0+ks_1\le\lambda/\alpha,\\
 \EE W_k^4
 &\le3t_k^2+\epsilon_0^2t_k
 \le C_4\lambda(\lambda+\epsilon_0^2)=:M_4^*,
 \end{aligned}
\end{equation}
where, for example, \(C_4=3\alpha^{-2}+\alpha^{-1}\).

Put \(z_k=(k-np)/\sqrt{nv}\).
For each fixed \(M>0\), all integers with \(|z_k|\le M\)
belong to \([0,n]\) once \(n\) is sufficiently large:
their distance from \(np\) is at most \(M\sqrt n/2\),
whereas \(np,nq\ge\alpha n\).
For such integers,
\begin{equation}\label{eq:cluster-central-variance}
 \begin{aligned}
 |t_k-\lambda|
 &=|k-np|\,|s_1-s_0|
 \le M\sqrt{nv}\,\frac{s}{\alpha}
 \le\frac{M}{2\alpha}\frac{\lambda}{\sqrt n}.
 \end{aligned}
\end{equation}
Thus, if \(\lambda>0\), increasing \(n\) by a bound depending
only on \(I_*,M\) ensures \(\lambda/2\le t_k\le3\lambda/2\).
Hölder's inequality, with exponents \(3/2\) and \(3\), gives
\[
 t_k=\EE\bigl[|W_k|^{2/3}|W_k|^{4/3}\bigr]
 \le(\EE|W_k|)^{2/3}(\EE W_k^4)^{1/3}.
\]
Since \(t_k>0\), rearrangement and
\eqref{eq:cluster-fourth-moment} yield
\begin{equation}\label{eq:cluster-absolute-loss}
 \EE|W_k|
 \ge\frac{t_k^{3/2}}{(\EE W_k^4)^{1/2}}
 \ge\frac{t_k}{\sqrt{3t_k+\epsilon_0^2}}
 \ge\frac{c_M\lambda}{\sqrt{C\lambda+\epsilon_0^2}}.
\end{equation}
The sample-size requirement in this argument does not depend
on how small the positive number \(\lambda\) is.

\textbf{Separate the contribution from the nearest binomial count.}
Let \(F\) be the CDF of \(\sum_iY_i\).
Every real threshold can be written as \(x=k+u\), with
\(k\in\mathbb Z\) and \(|u|\le1/2\).
For \(0\le k\le n\), conditioning on \(K\) gives
\begin{align}
 F(k+u)&=F_B(k-1)+b_k\PP(W_k\le u)+e^+_{k,u},
 \label{eq:cluster-right-block}\\
 F((k+u)-)&=F_B(k-1)+b_k\PP(W_k<u)+e^-_{k,u}.
 \label{eq:cluster-left-block}
\end{align}
More explicitly,
\[
 \begin{aligned}
 e^+_{k,u}
 &=\sum_{j<k}b_j\{\PP(W_j\le k+u-j)-1\}
   +\sum_{j>k}b_j\PP(W_j\le k+u-j),\\
 e^-_{k,u}
 &=\sum_{j<k}b_j\{\PP(W_j<k+u-j)-1\}
   +\sum_{j>k}b_j\PP(W_j<k+u-j),
 \end{aligned}
\]
where both sums are restricted to \(0\le j\le n\).
For \(j<k\), exclusion of that contribution requires a noise
of magnitude at least \(k-j-1/2\).
For \(j>k\), inclusion requires a noise of magnitude at least
\(j-k-1/2\).
Markov's inequality, \eqref{eq:binomial-mass-bound}, and
\eqref{eq:cluster-fourth-moment} therefore give
\begin{equation}\label{eq:cluster-leakage}
 \begin{aligned}
 |e^\pm_{k,u}|
 &\le\sum_{\substack{0\le j\le n\\j\ne k}}
       \frac{b_j\EE W_j^4}{(|j-k|-1/2)^4}\\
 &\le\frac{CM_4^*}{\sqrt n}.
 \end{aligned}
\end{equation}
The constant is uniform in \(k,n\), since
\[
 \sum_{\substack{j\in\mathbb Z\\j\ne k}}
       (|j-k|-1/2)^{-4}
 =2\sum_{m=1}^{\infty}(m-1/2)^{-4}
 \le32\sum_{m=1}^{\infty}m^{-4}<\infty.
\]
The use of weak inequalities in the tail bound includes
atoms at the threshold and either choice of \(k\) at a
half-integer.
The same argument applies when \(k<0\) or \(k>n\), with
no term involving \(W_k\). In these cases it gives
\begin{equation}\label{eq:cluster-outside-blocks}
 \begin{aligned}
 0\le F((k+u)-)\le F(k+u)
 &\le CM_4^*/\sqrt n &&(k<0),\\
 1-CM_4^*/\sqrt n\le F((k+u)-)\le F(k+u)
 &\le1 &&(k>n).
 \end{aligned}
\end{equation}

\textbf{Apply the one-sided loss estimate in a central interval.}
For every integer \(k\), define
\[
 G_k(u)=\Phi\left(\frac{k-np+u}{\sqrt{n(v+s)}}\right).
\]
For \(|z_k|\le M\) and sufficiently large \(n\), we have
\(0\le k\le n\) and \(b_k>0\), and may define
\[
 f_k(u)=\frac{G_k(u)-G_k(0)}{b_k}.
\]
This function is continuously differentiable and \(f_k(0)=0\).
To verify its derivative bounds, put
\[
 w_{k,u}=\sqrt{\frac v{v+s}}\,z_k
             +\frac{u}{\sqrt{n(v+s)}}.
\]
Direct differentiation gives
\begin{equation}\label{eq:cluster-derivative-ratio}
 f_k'(u)
 =\sqrt{\frac v{v+s}}\,
   \frac{\phi(w_{k,u})}{\phi(z_k)}\,
   \frac{\phi(z_k)}{\sqrt{nv}\,b_k}.
\end{equation}
For \(|z_k|\le M\), \(|u|\le3/5\), and
\(s\le\lambda_0/n\le1/n\),
\[
 |w_{k,u}-z_k|
 \le\frac{Ms}{2v_*}+\frac{3}{5\sqrt{nv_*}}
 \le C_M n^{-1/2}.
\]
Here we used \(0\le1-(1+r)^{-1/2}\le r/2\) for \(r\ge0\).
Thus the first factor in \eqref{eq:cluster-derivative-ratio}
tends uniformly to one. The second does also, because
\(\phi(w_{k,u})/\phi(z_k)=
\exp((z_k^2-w_{k,u}^2)/2)\) and \(|z_k|\le M\).
The third factor tends uniformly to one by
\eqref{eq:binomial-local} and \(\phi(z_k)\ge\phi(M)>0\).
Consequently, for all sufficiently large \(n\),
\begin{equation}\label{eq:cluster-derivative-bounds}
 \frac34\le f_k'(u)\le\frac54,
 \qquad |z_k|\le M,\quad |u|\le3/5,
\end{equation}
uniformly over the conditional noise laws.
The local limit estimate is used only to obtain these
inequalities; it introduces no additive error into the
comparison that follows.

Require also \(C_4\lambda_0(\lambda_0+\epsilon_0^2)\le\kappa\),
where \(\kappa\) is from Lemma~\ref{lem:one-sided-loss}.
Then \(\EE W_k^4\le\kappa\) by
\eqref{eq:cluster-fourth-moment}.
Rearranging \eqref{eq:cluster-right-block} and
\eqref{eq:cluster-left-block}, and using
\(F_B(k)=F_B(k-1)+b_k\), gives the exact identities
\begin{equation}\label{eq:cluster-signed-blocks}
 \begin{aligned}
 F(k+u)-G_k(u)
 &=F_B(k)-G_k(0)
   -b_k\{\PP(W_k>u)+f_k(u)\}+e^+_{k,u},\\
 G_k(u)-F((k+u)-)
 &=G_k(0)-F_B(k-1)
   -b_k\{\PP(W_k<u)-f_k(u)\}-e^-_{k,u}.
 \end{aligned}
\end{equation}
Apply \eqref{eq:loss-upper} and \eqref{eq:loss-lower} to
\(W_k,f_k\), and then use \eqref{eq:cluster-leakage}.
For \(|u|\le1/2\), the resulting bounds are
\[
 \begin{aligned}
 F(k+u)-G_k(u)
 &\le F_B(k)-G_k(0)-c_0b_k\EE|W_k|
       +Cb_k\EE W_k^4+CM_4^*/\sqrt n,\\
 G_k(u)-F((k+u)-)
 &\le G_k(0)-F_B(k-1)-c_0b_k\EE|W_k|
       +Cb_k\EE W_k^4+CM_4^*/\sqrt n.
 \end{aligned}
\]

To express these inequalities in the normalization of \(R_n\),
put
\[
 A=\frac{(v+s)^{3/2}}{\rho},\qquad
 A_B=\frac{\sqrt v}{\tau},\qquad a_n=\sqrt n A.
\]
By \eqref{eq:raw-normalization}, \(a_n=\sqrt n/\beta(Y)\).
The moment bounds above imply \(0<c\le A,A_B\le C\).
For every integer \(k\), define
\begin{equation}\label{eq:cluster-normalized-branches}
 \begin{aligned}
 D^+_{k,u}&=a_n[F(k+u)-G_k(u)],&
 D^-_{k,u}&=a_n[G_k(u)-F((k+u)-)],\\
 \widetilde R^+_{n,k}&=a_n[F_B(k)-G_k(0)],&
 \widetilde R^-_{n,k}&=a_n[G_k(0)-F_B(k-1)].
 \end{aligned}
\end{equation}
The central local limit estimate also gives
\(0<c_M\le\sqrt n\,b_k\le C_M\) for \(|z_k|\le M\).
Multiplying the two preceding bounds by \(a_n\), and using
\eqref{eq:cluster-fourth-moment} and
\eqref{eq:cluster-absolute-loss}, therefore yields, for
\(0<\lambda\le\lambda_0\),
\begin{equation}\label{eq:cluster-central-comparison}
 D^\pm_{k,u}
 \le\widetilde R^\pm_{n,k}
   -\frac{c_M\lambda}{\sqrt{C\lambda+\epsilon_0^2}}
   +C_M\lambda(\lambda+\epsilon_0^2),
 \qquad |z_k|\le M,\quad |u|\le1/2.
\end{equation}
In particular, multiplying the leakage error by \(a_n\)
gives \(O(M_4^*)\), of the same order as the fourth-moment
error from Lemma~\ref{lem:one-sided-loss}.

\textbf{Control the change in variance and normalization.}
We compare \(\widetilde R^\pm_{n,k}\) with the binomial
quantities \(R^\pm_{n,k}(B_p)\) defined in
Lemma~\ref{lem:binomial-estimates}.
First, the identity
\[
 A-A_B
 =\frac{(v+s)^{3/2}-v^{3/2}}{\rho}
   +\frac{v^{3/2}(v\tau-\rho)}{\rho\,v\tau}
\]
and \eqref{eq:cluster-moment-bounds} show that
\begin{equation}\label{eq:cluster-prefactor-control}
 |A-A_B|\le Cs.
\end{equation}
Indeed, the denominators are bounded away from zero,
\(|\rho-v\tau|\le4s\), and the mean value theorem bounds
\(|(v+s)^{3/2}-v^{3/2}|\) by \(Cs\).

The Gaussian variance change is uniform over all thresholds.
For \(r\ge0\), putting \(w=z/\sqrt{1+r}\), we have
\[
 \left|\frac{\partial}{\partial r}
       \Phi\left(\frac z{\sqrt{1+r}}\right)\right|
 =\frac{|w|\phi(w)}{2(1+r)}
 \le\frac12\sup_{w\in\R}|w|\phi(w)<\infty.
\]
Integrating from \(0\) to \(s/v\), and using \(v\ge v_*\),
gives
\begin{equation}\label{eq:cluster-gaussian-variance}
 \sup_{z\in\R}
 \left|\Phi\left(\frac z{\sqrt{1+s/v}}\right)-\Phi(z)\right|
 \le Cs.
\end{equation}
By \eqref{eq:structural}, affine invariance, and
\eqref{eq:binomial-standardized-moments},
\[
 \Delta_n(B_p)
 \le\frac{0.4690\,\beta(B_p)}{\sqrt n}
 \le\frac C{\sqrt n},\qquad p\in I_*.
\]
Both \(|F_B(k)-\Phi(z_k)|\) and
\(|\Phi(z_k)-F_B(k-1)|\) are bounded by \(\Delta_n(B_p)\).
For the second quantity, take a limit of thresholds increasing
to \(k\), using \(F_B(k-)=F_B(k-1)\) and continuity of
the Gaussian CDF.

For example, the upper branch satisfies
\[
 \begin{aligned}
 \widetilde R^+_{n,k}-R^+_{n,k}(B_p)
 &=\sqrt n(A-A_B)[F_B(k)-\Phi(z_k)]\\
 &\quad+\sqrt n A[\Phi(z_k)-G_k(0)].
 \end{aligned}
\]
For the lower branch the corresponding identity has
\(\Phi(z_k)-F_B(k-1)\) in the first bracket and
\(G_k(0)-\Phi(z_k)\) in the second.
Since \(G_k(0)=\Phi(z_k/\sqrt{1+s/v})\),
\eqref{eq:cluster-prefactor-control} and
\eqref{eq:cluster-gaussian-variance} imply
\begin{equation}\label{eq:cluster-normalization-error}
 \bigl|\widetilde R^\pm_{n,k}-R^\pm_{n,k}(B_p)\bigr|
 \le Cs+C\sqrt n\,s
 \le C\lambda/\sqrt n,\qquad k\in\mathbb Z.
\end{equation}
Thus this comparison applies at every integer \(k\), not
only in the central interval.

For \(\lambda>0\), the ratios of the two positive errors
in \eqref{eq:cluster-central-comparison} and
\eqref{eq:cluster-normalization-error} to the loss term
are bounded by
\begin{equation}\label{eq:cluster-error-ratios}
 C_M(\lambda+\epsilon_0^2)^{3/2},
 \qquad
 \frac{C_M\sqrt{C\lambda+\epsilon_0^2}}{\sqrt n},
\end{equation}
respectively. For the first ratio, cancel \(\lambda\) and
use \(\sqrt{C\lambda+\epsilon_0^2}
\le C\sqrt{\lambda+\epsilon_0^2}\); for the second, cancel
\(\lambda\) directly.
For each fixed \(M\), reducing \(\lambda_0,\epsilon_0\)
and then increasing \(N_0\) makes each ratio at most \(1/4\),
uniformly for \(0<\lambda\le\lambda_0\) and \(n\ge N_0\).
Since \(R^\pm_{n,k}(B_p)\le R_n(B_p)\), at least half
the loss remains. Decreasing \(c_M\) to include this factor,
we obtain
\begin{equation}\label{eq:cluster-strict-central}
 D^\pm_{k,u}
 \le R_n(B_p)-\frac{c_M\lambda}
                         {\sqrt{C\lambda+\epsilon_0^2}},
 \qquad |z_k|\le M,\quad |u|\le1/2.
\end{equation}
The estimate \(R^-_{n,k}(B_p)\le R_n(B_p)\) uses the same
left-limit argument as above.
The restrictions just described will be imposed after
\(M\) is fixed below. No positive lower bound on
\(\lambda\) has been used.

\textbf{Bound the discrepancies outside the central interval.}
For \(0\le k\le n\), the block identities and
\eqref{eq:cluster-leakage} give
\[
 F_B(k-1)-\frac{CM_4^*}{\sqrt n}
 \le F((k+u)-)\le F(k+u)
 \le F_B(k)+\frac{CM_4^*}{\sqrt n}.
\]
Also, for \(|u|\le1/2\), the bound \(\phi\le\phi(0)\)
and \(\rho\ge v\tau\) imply
\[
 \begin{aligned}
 a_n|G_k(u)-G_k(0)|
 &\le\frac{\phi(0)(v+s)}{2\rho}\\
 &\le\frac{\phi(0)}{2\tau}
       +\frac{\phi(0)s}{2v\tau}
 \le\frac{\phi(0)}{2\tau}+Cs.
 \end{aligned}
\]
Consequently,
\[
 D^\pm_{k,u}
 \le\widetilde R^\pm_{n,k}
       +\frac{\phi(0)}{2\tau}+Cs+CM_4^*.
\]
Apply \eqref{eq:cluster-normalization-error} and the uniform
expansion \eqref{eq:binomial-envelopes}, both on the original
interval \(I_*\). This yields
\begin{equation}\label{eq:cluster-far-comparison}
 \begin{aligned}
 \max\{D^+_{k,u},D^-_{k,u}\}
 &\le\max\{A_p^+(z_k),A_p^-(z_k)\}
       +\frac{\phi(0)}{2\tau}\\
 &\quad+C\lambda(\lambda+\epsilon_0^2)
       +C\lambda/\sqrt n+\omega_n,
 \end{aligned}
\end{equation}
where
\[
 \omega_n:=
 \sup_{\substack{p\in I_*\\k\in\mathbb Z}}
 \max_{\theta\in\{+,-\}}
 |R^\theta_{n,k}(B_p)-A_p^\theta(z_k)|
 \longrightarrow0.
\]
In particular, \(\omega_n\) does not depend on the noise laws.
The term \(Cs\) has been included in \(C\lambda/\sqrt n\),
since \(s=\lambda/n\).

Because \(\tau\ge1/2\), the second term in
\eqref{eq:cluster-far-comparison} is at most \(\phi(0)\).
By \eqref{eq:ce}, \(\ce>\phi(0)\), so fix
\[
 \gamma=\frac{\ce-\phi(0)}8>0.
\]
Lemma~\ref{lem:binomial-estimates} gives uniform decay of
\(A_p^\pm(z)\) as \(|z|\to\infty\). Choose \(M>0\) so that
\[
 \sup_{p\in I_*}\sup_{|z|>M}
 \max\{|A_p^+(z)|,|A_p^-(z)|\}\le\gamma.
\]
With this \(M\) fixed, reduce \(\lambda_0,\epsilon_0\)
to satisfy the central-interval requirements above and
to make the fourth-moment term in
\eqref{eq:cluster-far-comparison} at most \(\gamma\).
Increase \(N_0\) so that \(C\lambda_0/\sqrt n+\omega_n\le\gamma\)
for every \(n\ge N_0\).
Then, for \(0\le k\le n\) and \(|z_k|>M\),
\[
 \max\{D^+_{k,u},D^-_{k,u}\}
 \le\phi(0)+3\gamma=\ce-5\gamma\le\ce-2\gamma.
\]

For completeness, consider \(k<0\) or \(k>n\).
If \(k<0\), then \(k+u\le-1/2\), so the argument of \(G_k(u)\)
is at most \(-c\sqrt n\).
If \(k>n\), then \(k+u\ge n+1/2\), so that argument is
at least \(c\sqrt n\).
Indeed, \(np,nq\ge\alpha n\) and \(v+s\le1/4+1\).
The Gaussian tail bound
\[
 1-\Phi(x)=\int_x^\infty\phi(t)\,dt
 \le\frac1x\int_x^\infty t\phi(t)\,dt
 =\frac{\phi(x)}x,\qquad x>0,
\]
therefore implies \(G_k(u)\le Ce^{-cn}\) when \(k<0\),
and \(1-G_k(u)\le Ce^{-cn}\) when \(k>n\).
Together with \eqref{eq:cluster-outside-blocks} and
\(a_n\le C\sqrt n\), this gives
\[
 \begin{array}{lll}
 D^+_{k,u}\le CM_4^*,&
 D^-_{k,u}\le C\sqrt n e^{-cn},& k<0,\\
 D^+_{k,u}\le C\sqrt n e^{-cn},&
 D^-_{k,u}\le CM_4^*,& k>n.
 \end{array}
\]
Reducing \(\lambda_0,\epsilon_0\) if necessary makes
\(CM_4^*\le\gamma\), and increasing \(N_0\) makes the
Gaussian tail term at most \(\gamma\).
These branches are therefore also at most \(\ce-2\gamma\).
Such reductions preserve every earlier requirement.

\textbf{Compare with the binomial maximum and take the supremum.}
At an integer \(k_n\) nearest to \(np\),
\eqref{eq:binomial-central-branch} gives, uniformly on \(I_*\),
\[
 R_n(B_p)\ge R^+_{n,k_n}(B_p)=g(p)+o(1),\qquad
 g(p)=\frac{\phi(0)(2-p)}{3[p^2+(1-p)^2]}.
\]
The function \(g\) is continuous and \(g(\pe)=\ce\).
Indeed, \eqref{eq:esseen-parameters} gives
\[
 2-\pe=\frac{\sqrt{10}}2,\qquad
 \pe^2+\qe^2=10-3\sqrt{10},
\]
so
\[
 g(\pe)
 =\frac{\phi(0)\sqrt{10}}{6(10-3\sqrt{10})}
 =\frac{\phi(0)(3+\sqrt{10})}{6}=\ce.
\]
Choose a compact interval \(I\subset I_*\), with \(\pe\)
in its interior, so that \(g(p)\ge\ce-\gamma/2\) on \(I\).
Increase \(N_0\) so that the absolute value of the uniform
remainder in \eqref{eq:binomial-central-branch} is at most
\(\gamma/2\). Then
\begin{equation}\label{eq:cluster-binomial-lower}
 R_n(B_p)\ge\ce-\gamma,\qquad p\in I,\quad n\ge N_0.
\end{equation}
Every branch treated outside the central interval is at
most \(\ce-2\gamma\), hence at most \(R_n(B_p)-\gamma\).
Together with \eqref{eq:cluster-strict-central}, this bounds
both \(D^+_{k,u}\) and \(D^-_{k,u}\), for every real threshold
\(x=k+u\), by \(R_n(B_p)\) minus the minimum of the two gaps.
The ordinary lower discrepancy is also bounded, because
\[
 a_n[G_k(u)-F(k+u)]
 \le a_n[G_k(u)-F((k+u)-)]=D^-_{k,u}.
\]
Taking the supremum over \(x\), and using
\eqref{eq:raw-normalization}, therefore gives
\begin{equation}\label{eq:cluster-uniform-gap}
 R_n(Y)\le R_n(B_p)-
 \min\left\{
   \frac{c_M\lambda}{\sqrt{C\lambda+\epsilon_0^2}},\,
   \gamma
 \right\},\qquad 0<\lambda\le\lambda_0.
\end{equation}
For each positive \(\lambda\), this minimum is positive
and independent of the threshold, so taking the supremum
preserves strict inequality.

All choices can be made in the following order:
fix \(I_*\) and \(\gamma\); choose \(M\); choose positive
\(\lambda_0,\epsilon_0\) small enough for the fourth-moment,
central-error, and far-threshold bounds; choose \(I\);
and finally choose \(N_0\) large enough for all the
sample-size requirements.
These requirements depend only on the fixed parameters,
not on the particular conditional laws or a positive
lower bound for \(\lambda\).
If \(\lambda=0\), then \(\EE U^2=0\), so \(U=0\)
almost surely and \(R_n(Y)=R_n(B_p)\).
Finally, Schulz's two-point bound
\eqref{eq:schulz-two-point} gives \(R_n(B_p)<\ce\) for
every \(n\ge1\). This proves all assertions of the lemma.
\end{proof}

We next prepare for the case in which the accumulated noise
variance \(n\EE U^2\) is bounded below by a fixed positive
constant.
The following lemma gives a uniform lower bound
\(c_*/\sqrt n\) for the probability of an interval
\((x+a,x+b)\), where \(a<b\) are fixed and
\(|x-np|\le R\sqrt n\).
Lemma~\ref{lem:accumulated-cluster-variance} will use this
bound to strengthen the inequalities comparing the
distribution function of the original sum with that of
the sum after adding an independent uniform variable.
This will yield the strict asymptotic bound
\eqref{eq:macroscopic-gap}.

\begin{lemma}[Uniform local mass]\label{lem:two-cluster-local-mass}
Fix a compact interval $I\subset(0,1)$, $\lambda_*>0$, $R\ge0$,
and real numbers $a<b$. There exist $\epsilon_*>0$, $c_*>0$, and
$N_*<\infty$ with the following property. Let
\[
 B\sim\mathrm{Bernoulli}(p),\qquad p\in I,\qquad
 \EE(U\mid B)=0,\qquad |U|\le\epsilon_*,\qquad Y=B+U.
\]
For every $n\ge N_*$ satisfying $L:=n\EE U^2\ge\lambda_*$,
and independent copies $Y_1,\ldots,Y_n$ of $Y$, one has
\[
 \inf_{|x-np|\le R\sqrt n}
 \PP\left\{\sum_{i=1}^nY_i\in(x+a,x+b)\right\}
 \ge\frac{c_*}{\sqrt n}.
\]
\end{lemma}

\begin{proof}
\textbf{Represent the sum conditional on the binomial count.}
Put
\[
 \alpha=\min_{p\in I}\min(p,1-p)>0,\qquad
 v_I=\min_{p\in I}p(1-p)>0,
\]
and initially require \(0<\epsilon_*\le1\).
Let \((B_i,U_i)\), \(1\le i\le n\), be independent copies
of \((B,U)\), and write
\[
 S_n=\sum_{i=1}^n(B_i+U_i),\qquad
 K=\sum_{i=1}^nB_i,\qquad
 s_b=\EE(U^2\mid B=b),\qquad
 s=(1-p)s_0+ps_1=L/n.
\]
Since \(\EE(U\mid B)=0\), we have
\(\EE U=0\) and \(\EE S_n=np\).
Also, \(p,1-p\ge\alpha\) imply
\begin{equation}\label{eq:local-mass-conditional-variances}
 s\ge\alpha(s_0+s_1),\qquad
 |s_1-s_0|\le s_0+s_1\le s/\alpha.
\end{equation}

For \(0\le k\le n\), conditional on \(K=k\), the variable
\(S_n\) has the law of \(k+W_k\), with the same representation
as in Lemma~\ref{lem:small-cluster-variance}.
Namely, \(W_k\) is a sum of \(n-k\) independent variables
with law \(U\mid B=0\) and \(k\) independent variables
with law \(U\mid B=1\).
To see this, first condition on the full vector of labels
\((B_1,\ldots,B_n)\). The conditional noise variables are
then independent, and their sum has the same convolution
law for every label vector containing \(k\) ones.
Conditioning only on \(K=k\) therefore gives that law.
Each conditional summand is centered and bounded in
absolute value by \(\epsilon_*\). Thus
\begin{equation}\label{eq:local-mass-noise-variance}
 \EE W_k=0,\qquad
 t_k:=\EE W_k^2=(n-k)s_0+ks_1
      =L+(k-np)(s_1-s_0).
\end{equation}
The last identity follows by substituting
\(L=n[(1-p)s_0+ps_1]\).

\textbf{Select counts with uniformly positive binomial masses.}
Fix \(x\) with \(|x-np|\le R\sqrt n\), and define
\begin{equation}\label{eq:local-mass-selected-counts}
 r=\max(1,\sqrt L),\qquad
 \mathcal K_x=\{k\in\mathbb Z:|k-x|\le r\}.
\end{equation}
The assumption \(|U|\le\epsilon_*\le1\) gives
\[
 L=n\EE U^2\le n\epsilon_*^2\le n.
\]
For \(n\ge1\), it follows that \(r\le1+\sqrt n\le2\sqrt n\).
Consequently, every \(k\in\mathcal K_x\) satisfies
\begin{equation}\label{eq:local-mass-count-location}
 |k-np|\le|k-x|+|x-np|\le(R+2)\sqrt n.
\end{equation}
Since \(np,n(1-p)\ge\alpha n\), all these integers lie
in \([0,n]\) once \((R+2)\sqrt n\le\alpha n\).
Choose \(N_*\) large enough that this holds for all
\(n\ge N_*\).

We may now apply \eqref{eq:local-mass-noise-variance} to
every selected count. Using
\eqref{eq:local-mass-conditional-variances},
\eqref{eq:local-mass-count-location}, and \(s=L/n\), we obtain
\[
 \left|\frac{t_k}{L}-1\right|
 =\frac{|k-np|\,|s_1-s_0|}{L}
 \le\frac{(R+2)\sqrt n}{L}\frac{L}{\alpha n}
 =\frac{R+2}{\alpha\sqrt n}.
\]
Increase \(N_*\) so that the last bound is at most \(1/2\).
Then
\begin{equation}\label{eq:local-mass-variance-comparison}
 \frac L2\le t_k\le\frac{3L}{2},
 \qquad k\in\mathcal K_x,\quad n\ge N_*.
\end{equation}
In particular, \(t_k>0\), because \(L\ge\lambda_*>0\).
This comparison is uniform without any fixed upper bound on \(L\).

To bound the binomial masses, put \(v=p(1-p)\) and
\[
 z_k=\frac{k-np}{\sqrt{nv}},\qquad
 M=\frac{R+2}{\sqrt{v_I}}.
\]
Equation~\eqref{eq:local-mass-count-location} implies
\(|z_k|\le M\) for all selected counts.
The local limit estimate \eqref{eq:binomial-local}, together
with \(\phi(z_k)\ge\phi(M)>0\), gives
\[
 \frac{\sqrt{nv}\,\PP(K=k)}{\phi(z_k)}
 \longrightarrow1
\]
uniformly over \(p\in I\) and these integers.
Indeed, multiplication of the uniform \(o(n^{-1/2})\)
remainder in \eqref{eq:binomial-local} by
\(\sqrt{nv}/\phi(z_k)\) gives \(o(1)\), since
\(\sqrt v\le1/2\) and \(\phi(z_k)\ge\phi(M)\).
Increasing \(N_*\) makes the displayed ratio at least \(1/2\).
Hence
\begin{equation}\label{eq:local-mass-binomial-lower}
 \PP(K=k)
 \ge\frac{\phi(z_k)}{2\sqrt{nv}}
 \ge\frac{\phi(M)}{\sqrt n}
 =:\frac{c_R}{\sqrt n},
 \qquad k\in\mathcal K_x.
\end{equation}
Here \(c_R=\phi(M)>0\) depends only on \(I,R\).

\textbf{Bound the conditional interval probability by normal approximation.}
Let \(\xi_{k,1},\ldots,\xi_{k,n}\) be the independent
centered summands representing \(W_k\).
Since \(|\xi_{k,i}|\le\epsilon_*\),
\[
 \sum_{i=1}^n\EE|\xi_{k,i}|^3
 \le\epsilon_*\sum_{i=1}^n\EE\xi_{k,i}^2
 =\epsilon_*t_k.
\]
Define \(H_k(z)=\PP(W_k/\sqrt{t_k}\le z)\).
The independent-summand Berry--Esseen inequality
\eqref{eq:independent-berry-esseen} therefore gives
\begin{equation}\label{eq:local-mass-conditional-be}
 \sup_{z\in\R}|H_k(z)-\Phi(z)|
 \le C_{\mathrm{ind}}\frac{\epsilon_*t_k}{t_k^{3/2}}
 =\frac{C_{\mathrm{ind}}\epsilon_*}{\sqrt{t_k}}
 =:\delta_k,
\end{equation}
where \(C_{\mathrm{ind}}=1\).
All its hypotheses hold: the summands are independent and
centered, have finite third absolute moments, and have
total variance \(t_k>0\).
For any \(z\), letting \(z_m\uparrow z\) in
\eqref{eq:local-mass-conditional-be} and using continuity
of \(\Phi\) also gives
\begin{equation}\label{eq:local-mass-left-limit-be}
 |H_k(z-)-\Phi(z)|\le\delta_k.
\end{equation}

The standardized endpoints of the conditional interval are
\[
 \ell_k=\frac{x+a-k}{\sqrt{t_k}},\qquad
 u_k=\frac{x+b-k}{\sqrt{t_k}}.
\]
For \(c\in\{a,b\}\), \eqref{eq:local-mass-selected-counts}
and \eqref{eq:local-mass-variance-comparison} imply
\[
 \begin{aligned}
 \frac{|x+c-k|}{\sqrt{t_k}}
 &\le\frac{r+|c|}{\sqrt{L/2}}\\
 &=\sqrt2\left[\max(1,L^{-1/2})+|c|L^{-1/2}\right]\\
 &\le\sqrt2\left[
       \max(1,\lambda_*^{-1/2})
       +\max(|a|,|b|)\lambda_*^{-1/2}
       \right]=:A.
 \end{aligned}
\]
Thus \(\ell_k,u_k\in[-A,A]\), where \(A\) depends only
on \(\lambda_*,a,b\).
Since \(\phi(z)\ge\phi(A)\) on this interval and
\(u_k-\ell_k=(b-a)/\sqrt{t_k}\), we have
\begin{equation}\label{eq:local-mass-normal-interval}
 \begin{aligned}
 \Phi(u_k)-\Phi(\ell_k)
 &=\int_{\ell_k}^{u_k}\phi(z)\,dz\\
 &\ge\frac{(b-a)\phi(A)}{\sqrt{t_k}}
 \ge\frac{c_1}{\sqrt L},\\
 c_1&:=(b-a)\sqrt{2/3}\,\phi(A)>0.
 \end{aligned}
\end{equation}
The last inequality uses \(t_k\le3L/2\).

Because the target interval is open at both endpoints,
its exact conditional probability is
\[
 \PP\{W_k\in(x+a-k,x+b-k)\}
 =H_k(u_k-)-H_k(\ell_k).
\]
Use \eqref{eq:local-mass-left-limit-be} at \(u_k\) and
\eqref{eq:local-mass-conditional-be} at \(\ell_k\).
Together with \eqref{eq:local-mass-normal-interval} and
\(t_k\ge L/2\), these estimates yield
\begin{equation}\label{eq:local-mass-conditional-lower}
 \begin{aligned}
 \PP\{W_k\in(x+a-k,x+b-k)\}
 &\ge\Phi(u_k)-\Phi(\ell_k)-2\delta_k\\
 &\ge\frac{c_1}{\sqrt L}
       -\frac{2C_{\mathrm{ind}}\epsilon_*}{\sqrt{t_k}}\\
 &\ge\frac{c_1-2\sqrt2 C_{\mathrm{ind}}\epsilon_*}{\sqrt L}.
 \end{aligned}
\end{equation}
This calculation allows atoms at either endpoint.
The constant \(c_1\) has already been fixed independently
of \(\epsilon_*\). We may therefore choose
\[
 0<\epsilon_*
 \le\min\left\{1,\frac{c_1}{4\sqrt2 C_{\mathrm{ind}}}\right\}.
\]
Then \eqref{eq:local-mass-conditional-lower} gives
\begin{equation}\label{eq:local-mass-positive-conditional}
 \PP\{W_k\in(x+a-k,x+b-k)\}
 \ge\frac{c_1}{2\sqrt L},
 \qquad k\in\mathcal K_x.
\end{equation}
This choice preserves all preceding estimates, which used
only \(\epsilon_*\le1\).

\textbf{Sum over the selected counts.}
The interval \([x-r,x+r]\) contains
\(\lfloor x+r\rfloor-\lceil x-r\rceil+1\) integers.
Using \(\lfloor y\rfloor\ge y-1\),
\(\lceil y\rceil\le y+1\), and \(r\ge1\), we obtain
\begin{equation}\label{eq:local-mass-number-of-counts}
 \#\mathcal K_x
 \ge2r-1\ge r\ge\sqrt L.
\end{equation}
Conditioning on \(K\), retaining only \(k\in\mathcal K_x\),
and using \eqref{eq:local-mass-binomial-lower} and
\eqref{eq:local-mass-positive-conditional}, we conclude that
\[
 \begin{aligned}
 \PP\{S_n\in(x+a,x+b)\}
 &=\sum_{k=0}^n\PP(K=k)\,
       \PP\{W_k\in(x+a-k,x+b-k)\}\\
 &\ge\sum_{k\in\mathcal K_x}
       \frac{c_R}{\sqrt n}\frac{c_1}{2\sqrt L}\\
 &\ge\frac{c_Rc_1}{2\sqrt n}.
 \end{aligned}
\]
The last inequality uses \eqref{eq:local-mass-number-of-counts}.
Set \(c_*=c_Rc_1/2>0\).
Every estimate is uniform over \(|x-np|\le R\sqrt n\),
so taking the infimum over these \(x\) proves the assertion.
The constants \(\epsilon_*,c_*,N_*\) depend only on
\(I,\lambda_*,R,a,b\), not on the particular noise law or
the value of \(L\ge\lambda_*\).
\end{proof}

We now combine the local probability bound from
Lemma~\ref{lem:two-cluster-local-mass} with the smoothing
expansion to obtain a strict asymptotic bound when the
accumulated noise variance stays bounded away from zero.
At either Esseen parameter, this will show that
\(R_{n_j}(Y_j)\to\ce\) requires \(n_j\EE U_j^2\to0\).
In the proof of Proposition~\ref{prop:clusters}, this
conclusion will allow us to apply
Lemma~\ref{lem:small-cluster-variance}.

\begin{lemma}\label{lem:accumulated-cluster-variance}
Let $p_j\to p\in(0,1)$, $n_j\to\infty$, and
\[
 Y_j=B_j+U_j,\qquad B_j\sim\mathrm{Bernoulli}(p_j),\qquad
 \EE(U_j\mid B_j)=0,\qquad |U_j|\le\epsilon_j\to0.
\]
Write $s_j=\EE U_j^2$, $v=p(1-p)$ and $d=1-2p$.
If $\liminf_j n_js_j>0$, then for some $\delta>0$,
\begin{equation}\label{eq:macroscopic-gap}
 \limsup_j\sqrt{n_j}\,\Delta_{n_j}(Y_j)
 \le\frac{\phi(0)(3+|d|)}{6\sqrt v}-\delta.
\end{equation}
Consequently, if $p\in\{\pe,1-\pe\}$ and
$R_{n_j}(Y_j)\to\ce$, then $n_js_j\to0$.
\end{lemma}

\begin{proof}
\textbf{Obtain a uniform lower bound for local interval probabilities.}
First assume \(\liminf_j n_js_j>0\).
Choose \(\lambda_*>0\) so that \(n_js_j\ge\lambda_*\)
for all sufficiently large \(j\), and fix a compact interval
\(I\subset(0,1)\) containing all sufficiently late \(p_j\).
Let \(Y_{j,1},\ldots,Y_{j,n_j}\) be independent copies of \(Y_j\).

Fix a raw radius \(T\ge0\) and real numbers \(a<b\).
Lemma~\ref{lem:two-cluster-local-mass} supplies constants
\(\epsilon_*,c_*,N_*>0\) for these fixed parameters,
\(I\), and \(\lambda_*\).
Its assumptions hold for all sufficiently large \(j\):
\[
 p_j\in I,\qquad n_j\ge N_*,\qquad
 |U_j|\le\epsilon_j\le\epsilon_*,
 \qquad n_j\EE U_j^2=n_js_j\ge\lambda_*.
\]
Here we used \(p_j\to p\), \(n_j\to\infty\), and
\(\epsilon_j\to0\).
Consequently,
\begin{equation}\label{eq:cluster-local-mass}
 \inf_{|x-n_jp_j|\le T\sqrt{n_j}}
 \PP\left\{\sum_{i=1}^{n_j}Y_{j,i}\in(x+a,x+b)\right\}
 \ge\frac{c_*}{\sqrt{n_j}}
\end{equation}
eventually.
The constant may depend on \(I,\lambda_*,T,a,b\), but not
on \(j\) or \(x\). In particular, no upper bound on
\(n_js_j\) is required.

\textbf{Apply the smoothing expansion to the unstandardized sums.}
Write
\[
 S_j=\sum_{i=1}^{n_j}Y_{j,i},\qquad F_j(x)=\PP(S_j\le x),
\]
and let \(J_j\) be the CDF of \(S_j+H\), where
\(H\sim\mathrm{Unif}[-1/2,1/2]\) is independent of the row.
Set
\[
 v_j=p_j(1-p_j),\qquad
 \sigma_j^2=v_j+s_j,\qquad
 z=z_j(x)=\frac{x-n_jp_j}{\sigma_j\sqrt{n_j}},\qquad
 \gamma_j=\frac{\EE(Y_j-p_j)^3}{\sigma_j^3}.
\]
We use these quantities for sufficiently large \(j\), when
\(p_j\in I\) and hence \(\sigma_j^2\ge v_j>0\).

Let \(Q_j\) be the law of \(Y_j\), and let
\(Q=\mathrm{Bernoulli}(p)\).
We verify the hypotheses of
Proposition~\ref{prop:jitter-width}\textup{(ii)}.
The given coupling \(Y_j=B_j+U_j\) implies
\[
 W_3\bigl(Q_j,\mathrm{Bernoulli}(p_j)\bigr)
 \le(\EE|U_j|^3)^{1/3}\le\epsilon_j.
\]
For a uniform variable \(V\) on \([0,1]\), the indicators
\(\ind_{\{V\le p_j\}}\) and \(\ind_{\{V\le p\}}\) have
Bernoulli parameters \(p_j\) and \(p\), and
\[
 \EE\left|\ind_{\{V\le p_j\}}-\ind_{\{V\le p\}}\right|^3
 =|p_j-p|.
\]
The triangle inequality for \(W_3\) therefore gives
\begin{equation}\label{eq:cluster-raw-wasserstein}
 W_3(Q_j,Q)\le\epsilon_j+|p_j-p|^{1/3}\longrightarrow0.
\end{equation}

Conditional centering gives
\[
 \EE U_j=\EE[\EE(U_j\mid B_j)]=0,\qquad
 \EE(B_jU_j)=\EE[B_j\EE(U_j\mid B_j)]=0.
\]
Thus
\[
 \EE Y_j=p_j,\qquad
 \operatorname{Var}(Y_j)
 =\operatorname{Var}(B_j)+\EE U_j^2
     +2\EE[(B_j-p_j)U_j]
 =v_j+s_j=\sigma_j^2.
\]
Moreover,
\begin{equation}\label{eq:cluster-variance-limit}
 0\le s_j\le\epsilon_j^2\longrightarrow0,\qquad
 v_j\longrightarrow v>0,\qquad
 \sigma_j\longrightarrow\sqrt v.
\end{equation}
Each \(Q_j\) has finite third absolute moment, since
\(|Y_j|\le1+\epsilon_j\).
The limiting law \(Q\) has positive variance \(v\) and
maximal lattice span one, because its two atoms are \(0\)
and \(1\), both with positive mass.
We therefore apply \eqref{eq:jitter-raw-expansion} with
the raw jitter widths equal to one.
Evaluating it at \(z=z_j(x)\) yields
\begin{equation}\label{eq:cluster-jitter-expansion}
 J_j(x)=\Phi(z)+\frac{\gamma_j}{6\sqrt{n_j}}
                    (1-z^2)\phi(z)+o(n_j^{-1/2}),
\end{equation}
uniformly over \(x\in\R\).
The uniformity follows because \(x\mapsto z_j(x)\) is a
bijection of \(\R\).

\textbf{Use the local mass to strengthen the CDF comparison.}
The CDF of \(H\) is
\[
 F_H(t)=
 \begin{cases}
 0,&t\le-1/2,\\
 t+1/2,&-1/2<t<1/2,\\
 1,&t\ge1/2.
 \end{cases}
\]
Independence gives \(J_j(y)=\EE F_H(y-S_j)\).
At \(y=x+1/2\) and \(y=x-1/2\), the displayed formula gives
\[
 \begin{aligned}
 J_j(x+1/2)
 &=F_j(x)+
   \EE[(x+1-S_j)\ind_{\{x<S_j<x+1\}}],\\
 J_j(x-1/2)
 &=F_j(x-1)+
   \EE[(x-S_j)\ind_{\{x-1<S_j<x\}}].
 \end{aligned}
\]
Since
\[
 F_j(x)-F_j(x-1)
 =\PP(x-1<S_j<x)+\PP(S_j=x),
\]
we obtain the exact identities
\begin{equation}\label{eq:cluster-jitter-cdf-identities}
 \begin{aligned}
 J_j(x+1/2)-F_j(x)
 &=\EE[(x+1-S_j)\ind_{\{x<S_j<x+1\}}],\\
 F_j(x)-J_j(x-1/2)
 &=\EE[(S_j-x+1)\ind_{\{x-1<S_j<x\}}]
      +\PP(S_j=x).
 \end{aligned}
\end{equation}
In particular, both right-hand sides are nonnegative, so
\begin{equation}\label{eq:cluster-global-jitter-sandwich}
 J_j(x-1/2)\le F_j(x)\le J_j(x+1/2),
 \qquad x\in\R.
\end{equation}
The atom term in the second identity is retained, so these
formulas do not require continuity of \(F_j\).

On \(x+1/4<S_j<x+1/2\), the first weight
\(x+1-S_j\) is at least \(1/2\).
On \(x-1/2<S_j<x-1/4\), the second weight
\(S_j-x+1\) is at least \(1/2\).
Hence \eqref{eq:cluster-jitter-cdf-identities} implies
\[
 \begin{aligned}
 J_j(x+1/2)-F_j(x)
 &\ge\frac12\PP\{S_j\in(x+1/4,x+1/2)\},\\
 F_j(x)-J_j(x-1/2)
 &\ge\frac12\PP\{S_j\in(x-1/2,x-1/4)\}.
 \end{aligned}
\]

Fix any standardized radius \(R>0\).
By \eqref{eq:cluster-variance-limit}, there are \(j_0\)
and a finite constant \(\Sigma>0\) such that
\(\sigma_j\le\Sigma\) for \(j\ge j_0\).
If \(|z_j(x)|\le R\), then
\[
 |x-n_jp_j|=\sigma_j\sqrt{n_j}|z_j(x)|
 \le R\Sigma\sqrt{n_j}.
\]
Apply \eqref{eq:cluster-local-mass} with \(T=R\Sigma\),
first with \((a,b)=(1/4,1/2)\) and then with
\((a,b)=(-1/2,-1/4)\).
Let the resulting positive constants be \(c_R^+\) and
\(c_R^-\), and set \(c_R=\frac12\min(c_R^+,c_R^-)>0\).
The preceding probability bounds give, eventually,
\begin{equation}\label{eq:cluster-strict-jitter-sandwich}
 J_j(x-1/2)+\frac{c_R}{\sqrt{n_j}}
 \le F_j(x)
 \le J_j(x+1/2)-\frac{c_R}{\sqrt{n_j}},
 \qquad |z_j(x)|\le R.
\end{equation}
This conclusion is available for every fixed \(R\);
the radius used for the final bound will be chosen below.

\textbf{Expand the shifted comparison functions with a uniform error.}
We first verify the convergence of the signed third moments.
Conditional centering implies
\[
 \EE[(B_j-p_j)^2U_j]
 =\EE[(B_j-p_j)^2\EE(U_j\mid B_j)]=0.
\]
Expanding the cube therefore gives
\[
 \begin{aligned}
 \EE(Y_j-p_j)^3
 &=\EE(B_j-p_j)^3
      +3\EE[(B_j-p_j)U_j^2]+\EE U_j^3\\
 &=v_j(1-2p_j)
      +3\EE[(B_j-p_j)U_j^2]+\EE U_j^3.
 \end{aligned}
\]
Here the Bernoulli third moment is
\[
 \EE(B_j-p_j)^3
 =p_j(1-p_j)^3-(1-p_j)p_j^3
 =v_j(1-2p_j).
\]
Since \(|B_j-p_j|\le1\) and \(|U_j|\le\epsilon_j\),
\[
 \left|\EE(Y_j-p_j)^3-v_j(1-2p_j)\right|
 \le3s_j+\epsilon_js_j\longrightarrow0.
\]
Together with \eqref{eq:cluster-variance-limit}, this proves
\begin{equation}\label{eq:cluster-skewness-limit}
 \gamma_j\longrightarrow\frac d{\sqrt v}.
\end{equation}
Thus \(\gamma_j\) is bounded, and \(\sigma_j\) is bounded
away from zero, for all sufficiently large \(j\).

Put \(\psi(z)=(1-z^2)\phi(z)\), and write the uniform
remainder in \eqref{eq:cluster-jitter-expansion} as
\[
 r_j:=
 \sup_{x\in\R}\left|
 J_j(x)-\Phi(z_j(x))
       -\frac{\gamma_j}{6\sqrt{n_j}}\psi(z_j(x))
 \right|=o(n_j^{-1/2}).
\]
Let
\[
 a_j=\frac1{2\sigma_j\sqrt{n_j}},\qquad
 M_1=\|\phi'\|_\infty,\qquad M_2=\|\psi'\|_\infty.
\]
These derivative bounds are finite because
\(\phi'(z)=-z\phi(z)\) and
\(\psi'(z)=(z^3-3z)\phi(z)\).
For \(\theta\in\{-1,1\}\), the standardized threshold
at \(x+\theta/2\) is \(z+\theta a_j\).
Taylor's theorem and the mean value theorem give
\[
 \begin{aligned}
 |\Phi(z+\theta a_j)-\Phi(z)-\theta a_j\phi(z)|
 &\le M_1a_j^2/2,\\
 |\psi(z+\theta a_j)-\psi(z)|
 &\le M_2a_j.
 \end{aligned}
\]
Substitution in \eqref{eq:cluster-jitter-expansion} yields
\begin{equation}\label{eq:cluster-shifted-jitter-expansion}
 J_j(x+\theta/2)
 =\Phi(z)+\frac1{\sqrt{n_j}}
   \left[\frac{\theta}{2\sigma_j}\phi(z)
         +\frac{\gamma_j}{6}\psi(z)\right]
   +E_{j,\theta}(x),
\end{equation}
where, uniformly over \(x\in\R\) and \(\theta\in\{-1,1\}\),
\[
 \begin{aligned}
 |E_{j,\theta}(x)|
 &\le r_j+\frac{M_1}{8\sigma_j^2n_j}
             +\frac{|\gamma_j|M_2}{12\sigma_jn_j}\\
 &\le r_j+\frac{K}{n_j}
 \end{aligned}
\]
for some finite \(K\) independent of \(j,x,\theta\), after
discarding finitely many indices.
Define
\begin{equation}\label{eq:cluster-shifted-jitter-error}
 \eta_j=\sqrt{n_j}\,r_j+\frac K{\sqrt{n_j}}
 \longrightarrow0.
\end{equation}

Use \(\theta=1\) for the upper bound in
\eqref{eq:cluster-strict-jitter-sandwich}, and
\(\theta=-1\) for the lower bound.
Equations \eqref{eq:cluster-shifted-jitter-expansion} and
\eqref{eq:cluster-shifted-jitter-error} then give
\begin{equation}\label{eq:cluster-improved-envelopes}
 \begin{aligned}
 \sqrt{n_j}[F_j(x)-\Phi(z)]
 &\le\frac{\phi(z)}{2\sigma_j}
       +\frac{\gamma_j}{6}\psi(z)-c_R+\eta_j,\\
 \sqrt{n_j}[\Phi(z)-F_j(x)]
 &\le\frac{\phi(z)}{2\sigma_j}
       -\frac{\gamma_j}{6}\psi(z)-c_R+\eta_j,
 \qquad |z|\le R.
 \end{aligned}
\end{equation}
Using the global comparison
\eqref{eq:cluster-global-jitter-sandwich} instead gives
the same two bounds without \(-c_R\), for every real \(x\).
Taking the maximum of these two signed bounds gives
\begin{equation}\label{eq:cluster-global-envelope}
 \sqrt{n_j}|F_j(x)-\Phi(z)|
 \le\mathcal E_j(z)+\eta_j,\qquad
 \mathcal E_j(z)=\frac{\phi(z)}{2\sigma_j}
                  +\frac{|\gamma_j|}{6}|\psi(z)|.
\end{equation}

\textbf{Combine the central improvement with the bound at distant thresholds.}
By \eqref{eq:jitter-envelope-gaussian-suprema},
\(\|\phi\|_\infty=\|\psi\|_\infty=\phi(0)\).
Hence
\[
 \mathcal E_j(z)\le
 L_j:=\phi(0)\left[\frac1{2\sigma_j}+\frac{|\gamma_j|}{6}\right].
\]
Equations \eqref{eq:cluster-variance-limit} and
\eqref{eq:cluster-skewness-limit} imply
\begin{equation}\label{eq:cluster-envelope-limit}
 L_j\longrightarrow
 L_*:=\frac{\phi(0)(3+|d|)}{6\sqrt v}>0.
\end{equation}
It follows from \eqref{eq:cluster-improved-envelopes} that,
for each fixed \(R\) and all sufficiently large \(j\),
\begin{equation}\label{eq:cluster-central-deficit}
 \sup_{|z_j(x)|\le R}
 \sqrt{n_j}|F_j(x)-\Phi(z_j(x))|
 \le L_j-c_R+\eta_j.
\end{equation}

For distant thresholds, define
\[
 \mathcal E(z)
 =\frac{\phi(z)}{2\sqrt v}
       +\frac{|d|}{6\sqrt v}|\psi(z)|.
\]
The same Gaussian suprema show that
\[
 \|\mathcal E_j-\mathcal E\|_\infty
 \le\phi(0)\left[
    \left|\frac1{2\sigma_j}-\frac1{2\sqrt v}\right|
    +\frac16\left||\gamma_j|-\frac{|d|}{\sqrt v}\right|
 \right]\longrightarrow0.
\]
Also, \(\phi(z)\to0\) and \(\psi(z)\to0\) as
\(|z|\to\infty\), so \(\mathcal E(z)\to0\).
We now choose a fixed \(R>0\) large enough that
\[
 \sup_{|z|>R}\mathcal E(z)\le L_*/2.
\]
For this chosen \(R\), take the positive constant \(c_R\)
in \eqref{eq:cluster-strict-jitter-sandwich}.
Equation~\eqref{eq:cluster-global-envelope} gives
\begin{equation}\label{eq:cluster-distant-deficit}
 \sup_{|z_j(x)|>R}
 \sqrt{n_j}|F_j(x)-\Phi(z_j(x))|
 \le\frac{L_*}{2}
       +\|\mathcal E_j-\mathcal E\|_\infty+\eta_j.
\end{equation}
The choice of \(R\) uses only the limiting function
\(\mathcal E\); the local mass argument then supplies
\(c_R>0\) for that fixed radius.

By \eqref{eq:raw-normalization},
\[
 \Delta_{n_j}(Y_j)
 =\sup_{x\in\R}|F_j(x)-\Phi(z_j(x))|.
\]
Splitting this supremum into \(|z_j(x)|\le R\) and
\(|z_j(x)|>R\), and applying
\eqref{eq:cluster-central-deficit} and
\eqref{eq:cluster-distant-deficit}, we obtain
\[
 \begin{aligned}
 \limsup_{j\to\infty}\sqrt{n_j}\,\Delta_{n_j}(Y_j)
 &\le\max\{L_*-c_R,L_*/2\}\\
 &=L_*-\min\{c_R,L_*/2\}.
 \end{aligned}
\]
Thus \eqref{eq:macroscopic-gap} holds with
\(\delta=\min\{c_R,L_*/2\}>0\).

\textbf{Deduce vanishing accumulated variance at the Esseen parameters.}
We first verify convergence of the normalizing third
absolute moments.
Since \(|B_j-p_j|\le1\) and \(|U_j|\le\epsilon_j\),
the mean value theorem for \(t\mapsto|t|^3\) gives
\[
 \bigl||Y_j-p_j|^3-|B_j-p_j|^3\bigr|
 \le3(1+\epsilon_j)^2|U_j|.
\]
Indeed, the derivative has absolute value \(3t^2\), and
the interval between \(B_j-p_j\) and \(Y_j-p_j\) lies
in \([-1-\epsilon_j,1+\epsilon_j]\).
Taking expectations and using the Bernoulli absolute
third moment yields
\[
 \begin{aligned}
 \left|\EE|Y_j-p_j|^3
       -v_j[p_j^2+(1-p_j)^2]\right|
 &\le3(1+\epsilon_j)^2\EE|U_j|\\
 &\le3(1+\epsilon_j)^2\epsilon_j\longrightarrow0.
 \end{aligned}
\]
Dividing by \(\sigma_j^3\to v^{3/2}>0\), we conclude that
\begin{equation}\label{eq:cluster-absolute-third-limit}
 \beta(Y_j)\longrightarrow
 \beta(B_p)=\frac{p^2+(1-p)^2}{\sqrt v}>0.
\end{equation}

If \(p\in\{\pe,1-\pe\}\), then
\eqref{eq:esseen-parameters} gives
\[
 p^2+(1-p)^2=10-3\sqrt{10},\qquad
 |d|=|1-2p|=\sqrt{10}-3.
\]
Therefore
\begin{equation}\label{eq:cluster-esseen-normalization}
 \frac{L_*}{\beta(B_p)}
 =\frac{\phi(0)\sqrt{10}}{6(10-3\sqrt{10})}
 =\frac{\phi(0)(3+\sqrt{10})}{6}
 =\ce,
\end{equation}
where the last equality is \eqref{eq:ce}.

Now suppose \(R_{n_j}(Y_j)\to\ce\) but \(n_js_j\not\to0\).
Since \(n_js_j\ge0\), there are a number \(\lambda_*>0\)
and a subsequence on which \(n_js_j\ge\lambda_*\).
All the original convergence assumptions remain valid
along this subsequence.
The first part of the proof therefore gives
\[
 \limsup\sqrt{n_j}\,\Delta_{n_j}(Y_j)\le L_*-\delta
\]
along it, for some \(\delta>0\).
On the other hand, \eqref{eq:raw-normalization},
\eqref{eq:cluster-absolute-third-limit}, and
\eqref{eq:cluster-esseen-normalization} imply along the
same subsequence that
\[
 \sqrt{n_j}\,\Delta_{n_j}(Y_j)
 =\beta(Y_j)R_{n_j}(Y_j)
 \longrightarrow\beta(B_p)\ce=L_*.
\]
This contradiction proves \(n_js_j\to0\).
\end{proof}

\begin{proof}[Proof of Proposition~\ref{prop:clusters}]
\textbf{Choose a sequence of counterexamples and center each cluster.}
We first prove the assertion with upper-cluster mass close
to \(\pe\).
Suppose that no pair \((\eta,N)\) has the asserted property.
Taking \(\eta=1/j\) and \(N=j\), for every sufficiently
large integer \(j>4\) we may choose a law \(Q_j\) and an
integer \(n_j\ge j\) such that
\begin{equation}\label{eq:clusters-counterexample-sequence}
 \begin{gathered}
 \supp Q_j\subset[-1/j,1/j]\cup[1-1/j,1+1/j],\\
 \left|Q_j([1-1/j,1+1/j])-\pe\right|<1/j,
 \qquad R_{n_j}(Q_j)>\ce.
 \end{gathered}
\end{equation}
In particular, \(n_j\to\infty\).

Let \(X_j\sim Q_j\), and define
\[
 B_j=\ind_{\{X_j\in[1-1/j,1+1/j]\}},\qquad
 p_j=\PP(B_j=1).
\]
Then \(|p_j-\pe|<1/j\), so \(p_j\to\pe\in(0,1)\).
Thus both cluster probabilities are positive for all
sufficiently large \(j\).
The intervals are disjoint, and \(X_j\) is bounded, so
\(Q_j\) has positive variance and finite third absolute
moment.
The conditional means
\[
 m_{bj}=\EE(X_j\mid B_j=b),\qquad b=0,1,
\]
are therefore well defined for these indices.
Since a mean of a variable supported in an interval belongs
to that interval, we have
\[
 -1/j\le m_{0j}\le1/j,\qquad
 1-1/j\le m_{1j}\le1+1/j.
\]
Consequently, with \(d_j=m_{1j}-m_{0j}\),
\begin{equation}\label{eq:clusters-mean-separation}
 1-2/j\le d_j\le1+2/j,\qquad d_j>0.
\end{equation}

Define the transformed variable and its conditional noise by
\begin{equation}\label{eq:clusters-centered-transformation}
 Y_j=\frac{X_j-m_{0j}}{d_j},\qquad
 U_j=\frac{X_j-m_{B_j,j}}{d_j}.
\end{equation}
Here \(m_{B_j,j}=m_{0j}+d_jB_j\), so
\[
 Y_j=B_j+U_j,\qquad
 \EE(U_j\mid B_j=b)
 =\frac{\EE(X_j\mid B_j=b)-m_{bj}}{d_j}=0.
\]
On each cluster, \(X_j\) and its conditional mean belong
to the same interval of length \(2/j\).
Hence \eqref{eq:clusters-mean-separation} gives
\begin{equation}\label{eq:clusters-noise-bound}
 |U_j|\le\frac{2/j}{1-2/j}=:\epsilon_j
 \longrightarrow0.
\end{equation}
By affine invariance in \eqref{eq:raw-normalization},
\begin{equation}\label{eq:clusters-transformed-violation}
 R_{n_j}(Y_j)=R_{n_j}(Q_j)>\ce.
\end{equation}

\textbf{Verify convergence of the standardized laws.}
Put
\[
 s_j=\EE U_j^2,\qquad
 \sigma_j^2=p_j(1-p_j)+s_j.
\]
Conditional centering gives
\[
 \EE U_j=0,\qquad
 \EE[(B_j-p_j)U_j]=0,
\]
and therefore \(\EE Y_j=p_j\) and
\(\operatorname{Var}(Y_j)=\sigma_j^2\).
Moreover,
\begin{equation}\label{eq:clusters-variance-convergence}
 0\le s_j\le\epsilon_j^2\longrightarrow0,\qquad
 \sigma_j^2\longrightarrow\pe\qe=\se^2>0.
\end{equation}
Let \(P_j\) be the law of
\[
 Z_j=\frac{Y_j-p_j}{\sigma_j}.
\]
Then \(P_j\in\Pthree\): the variable \(Z_j\) has mean zero,
variance one, and finite third absolute moment because
\(Y_j\) is bounded.

We now verify \(P_j\to\Pe\) in \(W_3\).
Let \(\widehat P_j\) be the law of
\((B_j-p_j)/\sigma_j\).
The given coupling \(Y_j=B_j+U_j\) gives
\[
 W_3(P_j,\widehat P_j)
 \le\frac{(\EE|U_j|^3)^{1/3}}{\sigma_j}
 \le\frac{\epsilon_j}{\sigma_j}.
\]
To compare \(\widehat P_j\) with \(\Pe\), take
\(V\sim\mathrm{Unif}[0,1]\), and put
\[
 B_j^*=\ind_{\{V\le p_j\}},\qquad
 B_{\mathrm E}^*=\ind_{\{V\le\pe\}}.
\]
Then
\[
 \EE|B_j^*-B_{\mathrm E}^*|^3=|p_j-\pe|,
\]
and the law of \((B_{\mathrm E}^*-\pe)/\se\) is
\(\Pe\), by \eqref{eq:esseen-law}.
The identity
\[
 \begin{aligned}
 \frac{B_j^*-p_j}{\sigma_j}
       -\frac{B_{\mathrm E}^*-\pe}{\se}
 &=
 \frac{B_j^*-B_{\mathrm E}^*-(p_j-\pe)}{\sigma_j}\\
 &\quad+
 \left(\frac1{\sigma_j}-\frac1\se\right)
       (B_{\mathrm E}^*-\pe)
 \end{aligned}
\]
and the triangle inequality for the \(L^3\) norm yield
\[
 W_3(\widehat P_j,\Pe)
 \le
 \frac{|p_j-\pe|^{1/3}+|p_j-\pe|}{\sigma_j}
 +\left|\frac1{\sigma_j}-\frac1\se\right|.
\]
Here we used \(|B_{\mathrm E}^*-\pe|\le1\).
Combining the two Wasserstein bounds gives
\begin{equation}\label{eq:clusters-standardized-convergence}
 \begin{aligned}
 W_3(P_j,\Pe)
 &\le
 \frac{\epsilon_j+|p_j-\pe|^{1/3}+|p_j-\pe|}{\sigma_j}
 +\left|\frac1{\sigma_j}-\frac1\se\right|
 \longrightarrow0.
 \end{aligned}
\end{equation}
Every term tends to zero by \(p_j\to\pe\),
\eqref{eq:clusters-noise-bound}, and
\eqref{eq:clusters-variance-convergence}.
In particular, convergence of third absolute moments in
\(W_3\) implies
\begin{equation}\label{eq:clusters-third-moment-convergence}
 \beta(Y_j)=\beta(P_j)\longrightarrow\beta(\Pe)=\be>0.
\end{equation}

\textbf{Identify the limiting error and apply the small-variance comparison.}
We can now apply Corollary~\ref{cor:jitter-envelopes}
to \(P_j\to\Pe\) in \(W_3\) and \(n_j\to\infty\).
For the limiting law, the third moment is \(\ke\) and
the maximal lattice span is \(h_{\mathrm E}\).
Equation~\eqref{eq:triangular-upper} gives
\begin{equation}\label{eq:clusters-asymptotic-upper}
 \limsup_{j\to\infty}
       \sqrt{n_j}\,\Delta_{n_j}(P_j)
 \le\frac{\phi(0)}6\bigl(|\ke|+3h_{\mathrm E}\bigr)
 =\ce\be.
\end{equation}
The equality uses the equality case of
\eqref{eq:esseen-moment} for \(\Pe\) and
\(\ce=\phi(0)(3+\sqrt{10})/6\) from \eqref{eq:ce}.
By affine invariance and
\eqref{eq:clusters-third-moment-convergence},
\[
 R_{n_j}(Y_j)
 =\frac{\sqrt{n_j}\,\Delta_{n_j}(P_j)}{\beta(P_j)}.
\]
To make the division of the asymptotic bound explicit, fix
\(0<\varepsilon<\be\).
Equations \eqref{eq:clusters-asymptotic-upper} and
\eqref{eq:clusters-third-moment-convergence} ensure that,
for all sufficiently large \(j\),
\[
 R_{n_j}(Y_j)
 \le\frac{\ce\be+\varepsilon}{\be-\varepsilon}.
\]
First take the upper limit in \(j\), and then let
\(\varepsilon\downarrow0\), to obtain
\(\limsup_jR_{n_j}(Y_j)\le\ce\).
Since every term is greater than \(\ce\) by
\eqref{eq:clusters-transformed-violation}, it follows that
\begin{equation}\label{eq:clusters-error-convergence}
 R_{n_j}(Y_j)\longrightarrow\ce.
\end{equation}

The representation \eqref{eq:clusters-centered-transformation}
and the bounds above satisfy the hypotheses of
Lemma~\ref{lem:accumulated-cluster-variance}:
\[
 p_j\to\pe,\qquad n_j\to\infty,\qquad
 \EE(U_j\mid B_j)=0,\qquad |U_j|\le\epsilon_j\to0.
\]
Its final assertion, together with
\eqref{eq:clusters-error-convergence}, therefore gives
\begin{equation}\label{eq:clusters-accumulated-variance}
 n_js_j=n_j\EE U_j^2\longrightarrow0.
\end{equation}

Let \(I_0,\epsilon_0,\lambda_0,N_0\) denote the fixed
interval and constants supplied by
Lemma~\ref{lem:small-cluster-variance}.
Because \(\pe\) lies in the interior of \(I_0\),
\(p_j\to\pe\) implies \(p_j\in I_0\) eventually.
Equations \eqref{eq:clusters-noise-bound} and
\eqref{eq:clusters-accumulated-variance}, together with
\(n_j\to\infty\), also give, eventually,
\[
 |U_j|\le\epsilon_j\le\epsilon_0,\qquad
 n_js_j\le\lambda_0,\qquad n_j\ge N_0.
\]
Thus every hypothesis of that lemma holds, and it yields
\[
 R_{n_j}(Y_j)\le R_{n_j}(B_j)<\ce.
\]
This contradicts \eqref{eq:clusters-transformed-violation}.
The contradiction proves the existence of fixed
\(\eta\in(0,1/4)\) and \(N<\infty\) for the assertion
at \(\pe\).

\textbf{Transfer the same constants to the reflected parameter.}
Fix the \(\eta,N\) just obtained, and suppose \(Q\) satisfies
the support condition of the proposition and
\[
 |Q([1-\eta,1+\eta])-(1-\pe)|<\eta.
\]
For \(X\sim Q\), let \(Q'\) be the law of \(1-X\).
The map \(x\mapsto1-x\) exchanges the two intervals, so
\[
 \supp Q'\subset[-\eta,\eta]\cup[1-\eta,1+\eta].
\]
Since these intervals are disjoint and carry all the mass,
\[
 \begin{aligned}
 Q'([1-\eta,1+\eta])
 &=Q([-\eta,\eta])
 =1-Q([1-\eta,1+\eta]),\\
 |Q'([1-\eta,1+\eta])-\pe|
 &=|Q([1-\eta,1+\eta])-(1-\pe)|<\eta.
 \end{aligned}
\]
The assertion already proved at \(\pe\) gives
\(R_n(Q')\le\ce\) for every integer \(n\ge N\).
Affine invariance gives \(R_n(Q)=R_n(Q')\), proving the
reflected assertion with the same \(\eta,N\).
\end{proof}

% ===== extremizers =====
\section{Extremizers and uniform support bounds}
\label{sec:extremizers}

We use the extremal constants \(C_n\) from \eqref{eq:cn}, the moment
cutoff \(B_*\) from \eqref{eq:moment-cutoff}, and the convergence in
\eqref{eq:cn-limit}. Write \(C_0=0.4690\), so that \(C_n\le C_0\) for
all \(n\).

Extremal problems for the distribution function of an iid sum under
moment and range constraints were studied by Hoeffding and
Shrikhande~\cite{HoeffdingShrikhande1955}.
For finite-support and few-point reduction principles in related
moment problems, see also \cite[Section~3]{MattnerShevtsova2019}.
Below we differentiate the normalized signed discrepancy at a fixed
active threshold to obtain contact identities for maximizing laws.

We first show that, whenever \(C_n>\ce\), there is a law \(P\)
and a threshold \(t\) at which the normalized signed discrepancy
equals \(C_n\). This gives the maximizing pair needed for the
perturbation argument in Lemma~\ref{lem:influence}.
The moment cutoff \eqref{eq:moment-cutoff} allows us to pass
to a weak limit of a maximizing sequence, and maximality then
shows that this limit attains the same value.

\begin{lemma}[Attainment above the Esseen constant]
\label{lem:attainment}
If \(C_n>\ce\), there are \(P\in\Pthree\) and \(t\in\R\) such that
\begin{equation}
\label{eq:active-branch}
\frac{\sqrt n}{\beta(P)}
\left\{F_{n,P}(t)-\Phi(t/\sqrt n)\right\}
=R_n(P)=C_n,
\qquad
1\le\beta(P)<B_*.
\end{equation}
\end{lemma}

\begin{proof}
Fix \(n\) such that \(C_n>\ce\).

\textbf{Choose a maximizing sequence and identify its limiting moments.}
By the signed-supremum representation in \eqref{eq:cn}, we can
choose \(P_k\in\Pthree\) and \(t_k\in\R\) such that, with
\(\beta_k=\beta(P_k)\),
\begin{equation}\label{eq:attainment-maximizing-sequence}
 r_k:=
 \frac{\sqrt n}{\beta_k}
 \left\{F_{n,P_k}(t_k)-\Phi(t_k/\sqrt n)\right\}
 >\ce,
 \qquad r_k\longrightarrow C_n.
\end{equation}
For every \(Q\in\Pthree\), H\"older's inequality gives
\begin{equation}\label{eq:attainment-moment-lower}
 1=\int x^2\,Q(dx)
 \le \left(\int |x|^3\,Q(dx)\right)^{2/3},
 \qquad \beta(Q)\ge1.
\end{equation}
Also, \(R_n(P_k)\ge r_k>\ce\), so the moment cutoff
\eqref{eq:moment-cutoff} applies to every \(P_k\). Consequently,
\begin{equation}\label{eq:attainment-moment-bounds}
 1\le\beta_k<B_*.
\end{equation}

Since \(\int x^2\,P_k(dx)=1\), for every \(A>0\) we have
\[
 \sup_k P_k(\{|x|>A\})\le \frac1{A^2}.
\]
Thus \((P_k)\) is tight. By
\cite[Theorem~5.1]{Billingsley1999}, it has a subsequence
converging weakly to a probability law \(P\).
The boundedness in \eqref{eq:attainment-moment-bounds} permits
a further subsequence along which \(\beta_k\) converges.
Retaining the same indices, we therefore have
\begin{equation}\label{eq:attainment-weak-limit}
 P_k\Longrightarrow P,
 \qquad \beta_k\longrightarrow b\in[1,B_*].
\end{equation}

We check explicitly that the first two moments are preserved.
For \(A>0\), \eqref{eq:attainment-moment-bounds} gives
\begin{equation}\label{eq:attainment-moment-tails}
 \begin{aligned}
 \int_{|x|>A}x^2\,P_k(dx)
 &\le\frac1A\int_{|x|>A}|x|^3\,P_k(dx)
 \le\frac{B_*}{A},\\
 \int_{|x|>A}|x|\,P_k(dx)
 &\le\frac1{A^2}\int_{|x|>A}|x|^3\,P_k(dx)
 \le\frac{B_*}{A^2}.
 \end{aligned}
\end{equation}
Define the bounded continuous functions
\[
 q_A(x)=\min\{x^2,A^2\},
 \qquad h_A(x)=\max\{-A,\min\{x,A\}\}.
\]
Because
\[
 0\le x^2-q_A(x)\le x^2\mathbf1_{\{|x|>A\}},
 \qquad
 |x-h_A(x)|\le |x|\mathbf1_{\{|x|>A\}},
\]
the moment conditions on \(P_k\) and
\eqref{eq:attainment-moment-tails} imply
\[
 1-\frac{B_*}{A}\le\int q_A\,dP_k\le1,
 \qquad
 \left|\int h_A\,dP_k\right|\le\frac{B_*}{A^2}.
\]
For each fixed \(A\), weak convergence in
\eqref{eq:attainment-weak-limit} allows us to pass to the limit
in these integrals:
\begin{equation}\label{eq:attainment-truncated-moments}
 1-\frac{B_*}{A}\le\int q_A\,dP\le1,
 \qquad
 \left|\int h_A\,dP\right|\le\frac{B_*}{A^2}.
\end{equation}
As \(A\uparrow\infty\), \(q_A(x)\uparrow x^2\).
The monotone convergence theorem and the first inequality in
\eqref{eq:attainment-truncated-moments} give
\(\int x^2\,P(dx)=1\).
In particular, the Cauchy--Schwarz inequality gives
\(\int |x|\,P(dx)\le1\).
Since \(h_A(x)\to x\) and \(|h_A(x)|\le|x|\), dominated
convergence and the second inequality in
\eqref{eq:attainment-truncated-moments} now give
\(\int x\,P(dx)=0\).

To control the third absolute moment, apply weak convergence
to the bounded continuous function \(\min\{|x|^3,A^3\}\):
\[
 \int \min\{|x|^3,A^3\}\,P(dx)
 =\lim_{k\to\infty}
   \int \min\{|x|^3,A^3\}\,P_k(dx)
 \le\lim_{k\to\infty}\beta_k=b.
\]
Letting \(A\uparrow\infty\) and using monotone convergence
shows that \(P\in\Pthree\) and
\begin{equation}\label{eq:attainment-limit-moment}
 1\le\beta(P)\le b.
\end{equation}
Here the lower bound follows from
\eqref{eq:attainment-moment-lower}.

\textbf{Show that the thresholds remain bounded for this fixed \(n\).}
By \eqref{eq:attainment-maximizing-sequence} and
\eqref{eq:attainment-moment-bounds}, the positive discrepancies
satisfy
\begin{equation}\label{eq:attainment-discrepancy-lower}
 D_k:=F_{n,P_k}(t_k)-\Phi(t_k/\sqrt n)
 =\frac{r_k\beta_k}{\sqrt n}
 >\frac{\ce}{\sqrt n}.
\end{equation}
If \(t_k\ge0\), then
\[
 D_k\le1-\Phi(t_k/\sqrt n).
\]
Choose \(T_+>0\) so that
\(1-\Phi(T_+/\sqrt n)<\ce/\sqrt n\), which is possible
because the Gaussian upper tail tends to zero.
The last inequality and \eqref{eq:attainment-discrepancy-lower}
then imply \(t_k<T_+\) whenever \(t_k\ge0\).

For \(t_k<0\), let \(S_k\) be the sum of \(n\) independent
variables with common law \(P_k\). Then
\(\EE S_k=0\) and \(\EE S_k^2=n\), so Chebyshev's inequality gives
\[
 D_k\le F_{n,P_k}(t_k)
 =\PP\{S_k\le t_k\}
 \le\PP\{|S_k|\ge|t_k|\}
 \le\frac{n}{t_k^2}.
\]
Together with \eqref{eq:attainment-discrepancy-lower}, this gives
\(|t_k|<n^{3/4}/\sqrt{\ce}\) whenever \(t_k<0\).
Thus \((t_k)\) is bounded. Passing to a further subsequence,
we may assume
\begin{equation}\label{eq:attainment-threshold-limit}
 t_k\longrightarrow t\in\R.
\end{equation}

\textbf{Pass to the limit at the moving thresholds.}
Since \(n\) is fixed, weak convergence of \(P_k\) implies
weak convergence of the product laws \(P_k^{\otimes n}\) to
\(P^{\otimes n}\); see
\cite[Theorem~2.8(ii)]{Billingsley1999}.
Applying the continuous sum map
\((x_1,\ldots,x_n)\mapsto x_1+\cdots+x_n\) gives
\(P_k^{*n}\Longrightarrow P^{*n}\), by
\cite[Theorem~2.7]{Billingsley1999}.
Fix \(\varepsilon>0\). By
\eqref{eq:attainment-threshold-limit}, \(t_k\le t+\varepsilon\)
for all sufficiently large \(k\). Hence the closed-set
inequality in the Portmanteau theorem
\cite[Theorem~2.1(iii)]{Billingsley1999} yields
\[
 \begin{aligned}
 \limsup_{k\to\infty}F_{n,P_k}(t_k)
 &\le \limsup_{k\to\infty}
       P_k^{*n}((-\infty,t+\varepsilon])\\
 &\le P^{*n}((-\infty,t+\varepsilon])
 =F_{n,P}(t+\varepsilon).
 \end{aligned}
\]
Letting \(\varepsilon\downarrow0\) and using the right
continuity of \(F_{n,P}\), we obtain
\begin{equation}\label{eq:attainment-moving-threshold}
 F_{n,P}(t)\ge
 \limsup_{k\to\infty}F_{n,P_k}(t_k).
\end{equation}
On the other hand,
\eqref{eq:attainment-maximizing-sequence},
\eqref{eq:attainment-weak-limit}, and
\eqref{eq:attainment-threshold-limit}, together with the
continuity of \(\Phi\), give
\[
 F_{n,P_k}(t_k)
 =\Phi(t_k/\sqrt n)+\frac{r_k\beta_k}{\sqrt n}
 \longrightarrow
 \Phi(t/\sqrt n)+\frac{C_n b}{\sqrt n}.
\]
Substitution into \eqref{eq:attainment-moving-threshold}
therefore gives
\begin{equation}\label{eq:attainment-limit-discrepancy}
 F_{n,P}(t)-\Phi(t/\sqrt n)
 \ge\frac{C_n b}{\sqrt n}.
\end{equation}

\textbf{Use maximality to obtain equality.}
Because \(P\in\Pthree\), the definition \eqref{eq:cn},
the discrepancy bound \eqref{eq:attainment-limit-discrepancy},
and the moment bound \eqref{eq:attainment-limit-moment} give,
in this order,
\begin{equation}\label{eq:attainment-equality-chain}
 C_n\ge
 \frac{\sqrt n\{F_{n,P}(t)-\Phi(t/\sqrt n)\}}{\beta(P)}
 \ge C_n\frac{b}{\beta(P)}
 \ge C_n.
\end{equation}
All inequalities in \eqref{eq:attainment-equality-chain}
must therefore be equalities. Since \(C_n>0\), this implies
\(\beta(P)=b\), as well as
\[
 \frac{\sqrt n\{F_{n,P}(t)-\Phi(t/\sqrt n)\}}{\beta(P)}
 =C_n.
\]
Finally, the definition of \(R_n(P)\) as a supremum of
absolute discrepancies and \eqref{eq:cn} imply
\[
 C_n=
 \frac{\sqrt n\{F_{n,P}(t)-\Phi(t/\sqrt n)\}}{\beta(P)}
 \le R_n(P)\le C_n.
\]
Thus \(R_n(P)=C_n>\ce\). Applying
\eqref{eq:moment-cutoff} once more gives \(\beta(P)<B_*\).
Together with \eqref{eq:attainment-limit-moment}, this proves
\eqref{eq:active-branch}.
\end{proof}

To obtain the uniform support bound in
Proposition~\ref{prop:bounded-extremizers}, we will need to
choose sample sizes carefully. Under the assumption that
\(C_n>\ce\) for arbitrarily large \(n\), we seek indices
\(n_j\uparrow\infty\) such that
\[
 C_{n_j}>\ce,
 \qquad
 n_j(C_{n_j-1}-C_{n_j})_+\longrightarrow0.
\]
The second condition controls the term involving
\(n(C_{n-1}-C_n)_+\) in \eqref{eq:influence-bound}.
The following elementary lemma supplies these indices when
applied to \(a_n=C_n-\ce\), using \eqref{eq:cn-limit}.
No monotonicity of \((C_n)\) is required.

\begin{lemma}[Selection by harmonic increments]
\label{lem:selection}
Suppose that \(a_n\to0\) and that \(a_n>0\) for arbitrarily large
integers \(n\). There are integers \(n_j\uparrow\infty\) such that
\begin{equation}
\label{eq:harmonic-selection}
a_{n_j}>0,
\qquad
n_j(a_{n_j-1}-a_{n_j})_+\longrightarrow0.
\end{equation}
\end{lemma}

\begin{proof}
\textbf{Find arbitrarily late positive terms with a small preceding decrease.}
We first prove that, for every \(\varepsilon>0\) and every
integer \(N\ge2\), there is an integer \(n\ge N\) such that
\begin{equation}\label{eq:selection-good-index}
 a_n>0,
 \qquad n(a_{n-1}-a_n)_+<\varepsilon.
\end{equation}
Suppose otherwise. Then there are \(\varepsilon_0>0\) and
an integer \(N_0\ge2\) such that
\[
 n(a_{n-1}-a_n)_+\ge\varepsilon_0
 \qquad
 \text{whenever } n\ge N_0 \text{ and } a_n>0.
\]
Since \(\varepsilon_0/n>0\), the positive part in this
inequality must be strictly positive and therefore equals
\(a_{n-1}-a_n\). Consequently,
\begin{equation}\label{eq:selection-forced-decrease}
 a_n>0
 \quad\Longrightarrow\quad
 a_{n-1}\ge a_n+\frac{\varepsilon_0}{n}>0,
 \qquad n\ge N_0.
\end{equation}

We now show that every term with index at least \(N_0\)
is positive. Fix an integer \(m\ge N_0\).
By the hypothesis of arbitrarily late positive terms,
there is an integer \(r\ge m\) with \(a_r>0\).
If \(r>m\), apply \eqref{eq:selection-forced-decrease}
successively at the indices \(r,r-1,\ldots,m+1\).
At each step, positivity of the current term implies
positivity of the preceding one, so this gives \(a_m>0\).
If \(r=m\), the same conclusion already holds.
Since \(m\ge N_0\) was arbitrary, all terms from \(N_0\)
onward are positive.

For any integer \(m>N_0\), we may therefore apply
\eqref{eq:selection-forced-decrease} at every index
\(k=N_0+1,\ldots,m\) and sum the resulting inequalities:
\[
 a_{N_0}-a_m
 =\sum_{k=N_0+1}^{m}(a_{k-1}-a_k)
 \ge\varepsilon_0\sum_{k=N_0+1}^{m}\frac1k.
\]
It follows that
\begin{equation}\label{eq:selection-harmonic-contradiction}
 0<a_m\le
 a_{N_0}-\varepsilon_0\sum_{k=N_0+1}^{m}\frac1k.
\end{equation}
The harmonic sum tends to \(+\infty\) as \(m\to\infty\),
so the right-hand side of
\eqref{eq:selection-harmonic-contradiction} is negative
for all sufficiently large \(m\). This is a contradiction,
which proves \eqref{eq:selection-good-index}.

\textbf{Choose the subsequence recursively.}
Set \(n_0=1\).
For each \(j\ge1\), apply \eqref{eq:selection-good-index}
with
\[
 \varepsilon=\frac1j,
 \qquad N=n_{j-1}+1\ge2.
\]
This gives an integer \(n_j>n_{j-1}\) satisfying
\[
 a_{n_j}>0,
 \qquad
 0\le n_j(a_{n_j-1}-a_{n_j})_+<\frac1j.
\]
The strictly increasing sequence of integers \(n_j\)
tends to infinity, and the last inequality implies
\(n_j(a_{n_j-1}-a_{n_j})_+\to0\).
These are precisely the conclusions in
\eqref{eq:harmonic-selection}.
\end{proof}

We now derive a first-order condition for a maximizing pair.
We add a small mass at an arbitrary point \(y\), keeping the
raw threshold \(t\) fixed and allowing the mean and variance
to change. These changes are included in the normalization
\eqref{eq:raw-normalization}.
We will show that the resulting derivative is nonpositive
everywhere and vanishes at every point of the support.
Proposition~\ref{prop:bounded-extremizers} will use this
vanishing condition to exclude distant support points.
The resulting contact equation will also be used in
Lemma~\ref{lem:support-separation} to compare distinct
support points.

\begin{lemma}[Contamination influence and support contact]
\label{lem:influence}
Let \(n\ge2\), \(P\in\Pthree\), and \(t\in\R\) satisfy
\begin{equation}\label{eq:influence-active-pair}
 \frac{\sqrt n}{\beta(P)}
 \left\{F_{n,P}(t)-\Phi(t/\sqrt n)\right\}
 =R_n(P)=C_n=:R>0.
\end{equation} Put
\[
z=\frac{t}{\sqrt n},\qquad
\beta=\beta(P),\qquad M=\int x|x|\,P(dx),\qquad G=F_{n-1,P}.
\]
For \(y\in\R\), contaminate the raw law by
\(P_e=(1-e)P+e\delta_y\), and evaluate its normalized signed discrepancy
at the same raw threshold \(t\). Its right derivative \(I(y)\) at
\(e=0\) exists and satisfies
\begin{align}
\beta I(y)={}&n^{3/2}\{G(t-y)-F_{n,P}(t)\}
+ny\phi(z)+\frac{\sqrt n}{2}z\phi(z)(y^2-1)\notag\\
&+\frac32R\beta(y^2-1)
-R\{|y|^3-\beta-3yM\}.\label{eq:influence}
\end{align}
Moreover,
\begin{equation}
\label{eq:contact}
I(y)\le0\quad(y\in\R),
\qquad I(x)=0\quad(x\in\supp P).
\end{equation}
In particular, every \(x\in\supp P\) satisfies the exact contact equation
\begin{align}
n\{G(t-x)-F_{n,P}(t)\}={}&-\sqrt n\,\phi(z)x
-\frac{z\phi(z)}2(x^2-1)\notag\\
&+\frac{R}{\sqrt n}
\left\{|x|^3-\beta-3Mx-\frac32\beta(x^2-1)\right\}.
\label{eq:contact-equation}
\end{align}
\end{lemma}

\begin{proof}
\textbf{Keep the raw threshold fixed and use maximality.}
Fix \(y\in\R\). For \(0\le e<1\), let \(m_e,v_e,\rho_e\)
denote the mean, variance, and third absolute central moment
of \(P_e=(1-e)P+e\delta_y\).
The first two moment conditions on \(P\) give
\begin{equation}\label{eq:influence-raw-moments}
 \begin{aligned}
 m_e&=ey,\\
 v_e&=(1-e)\int x^2\,P(dx)+ey^2-m_e^2\\
    &=1+e(y^2-1)-e^2y^2
      =(1-e)(1+ey^2)>0.
 \end{aligned}
\end{equation}
Also,
\[
 \int |x|^3\,P_e(dx)=(1-e)\beta+e|y|^3<\infty.
\]
The inequality
\(|x-m_e|^3\le4(|x|^3+|m_e|^3)\) therefore shows that
\(\rho_e<\infty\). Moreover, \(\rho_e=0\) would force
\(P_e=\delta_{m_e}\) and hence \(v_e=0\).
Thus \(\rho_e>0\).

Write
\[
 A_e=P_e^{*n}((-\infty,t]),\qquad
 z_e=\frac{t-nm_e}{\sqrt{nv_e}},\qquad
 \beta_e=\frac{\rho_e}{v_e^{3/2}}.
\]
Here \(t\) is fixed throughout the perturbation.
By the normalization \eqref{eq:raw-normalization}, the
normalized signed discrepancy at this threshold is
\begin{equation}\label{eq:influence-perturbed-discrepancy}
 D_e:=A_e-\Phi(z_e),
 \qquad
 r(e):=\frac{\sqrt n\,D_e}{\beta_e}.
\end{equation}
Let \(\widetilde P_e\) be the law of
\((X_e-m_e)/\sqrt{v_e}\), where \(X_e\) has law \(P_e\).
Then \(\widetilde P_e\in\Pthree\),
\(\beta(\widetilde P_e)=\beta_e\), and
\[
 A_e=F_{n,\widetilde P_e}(\sqrt n\,z_e).
\]
Thus \(r(e)\) is one of the signed ratios in the supremum
\eqref{eq:cn}. Since \(m_0=0\), \(v_0=1\), and
\(\rho_0=\beta\), the active-pair identity
\eqref{eq:influence-active-pair} gives
\begin{equation}\label{eq:influence-maximality}
 r(e)\le C_n=R=r(0),\qquad 0\le e<1.
\end{equation}

\textbf{Differentiate the moment terms.}
All derivatives at \(e=0\) below are right derivatives.
The exact formulas \eqref{eq:influence-raw-moments} give
\(m'_0=y\) and \(v'_0=y^2-1\).
Since \(m_e=ey\), the third absolute central moment is
\begin{equation}\label{eq:influence-central-moment}
 \rho_e
 =(1-e)\int |x-ey|^3\,P(dx)
   +e(1-e)^3|y|^3.
\end{equation}
For each \(x\), the integrand is differentiable in \(e\), with
\[
 \frac{\partial}{\partial e}|x-ey|^3
 =-3y(x-ey)|x-ey|.
\]
For \(0\le e\le1/2\), its absolute value is bounded by
\[
 3|y|(|x|+|y|)^2
 \le6|y|x^2+6|y|^3.
\]
This bound is integrable under \(P\), since
\(\int x^2\,P(dx)=1\).
By the mean value theorem, the same bound controls the
difference quotients
\((|x-ey|^3-|x|^3)/e\) for \(0<e\le1/2\).
Dominated convergence therefore gives
\[
 \left.\frac{d}{de}\int |x-ey|^3\,P(dx)\right|_{e=0+}
 =-3y\int x|x|\,P(dx)=-3yM.
\]
Differentiating \eqref{eq:influence-central-moment},
including both mixing weights, now yields
\begin{equation}\label{eq:influence-moment-derivatives}
 m'_0=y,\qquad
 v'_0=y^2-1,\qquad
 \rho'_0=|y|^3-\beta-3yM.
\end{equation}
Because \(\beta_e=\rho_e v_e^{-3/2}\) and
\((v_0,\rho_0)=(1,\beta)\), the product and chain rules give
\begin{equation}\label{eq:influence-beta-derivative}
 \begin{aligned}
 \beta'_0
 &=\rho'_0-\frac32\beta v'_0\\
 &=|y|^3-\beta-3yM-\frac32\beta(y^2-1).
 \end{aligned}
\end{equation}

\textbf{Differentiate the convolution and the normal term.}
Expanding the convolution of the mixture gives the exact identity
\begin{equation}\label{eq:influence-mixture-expansion}
 A_e
 =\sum_{\ell=0}^{n}
   \binom n\ell(1-e)^{n-\ell}e^\ell
   P^{*(n-\ell)}((-\infty,t-\ell y]).
\end{equation}
Here \(P^{*0}:=\delta_0\).
The term indexed by \(\ell\) corresponds to choosing
\(\ell\) of the \(n\) summands from the point mass at \(y\).
All distribution-function values on the right side of
\eqref{eq:influence-mixture-expansion} are independent of \(e\).
The term \(\ell=0\) has derivative \(-nF_{n,P}(t)\)
at zero, the term \(\ell=1\) has derivative \(nG(t-y)\),
and every term with \(\ell\ge2\) has derivative zero.
Consequently,
\begin{equation}\label{eq:influence-convolution-derivative}
 A'_0=n\{G(t-y)-F_{n,P}(t)\}.
\end{equation}
Thus the convolution derivative is obtained by differentiating
the mixing weights at the fixed raw threshold \(t\).

Next, write \(z_e=(t-nm_e)v_e^{-1/2}/\sqrt n\).
Using \(z_0=z\) and
\eqref{eq:influence-moment-derivatives}, we obtain
\begin{equation}\label{eq:influence-normal-argument-derivative}
 \begin{aligned}
 z'_0
 &=-\sqrt n\,m'_0-\frac{t}{2\sqrt n}\,v'_0\\
 &=-\sqrt n\,y-\frac z2(y^2-1).
 \end{aligned}
\end{equation}
Since \(\Phi'=\phi\), equations
\eqref{eq:influence-convolution-derivative} and
\eqref{eq:influence-normal-argument-derivative} give
\begin{equation}\label{eq:influence-discrepancy-derivative}
 \begin{aligned}
 D'_0
 &=A'_0-\phi(z)z'_0\\
 &=n\{G(t-y)-F_{n,P}(t)\}
   +\sqrt n\,y\phi(z)
   +\frac{z\phi(z)}2(y^2-1).
 \end{aligned}
\end{equation}

Both \(D_e\) and \(\beta_e\) have right derivatives at zero,
and \(\beta_0=\beta>0\). Hence the right derivative
\(I(y):=r'(0+)\) exists. Differentiating
\eqref{eq:influence-perturbed-discrepancy} and using
\(\sqrt nD_0/\beta=R\), we find
\begin{equation}\label{eq:influence-quotient-derivative}
 I(y)
 =\frac{\sqrt n}{\beta}D'_0
  -\frac{\sqrt nD_0}{\beta^2}\beta'_0,
 \qquad
 \beta I(y)=\sqrt nD'_0-R\beta'_0.
\end{equation}
Substituting \eqref{eq:influence-beta-derivative} and
\eqref{eq:influence-discrepancy-derivative} into
\eqref{eq:influence-quotient-derivative} gives
\eqref{eq:influence}.
Moreover, \eqref{eq:influence-maximality} implies
\[
 I(y)=\lim_{e\downarrow0}\frac{r(e)-r(0)}e\le0.
\]
Since \(y\in\R\) was arbitrary, this inequality holds
at every \(y\).

\textbf{Average the derivative under the maximizing law.}
For fixed \(n,P,t\), formula \eqref{eq:influence} gives
\[
 |I(y)|\le C_{n,P,t}(1+|y|^3),\qquad y\in\R,
\]
for a finite constant \(C_{n,P,t}\).
Indeed, the distribution-function values lie in \([0,1]\),
and the remaining terms grow at most cubically in \(|y|\).
Thus \(I\) is integrable under \(P\).

By the definition of \(G=F_{n-1,P}\) and the convolution
identity, Tonelli's theorem gives
\begin{equation}\label{eq:influence-convolution-average}
 \begin{aligned}
 \int G(t-y)\,P(dy)
 &=\iint \mathbf1_{\{s+y\le t\}}\,
          P^{*(n-1)}(ds)\,P(dy)\\
 &=P^{*n}((-\infty,t])=F_{n,P}(t).
 \end{aligned}
\end{equation}
The moment conditions on \(P\) also give
\[
 \int y\,P(dy)=0,\qquad
 \int(y^2-1)\,P(dy)=0,
\]
and
\[
 \int\{|y|^3-\beta-3yM\}\,P(dy)
 =\beta-\beta-3M\cdot0=0.
\]
Integrating \eqref{eq:influence} and using these identities
and \eqref{eq:influence-convolution-average}, we obtain
\begin{equation}\label{eq:influence-zero-average}
 \int I(y)\,P(dy)=0.
\end{equation}
Since \(-I\ge0\), its zero integral implies
\(I=0\) \(P\)-almost surely.

\textbf{Extend the identity to every support point.}
We first verify upper semicontinuity.
If \(y_k\to y\), then for every \(\delta>0\) we have
\(t-y_k\le t-y+\delta\) for all sufficiently large \(k\).
Monotonicity of \(G\) therefore gives
\[
 \limsup_{k\to\infty}G(t-y_k)\le G(t-y+\delta).
\]
Letting \(\delta\downarrow0\) and using the right continuity
of \(G\), we obtain
\[
 \limsup_{k\to\infty}G(t-y_k)\le G(t-y).
\]
In \eqref{eq:influence}, the coefficient of \(G(t-y)\)
is positive and all other terms are continuous in \(y\).
Since \(\beta>0\), it follows that
\[
 \limsup_{k\to\infty}I(y_k)\le I(y).
\]
Thus \(I\) is upper semicontinuous.

Because \(I\le0\) everywhere, its zero set can be written as
\[
 E:=\{y:I(y)=0\}=\{y:I(y)\ge0\}.
\]
Upper semicontinuity makes \(E\) closed.
The almost-sure equality proved above gives \(P(E)=1\).
If a point \(x\in\supp P\) did not belong to \(E\),
the open set \(\R\setminus E\) would contain an open interval
around \(x\). That interval would have zero \(P\)-measure
because \(P(E)=1\), contradicting the defining property
of a support point. Therefore \(I(x)=0\) for every
\(x\in\supp P\), which completes \eqref{eq:contact}.

Finally, substitute \(I(x)=0\) into \eqref{eq:influence}
and divide by \(\sqrt n\):
\[
 \begin{aligned}
 0={}&n\{G(t-x)-F_{n,P}(t)\}
      +\sqrt n\,\phi(z)x
      +\frac{z\phi(z)}2(x^2-1)\\
    &+\frac{3R\beta}{2\sqrt n}(x^2-1)
      -\frac{R}{\sqrt n}
       \{|x|^3-\beta-3Mx\}.
 \end{aligned}
\]
Rearranging this equality gives
\eqref{eq:contact-equation}.
\end{proof}

We next estimate the Gaussian terms that arise when \(G(t-y)\)
in \eqref{eq:influence} is compared with
\(\Phi((t-y)/\sqrt{n-1})\).
With \(z=t/\sqrt n\), these terms form \(H_n(z,y)\), defined
below. The global bound will be used in
Proposition~\ref{prop:bounded-extremizers} to obtain a common
bounded interval containing the supports of the selected
maximizing laws.
The expansion for bounded \(y\) will be used in
Lemma~\ref{lem:support-interval} to locate limits of support
points.

\begin{lemma}[Uniform Gaussian cancellation]
\label{lem:gaussian-expansion}
For \(n\ge2\) and \(y,z\in\R\), define
\begin{align}
H_n(z,y):={}&n^{3/2}
\left\{\Phi\left(\frac{z-y/\sqrt n}{\sqrt{1-1/n}}\right)-\Phi(z)\right\}
+ny\phi(z)\notag\\
&+\frac{\sqrt n}{2}z\phi(z)(y^2-1).\label{eq:gaussian-H}
\end{align}
There is an absolute constant \(A\) such that
\begin{equation}
\label{eq:gaussian-global}
|H_n(z,y)|\le\frac{\phi(0)}6|y|^3+A(1+|y|)
\qquad(n\ge2,\ y,z\in\R).
\end{equation}
For every \(L<\infty\), there is a constant \(A_L<\infty\) such that
\begin{equation}
\label{eq:gaussian-remainder}
\left|H_n(z,y)-\phi''(z)\left(\frac y2-\frac{y^3}{6}\right)\right|
\le\frac{A_L}{\sqrt n}
\qquad(n\ge2,\ |y|\le L,\ z\in\R).
\end{equation}
Consequently, if \(n_j\to\infty\) and \(z_j\to0\), then
\begin{equation}
\label{eq:gaussian-compact}
H_{n_j}(z_j,y)\longrightarrow\frac{\phi(0)}6(y^3-3y)
\quad\text{uniformly for }|y|\le L.
\end{equation}
\end{lemma}

\begin{proof}
Put \(u=z-y/\sqrt n\).
We separate the change from \(z\) to
\(u/\sqrt{1-1/n}\) into the translation from \(z\) to \(u\)
and the change of scale from \(u\) to \(u/\sqrt{1-1/n}\).

\textbf{Bound the translation remainder.}
The derivatives of the normal density satisfy
\[
 \phi'(x)=-x\phi(x),\qquad
 \phi''(x)=(x^2-1)\phi(x),\qquad
 \phi'''(x)=x(3-x^2)\phi(x).
\]
Thus the only critical points of \(\phi''\) are
\(0,\sqrt3,-\sqrt3\), and \(\phi''(x)\to0\) as
\(|x|\to\infty\).
At these points, its absolute values are respectively
\(\phi(0)\) and \(2e^{-3/2}\phi(0)<\phi(0)\).
Consequently,
\begin{equation}\label{eq:gaussian-derivative-bounds}
 \|\phi''\|_\infty=\phi(0),
 \qquad M_3:=\|\phi'''\|_\infty<\infty.
\end{equation}
The finiteness of \(M_3\) follows from the displayed formula
for \(\phi'''\), since a polynomial times \(\phi\) is bounded.

Since \(\Phi'''=\phi''\), the second-order Taylor formula
at \(z\), with increment \(u-z=-y/\sqrt n\), gives
\[
 \Phi(u)=\Phi(z)-\frac{y}{\sqrt n}\phi(z)
         +\frac{y^2}{2n}\phi'(z)+R_{2,n}(z,y),
 \qquad
 |R_{2,n}(z,y)|
 \le\frac{\phi(0)|y|^3}{6n^{3/2}}.
\]
After multiplying by \(n^{3/2}\), we cancel the terms
linear and quadratic in \(y\) and write
\begin{equation}\label{eq:gaussian-translation-term}
 T_n(z,y):=
 n^{3/2}\{\Phi(u)-\Phi(z)\}
 +ny\phi(z)+\frac{\sqrt n}{2}zy^2\phi(z).
\end{equation}
Using \(\phi'(z)=-z\phi(z)\) in the Taylor formula shows that
\(T_n(z,y)=n^{3/2}R_{2,n}(z,y)\). Hence
\begin{equation}\label{eq:translation-bound}
 |T_n(z,y)|\le\frac{\phi(0)}6|y|^3
 \qquad(n\ge2,\ y,z\in\R).
\end{equation}

\textbf{Control the change of scale and obtain the global bound.}
For fixed \(u\in\R\), set
\[
 f_u(\lambda)=\Phi\left(\frac{u}{\sqrt{1-\lambda}}\right),
 \qquad
 r=\frac{u}{\sqrt{1-\lambda}},
 \qquad 0\le\lambda\le\frac12.
\]
Since \(dr/d\lambda=r/[2(1-\lambda)]\), differentiation gives
\begin{equation}\label{eq:gaussian-scale-derivatives}
 \begin{aligned}
 f_u'(\lambda)
 &=\frac{r\phi(r)}{2(1-\lambda)},\\
 f_u''(\lambda)
 &=\frac{[\phi(r)+r\phi'(r)]r}{4(1-\lambda)^2}
   +\frac{r\phi(r)}{2(1-\lambda)^2}\\
 &=\frac{r(3-r^2)\phi(r)}{4(1-\lambda)^2}
  =\frac{\phi'''(r)}{4(1-\lambda)^2}.
 \end{aligned}
\end{equation}
For \(0\le\lambda\le1/2\), the denominator
\(4(1-\lambda)^2\) is at least one. Therefore
\eqref{eq:gaussian-derivative-bounds} implies
\[
 |f_u''(\lambda)|\le M_3
 \qquad(u\in\R,\ 0\le\lambda\le1/2).
\]
As \(1/n\le1/2\), the first-order Taylor formula for \(f_u\)
at zero yields
\begin{equation}\label{eq:gaussian-scale-expansion}
 \Phi\left(\frac{u}{\sqrt{1-1/n}}\right)-\Phi(u)
 =\frac{u\phi(u)}{2n}+E_n(u),
 \qquad |E_n(u)|\le\frac{M_3}{2n^2}.
\end{equation}
This remainder bound is uniform over all \(u\in\R\).

Let \(g(x)=x\phi(x)\). Combining the definition
\eqref{eq:gaussian-H} with
\eqref{eq:gaussian-translation-term} and
\eqref{eq:gaussian-scale-expansion}, we obtain
\begin{equation}\label{eq:gaussian-decomposition}
 \begin{aligned}
 H_n(z,y)
 &=T_n(z,y)
   +n^{3/2}\left\{
     \Phi\left(\frac{u}{\sqrt{1-1/n}}\right)-\Phi(u)
    \right\}
   -\frac{\sqrt n}{2}z\phi(z)\\
 &=T_n(z,y)
   +\frac{\sqrt n}{2}\{g(u)-g(z)\}
   +n^{3/2}E_n(u).
 \end{aligned}
\end{equation}
The derivative
\begin{equation}\label{eq:gaussian-g-derivatives}
 g'(x)=(1-x^2)\phi(x)=-\phi''(x),
 \qquad g''(x)=-\phi'''(x)
\end{equation}
therefore satisfies \(\|g'\|_\infty=\phi(0)\).
The mean value theorem gives
\[
 \frac{\sqrt n}{2}|g(u)-g(z)|
 \le\frac{\sqrt n}{2}\phi(0)|u-z|
 =\frac{\phi(0)}2|y|.
\]
Applying this estimate, \eqref{eq:translation-bound}, and
\eqref{eq:gaussian-scale-expansion} to
\eqref{eq:gaussian-decomposition}, we find
\begin{equation}\label{eq:gaussian-global-components}
 |H_n(z,y)|
 \le\frac{\phi(0)}6|y|^3
     +\frac{\phi(0)}2|y|
     +\frac{M_3}{2\sqrt n}.
\end{equation}
In particular, \eqref{eq:gaussian-global} holds with the
absolute constant \(A=(\phi(0)+M_3)/2\).

\textbf{Retain the next terms for bounded \(y\).}
Fix \(0\le L<\infty\) and suppose \(|y|\le L\).
Since \(\Phi^{(4)}=\phi'''\), the third-order Taylor formula
at \(z\) gives
\[
 \begin{aligned}
 \Phi(u)
 &=\Phi(z)-\frac{y}{\sqrt n}\phi(z)
   +\frac{y^2}{2n}\phi'(z)
   -\frac{y^3}{6n^{3/2}}\phi''(z)+R_{3,n}(z,y),\\
 |R_{3,n}(z,y)|
 &\le\frac{M_3|y|^4}{24n^2}.
 \end{aligned}
\]
Substituting this formula into
\eqref{eq:gaussian-translation-term} cancels the terms of
degrees one and two in \(y\), and yields
\begin{equation}\label{eq:gaussian-translation-expansion}
 T_n(z,y)=-\frac{y^3}{6}\phi''(z)+E_{T,n}(z,y),
 \qquad
 |E_{T,n}(z,y)|\le\frac{M_3|y|^4}{24\sqrt n}.
\end{equation}
Here \(E_{T,n}=n^{3/2}R_{3,n}\).

By \eqref{eq:gaussian-g-derivatives},
\(\|g''\|_\infty=M_3\).
The first-order Taylor formula for \(g\) gives
\[
 g(u)=g(z)+(u-z)g'(z)+R_{g,n}(z,y),
 \qquad
 |R_{g,n}(z,y)|
 \le\frac{M_3}{2}|u-z|^2
 =\frac{M_3y^2}{2n}.
\]
Multiplying by \(\sqrt n/2\) and using
\(u-z=-y/\sqrt n\), we obtain
\begin{equation}\label{eq:gaussian-g-expansion}
 \frac{\sqrt n}{2}\{g(u)-g(z)\}
 =-\frac y2g'(z)+E_{g,n}(z,y),
 \qquad
 |E_{g,n}(z,y)|\le\frac{M_3y^2}{4\sqrt n}.
\end{equation}
Substituting \eqref{eq:gaussian-translation-expansion} and
\eqref{eq:gaussian-g-expansion} into
\eqref{eq:gaussian-decomposition}, and using
\(g'(z)=-\phi''(z)\), gives
\[
 \begin{aligned}
 H_n(z,y)
 &=\phi''(z)\left(\frac y2-\frac{y^3}{6}\right)
   +E_{T,n}(z,y)+E_{g,n}(z,y)+n^{3/2}E_n(u).
 \end{aligned}
\]
The three remainder bounds now imply
\begin{equation}\label{eq:gaussian-compact-error-bound}
 \left|
 H_n(z,y)-\phi''(z)\left(\frac y2-\frac{y^3}{6}\right)
 \right|
 \le
 \frac{M_3}{\sqrt n}
 \left(\frac{|y|^4}{24}+\frac{y^2}{4}+\frac12\right).
\end{equation}
All derivative bounds used above hold on the whole real line,
so \eqref{eq:gaussian-compact-error-bound} is uniform in \(z\).
Taking
\[
 A_L=M_3\left(\frac{L^4}{24}+\frac{L^2}{4}+\frac12\right)
\]
proves \eqref{eq:gaussian-remainder}.

\textbf{Pass to the limit when \(z_j\to0\).}
Since \(\phi''(0)=-\phi(0)\),
\[
 \phi''(0)\left(\frac y2-\frac{y^3}{6}\right)
 =\frac{\phi(0)}6(y^3-3y).
\]
For all sufficiently large \(j\), we have \(n_j\ge2\).
Applying \eqref{eq:gaussian-remainder} at \((n,z)=(n_j,z_j)\)
therefore gives
\[
 \begin{aligned}
 &\sup_{|y|\le L}
 \left|
 H_{n_j}(z_j,y)-\frac{\phi(0)}6(y^3-3y)
 \right|\\
 &\qquad\le
 \frac{A_L}{\sqrt{n_j}}
 +\left(\frac L2+\frac{L^3}{6}\right)
  |\phi''(z_j)-\phi''(0)|.
 \end{aligned}
\]
The first term tends to zero because \(n_j\to\infty\).
The second tends to zero because \(z_j\to0\) and
\(\phi''\) is continuous.
This proves the uniform convergence in
\eqref{eq:gaussian-compact}.
\end{proof}

We now combine the preceding lemmas to construct the sequence
of maximizing laws used in Section~\ref{sec:completion}.
Our goal is to place their entire supports in one fixed
bounded interval.
Comparing the convolution term in \eqref{eq:influence} with
its normal approximation, and applying
Lemma~\ref{lem:gaussian-expansion}, gives an upper bound
with a negative cubic term and an additional term involving
\(\beta n(C_{n-1}-C_n)_+\).
Lemma~\ref{lem:selection}, together with the moment bound
\(\beta<B_*\) from Lemma~\ref{lem:attainment}, makes this
additional term tend to zero along the selected sample sizes.
The upper bound is then negative for all sufficiently large
\(|y|\), uniformly along a tail of the sequence.
The support contact condition \eqref{eq:contact} excludes
such points from the supports.

\begin{proposition}[A uniformly bounded sequence of extremizers]
\label{prop:bounded-extremizers}
Suppose that \(C_n>\ce\) for arbitrarily large integers \(n\). Then
there are integers \(n_j\uparrow\infty\), laws \(P_j\in\Pthree\),
thresholds \(t_j\in\R\), and a finite constant \(L\) such that, for
every \(j\),
\begin{gather}
\label{eq:selected-extremizers}
\frac{\sqrt{n_j}}{\beta_j}
\left\{F_{n_j,P_j}(t_j)-\Phi(t_j/\sqrt{n_j})\right\}
=R_{n_j}(P_j)=C_{n_j}>\ce,
\qquad 1\le\beta_j:=\beta(P_j)<B_*,\\
\label{eq:selected-increments}
n_j(C_{n_j-1}-C_{n_j})_+\longrightarrow0,\\
\label{eq:bounded-extremizer-support}
\supp P_j\subset[-L,L].
\end{gather}
The constant \(L\) may be chosen absolute.
\end{proposition}

\begin{proof}
\textbf{Select sample sizes and maximizing pairs.}
Set \(a_n=C_n-\ce\).
By \eqref{eq:cn-limit}, \(a_n\to0\), and the hypothesis
of the proposition says that \(a_n>0\) for arbitrarily
large \(n\).
Lemma~\ref{lem:selection} therefore gives integers
\(n_j\uparrow\infty\) such that
\[
 C_{n_j}>\ce,\qquad
 n_j(C_{n_j-1}-C_{n_j})_+\longrightarrow0,
\]
because \(a_{n_j-1}-a_{n_j}=C_{n_j-1}-C_{n_j}\).
This proves \eqref{eq:selected-increments}.
Deleting any initial indices below two, we may assume
\(n_j\ge2\).

For each selected \(n_j\), apply Lemma~\ref{lem:attainment}
to choose \(P_j\in\Pthree\) and \(t_j\in\R\) satisfying
\eqref{eq:selected-extremizers}.
In particular, with \(\beta_j=\beta(P_j)\), we have
\(1\le\beta_j<B_*\).
Let \(I_j\) be the contamination derivative from
Lemma~\ref{lem:influence} for this maximizing pair.
It remains to obtain a common absolute bound on the supports.

\textbf{Compare the convolution term with its normal approximation.}
Fix a selected index and suppress the subscript \(j\).
Thus \(n=n_j\), \(P=P_j\), \(t=t_j\), and \(R=C_n\).
Retain the notation
\[
 \beta=\beta(P),\qquad z=\frac{t}{\sqrt n},\qquad
 M=\int x|x|\,P(dx),\qquad G=F_{n-1,P},
\]
and write \(I=I_j\).
For arbitrary \(y\in\R\), put
\begin{equation}\label{eq:bounded-normal-argument}
 w=\frac{t-y}{\sqrt{n-1}}
  =\frac{z-y/\sqrt n}{\sqrt{1-1/n}}.
\end{equation}
By the definition of \(R_{n-1}(P)\) in
\eqref{eq:raw-normalization} and the supremum in
\eqref{eq:cn},
\[
 \frac{\sqrt{n-1}}{\beta}
 \{G(t-y)-\Phi(w)\}
 \le R_{n-1}(P)\le C_{n-1}.
\]
Consequently,
\begin{equation}\label{eq:bounded-normal-comparison}
 G(t-y)-\Phi(w)
 \le\frac{\beta C_{n-1}}{\sqrt{n-1}}.
\end{equation}
The maximizing identity \eqref{eq:selected-extremizers} gives
\[
 F_{n,P}(t)-\Phi(z)=\frac{R\beta}{\sqrt n}.
\]
Subtracting this identity from
\eqref{eq:bounded-normal-comparison}, with the Gaussian
difference kept explicit, yields
\begin{equation}\label{eq:bounded-convolution-comparison}
 \begin{aligned}
 G(t-y)-F_{n,P}(t)
 &=[G(t-y)-\Phi(w)]+[\Phi(w)-\Phi(z)]\\
 &\quad-[F_{n,P}(t)-\Phi(z)]\\
 &\le\Phi(w)-\Phi(z)
     +\frac{\beta C_{n-1}}{\sqrt{n-1}}
     -\frac{R\beta}{\sqrt n}.
 \end{aligned}
\end{equation}
Multiply \eqref{eq:bounded-convolution-comparison} by
\(n^{3/2}\) and substitute it into \eqref{eq:influence}.
By \eqref{eq:bounded-normal-argument}, the Gaussian terms
are exactly \(H_n(z,y)\) from \eqref{eq:gaussian-H}.
Thus
\begin{equation}\label{eq:bounded-influence-comparison}
 \begin{aligned}
 \beta I(y)\le{}&
 H_n(z,y)
 +\beta\left\{
   \frac{n^{3/2}}{\sqrt{n-1}}C_{n-1}-nR
  \right\}\\
 &+\frac32R\beta(y^2-1)
  -R\{|y|^3-\beta-3yM\}.
 \end{aligned}
\end{equation}

\textbf{Bound the increment term and the remaining polynomial terms.}
For \(n\ge2\), rationalizing the square-root difference gives
\begin{equation}\label{eq:bounded-scale-factor}
 \begin{aligned}
 0\le n\left(\sqrt{\frac n{n-1}}-1\right)
 &=\frac{n}{(n-1)(\sqrt{n/(n-1)}+1)}\\
 &\le\frac{n}{2(n-1)}\le1.
 \end{aligned}
\end{equation}
Recall \(C_0=0.4690\).
By \eqref{eq:cn} and \eqref{eq:cn-limit},
\(0\le C_{n-1}\le C_0\) and \(0<R\le C_0\).
Using \eqref{eq:bounded-scale-factor}, we therefore obtain
\begin{equation}\label{eq:bounded-extremal-increment}
 \begin{aligned}
 \frac{n^{3/2}}{\sqrt{n-1}}C_{n-1}-nR
 &=n(C_{n-1}-R)
   +n\left(\sqrt{\frac n{n-1}}-1\right)C_{n-1}\\
 &\le n(C_{n-1}-R)_++C_0.
 \end{aligned}
\end{equation}
Also,
\[
 |M|\le\int |x|^2\,P(dx)=1,
 \qquad 1\le\beta<B_*.
\]
Apply the global Gaussian bound \eqref{eq:gaussian-global}
and \eqref{eq:bounded-extremal-increment} to
\eqref{eq:bounded-influence-comparison}.
Expanding the last two terms in that inequality gives
\begin{equation}\label{eq:bounded-polynomial-terms}
 \begin{aligned}
 \beta I(y)\le{}&
 -\left(R-\frac{\phi(0)}6\right)|y|^3
 +\beta n(C_{n-1}-R)_+\\
 &+A(1+|y|)+\beta C_0
   +\frac32R\beta y^2-\frac12R\beta+3RMy\\
 \le{}&
 -\left(R-\frac{\phi(0)}6\right)|y|^3
 +\beta n(C_{n-1}-R)_+\\
 &+A(1+|y|)+C_0B_*
   +\frac32C_0B_*y^2+3C_0|y|.
 \end{aligned}
\end{equation}
The second inequality uses \(|M|\le1\), \(\beta<B_*\),
and \(0<R\le C_0\), and discards the nonpositive term
\(-R\beta/2\).

Since \(1\le1+y^2\) and \(|y|\le(1+y^2)/2\), we have
\[
 \begin{aligned}
 A(1+|y|)&\le\frac32A(1+y^2),\\
 C_0B_*+\frac32C_0B_*y^2
 &\le\frac32C_0B_*(1+y^2),\\
 3C_0|y|&\le\frac32C_0(1+y^2).
 \end{aligned}
\]
We may therefore take
\begin{equation}\label{eq:bounded-polynomial-constant}
 A_1=\frac32(A+C_0+C_0B_*).
\end{equation}
Here \(A\) is the absolute constant in
\eqref{eq:gaussian-global}; hence \(A_1\) is independent of
\(n,P,t,y\) and of the selected sequence.
Equation \eqref{eq:bounded-polynomial-terms} now gives
\begin{equation}\label{eq:influence-bound}
 \beta I(y)\le
 -\left(R-\frac{\phi(0)}6\right)|y|^3
 +A_1(1+y^2)+\beta n(C_{n-1}-R)_+
 \qquad(y\in\R).
\end{equation}

\textbf{Use the selected increments to exclude distant support points.}
Restore the subscript \(j\), and write \(R_j=C_{n_j}\).
By \eqref{eq:ce}, the constant
\begin{equation}\label{eq:bounded-cubic-gap}
 \gamma:=\ce-\frac{\phi(0)}6
 =\frac{(2+\sqrt{10})\phi(0)}6
 >0
\end{equation}
is absolute.
Since \(R_j>\ce\), we have
\(R_j-\phi(0)/6\ge\gamma\).
Moreover, \eqref{eq:selected-extremizers} and
\eqref{eq:selected-increments} imply
\begin{equation}\label{eq:bounded-uniform-error}
 \begin{aligned}
 \varepsilon_j
 &:=\beta_j n_j(C_{n_j-1}-C_{n_j})_+,\\
 0\le\varepsilon_j
 &\le B_* n_j(C_{n_j-1}-C_{n_j})_+
 \longrightarrow0.
 \end{aligned}
\end{equation}
Applying \eqref{eq:influence-bound} to each selected pair
therefore gives
\begin{equation}\label{eq:bounded-selected-influence}
 \beta_jI_j(y)
 \le-\gamma|y|^3+A_1(1+y^2)+\varepsilon_j,
 \qquad y\in\R.
\end{equation}
The error \(\varepsilon_j\) does not depend on \(y\).
Thus its convergence to zero controls the bound
simultaneously for all \(y\in\R\).

Choose
\begin{equation}\label{eq:bounded-support-radius}
 L=\max\left\{1,\frac{2A_1+1}{\gamma}\right\}.
\end{equation}
This is an absolute constant because \(A_1\) and \(\gamma\)
are absolute.
For every \(r>L\), we have \(r>1\) and
\(\gamma r>2A_1+1\), so
\begin{equation}\label{eq:bounded-cubic-domination}
 A_1(1+r^2)+1
 \le(2A_1+1)r^2
 <\gamma r^3.
\end{equation}
By \eqref{eq:bounded-uniform-error}, there is an index
\(j_0\) such that \(\varepsilon_j\le1\) for every \(j\ge j_0\).
For all such \(j\), equations
\eqref{eq:bounded-selected-influence} and
\eqref{eq:bounded-cubic-domination} imply
\[
 \beta_jI_j(y)
 \le-\gamma|y|^3+A_1(1+y^2)+1<0
 \qquad(|y|>L).
\]
Since \(\beta_j>0\), this gives \(I_j(y)<0\) whenever
\(|y|>L\).
On the other hand, the support contact condition
\eqref{eq:contact} gives \(I_j(x)=0\) at every
\(x\in\supp P_j\).
No support point can therefore satisfy \(|x|>L\), and
\[
 \supp P_j\subset[-L,L]\qquad(j\ge j_0).
\]
Discard the finitely many terms with \(j<j_0\) and reindex
the remaining sequence.
The maximizing identities \eqref{eq:selected-extremizers}
continue to hold, and taking a tail preserves both
\(n_j\uparrow\infty\) and \eqref{eq:selected-increments}.
The support inclusion now holds for every index of the
reindexed sequence, proving
\eqref{eq:bounded-extremizer-support}.
\end{proof}

% ===== completion =====
\section{Separation of support points and proof of the main theorem}\label{sec:completion}

In this section we work with the maximizing pairs supplied by
Proposition~\ref{prop:bounded-extremizers}.
Thus \(P_j\in\Pthree\), \(n_j\uparrow\infty\), and, with
\[
 \beta_j=\beta(P_j),\qquad
 z_j=\frac{t_j}{\sqrt{n_j}},\qquad
 M_j=\int x|x|\,P_j(dx),
\]
we have
\begin{equation}\label{eq:selected-sequence}
\begin{gathered}
 \frac{\sqrt{n_j}}{\beta_j}
 \{F_{n_j,P_j}(t_j)-\Phi(z_j)\}
 =R_{n_j}(P_j)=C_{n_j}>\ce,\\
 n_j(C_{n_j-1}-C_{n_j})_+\longrightarrow0,\\
 \supp P_j\subset[-L,L],\qquad 1\le\beta_j<B_*,
\end{gathered}
\end{equation}
where \(L<\infty\) is fixed.
In particular, the signed discrepancy at each \(t_j\)
is positive and attains the maximal normalized error.

We first identify a subsequential limit of these laws.
The common support bound gives a convergent subsequence,
and the upper bound from Corollary~\ref{cor:jitter-envelopes}
forces equality in Esseen's moment inequality.
Its equality characterization leaves two possible limits:
\(\Pe\) and its reflection.
We then use the positive signed maximizing identity in
\eqref{eq:selected-sequence} to exclude the reflected law
and show that \(z_j\to0\).
These limits provide the moment convergence and threshold
location used in Lemmas~\ref{lem:support-interval}
and~\ref{lem:support-separation}.

\begin{lemma}\label{lem:limit-extremizer}
After passing to a subsequence in \eqref{eq:selected-sequence},
\[
 P_j\longrightarrow\Pe\quad\text{in }W_3,
 \qquad z_j\longrightarrow0.
\]
\end{lemma}

\begin{proof}
\textbf{Pass to a limit and record convergence of the moments.}
The common support bound in \eqref{eq:selected-sequence}
implies that \((P_j)\) is tight.
The compactness criterion recalled in
Subsection~\ref{subsec:notation-convergence} therefore gives
a weakly convergent subsequence.
Retaining the same indices, write \(P_j\Longrightarrow P\).
Since \([-L,L]\) is closed, the Portmanteau theorem gives
\[
 P([-L,L])\ge\limsup_{j\to\infty}P_j([-L,L])=1.
\]
Thus \(P\) is also supported in \([-L,L]\).

To pass to the moments, define the bounded continuous function
\[
 \chi_L(x)=\max\{-L,\min\{x,L\}\}.
\]
For each of
\[
 f(x)=x,\quad x^2,\quad |x|^3,\quad x^3,\quad x|x|,
\]
the function \(f\circ\chi_L\) is bounded and continuous
on \(\R\), and agrees with \(f\) on the supports of
\(P_j\) and \(P\).
Weak convergence therefore gives
\[
 \int f(x)\,P_j(dx)
 =\int f(\chi_L(x))\,P_j(dx)
 \longrightarrow
 \int f(\chi_L(x))\,P(dx)
 =\int f(x)\,P(dx).
\]
In particular,
\[
 \int x\,P(dx)=0,\qquad \int x^2\,P(dx)=1.
\]
Also, \(\int |x|^3P(dx)\le L^3<\infty\), so
\(P\in\Pthree\).
Set
\[
 \beta=\beta(P),\qquad
 \kappa=\kappa(P),\qquad
 h=h(P),\qquad M=\int x|x|\,P(dx),
 \qquad \mu_j=\int x^3\,P_j(dx).
\]
The moment convergences just proved give
\begin{equation}\label{eq:limit-moment-convergence}
 \beta_j\longrightarrow\beta,\qquad
 \mu_j\longrightarrow\kappa,\qquad
 M_j\longrightarrow M.
\end{equation}
Weak convergence together with convergence of the third
absolute moments implies
\begin{equation}\label{eq:limit-wasserstein-convergence}
 P_j\longrightarrow P\quad\text{in }W_3,
\end{equation}
by the characterization recalled in
Subsection~\ref{subsec:notation-convergence}.

\textbf{Use the signed upper bound to force equality in Esseen's inequality.}
Apply Corollary~\ref{cor:jitter-envelopes} to
\eqref{eq:limit-wasserstein-convergence} and \(n_j\to\infty\).
The distribution function \(F_j\) in that corollary is
for the sum divided by \(\sqrt{n_j}\), so here
\[
 F_j(z)=F_{n_j,P_j}(\sqrt{n_j}z).
\]
Its upper bound \eqref{eq:jitter-upper}, multiplied by
\(\sqrt{n_j}\), gives a sequence \(r_j\ge0\) with \(r_j\to0\)
such that
\begin{equation}\label{eq:limit-row-envelope}
 \sqrt{n_j}
 \{F_{n_j,P_j}(\sqrt{n_j}z)-\Phi(z)\}
 \le
 \phi(z)\left[\frac h2+\frac{\mu_j}{6}(1-z^2)\right]+r_j,
 \qquad z\in\R.
\end{equation}
Here the same \(r_j\) works for every \(z\), because the
remainder in \eqref{eq:jitter-upper} is uniform.

We replace \(\mu_j\) by its limit \(\kappa\).
The Gaussian bounds established in
\eqref{eq:jitter-envelope-gaussian-suprema} are
\[
 \phi(z)\le\phi_0,\qquad
 |1-z^2|\phi(z)\le\phi_0,\qquad z\in\R.
\]
Consequently,
\[
 \sup_{z\in\R}
 \left|\frac{\mu_j-\kappa}{6}(1-z^2)\phi(z)\right|
 \le\frac{\phi_0}{6}|\mu_j-\kappa|
 \longrightarrow0
\]
by \eqref{eq:limit-moment-convergence}.
Define
\begin{equation}\label{eq:limit-envelope-error}
 \varepsilon_j:=r_j+\frac{\phi_0}{6}|\mu_j-\kappa|
 \longrightarrow0.
\end{equation}
Then \eqref{eq:limit-row-envelope} implies
\begin{equation}\label{eq:limit-upper-envelope}
 \sqrt{n_j}
 \{F_{n_j,P_j}(\sqrt{n_j}z)-\Phi(z)\}
 \le
 \phi(z)\left[\frac h2+\frac{\kappa}{6}(1-z^2)\right]
 +\varepsilon_j,
 \qquad z\in\R.
\end{equation}
Because \(\varepsilon_j\) is independent of \(z\), this bound
can be evaluated at \(z=z_j\) without any prior bound on
the sequence \((z_j)\).
The positive signed maximizing identity in
\eqref{eq:selected-sequence} therefore gives
\begin{equation}\label{eq:limit-active-envelope}
 C_{n_j}\beta_j
 \le
 \phi(z_j)\left[\frac h2+\frac{\kappa}{6}(1-z_j^2)\right]
 +\varepsilon_j.
\end{equation}
Moreover, \eqref{eq:cn-limit} and
\eqref{eq:limit-moment-convergence} imply
\begin{equation}\label{eq:limit-maximal-error}
 C_{n_j}\beta_j\longrightarrow\ce\beta.
\end{equation}

Since \(h\ge0\), the Gaussian bounds above give
\[
 \begin{aligned}
 \phi(z)\left[\frac h2+\frac{\kappa}{6}(1-z^2)\right]
 &\le \frac h2\phi(z)
      +\frac{|\kappa|}{6}|1-z^2|\phi(z)\\
 &\le\frac{\phi_0}{6}(3h+|\kappa|).
 \end{aligned}
\]
Using this estimate in \eqref{eq:limit-active-envelope},
and taking the limit with
\eqref{eq:limit-envelope-error} and
\eqref{eq:limit-maximal-error}, yields
\[
 \ce\beta\le\frac{\phi_0}{6}(3h+|\kappa|).
\]
On the other hand, Esseen's moment inequality
\eqref{eq:esseen-moment} applies to \(P\in\Pthree\).
Together with the value of \(\ce\) in \eqref{eq:ce},
it gives
\begin{equation}\label{eq:limit-esseen-equality}
 \ce\beta
 \le\frac{\phi_0}{6}(3h+|\kappa|)
 \le\frac{\phi_0}{6}(3+\sqrt{10})\beta
 =\ce\beta.
\end{equation}
The endpoints coincide, so both inequalities are equalities.
In particular, equality holds in \eqref{eq:esseen-moment}.
The equality characterization stated with that inequality
therefore implies that \(P\) is either \(\Pe\) or its reflection.

\textbf{Use the positive signed discrepancy to exclude the reflected law.}
Reflection preserves the maximal lattice span and the
third absolute moment, and reverses the sign of the third
moment. Thus both candidates have
\[
 h=h_{\mathrm E},\qquad \beta=\be,
\]
whereas \(\kappa=\ke\) for \(\Pe\) and
\(\kappa=-\ke\) for its reflection.
Put
\[
 A=\frac{h_{\mathrm E}}2,\qquad B=\frac{\ke}{6}.
\]
The parameters in \eqref{eq:esseen-parameters} and
the identities in \eqref{eq:esseen-identities} give
\[
 \ke=(\qe-\pe)h_{\mathrm E},\qquad
 0<\qe-\pe=\sqrt{10}-3<1.
\]
Hence \(h_{\mathrm E}>\ke>0\), and therefore \(A>3B>0\).
The corresponding upper envelopes in
\eqref{eq:limit-upper-envelope} are
\begin{equation}\label{eq:limit-candidate-envelopes}
 E_+(z)=\phi(z)(A+B-Bz^2),\qquad
 E_-(z)=\phi(z)(A-B+Bz^2),
\end{equation}
for \(\Pe\) and its reflection, respectively.
Since \(\Pe\) is an equality case of
\eqref{eq:esseen-moment} and \(\ke>0\), we also have
\begin{equation}\label{eq:limit-envelope-height}
 \phi_0(A+B)
 =\frac{\phi_0}{6}(3h_{\mathrm E}+\ke)
 =\ce\be.
\end{equation}

To bound the envelope for the reflected law, differentiate
using \(\phi'(z)=-z\phi(z)\):
\[
 \begin{aligned}
 E_-'(z)
 &=-z\phi(z)(A-B+Bz^2)+2Bz\phi(z)\\
 &=z\phi(z)(3B-A-Bz^2).
 \end{aligned}
\]
For \(z>0\), this derivative is strictly negative because
\(A>3B\).
Also, \(E_-\) is even.
It follows that
\begin{equation}\label{eq:limit-reflected-envelope-gap}
 \sup_{z\in\R}E_-(z)
 =E_-(0)=\phi_0(A-B)
 <\phi_0(A+B)=\ce\be,
\end{equation}
where the strict inequality uses \(B>0\).
If \(P\) were the reflected law, then
\eqref{eq:limit-active-envelope} would give
\[
 C_{n_j}\beta_j
 \le E_-(z_j)+\varepsilon_j
 \le\phi_0(A-B)+\varepsilon_j.
\]
Taking the limit using \eqref{eq:limit-maximal-error},
\(\beta=\be\), and \eqref{eq:limit-envelope-error} would imply
\[
 \ce\be\le\phi_0(A-B),
\]
contradicting \eqref{eq:limit-reflected-envelope-gap}.
Therefore \(P=\Pe\).
Together with \eqref{eq:limit-wasserstein-convergence},
this proves \(P_j\to\Pe\) in \(W_3\).

\textbf{Locate the standardized maximizing thresholds.}
Now the upper envelope in
\eqref{eq:limit-active-envelope} is \(E_+\).
Since \(B>0\) and \(\phi(z)>0\),
\[
 E_+(z)
 =\phi(z)(A+B-Bz^2)
 \le(A+B)\phi(z)
 =E_+(0)e^{-z^2/2}.
\]
At \(z=0\), equality holds, whereas for \(z\ne0\)
the last expression is strictly less than \(E_+(0)\).
Thus \(E_+\) has its unique global maximum at zero.
More explicitly, for every \(\delta>0\),
\begin{equation}\label{eq:limit-positive-envelope-gap}
 \sup_{|z|\ge\delta}E_+(z)
 \le e^{-\delta^2/2}E_+(0)
 <E_+(0)=\ce\be.
\end{equation}
If \(z_j\) did not converge to zero, there would be
\(\delta>0\) and an infinite subsequence on which
\(|z_j|\ge\delta\).
Along that subsequence,
\eqref{eq:limit-active-envelope} and
\eqref{eq:limit-positive-envelope-gap} would give
\[
 C_{n_j}\beta_j
 \le e^{-\delta^2/2}\ce\be+\varepsilon_j.
\]
Taking the limit, using
\eqref{eq:limit-maximal-error} with \(\beta=\be\)
and \eqref{eq:limit-envelope-error}, yields
\[
 \ce\be\le e^{-\delta^2/2}\ce\be<\ce\be,
\]
a contradiction. Hence \(z_j\to0\).

Finally, substituting \(P=\Pe\) into
\eqref{eq:limit-moment-convergence} and using
\eqref{eq:esseen-identities}, we record
\begin{equation}\label{eq:limit-esseen-moments}
 \beta_j\longrightarrow\be,\qquad
 \mu_j\longrightarrow\ke,\qquad
 M_j\longrightarrow M_{\mathrm E}=\qe-\pe.
\end{equation}
\end{proof}

The convergence in Lemma~\ref{lem:limit-extremizer} does not
by itself exclude additional support points away from the
two Esseen atoms.
We now return to the support contact condition
\eqref{eq:contact} from Lemma~\ref{lem:influence}.
For the contamination derivative \(I_j\) associated with
\((P_j,t_j)\), it gives \(I_j(y_j)=0\) at every
point \(y_j\in\supp P_j\).
Combining this identity with the compact Gaussian expansion
in Lemma~\ref{lem:gaussian-expansion} and the selected increment
condition in \eqref{eq:selected-sequence}, we obtain a
polynomial inequality for every limit of support points.
Solving this inequality gives the location bound below.

\begin{lemma}\label{lem:support-interval}
Under \eqref{eq:selected-sequence} and the subsequence of Lemma~\ref{lem:limit-extremizer}, every convergent sequence \(y_j\in\supp(P_j)\), with limit \(y\), satisfies
\begin{equation}\label{eq:limiting-support-interval}
 y\in\left[-\frac{\sqrt{10}}2\aE,
                \frac{\sqrt{10}}2\bE\right].
\end{equation}
\end{lemma}
\begin{proof}
\textbf{Apply the contact identity at the selected support points.}
Fix \(y_j\in\supp P_j\) with \(y_j\to y\).
By \eqref{eq:selected-sequence}, \(|y_j|\le L\), and hence
\(|y|\le L\).
Since \(n_j\to\infty\), we may discard finitely many terms
and assume \(n_j\ge2\).
Let \(I_j\) be the contamination derivative from
Lemma~\ref{lem:influence} for the maximizing pair
\((P_j,t_j)\).
The support contact condition \eqref{eq:contact} gives
\begin{equation}\label{eq:support-point-contact}
 I_j(y_j)=0.
\end{equation}
Put
\[
 G_j=F_{n_j-1,P_j},\qquad
 w_j=\frac{t_j-y_j}{\sqrt{n_j-1}}
     =\frac{z_j-y_j/\sqrt{n_j}}{\sqrt{1-1/n_j}}.
\]
Since \(P_j\in\Pthree\), the definition of \(C_{n_j-1}\)
in \eqref{eq:cn} implies
\[
 \frac{\sqrt{n_j-1}}{\beta_j}
 \{G_j(t_j-y_j)-\Phi(w_j)\}
 \le R_{n_j-1}(P_j)\le C_{n_j-1}.
\]
Thus
\begin{equation}\label{eq:support-convolution-upper}
 G_j(t_j-y_j)-\Phi(w_j)
 \le\frac{C_{n_j-1}\beta_j}{\sqrt{n_j-1}}.
\end{equation}
The signed maximizing identity in \eqref{eq:selected-sequence}
also gives
\[
 F_{n_j,P_j}(t_j)-\Phi(z_j)
 =\frac{C_{n_j}\beta_j}{\sqrt{n_j}}.
\]
Subtracting this identity from
\eqref{eq:support-convolution-upper}, while retaining the
Gaussian difference, yields
\[
 \begin{aligned}
 G_j(t_j-y_j)-F_{n_j,P_j}(t_j)
 &=[G_j(t_j-y_j)-\Phi(w_j)]
   +[\Phi(w_j)-\Phi(z_j)]\\
 &\quad-[F_{n_j,P_j}(t_j)-\Phi(z_j)]\\
 &\le\Phi(w_j)-\Phi(z_j)
   +\frac{C_{n_j-1}\beta_j}{\sqrt{n_j-1}}
   -\frac{C_{n_j}\beta_j}{\sqrt{n_j}}.
 \end{aligned}
\]
Insert this bound into \eqref{eq:influence}, with
\(R=C_{n_j}\), and use \eqref{eq:support-point-contact}.
The Gaussian terms combine into \(H_{n_j}(z_j,y_j)\)
by \eqref{eq:gaussian-H}, so
\begin{equation}\label{eq:support-contact-upper}
 \begin{aligned}
 0=\beta_j I_j(y_j)\le{}&
 H_{n_j}(z_j,y_j)
 +\beta_j\left[
   \frac{n_j^{3/2}C_{n_j-1}}{\sqrt{n_j-1}}
   -n_jC_{n_j}\right]\\
 &+\frac32 C_{n_j}\beta_j(y_j^2-1)
   -C_{n_j}\{|y_j|^3-\beta_j-3y_jM_j\}.
 \end{aligned}
\end{equation}

\textbf{Evaluate the Gaussian term at the moving support points.}
Lemma~\ref{lem:limit-extremizer} gives \(z_j\to0\).
Together with \(n_j\to\infty\), this allows us to apply
the uniform convergence in \eqref{eq:gaussian-compact}.
Because \(|y_j|\le L\), we have
\[
 \begin{aligned}
 &\left|H_{n_j}(z_j,y_j)-\frac{\phi_0}{6}(y^3-3y)\right|\\
 &\quad\le
 \sup_{|u|\le L}
 \left|H_{n_j}(z_j,u)-\frac{\phi_0}{6}(u^3-3u)\right|\\
 &\qquad+\frac{\phi_0}{6}
       |(y_j^3-3y_j)-(y^3-3y)|
 \longrightarrow0.
 \end{aligned}
\]
The first term tends to zero by
\eqref{eq:gaussian-compact}; the second tends to zero
because \(y_j\to y\).
In the notation of \eqref{eq:gaussian-H}, this gives
\begin{equation}\label{eq:sharp-gaussian-limit}
 \begin{aligned}
 H_{n_j}(z_j,y_j)
 ={}&n_j^{3/2}[\Phi(w_j)-\Phi(z_j)]
      +n_jy_j\phi(z_j)\\
 &+\frac{\sqrt{n_j}}2z_j\phi(z_j)(y_j^2-1)\\
 &\longrightarrow\frac{\phi_0}{6}(y^3-3y).
 \end{aligned}
\end{equation}

\textbf{Retain the constant cancellation and control the selected increment.}
For \(n\ge2\), rationalizing the square-root difference gives
\[
 n\left(\sqrt{\frac n{n-1}}-1\right)
 =\frac{n}{(n-1)(\sqrt{n/(n-1)}+1)}
 \longrightarrow\frac12.
\]
Define
\[
 \theta_n
 :=n\left(\sqrt{\frac n{n-1}}-1\right)-\frac12,
 \qquad \theta_n\longrightarrow0.
\]
The identity
\[
 \frac{n_j^{3/2}}{\sqrt{n_j-1}}
 =n_j+\frac12+\theta_{n_j}
\]
therefore yields
\begin{equation}\label{eq:recursive-constant-expansion}
 \begin{aligned}
 &\beta_j\left[
   \frac{n_j^{3/2}C_{n_j-1}}{\sqrt{n_j-1}}
   -n_jC_{n_j}\right]\\
 &\qquad=
 \beta_j n_j(C_{n_j-1}-C_{n_j})
 +\frac12\beta_j C_{n_j-1}+\eta_j,
 \end{aligned}
\end{equation}
where
\[
 \eta_j:=\beta_j C_{n_j-1}\theta_{n_j}\longrightarrow0.
\]
Indeed, \eqref{eq:selected-sequence} gives \(\beta_j<B_*\),
and \eqref{eq:cn} and \eqref{eq:cn-limit} give
\(0\le C_{n_j-1}\le C_0=0.4690\), so
\[
 |\eta_j|\le B_*C_0|\theta_{n_j}|\longrightarrow0.
\]

The last two terms of \eqref{eq:support-contact-upper}
expand as
\[
 \begin{aligned}
 &\frac32 C_{n_j}\beta_j(y_j^2-1)
  -C_{n_j}\{|y_j|^3-\beta_j-3y_jM_j\}\\
 &\qquad=
 -C_{n_j}|y_j|^3+\frac32C_{n_j}\beta_j y_j^2
 +3C_{n_j}M_jy_j-\frac12C_{n_j}\beta_j.
 \end{aligned}
\]
Consequently, the constant
\(\beta_j C_{n_j-1}/2\) in
\eqref{eq:recursive-constant-expansion} combines with
\(-\beta_j C_{n_j}/2\) to give
\[
 \frac12\beta_j(C_{n_j-1}-C_{n_j})\longrightarrow0.
\]
Here \(C_{n_j-1}\to\ce\) and \(C_{n_j}\to\ce\) by
\eqref{eq:cn-limit}, since both \(n_j-1\) and \(n_j\)
tend to infinity, and \(\beta_j\) is bounded.
Put
\begin{equation}\label{eq:support-cancellation-error}
 r_j:=\eta_j+\frac12\beta_j(C_{n_j-1}-C_{n_j})
 \longrightarrow0.
\end{equation}
Substituting these identities into
\eqref{eq:support-contact-upper} gives
\begin{equation}\label{eq:support-contact-expanded}
 \begin{aligned}
 0\le{}&
 H_{n_j}(z_j,y_j)-C_{n_j}|y_j|^3
 +\frac32C_{n_j}\beta_j y_j^2
 +3C_{n_j}M_jy_j\\
 &+\beta_j n_j(C_{n_j-1}-C_{n_j})+r_j.
 \end{aligned}
\end{equation}

We control the remaining signed increment from above.
By \eqref{eq:selected-sequence},
\begin{equation}\label{eq:support-increment-upper}
 \begin{aligned}
 \beta_j n_j(C_{n_j-1}-C_{n_j})
 &\le\varepsilon_j
 :=\beta_j n_j(C_{n_j-1}-C_{n_j})_+,\\
 0\le\varepsilon_j
 &\le B_*n_j(C_{n_j-1}-C_{n_j})_+
 \longrightarrow0.
 \end{aligned}
\end{equation}
Thus \eqref{eq:support-contact-expanded} implies
\begin{equation}\label{eq:support-contact-prelimit}
 \begin{aligned}
 0\le{}&
 H_{n_j}(z_j,y_j)-C_{n_j}|y_j|^3
 +\frac32C_{n_j}\beta_j y_j^2
 +3C_{n_j}M_jy_j\\
 &+\varepsilon_j+r_j.
 \end{aligned}
\end{equation}
We now use \(y_j\to y\), \(C_{n_j}\to\ce\) from
\eqref{eq:cn-limit}, and
\[
 \beta_j\to\be,\qquad M_j\to M_{\mathrm E}=\qe-\pe
\]
from \eqref{eq:limit-esseen-moments}.
Together with \eqref{eq:sharp-gaussian-limit},
\eqref{eq:support-cancellation-error}, and
\eqref{eq:support-increment-upper}, these convergences
allow us to take the limit in
\eqref{eq:support-contact-prelimit}.
We obtain
\begin{equation}\label{eq:limiting-contact-inequality}
 0\le
 \frac{\phi_0}{6}(y^3-3y)-\ce|y|^3
 +\frac32\ce\be y^2+3\ce M_{\mathrm E}y,
 \qquad M_{\mathrm E}=\qe-\pe.
\end{equation}

\textbf{Solve the limiting polynomial inequality.}
Set \(c_*=3+\sqrt{10}=6\ce/\phi_0\), using \eqref{eq:ce}.
By \eqref{eq:esseen-identities},
\begin{equation}\label{eq:support-linear-cancellation}
 c_*M_{\mathrm E}
 =(3+\sqrt{10})(\sqrt{10}-3)=1.
\end{equation}
Multiply \eqref{eq:limiting-contact-inequality} by
\(6/\phi_0>0\) and use
\eqref{eq:support-linear-cancellation}:
\[
 \begin{aligned}
 0
 &\le y^3-3y-c_*|y|^3
      +\frac32c_*\be y^2+3c_*M_{\mathrm E}y\\
 &=y^3-c_*|y|^3+\frac32c_*\be y^2.
 \end{aligned}
\]
Rearranging gives
\begin{equation}\label{eq:contact-polynomial}
 c_*|y|^3-y^3-\frac32c_*\be y^2\le0.
\end{equation}

To express the endpoints in terms of the Esseen atoms,
the parameters in \eqref{eq:esseen-parameters} give
\[
 \pe^2+\qe^2
 =\frac{(4-\sqrt{10})^2+(\sqrt{10}-2)^2}{4}
 =10-3\sqrt{10}.
\]
Using \eqref{eq:esseen-identities}, it follows that
\begin{equation}\label{eq:support-third-moment-identity}
 c_*\be
 =\frac{(3+\sqrt{10})(10-3\sqrt{10})}{\se}
 =\frac{\sqrt{10}}{\se}.
\end{equation}
Also,
\begin{equation}\label{eq:support-parameter-identities}
 \begin{aligned}
 \frac3{c_*+1}
 &=\frac3{4+\sqrt{10}}
   =\frac{4-\sqrt{10}}2=\pe,\\
 \frac3{c_*-1}
 &=\frac3{2+\sqrt{10}}
   =\frac{\sqrt{10}-2}2=\qe.
 \end{aligned}
\end{equation}

If \(y>0\), then \(|y|^3=y^3\), and
\eqref{eq:contact-polynomial} becomes
\[
 y^2\left[(c_*-1)y-\frac32c_*\be\right]\le0.
\]
Dividing by \(y^2>0\) and then by \(c_*-1>0\) gives
\[
 y\le\frac{3c_*\be}{2(c_*-1)}
   =\frac{\sqrt{10}}{2\se}\frac3{c_*-1}
   =\frac{\sqrt{10}}2\,\frac{\qe}{\se}
   =\frac{\sqrt{10}}2\bE.
\]
These equalities use
\eqref{eq:support-third-moment-identity},
\eqref{eq:support-parameter-identities}, and
\(\bE=\qe/\se\) from \eqref{eq:esseen-law}.

If \(y<0\), put \(s=-y>0\).
Then \(|y|^3=s^3\), \(y^3=-s^3\), and \(y^2=s^2\),
so \eqref{eq:contact-polynomial} becomes
\[
 s^2\left[(c_*+1)s-\frac32c_*\be\right]\le0.
\]
Dividing by \(s^2>0\) and \(c_*+1>0\) gives
\[
 -y=s\le\frac{3c_*\be}{2(c_*+1)}
   =\frac{\sqrt{10}}{2\se}\frac3{c_*+1}
   =\frac{\sqrt{10}}2\,\frac{\pe}{\se}
   =\frac{\sqrt{10}}2\aE,
\]
using the same two identities and
\(\aE=\pe/\se\) from \eqref{eq:esseen-law}.
Equivalently, \(y\ge-\sqrt{10}\,\aE/2\).
Finally, \(y=0\) also belongs to the claimed interval,
because \(\aE,\bE>0\).
Combining the three cases proves
\eqref{eq:limiting-support-interval}.
\end{proof}

The interval in \eqref{eq:limiting-support-interval} is larger
than the convex hull of the two Esseen atoms.
We now show that any two distinct limits of support points
must be at least \(h_{\mathrm E}\) apart.
For support points \(x_j\to x\) and \(y_j\to y\), write
\(G_j=F_{n_j-1,P_j}\) and compare the increment
\[
 G_j(t_j-x_j)-G_j(t_j-y_j).
\]
If \(0<y-x<h_{\mathrm E}\), the smoothing expansion forces
this increment, multiplied by \(\sqrt{n_j-1}\), to tend
to zero.
The contact equations at \(x_j\) and \(y_j\), however,
give the strictly positive limit \(\phi_0(y-x)\).
This contradiction yields the separation property below.

\begin{lemma}[Separation of support limits]\label{lem:support-separation}
Along the subsequence of Lemma~\ref{lem:limit-extremizer}, suppose \(x_j,y_j\in\supp(P_j)\), \(x_j\to x\) and \(y_j\to y\). Then either \(x=y\) or
\[
 |x-y|\ge h_{\mathrm E}=1/\se.
\]
\end{lemma}

\begin{proof}
Put
\[
 m_j=n_j-1,\qquad G_j=F_{m_j,P_j},\qquad
 h=h_{\mathrm E}=\frac1{\se}>0.
\]
By Lemma~\ref{lem:limit-extremizer}, \(z_j\to0\).
Since \(n_j\to\infty\), we may discard finitely many terms
and assume \(n_j\ge2\) and \(|z_j|\le1\).

\textbf{Expand the contact equation uniformly over the support.}
For every \(v\in\supp P_j\), divide the contact equation
\eqref{eq:contact-equation} by \(n_j\), with \(R=C_{n_j}\).
This gives the exact identity
\begin{equation}\label{eq:separation-contact-cdf}
 G_j(t_j-v)
 =F_{n_j,P_j}(t_j)-\frac{\phi(z_j)v}{\sqrt{n_j}}
  +\mathcal R_j(v),
\end{equation}
where
\begin{equation}\label{eq:separation-contact-remainder}
 \begin{aligned}
 \mathcal R_j(v)
 :={}&-\frac{z_j\phi(z_j)}{2n_j}(v^2-1)\\
 &+\frac{C_{n_j}}{n_j^{3/2}}
 \left\{|v|^3-\beta_j-3M_jv
                -\frac32\beta_j(v^2-1)\right\}.
 \end{aligned}
\end{equation}
The common support and moment bounds in
\eqref{eq:selected-sequence}, together with \eqref{eq:cn-limit},
give
\[
 |v|\le L,\qquad \beta_j<B_*,
 \qquad C_{n_j}\le C_0=0.4690.
\]
Also,
\[
 |M_j|
 \le\int |x|^2\,P_j(dx)=1.
\]
Let
\[
 D_1:=\|\phi'\|_\infty
     =\sup_{z\in\R}|z\phi(z)|<\infty.
\]
The equality uses \(\phi'(z)=-z\phi(z)\), and finiteness
follows because \(z\phi(z)\) is bounded.
Applying these bounds to \eqref{eq:separation-contact-remainder},
we obtain
\begin{equation}\label{eq:separation-contact-remainder-bound}
 \begin{aligned}
 \sup_{v\in\supp P_j}|\mathcal R_j(v)|
 &\le\frac{D_1(L^2+1)}{2n_j}\\
 &\quad+\frac{C_0}{n_j^{3/2}}
 \left[L^3+B_*+3L+\frac32B_*(L^2+1)\right]\\
 &\le\frac{K_1}{n_j},
 \end{aligned}
\end{equation}
where the constant
\[
 K_1=\frac{D_1(L^2+1)}2
 +C_0\left[L^3+B_*+3L+\frac32B_*(L^2+1)\right]
\]
is independent of \(j\) and \(v\).
The last inequality uses \(n_j^{-3/2}\le n_j^{-1}\).

\textbf{Compare with the normal distribution at the same thresholds.}
For \(|v|\le L\), define
\[
 w_j(v)=\frac{t_j-v}{\sqrt{m_j}}
       =a_j\left(z_j-\frac{v}{\sqrt{n_j}}\right),
 \qquad
 a_j=\sqrt{\frac{n_j}{n_j-1}}.
\]
Equation \eqref{eq:bounded-scale-factor} gives
\(0\le a_j-1\le n_j^{-1}\).
Writing
\[
 e_j(v):=(a_j-1)\left(z_j-\frac{v}{\sqrt{n_j}}\right),
\]
we have
\[
 w_j(v)-z_j=-\frac{v}{\sqrt{n_j}}+e_j(v).
\]
Since \(|z_j|\le1\) and \(|v|\le L\),
\begin{equation}\label{eq:separation-normal-argument}
 |e_j(v)|\le\frac{1+L}{n_j},\qquad
 |w_j(v)-z_j|
 \le\frac{L}{\sqrt{n_j}}+\frac{1+L}{n_j}
 \le\frac{1+2L}{\sqrt{n_j}}.
\end{equation}
Taylor's theorem, using \(\Phi'=\phi\) and
\(\|\Phi''\|_\infty=D_1\), gives
\[
 |\Phi(w_j(v))-\Phi(z_j)-\phi(z_j)(w_j(v)-z_j)|
 \le\frac{D_1}{2}|w_j(v)-z_j|^2.
\]
Substituting \eqref{eq:separation-normal-argument}, and using
\(\phi(z_j)\le\phi_0\), therefore yields
\begin{equation}\label{eq:separation-normal-cdf}
 \Phi(w_j(v))
 =\Phi(z_j)-\frac{\phi(z_j)v}{\sqrt{n_j}}
  +\mathcal N_j(v),
 \qquad
 \sup_{|v|\le L}|\mathcal N_j(v)|\le\frac{K_2}{n_j},
\end{equation}
where
\[
 K_2=\phi_0(1+L)+\frac{D_1}{2}(1+2L)^2.
\]
Indeed, the error is bounded by
\(\phi_0|e_j(v)|+D_1|w_j(v)-z_j|^2/2\).

Now choose an arbitrary sequence \(v_j\in\supp P_j\), and put
\begin{equation}\label{eq:separation-contact-threshold}
 u_j=t_j-v_j,\qquad w_j=\frac{u_j}{\sqrt{m_j}}=w_j(v_j).
\end{equation}
Subtracting \eqref{eq:separation-normal-cdf} from
\eqref{eq:separation-contact-cdf} at \(v=v_j\) cancels their
terms linear in \(v_j\).
Using the signed maximizing identity in
\eqref{eq:selected-sequence}, we obtain
\[
 \begin{aligned}
 \sqrt{m_j}\{G_j(u_j)-\Phi(w_j)\}
 &=
 \sqrt{\frac{m_j}{n_j}}\,C_{n_j}\beta_j\\
 &\quad+\sqrt{m_j}\{\mathcal R_j(v_j)-\mathcal N_j(v_j)\}.
 \end{aligned}
\]
The last term has absolute value at most
\[
 (K_1+K_2)\frac{\sqrt{m_j}}{n_j}
 \le\frac{K_1+K_2}{\sqrt{n_j}}\longrightarrow0.
\]
Also, \(m_j/n_j=1-1/n_j\to1\),
\(C_{n_j}\to\ce\) by \eqref{eq:cn-limit}, and
\(\beta_j\to\be\) by \eqref{eq:limit-esseen-moments}.
Consequently,
\begin{equation}\label{eq:contact-saturation}
 \sqrt{m_j}\{G_j(u_j)-\Phi(w_j)\}
 \longrightarrow\ce\be.
\end{equation}
Moreover, \eqref{eq:separation-normal-argument} implies
\[
 |w_j|
 \le |z_j|+\frac{1+2L}{\sqrt{n_j}}
 \longrightarrow0.
\]
Neither conclusion requires convergence of the chosen
support points \(v_j\).

\textbf{Apply the smoothing expansion at the shifted contact threshold.}
Let \(X_{j,1},\ldots,X_{j,m_j}\) be independent with common
law \(P_j\), and write
\[
 T_j=\sum_{r=1}^{m_j}X_{j,r}.
\]
Let \(U\), independent of the row, be uniform on
\([-h/2,h/2]\), and define the distribution function
of the unnormalized smoothed sum by
\[
 J_j(a)=\PP\{T_j+U\le a\}.
\]
Lemma~\ref{lem:limit-extremizer} gives \(P_j\to\Pe\) in
\(W_3\), and \(m_j=n_j-1\to\infty\).
Since the limiting law has maximal lattice span \(h\),
Lemma~\ref{lem:jitter} applies with sample sizes \(m_j\).
Writing
\[
 \mu_j=\int x^3\,P_j(dx),
\]
its conclusion \eqref{eq:jitter-expansion} gives
\begin{equation}\label{eq:separation-jitter-expansion}
 J_j(\sqrt{m_j}z)
 =\Phi(z)+\frac{\mu_j}{6\sqrt{m_j}}(1-z^2)\phi(z)
   +\mathcal E_j(z),
 \qquad
 \sqrt{m_j}\|\mathcal E_j\|_\infty\longrightarrow0.
\end{equation}

Put
\[
 \alpha_j=\frac{h}{2\sqrt{m_j}},\qquad b_j=w_j+\alpha_j.
\]
Then \(\sqrt{m_j}b_j=u_j+h/2\), and \(b_j\to0\)
because \(w_j\to0\) and \(\alpha_j\to0\).
Evaluate \eqref{eq:separation-jitter-expansion} at \(z=b_j\),
subtract \(\Phi(w_j)\), and multiply by \(\sqrt{m_j}\):
\begin{equation}\label{eq:separation-shifted-jitter}
 \begin{aligned}
 &\sqrt{m_j}\{J_j(u_j+h/2)-\Phi(w_j)\}\\
 &\quad=
 \sqrt{m_j}\{\Phi(b_j)-\Phi(w_j)\}
 +\frac{\mu_j}{6}(1-b_j^2)\phi(b_j)
 +\sqrt{m_j}\mathcal E_j(b_j).
 \end{aligned}
\end{equation}
By the mean value theorem, there is \(\xi_j\) between
\(w_j\) and \(b_j\) such that
\[
 \sqrt{m_j}\{\Phi(b_j)-\Phi(w_j)\}
 =\sqrt{m_j}\alpha_j\phi(\xi_j)
 =\frac h2\phi(\xi_j)
 \longrightarrow\frac h2\phi_0.
\]
Here \(\xi_j\to0\), since both endpoints tend to zero.
Furthermore, \(\mu_j\to\ke\) by
\eqref{eq:limit-esseen-moments}, so
\[
 \frac{\mu_j}{6}(1-b_j^2)\phi(b_j)
 \longrightarrow\frac{\ke}{6}\phi_0.
\]
The last term of \eqref{eq:separation-shifted-jitter}
tends to zero by \eqref{eq:separation-jitter-expansion}.
Thus
\begin{equation}\label{eq:separation-smoothed-contact}
 \sqrt{m_j}\{J_j(u_j+h/2)-\Phi(w_j)\}
 \longrightarrow
 \phi_0\left(\frac h2+\frac{\ke}{6}\right)
 =\ce\be,
\end{equation}
where the final equality is \eqref{eq:limit-envelope-height}.
Subtracting \eqref{eq:contact-saturation} from
\eqref{eq:separation-smoothed-contact} gives
\begin{equation}\label{eq:jitter-saturation}
 \sqrt{m_j}\{J_j(u_j+h/2)-G_j(u_j)\}
 \longrightarrow0.
\end{equation}

\textbf{Convert the smoothing comparison into a bound on short increments.}
Since \(T_j\) has distribution function \(G_j\) and is
independent of \(U\), conditioning on \(U\) gives
\[
 \begin{aligned}
 J_j(u_j+h/2)
 &=\EE\,G_j(u_j+h/2-U)\\
 &=\frac1h\int_{-h/2}^{h/2}G_j(u_j+h/2-v)\,dv\\
 &=\frac1h\int_0^h G_j(u_j+s)\,ds.
 \end{aligned}
\]
The last equality uses the substitution \(s=h/2-v\).
Subtracting
\(G_j(u_j)=h^{-1}\int_0^hG_j(u_j)\,ds\) yields the exact identity
\begin{equation}\label{eq:jitter-contact-integral}
 J_j(u_j+h/2)-G_j(u_j)
 =\frac1h\int_0^h
       [G_j(u_j+s)-G_j(u_j)]\,ds.
\end{equation}
The integrand is nonnegative because \(s\ge0\) and
\(G_j\) is nondecreasing.

Let \(d_j\to d\in(0,h)\).
For all sufficiently large \(j\), we have \(0<d_j<h\).
Restricting the integral in \eqref{eq:jitter-contact-integral}
to \([d_j,h]\), and using monotonicity once more, gives
\[
 \begin{aligned}
 J_j(u_j+h/2)-G_j(u_j)
 &\ge\frac1h\int_{d_j}^h
       [G_j(u_j+s)-G_j(u_j)]\,ds\\
 &\ge\frac{h-d_j}{h}
       [G_j(u_j+d_j)-G_j(u_j)].
 \end{aligned}
\]
Since \(h-d_j>0\), it follows that
\[
 \begin{aligned}
 0
 &\le\sqrt{m_j}[G_j(u_j+d_j)-G_j(u_j)]\\
 &\le\frac{h}{h-d_j}
      \sqrt{m_j}[J_j(u_j+h/2)-G_j(u_j)]
 \longrightarrow0.
 \end{aligned}
\]
The convergence uses \eqref{eq:jitter-saturation} and
\(h/(h-d_j)\to h/(h-d)<\infty\).
We have therefore proved
\begin{equation}\label{eq:flat-contact-interval}
 \sqrt{m_j}[G_j(u_j+d_j)-G_j(u_j)]
 \longrightarrow0.
\end{equation}
This conclusion holds for every sequence
\(v_j\in\supp P_j\), with \(u_j=t_j-v_j\), and every
sequence \(d_j\to d\in(0,h)\).

\textbf{Rule out two distinct support limits less than \(h\) apart.}
Suppose that \(0<y-x<h\).
Set
\[
 v_j=y_j,\qquad u_j=t_j-y_j,\qquad d_j=y_j-x_j.
\]
Then \(d_j\to y-x\in(0,h)\), and
\(u_j+d_j=t_j-x_j\).
Applying \eqref{eq:flat-contact-interval} with these choices gives
\begin{equation}\label{eq:separation-zero-increment}
 \sqrt{m_j}[G_j(t_j-x_j)-G_j(t_j-y_j)]
 \longrightarrow0.
\end{equation}

On the other hand, subtract
\eqref{eq:separation-contact-cdf} at \(v=y_j\) from the
same identity at \(v=x_j\).
The terms \(F_{n_j,P_j}(t_j)\) cancel, leaving
\[
 \begin{aligned}
 &G_j(t_j-x_j)-G_j(t_j-y_j)\\
 &\quad=
 \frac{\phi(z_j)}{\sqrt{n_j}}(y_j-x_j)
 +\mathcal R_j(x_j)-\mathcal R_j(y_j).
 \end{aligned}
\]
After multiplication by \(\sqrt{m_j}\), this becomes
\begin{equation}\label{eq:separation-contact-difference}
 \begin{aligned}
 &\sqrt{m_j}[G_j(t_j-x_j)-G_j(t_j-y_j)]\\
 &\quad=
 \sqrt{\frac{m_j}{n_j}}\,\phi(z_j)(y_j-x_j)
 +\sqrt{m_j}\{\mathcal R_j(x_j)-\mathcal R_j(y_j)\}.
 \end{aligned}
\end{equation}
By \eqref{eq:separation-contact-remainder-bound},
the last term has absolute value at most
\[
 \frac{2K_1\sqrt{m_j}}{n_j}
 \le\frac{2K_1}{\sqrt{n_j}}\longrightarrow0.
\]
The first term tends to \(\phi_0(y-x)>0\), because
\(m_j/n_j\to1\), \(z_j\to0\), and \(y_j-x_j\to y-x\).
Thus \eqref{eq:separation-contact-difference} contradicts
\eqref{eq:separation-zero-increment}.
This excludes \(0<y-x<h\).
Interchanging the two support sequences excludes
\(0<x-y<h\).
Hence either \(x=y\) or \(|x-y|\ge h=h_{\mathrm E}\),
as asserted.
\end{proof}

We now combine the location bound in
Lemma~\ref{lem:support-interval} with the separation property
in Lemma~\ref{lem:support-separation}.
Together, they will show that the entire support of each
selected law approaches the two Esseen atoms.
After an affine change of coordinates, we verify the
support and probability conditions of
Proposition~\ref{prop:clusters} and obtain the final contradiction.

\begin{proof}[Proof of the existence assertion in Theorem~\ref{thm:main}]

\textbf{Choose a sequence of maximizing counterexamples.}
Suppose that no universal threshold exists.
Then, for every integer \(N\ge1\), there are an integer
\(n\ge N\) and a law \(P\in\Pthree\) such that
\(R_n(P)>\ce\).
By the definition \eqref{eq:cn}, this implies
\(C_n>\ce\) for arbitrarily large integers \(n\).
Proposition~\ref{prop:bounded-extremizers} therefore supplies
a sequence satisfying \eqref{eq:selected-sequence}.
Pass to the subsequence in Lemma~\ref{lem:limit-extremizer},
so that
\[
 P_j\longrightarrow\Pe\quad\text{in }W_3,
 \qquad z_j\longrightarrow0.
\]
The conditions in \eqref{eq:selected-sequence} and these
convergences remain valid after any further passage to a
subsequence.
Thus Lemmas~\ref{lem:support-interval}
and~\ref{lem:support-separation} also apply along any such
subsequence.

\textbf{Transform the two limiting atoms to \(0\) and \(1\).}
Let
\[
 T(x)=\pe+\se x,
\]
and let \(Q_j\) be the law of \(T(X_j)\) when \(X_j\sim P_j\).
The identities in \eqref{eq:esseen-law} give
\begin{equation}\label{eq:completion-affine-atoms}
 T(-\aE)=\pe-\se\aE=0,\qquad
 T(\bE)=\pe+\se\bE=1,\qquad
 \se h_{\mathrm E}=1.
\end{equation}
Because \(\se>0\), \(T\) is a continuous bijection with
continuous inverse.
For any open neighborhood \(V\) of \(T(x)\), its inverse
image \(T^{-1}(V)\) is an open neighborhood of \(x\), and
\(Q_j(V)=P_j(T^{-1}(V))\).
The definition of support therefore gives
\begin{equation}\label{eq:completion-affine-support}
 \supp Q_j=T(\supp P_j)
 \subset[\pe-\se L,\pe+\se L].
\end{equation}
In particular, these supports are nonempty compact sets
contained in one fixed bounded interval.

Consider any convergent sequence
\(v_{j_k}\in\supp Q_{j_k}\), with \(v_{j_k}\to v\).
By \eqref{eq:completion-affine-support}, the corresponding
points
\[
 x_k=\frac{v_{j_k}-\pe}{\se}\in\supp P_{j_k}
 \quad\text{satisfy}\quad
 x_k\longrightarrow\frac{v-\pe}{\se}.
\]
Lemma~\ref{lem:support-interval}, applied along these indices,
places the last limit in
\([-\sqrt{10}\,\aE/2,\sqrt{10}\,\bE/2]\).
Multiplication by \(\se>0\) and addition of \(\pe\)
therefore give
\begin{equation}\label{eq:raw-support-interval}
 v\in
 \left[\pe-\frac{\sqrt{10}}2\pe,
       \pe+\frac{\sqrt{10}}2\qe\right]
 =\left[-\pe\qe,\frac92-\sqrt{10}\right]
 \subset(-1,2).
\end{equation}
To check the endpoints, use
\(\pe=(4-\sqrt{10})/2\) and
\(\qe=(\sqrt{10}-2)/2\) from
\eqref{eq:esseen-parameters}:
\[
 \begin{aligned}
 \pe\left(1-\frac{\sqrt{10}}2\right)
 &=-\pe\qe,\\
 \pe+\frac{\sqrt{10}}2\qe
 &=\frac{4-\sqrt{10}}2
   +\frac{\sqrt{10}(\sqrt{10}-2)}4
 =\frac92-\sqrt{10}.
 \end{aligned}
\]
The lower endpoint is greater than \(-1\), since
\(0<\pe\qe\le1/4\).
The upper endpoint is less than \(2\), since
\(\sqrt{10}>5/2\).

The separation property also transforms explicitly.
Suppose
\[
 v_{j_k},\widetilde v_{j_k}\in\supp Q_{j_k},\qquad
 v_{j_k}\to v,\qquad \widetilde v_{j_k}\to\widetilde v.
\]
Their inverse images under \(T\) belong to \(\supp P_{j_k}\)
and converge to \((v-\pe)/\se\) and
\((\widetilde v-\pe)/\se\).
If \(v\ne\widetilde v\), Lemma~\ref{lem:support-separation}
and \eqref{eq:completion-affine-atoms} yield
\begin{equation}\label{eq:completion-transformed-separation}
 |v-\widetilde v|
 =\se\left|
       \frac{v-\pe}{\se}-\frac{\widetilde v-\pe}{\se}
      \right|
 \ge\se h_{\mathrm E}=1.
\end{equation}

\textbf{Find support points approaching both atoms and exclude other limits.}
Define
\[
 Q_{\mathrm E}=\qe\delta_0+\pe\delta_1.
\]
This is the image of \(\Pe\) under \(T\), by
\eqref{eq:completion-affine-atoms}.
For every bounded continuous function \(f\), the function
\(f\circ T\) is also bounded and continuous.
The weak convergence implied by \(P_j\to\Pe\) in \(W_3\)
therefore gives
\[
 \int f(v)\,Q_j(dv)
 =\int f(T(x))\,P_j(dx)
 \longrightarrow
 \int f(T(x))\,\Pe(dx)
 =\int f(v)\,Q_{\mathrm E}(dv).
\]
Hence
\begin{equation}\label{eq:completion-transformed-weak-limit}
 Q_j\Longrightarrow Q_{\mathrm E}.
\end{equation}

Fix \(a\in\{0,1\}\) and \(r\in(0,1/2)\).
The Portmanteau theorem, applied to the open interval
\((a-r,a+r)\), gives
\[
 \liminf_{j\to\infty}Q_j((a-r,a+r))
 \ge Q_{\mathrm E}((a-r,a+r))
 =Q_{\mathrm E}(\{a\})>0.
\]
Thus this interval has positive \(Q_j\)-probability for
all sufficiently large \(j\), and consequently meets
\(\supp Q_j\).
Since this holds for every such \(r\),
\[
 \dist(a,\supp Q_j)\longrightarrow0.
\]
Each support is nonempty and compact by
\eqref{eq:completion-affine-support}, so the distance from
\(a\) to \(\supp Q_j\) is attained.
Choose a point attaining it, for each \(a\) and \(j\).
This produces sequences
\begin{equation}\label{eq:completion-principal-support-points}
 v_j^{(0)},v_j^{(1)}\in\supp Q_j,\qquad
 v_j^{(0)}\longrightarrow0,\qquad
 v_j^{(1)}\longrightarrow1.
\end{equation}

Now let \(v_{j_k}\in\supp Q_{j_k}\) converge to \(v\)
along any subsequence.
By \eqref{eq:raw-support-interval}, \(v\in(-1,2)\).
Suppose \(v\notin\{0,1\}\).
There are three possibilities:
\[
 \begin{array}{lll}
 -1<v<0 &: &0<|v-0|=-v<1,\\
 0<v<1 &: &0<|v-0|=v<1,\\
 1<v<2 &: &0<|v-1|=v-1<1.
 \end{array}
\]
In the first two cases, compare \(v_{j_k}\) with
\(v_{j_k}^{(0)}\) from
\eqref{eq:completion-principal-support-points}.
Their distinct limits are \(v\) and \(0\), so
\eqref{eq:completion-transformed-separation} requires
\(|v-0|\ge1\), a contradiction.
In the third case, compare \(v_{j_k}\) with
\(v_{j_k}^{(1)}\).
The same separation property requires \(|v-1|\ge1\),
again a contradiction.
Therefore every such limit belongs to \(\{0,1\}\).

Returning to the original coordinates, if
\(x_k\in\supp P_{j_k}\) and \(x_k\to x\), then
\(T(x_k)\in\supp Q_{j_k}\) and \(T(x_k)\to T(x)\).
The preceding conclusion gives \(T(x)\in\{0,1\}\).
By \eqref{eq:completion-affine-atoms}, this means
\begin{equation}\label{eq:completion-only-support-limits}
 x_k\in\supp P_{j_k},\quad x_k\to x
 \quad\Longrightarrow\quad x\in\{-\aE,\bE\}.
\end{equation}

\textbf{Deduce uniform confinement of the entire support.}
We claim that
\begin{equation}\label{eq:uniform-confinement}
 \sup_{x\in\supp P_j}
       \dist(x,\{-\aE,\bE\})\longrightarrow0.
\end{equation}
If this failed, there would be \(\varepsilon>0\),
indices \(j_k\uparrow\infty\), and points
\(x_k\in\supp P_{j_k}\) such that
\[
 \dist(x_k,\{-\aE,\bE\})\ge\varepsilon
 \qquad\text{for every }k.
\]
The common support bound in \eqref{eq:selected-sequence}
places all \(x_k\) in \([-L,L]\).
After a further subsequence, \(x_k\) therefore converges
to some \(x\in[-L,L]\).
The distance function
\[
 x\longmapsto\dist(x,\{-\aE,\bE\})
 =\min\{|x+\aE|,|x-\bE|\}
\]
is continuous, so its value at this limit is at least
\(\varepsilon\).
This contradicts \eqref{eq:completion-only-support-limits}.
The contradiction proves \eqref{eq:uniform-confinement}.

By \eqref{eq:completion-affine-atoms},
\[
 \begin{aligned}
 \dist(T(x),\{0,1\})
 &=\min\{|\pe+\se x|,|\pe+\se x-1|\}\\
 &=\se\min\{|x+\aE|,|x-\bE|\}\\
 &=\se\dist(x,\{-\aE,\bE\}).
 \end{aligned}
\]
Together with \eqref{eq:completion-affine-support} and
\eqref{eq:uniform-confinement}, this gives
\begin{equation}\label{eq:completion-transformed-confinement}
 \sup_{v\in\supp Q_j}\dist(v,\{0,1\})
 =\se\sup_{x\in\supp P_j}\dist(x,\{-\aE,\bE\})
 \longrightarrow0.
\end{equation}

\textbf{Verify the local comparison conditions and obtain the contradiction.}
Fix the constants \(\eta\in(0,1/4)\) and \(N_{\mathrm{loc}}\)
provided by Proposition~\ref{prop:clusters}.
Equation \eqref{eq:completion-transformed-confinement}
implies that, for all sufficiently large \(j\),
\begin{equation}\label{eq:completion-cluster-support}
 \supp Q_j\subset[-\eta,\eta]\cup[1-\eta,1+\eta].
\end{equation}
Indeed, every point whose distance to \(\{0,1\}\) is at
most \(\eta\) lies in this union.

Let \(I=[1-\eta,1+\eta]\).
Since \(0<\eta<1/4\), this interval contains \(1\), excludes
\(0\), and its boundary consists of \(1-\eta\) and \(1+\eta\).
Neither boundary point is an atom of \(Q_{\mathrm E}\), so
\[
 Q_{\mathrm E}(\partial I)=0,\qquad
 Q_{\mathrm E}(I)=\pe.
\]
The continuity-set conclusion of the Portmanteau theorem,
applied to \eqref{eq:completion-transformed-weak-limit},
therefore gives
\begin{equation}\label{eq:completion-cluster-probability}
 Q_j([1-\eta,1+\eta])\longrightarrow\pe.
\end{equation}
In particular,
\(\bigl|Q_j([1-\eta,1+\eta])-\pe\bigr|<\eta\)
for all sufficiently large \(j\).
Also, \(n_j\ge N_{\mathrm{loc}}\) eventually, since
\(n_j\to\infty\).
Thus the support condition, probability condition, and
sample-size condition of Proposition~\ref{prop:clusters}
hold simultaneously along a tail of the sequence.

For completeness, we verify the normalization for \(Q_j\).
If \(X_j\sim P_j\) and \(Y_j=\pe+\se X_j\sim Q_j\), then
\[
 \EE Y_j=\pe,\qquad
 \operatorname{Var}(Y_j)=\se^2>0,\qquad
 \EE|Y_j-\pe|^3=\se^3\beta_j<\infty.
\]
Consequently, the definition \eqref{eq:raw-normalization}
gives
\[
 \beta(Q_j)=\frac{\se^3\beta_j}{(\se^2)^{3/2}}
 =\beta_j,
\]
where \(\se>0\).
For independent copies \(X_{j,r}\) with law \(P_j\), put
\(Y_{j,r}=\pe+\se X_{j,r}\).
Then the \(Y_{j,r}\) are independent with law \(Q_j\), and
\[
 \frac{\sum_{r=1}^{n_j}Y_{j,r}-n_j\pe}
      {\se\sqrt{n_j}}
 =\frac{\sum_{r=1}^{n_j}X_{j,r}}{\sqrt{n_j}}.
\]
The two standardized sums therefore have the same
distribution, so
\(\Delta_{n_j}(Q_j)=\Delta_{n_j}(P_j)\).
Together with \(\beta(Q_j)=\beta_j\), this proves
\begin{equation}\label{eq:completion-affine-error}
 R_{n_j}(Q_j)=R_{n_j}(P_j).
\end{equation}

Applying Proposition~\ref{prop:clusters} and then
\eqref{eq:completion-affine-error}, we obtain, for every
sufficiently large \(j\),
\[
 C_{n_j}=R_{n_j}(P_j)=R_{n_j}(Q_j)\le\ce.
\]
This contradicts \(C_{n_j}>\ce\) in
\eqref{eq:selected-sequence}.

Hence the set
\[
 \mathcal V=\{n\in\{1,2,\ldots\}:C_n>\ce\}
\]
is finite.
Choose
\[
 N=1+\max\bigl(\{0\}\cup\mathcal V\bigr).
\]
This definition also covers the case \(\mathcal V=\varnothing\).
For every integer \(n\ge N\) and every \(P\in\Pthree\),
\eqref{eq:cn} gives
\[
 R_n(P)\le C_n\le\ce.
\]
Equivalently, by \eqref{eq:raw-normalization},
\[
 \Delta_n(P)\le\ce\,\frac{\beta(P)}{\sqrt n}.
\]
The threshold \(N\) is independent of \(P\), since the
constants \(C_n\) already take the supremum over all
laws in \(\Pthree\).
This proves the existence assertion.
The explicit value in \eqref{eq:explicit-threshold}
is proved in Appendix~\ref{app:effective}.
\end{proof}

\section*{Acknowledgements}

We thank Prof. Lutz Mattner for his helpful suggestions on the title,
exposition, and references of this paper.
\newpage

% ===== explicit threshold =====
\appendix
\section{An explicit universal threshold}\label{app:effective}

We now prove the explicit threshold in
\eqref{eq:explicit-threshold}.
The argument starts from a violation at a single sample size
\(N\).
Lemma~\ref{lem:effective-selection} selects a violating
index \(n\in[\lceil N/2\rceil,N]\) with a controlled decrease
from \(C_{n-1}\) to \(C_n\), and shows that every
positive-discrepancy maximizing law at this index is
supported in \([-6,6]\).
Explicit smoothing estimates and stability of Esseen's moment
inequality then bound its \(W_1\) distance from \(\Pe\) in
Lemma~\ref{lem:effective-identification}.
After the affine transformation \(x\mapsto\pe+\se x\),
Proposition~\ref{prop:effective-confinement} places the
entire support near \(0\) and \(1\) and controls the cluster
probabilities.
The local comparison in
Proposition~\ref{prop:effective-clusters} then gives the
contradiction.

The variational inputs are the attainment result in
Lemma~\ref{lem:attainment} and the support contact identities
in Lemma~\ref{lem:influence}.
We combine these with explicit versions of the analytic
estimates used earlier.
The argument does not use the existence assertion proved
in Section~\ref{sec:completion}.

Recall \(c_*=3+\sqrt{10}\).
Throughout this appendix, set
\begin{equation}\label{eq:effective-threshold-parameters}
 \mathcal A=10^{14},\qquad
 \eta_*=e^{-2\mathcal A},\qquad
 N_*=\lceil e^{\mathcal A}\rceil,\qquad
 N_{\mathrm{conf}}=\lceil e^{1000\mathcal A}\rceil.
\end{equation}
The constants \(\eta_*\) and \(N_*\) specify the neighborhood
and sample-size threshold in the local comparison of
Proposition~\ref{prop:effective-clusters}.
The threshold \(N_{\mathrm{conf}}\) is used for the support
confinement in Proposition~\ref{prop:effective-confinement}.
Since selection may reduce the sample size from \(N\) to
\(\lceil N/2\rceil\), we start with \(N\ge2N_{\mathrm{conf}}\).
The resulting threshold is
\[
 2N_{\mathrm{conf}}
 =2\left\lceil\exp(10^{17})\right\rceil,
\]
because \(1000\mathcal A=10^{17}\).
The constants are chosen for convenience, without optimization.

For laws \(P,Q\) with finite first absolute moments, write
\[
 W_1(P,Q):=\inf\EE|X-Y|,
\]
where the infimum is over all couplings \((X,Y)\) with
marginal laws \(P\) and \(Q\).
We also use
\begin{equation}\label{eq:effective-increment}
 d_n:=n(C_{n-1}-C_n)_+,\qquad n\ge2.
\end{equation}

\subsection{Finite-sample selection and bounded support}

We begin with a finite-sample version of the selection and
support arguments in Section~\ref{sec:extremizers}.
Given \(C_N>\ce\), we will find a violating index between
\(\lceil N/2\rceil\) and \(N\) for which \(d_n\) is small.
The contact identity will then give the explicit support
bound \([-6,6]\).

We first record the numerical upper bounds for \(C_n\)
used in the selection argument.
For iid standardized summands, Shevtsova's estimate
\cite[Corollary~4.18, p.~303]{Shevtsova2012} gives
\[
 \Delta_n(P)
 \le\ce\frac{\beta(P)}{\sqrt n}
     +2.5786\frac{\beta(P)^2}{n}.
\]
This is the estimate used in \eqref{eq:asymptotic}
before enlarging the second coefficient to \(3\).
Multiplication by \(\sqrt n/\beta(P)\) gives
\[
 R_n(P)\le\ce+\frac{2.5786\,\beta(P)}{\sqrt n}.
\]
If \(R_n(P)>\ce\), the cutoff \eqref{eq:moment-cutoff}
implies \(\beta(P)<B_*<1.84\), and
\[
 2.5786\,\beta(P)<2.5786(1.84)=4.744624<4.75.
\]
If \(R_n(P)\le\ce\), the bound
\(R_n(P)\le\ce+4.75/\sqrt n\) already holds.
Taking the supremum over \(P\in\Pthree\), and using
\eqref{eq:cn-limit} for the first inequality below, yields
\begin{equation}\label{eq:effective-cn}
 C_n\le0.4690,\qquad
 C_n\le\ce+\frac{4.75}{\sqrt n},
 \qquad n\ge1.
\end{equation}

\begin{lemma}\label{lem:effective-selection}
If \(N\ge4\) is an integer and \(C_N>\ce\), there is an integer
\(n\in[\lceil N/2\rceil,N]\) such that \(C_n>\ce\) and
\[
 d_n<0.086,\qquad d_n\le\frac{10}{\sqrt{N-2}}.
\]
Every positive-discrepancy maximizing law at this index is
supported in \([-6,6]\).
\end{lemma}

\begin{proof}
\textbf{Select a violating index with a small preceding decrease.}
Set
\[
 m=\lceil N/2\rceil,\qquad a_j=C_j-\ce,\qquad
 H_N=\sum_{j=m}^{N}\frac1j,
\]
and put
\begin{equation}\label{eq:effective-selection-level}
 B_m=\min\left\{0.4690-\ce,\frac{4.75}{\sqrt{m-1}}\right\},
 \qquad D=\frac{B_m}{H_N}.
\end{equation}
Since \(N\ge4\), we have \(m\ge2\), so these quantities
are well defined.
Also, \(0.4690>\ce\), hence \(B_m>0\) and \(D>0\).
Equation \eqref{eq:effective-cn}, applied at \(m-1\), gives
\begin{equation}\label{eq:effective-selection-start-bound}
 a_{m-1}\le B_m.
\end{equation}

Suppose that no \(j\in\{m,\ldots,N\}\) with \(a_j>0\)
satisfies \(d_j\le D\).
Then each such positive term has \(d_j>D>0\).
Since \eqref{eq:effective-increment} gives
\[
 d_j=j(a_{j-1}-a_j)_+,
\]
we must have \(a_{j-1}-a_j>0\), and consequently
\begin{equation}\label{eq:effective-selection-backward}
 a_{j-1}>a_j+\frac Dj>0.
\end{equation}
Starting from \(a_N>0\), apply
\eqref{eq:effective-selection-backward} successively at
\(j=N,N-1,\ldots,m\).
Each step establishes positivity of the preceding term,
so the inequality is applicable at the next step.
Summing all these inequalities yields
\[
 a_{m-1}
 >a_N+D\sum_{j=m}^N\frac1j
 >DH_N=B_m,
\]
contradicting \eqref{eq:effective-selection-start-bound}.
Thus there is an index \(n\in\{m,\ldots,N\}\) such that
\(a_n>0\) and \(d_n\le D\).

\textbf{Obtain the two numerical bounds on the selected decrease.}
Since \(x\mapsto1/x\) is decreasing on \((0,\infty)\),
\[
 H_N
 \ge\sum_{j=m}^N\int_j^{j+1}\frac{dx}{x}
 =\log\frac{N+1}{m}
 \ge\log2.
\]
For the last inequality, if \(N\) is even then \(2m=N\),
and if \(N\) is odd then \(2m=N+1\).
In either case, \(N+1\ge2m\).
Also, \(m\ge N/2\), so
\[
 m-1\ge\frac{N-2}{2}>0.
\]
Using these bounds in \eqref{eq:effective-selection-level},
we obtain
\begin{equation}\label{eq:effective-selected-decrease}
 \begin{aligned}
 d_n
 &\le D\le\frac{0.4690-\ce}{\log2}<0.086,\\
 d_n
 &\le D\le\frac{4.75}{\sqrt{m-1}\log2}
 \le\frac{4.75\sqrt2}{\log2}\frac1{\sqrt{N-2}}
 <\frac{10}{\sqrt{N-2}}.
 \end{aligned}
\end{equation}
The numerical comparisons are
\[
 \frac{0.4690-\ce}{\log2}
 =0.08550538\ldots<0.086,\qquad
 \frac{4.75\sqrt2}{\log2}
 =9.69132474\ldots<10.
\]

\textbf{Make the global Gaussian bound explicit.}
We next verify that the function \(H_k\) defined in
\eqref{eq:gaussian-H} satisfies
\begin{equation}\label{eq:effective-H-global}
 |H_k(z,y)|
 \le\frac{\phi_0}{6}|y|^3+\frac{|y|}{5}+\frac12,
 \qquad k\ge2,\quad y,z\in\R.
\end{equation}
The estimate \eqref{eq:gaussian-global-components} from
Lemma~\ref{lem:gaussian-expansion} gives
\[
 |H_k(z,y)|
 \le\frac{\phi_0}{6}|y|^3
     +\frac{\phi_0}{2}|y|
     +\frac{M_3}{2\sqrt k},
 \qquad M_3=\|\phi'''\|_\infty.
\]
To bound \(M_3\), use
\(\phi'''(r)=r(3-r^2)\phi(r)\).
For \(r\ge0\),
\[
 \frac{d}{dr}\{r\phi(r)\}=(1-r^2)\phi(r),\qquad
 \frac{d}{dr}\{r^3\phi(r)\}=r^2(3-r^2)\phi(r).
\]
Both functions vanish at zero and tend to zero at infinity.
Their maxima therefore occur at \(r=1\) and \(r=\sqrt3\),
respectively, and symmetry gives
\[
 \sup_{r\in\R}|r|\phi(r)
 =\phi_0e^{-1/2}<\frac14,\qquad
 \sup_{r\in\R}|r|^3\phi(r)
 =3\sqrt3\,\phi_0e^{-3/2}<\frac12.
\]
It follows that
\[
 M_3
 \le3\sup_r|r|\phi(r)+\sup_r|r|^3\phi(r)
 <\frac34+\frac12=\frac54.
\]
Consequently, for \(k\ge2\),
\[
 \frac{M_3}{2\sqrt k}
 <\frac5{8\sqrt k}
 \le\frac5{8\sqrt2}<\frac12.
\]
Finally, \(\phi_0<2/5\) implies
\(\phi_0|y|/2\le |y|/5\).
These estimates prove \eqref{eq:effective-H-global}.

\textbf{Use the contact condition to exclude support points outside \([-6,6]\).}
We prove the support bound for any integer \(k\ge2\)
satisfying \(C_k>\ce\) and \(d_k\le1\).
Lemma~\ref{lem:attainment} ensures that a positive-discrepancy
maximizing pair exists at this index.
Let \((P,t)\) be any such pair, and write
\[
 R=R_k(P)=C_k,\qquad \beta=\beta(P),\qquad
 z=\frac{t}{\sqrt k},\qquad
 M=\int x|x|\,P(dx),\qquad G=F_{k-1,P}.
\]
In particular,
\[
 F_{k,P}(t)-\Phi(z)=\frac{R\beta}{\sqrt k}.
\]
Let \(I\) be the contamination derivative from
Lemma~\ref{lem:influence}.
For each \(y\in\R\), put
\[
 w=\frac{t-y}{\sqrt{k-1}}
  =\frac{z-y/\sqrt k}{\sqrt{1-1/k}}.
\]
The definition \eqref{eq:cn}, applied to the same law \(P\)
at sample size \(k-1\), gives
\[
 G(t-y)-\Phi(w)
 \le\frac{C_{k-1}\beta}{\sqrt{k-1}}.
\]
Subtract the maximizing identity above to obtain
\[
 G(t-y)-F_{k,P}(t)
 \le\Phi(w)-\Phi(z)
   +\frac{C_{k-1}\beta}{\sqrt{k-1}}
   -\frac{R\beta}{\sqrt k}.
\]
Substituting into \eqref{eq:influence} and collecting the
Gaussian terms by \eqref{eq:gaussian-H} yields
\[
 \begin{aligned}
 \beta I(y)\le{}&
 H_k(z,y)
 +\beta\left\{\frac{k^{3/2}}{\sqrt{k-1}}C_{k-1}-kR\right\}\\
 &+\frac32R\beta(y^2-1)
  -R\{|y|^3-\beta-3yM\}.
 \end{aligned}
\]
Since \(R=C_k\) and \(0\le C_{k-1}\le0.4690\) by
\eqref{eq:cn} and \eqref{eq:effective-cn},
the scale bound \eqref{eq:bounded-scale-factor} gives
\[
 \begin{aligned}
 \frac{k^{3/2}}{\sqrt{k-1}}C_{k-1}-kR
 &=k(C_{k-1}-C_k)
   +k\left(\sqrt{\frac{k}{k-1}}-1\right)C_{k-1}\\
 &\le d_k+0.4690.
 \end{aligned}
\]
Here the first term is at most \(d_k\) by
\eqref{eq:effective-increment}, and the coefficient of
\(C_{k-1}\) belongs to \([0,1]\).
Thus
\begin{equation}\label{eq:effective-influence-upper}
 \begin{aligned}
 \beta I(y)\le{}&
 H_k(z,y)+\beta(d_k+0.4690)\\
 &+\frac32R\beta(y^2-1)
  -R\{|y|^3-\beta-3yM\}.
 \end{aligned}
\end{equation}

The bounds \eqref{eq:effective-cn} and
\eqref{eq:moment-cutoff} imply
\[
 \frac25<\ce<R\le0.4690<\frac12,\qquad
 \beta<B_*<2,\qquad
 |M|\le\int x^2\,P(dx)=1.
\]
Put \(r=|y|\).
Using \eqref{eq:effective-H-global} in
\eqref{eq:effective-influence-upper}, expanding the last
two moment terms, and bounding \(My\le r\), we obtain
\begin{equation}\label{eq:effective-influence-polynomial-terms}
 \begin{aligned}
 \beta I(y)\le{}&
 -\left(R-\frac{\phi_0}{6}\right)r^3
 +\frac32R\beta r^2
 +\left(\frac15+3R\right)r\\
 &+\frac12+0.4690\beta-\frac12R\beta+\beta d_k.
 \end{aligned}
\end{equation}
The coefficients satisfy
\[
 R-\frac{\phi_0}{6}>\frac25-\frac1{15}=\frac13,\qquad
 \frac32R\beta<\frac32,\qquad
 \frac15+3R<\frac{17}{10}.
\]
Moreover, \(-R\beta/2\le0\), so
\[
 \frac12+0.4690\beta-\frac12R\beta
 \le\frac12+0.4690\beta
 <\frac12+0.938<\frac32.
\]
Since \(d_k\ge0\) and \(\beta<2\), we also have
\(\beta d_k\le2d_k\).
Consequently,
\begin{equation}\label{eq:effective-influence-polynomial}
 \beta I(y)
 \le-\frac{r^3}{3}+\frac{3r^2}{2}
      +\frac{17r}{10}+\frac32+2d_k.
\end{equation}

For \(d_k\le1\), the right side is bounded above by
\[
 q(r):=-\frac{r^3}{3}+\frac{3r^2}{2}
        +\frac{17r}{10}+\frac72.
\]
Direct calculation gives
\[
 q(6)=-72+54+\frac{51}{5}+\frac72=-\frac{43}{10}.
\]
For \(r\ge6\), its derivative satisfies
\[
 q'(r)=-r^2+3r+\frac{17}{10}
       =-r(r-3)+\frac{17}{10}
       \le-18+\frac{17}{10}
       =-\frac{163}{10}<0.
\]
Hence \(q(r)<0\) for every \(r\ge6\).
By \eqref{eq:effective-influence-polynomial} and \(\beta>0\),
this implies \(I(y)<0\) whenever \(|y|\ge6\).
The support contact condition \eqref{eq:contact} gives
\(I(y)=0\) at every \(y\in\supp P\), so no such point
can satisfy \(|y|\ge6\).
In particular, \(\supp P\subset[-6,6]\).

Finally, the selected index in
\eqref{eq:effective-selected-decrease} satisfies
\[
 n\ge m\ge2,\qquad C_n>\ce,\qquad d_n<0.086<1.
\]
The support argument therefore applies with \(k=n\).
Because the maximizing pair was arbitrary, it proves the
claimed support bound for every positive-discrepancy
maximizing law at that index.
\end{proof}

\subsection{Explicit smoothing and a near-lattice estimate}

We now make the Fourier estimates explicit for the bounded
laws obtained in Lemma~\ref{lem:effective-selection}.
The low-frequency estimate and the Gaussian tail bound below
will be used in the smoothing expansions of
Lemmas~\ref{lem:effective-cluster-jitter}
and~\ref{lem:effective-global-jitter}.
We then show that a violating law has small mean squared
distance from a suitable lattice.
This supplies the lattice used to construct the finite
approximation in Lemma~\ref{lem:effective-global-jitter}.

\textbf{Fix explicit constants in the smoothing inequality.}
In Lemma~\ref{lem:signed-smoothing}, one may take
\begin{equation}\label{eq:effective-smoothing-constants}
 C_1=\frac14,\qquad C_2=24.
\end{equation}
We verify these values using the kernel from that lemma.
The sum of two independent uniform variables on
\([-1/4,1/4]\) has density
\[
 q(v)=(2-4|v|)\mathbf 1_{\{|v|\le1/2\}}.
\]
The sum of four such variables therefore has density
\(q*q\) and characteristic function \(\sinc(t/4)^4\).
We have
\[
 (q*q)(0)=\int_{\R}q(v)^2\,dv
 =2\int_0^{1/2}(2-4v)^2\,dv=\frac43.
\]
The characteristic function \(\sinc(t/4)^4\) is integrable
because it is bounded near zero and is \(O(|t|^{-4})\)
at infinity.
Since \(q*q\) is continuous, Fourier inversion at zero gives
\[
 \int_{\R}\sinc(t/4)^4\,dt
 =2\pi(q*q)(0)=\frac{8\pi}{3}.
\]
Thus the normalized kernel is
\[
 K(x)=\frac3{8\pi}\sinc(x/4)^4.
\]
Using \(|\sin v|\le\min\{|v|,1\}\), we obtain, for \(x\ne0\),
\[
 K(x)\le\frac3{8\pi}\min\{1,256/x^4\}.
\]
If \(Z\) has density \(K\), splitting the first-moment
integral at \(4\) gives
\[
 \begin{aligned}
 \EE|Z|
 &\le\frac3{4\pi}
       \left(\int_0^4x\,dx
              +256\int_4^\infty x^{-3}\,dx\right)\\
 &=\frac3{4\pi}(8+8)=\frac{12}{\pi}<4.
 \end{aligned}
\]
Take \(a=16\).
Symmetry and Markov's inequality imply
\[
 \varepsilon=\PP(Z>16)
 =\frac12\PP(|Z|>16)
 \le\frac{\EE|Z|}{32}<\frac18.
\]
Symmetry also gives \(\EE Z=0\).
Using \(r_+=(r+|r|)/2\), we therefore have
\[
 C_a=\EE(16-Z)_+
 =\frac{16+\EE|16-Z|}{2}
 \le16+\frac12\EE|Z|<18.
\]
The constants obtained at the end of the proof of
Lemma~\ref{lem:signed-smoothing} consequently satisfy
\[
 \frac1{2\pi(1-2\varepsilon)}
 <\frac2{3\pi}<\frac14,\qquad
 \frac{C_a}{1-2\varepsilon}
 <\frac{18}{3/4}=24.
\]
Enlarging them gives \eqref{eq:effective-smoothing-constants}.

We next estimate the Fourier error on
\(|t|\le\sqrt k/2\).
The support bound gives a uniform fourth-moment bound,
which allows us to keep an explicit remainder in the
first-order expansion.

\begin{lemma}\label{lem:effective-low-frequency}
Let \(P\in\Pthree\) be supported in \([-6,6]\), with
\(\beta(P)<1.84\).
For \(X\sim P\), write
\[
 f(u)=\EE e^{iuX},\qquad \kappa=\EE X^3.
\]
For every integer \(k\ge2\) and every
\(|t|\le\sqrt k/2\),
\begin{equation}\label{eq:effective-fourier-error}
 \begin{aligned}
 &\left|f(t/\sqrt k)^k-
 e^{-t^2/2}\left(1+\frac{\kappa(it)^3}{6\sqrt k}\right)\right|\\
 &\qquad\le
 \frac{0.61t^4e^{-0.35t^2}+0.048t^6e^{-t^2/2}}{k}.
 \end{aligned}
\end{equation}
In particular,
\begin{equation}\label{eq:effective-fourier-integral}
 \int_{|t|\le\sqrt k/2}
 \frac{\left|f(t/\sqrt k)^k-
 e^{-t^2/2}\left(1+\frac{\kappa(it)^3}{6\sqrt k}\right)\right|}
 {|t|}\,dt
 \le\frac6k.
\end{equation}
The integrand at \(t=0\) is defined by continuity as zero.
\end{lemma}

\begin{proof}
\textbf{Expand the characteristic function and its logarithm.}
Write \(\beta=\beta(P)\).
Since \(|X|\le6\) almost surely,
\[
 \EE X^4\le6\EE|X|^3=6\beta<11.04,
 \qquad |\kappa|\le\beta<1.84.
\]
Apply the Taylor remainder bound
\eqref{eq:jitter-proof-taylor-remainder} with \(m=4\)
to \(uX\), and take expectations.
Using \(\EE X=0\) and \(\EE X^2=1\), we obtain
\begin{equation}\label{eq:effective-characteristic-expansion}
 f(u)=1-\frac{u^2}{2}+\frac{\kappa(iu)^3}{6}+r(u),
 \qquad
 |r(u)|\le\frac{\EE X^4}{24}u^4\le0.46u^4.
\end{equation}
The same Taylor bound with \(m=2\) gives separately
\[
 |f(u)-1|
 =\left|\EE(e^{iuX}-1-iuX)\right|
 \le\frac{u^2}{2}\EE X^2=\frac{u^2}{2}.
\]
For \(|u|\le1/2\), set \(w=f(u)-1\).
Then \(|w|\le1/8\), so \(f(u)\) lies in the disk
\(\{z:|z-1|\le1/8\}\), on which the principal logarithm
is well defined.
Its power series yields
\[
 \begin{aligned}
 |\log(1+w)-w|
 &\le\sum_{\ell=2}^{\infty}\frac{|w|^\ell}{\ell}\\
 &\le\frac{|w|^2}{2(1-|w|)}
 \le\frac{u^4}{8(1-1/8)}=\frac{u^4}{7}.
 \end{aligned}
\]
Combining this with
\eqref{eq:effective-characteristic-expansion}, and using
\(0.46+1/7<0.61\), gives
\begin{equation}\label{eq:effective-logarithmic-expansion}
 \left|\log f(u)+\frac{u^2}{2}
                   -\frac{\kappa(iu)^3}{6}\right|
 \le0.61u^4,\qquad |u|\le1/2.
\end{equation}

\textbf{Bound the modulus by a Gaussian function.}
Let \(X'\) be an independent copy of \(X\).
Independence and the first two moments of \(X\) give
\[
 \EE(X-X')^2=2,\qquad
 \EE(X-X')^4=2\EE X^4+6<28.08.
\]
Taking real parts in the fourth-order Taylor remainder
bound gives
\(\cos v\le1-v^2/2+v^4/24\) for every real \(v\).
Hence, for \(|u|\le1/2\),
\[
 \begin{aligned}
 |f(u)|^2
 &=\EE\cos(u(X-X'))\\
 &\le1-u^2+\frac{28.08}{24}u^4\\
 &=1-u^2+1.17u^4\\
 &\le1-\left(1-\frac{1.17}{4}\right)u^2
 \le1-0.7u^2.
 \end{aligned}
\]
Since \(1-s\le e^{-s}\) for \(s\ge0\), taking square
roots gives
\begin{equation}\label{eq:effective-characteristic-decay}
 |f(u)|\le e^{-0.35u^2},\qquad |u|\le1/2.
\end{equation}

\textbf{Compare the two exponentials and integrate the error.}
Fix \(|t|\le\sqrt k/2\), and put
\[
 u=\frac{t}{\sqrt k},\qquad
 A=k\log f(u),\qquad
 B=-\frac{ku^2}{2}+\frac{k\kappa(iu)^3}{6}.
\]
Equation \eqref{eq:effective-logarithmic-expansion} gives
\(|A-B|\le0.61ku^4\).
Also, \(\Re\log f(u)=\log|f(u)|\), so
\eqref{eq:effective-characteristic-decay} and the fact
that \(\kappa(iu)^3\) is purely imaginary imply
\[
 \Re A\le-0.35ku^2,\qquad
 \Re B=-\frac{ku^2}{2}.
\]
Consequently,
\[
 \Re\{(1-s)B+sA\}\le-0.35ku^2,\qquad 0\le s\le1.
\]
Integrating the derivative of
\(s\mapsto e^{(1-s)B+sA}\), we obtain
\begin{equation}\label{eq:effective-exponential-comparison}
 \begin{aligned}
 |e^A-e^B|
 &=\left|(A-B)\int_0^1e^{(1-s)B+sA}\,ds\right|\\
 &\le0.61ku^4e^{-0.35ku^2}
 =\frac{0.61t^4e^{-0.35t^2}}{k}.
 \end{aligned}
\end{equation}
Here \(e^A=f(u)^k\).
To compare \(e^B\) with the first-order expansion, set
\[
 v=-\frac{\kappa t^3}{6\sqrt k}\in\R.
\]
Then \(B=-t^2/2+iv\) and
\(iv=\kappa(it)^3/(6\sqrt k)\).
The Taylor bound
\eqref{eq:jitter-proof-taylor-remainder} with \(m=2\)
gives
\[
 \begin{aligned}
 &\left|e^B-e^{-t^2/2}
       \left(1+\frac{\kappa(it)^3}{6\sqrt k}\right)\right|\\
 &\qquad=e^{-t^2/2}|e^{iv}-1-iv|\\
 &\qquad\le e^{-t^2/2}\frac{v^2}{2}
 =\frac{\kappa^2t^6e^{-t^2/2}}{72k}
 \le\frac{0.048t^6e^{-t^2/2}}{k},
 \end{aligned}
\]
where \(\kappa^2/72<1.84^2/72<0.048\).
Adding this estimate to
\eqref{eq:effective-exponential-comparison} proves
\eqref{eq:effective-fourier-error}.

Finally, divide \eqref{eq:effective-fourier-error} by
\(|t|\) for \(t\ne0\), and enlarge the integration interval
to \(\R\).
The substitutions \(s=0.35t^2\) and \(s=t^2/2\) give
\[
 \int_{\R}|t|^3e^{-0.35t^2}\,dt=\frac1{(0.35)^2},
 \qquad
 \int_{\R}|t|^5e^{-t^2/2}\,dt=16.
\]
Thus the integral in
\eqref{eq:effective-fourier-integral} is at most
\[
 \frac1k\left(\frac{0.61}{(0.35)^2}+16(0.048)\right)
 <\frac6k.
\]
The bound in \eqref{eq:effective-fourier-error}, divided
by \(|t|\), also tends to zero as \(t\to0\), justifying
the stated value of the integrand there.
\end{proof}

\textbf{Prepare the comparison function and its Fourier tail.}
For \(k\ge2\) and \(|\kappa|<1.84\), define
\begin{equation}\label{eq:effective-edgeworth-comparison}
 \begin{aligned}
 \mathcal G_{k,\kappa}(x)
 &=\Phi(x)+\frac{\kappa}{6\sqrt k}(1-x^2)\phi(x),\\
 g_{k,\kappa}(x)
 &=\phi(x)\left[1+\frac{\kappa}{6\sqrt k}(x^3-3x)\right],\\
 b_k(t)
 &=e^{-t^2/2}\left[1+\frac{\kappa(it)^3}{6\sqrt k}\right].
 \end{aligned}
\end{equation}
Here \(g_{k,\kappa}=\mathcal G_{k,\kappa}'\), and \(b_k\)
is the Fourier transform of the signed measure
\(g_{k,\kappa}(x)\,dx\), by the calculation in
\eqref{eq:jitter-proof-comparison-density}--
\eqref{eq:jitter-proof-comparison-transform}.
The correction in \(g_{k,\kappa}\) is odd and integrable,
so this signed measure has total mass one.
Also,
\[
 |g_{k,\kappa}(x)|
 \le\left[1+\frac{1.84}{6}(|x|^3+3|x|)\right]\phi(x).
\]
The right side is integrable against \((1+|x|)\,dx\),
by the Gaussian moment calculation in
\eqref{eq:jitter-proof-density-moments}.
Thus the signed measure is finite and has finite first
absolute moment.
Since \(\mathcal G_{k,\kappa}(x)\to0\) as \(x\to-\infty\),
it is the primitive of this signed measure.

The estimate \(\|\phi'''\|_\infty<5/4\), proved in
Lemma~\ref{lem:effective-selection}, and
\((x^3-3x)\phi(x)=-\phi'''(x)\) give
\begin{equation}\label{eq:effective-edgeworth-density-bound}
 \|g_{k,\kappa}\|_\infty
 \le\phi_0+\frac{|\kappa|}{6\sqrt k}\|\phi'''\|_\infty
 <\frac25+\frac{1.84}{6}\frac54
 =\frac{47}{60}<1.
\end{equation}
It follows that \(\mathcal G_{k,\kappa}\) is globally
Lipschitz with constant at most one, as required for
Lemma~\ref{lem:signed-smoothing}.

We also need the explicit tail estimate
\begin{equation}\label{eq:effective-gaussian-tail}
 \int_{|t|\ge\sqrt k/2}\frac{|b_k(t)|}{|t|}\,dt
 \le5e^{-k/8},\qquad k\ge2.
\end{equation}
To prove it, set \(a=\sqrt k/2\).
By \eqref{eq:effective-edgeworth-comparison},
\[
 \int_{|t|\ge a}\frac{|b_k(t)|}{|t|}\,dt
 \le2\int_a^\infty\frac{e^{-t^2/2}}{t}\,dt
   +\frac{|\kappa|}{3\sqrt k}
        \int_a^\infty t^2e^{-t^2/2}\,dt.
\]
For \(t\ge a\), the inequality \(1/t\le t/a^2\) gives
\[
 2\int_a^\infty\frac{e^{-t^2/2}}{t}\,dt
 \le\frac2{a^2}\int_a^\infty te^{-t^2/2}\,dt
 =\frac2{a^2}e^{-a^2/2}.
\]
For the other term, integration by parts followed by
\(1\le t/a\) yields
\[
 \begin{aligned}
 \int_a^\infty t^2e^{-t^2/2}\,dt
 &=ae^{-a^2/2}+\int_a^\infty e^{-t^2/2}\,dt\\
 &\le ae^{-a^2/2}
      +\frac1a\int_a^\infty te^{-t^2/2}\,dt\\
 &=\left(a+\frac1a\right)e^{-a^2/2}.
 \end{aligned}
\]
Substituting \(a=\sqrt k/2\) and \(|\kappa|<1.84\)
therefore gives
\[
 \int_{|t|\ge\sqrt k/2}\frac{|b_k(t)|}{|t|}\,dt
 \le\left[
       \frac8k+1.84\left(\frac16+\frac2{3k}\right)
     \right]e^{-k/8}.
\]
For \(k\ge2\), the coefficient in brackets is at most
\(4+0.92=4.92<5\), proving
\eqref{eq:effective-gaussian-tail}.

We now use these estimates to locate a lattice near a
violating law.
For \(k\ge1024\), if \(|f(u)|\) were at most
\(1-\log k/k\) throughout \(1/2\le|u|\le1000\),
the smoothing inequality would give \(R_k(P)<\ce\).
A violation therefore forces a frequency \(u\) at which
\(|f(u)|\) is close to one.
At this frequency, the phases \(e^{iuX}\) have small
mean squared distance from a single point on the unit
circle, which will give the required lattice.

\begin{lemma}\label{lem:effective-near-lattice}
Suppose \(P\in\Pthree\) is supported in \([-6,6]\) and
\(R_k(P)>\ce\) for an integer \(k\ge1024\).
For \(X\sim P\), there are \(a\in\R\) and
\(h\in[\pi/500,4\pi]\) such that
\begin{equation}\label{eq:effective-near-lattice}
 \EE\dist(X,a+h\mathbb Z)^2
 \le2\pi^2\frac{\log k}{k}.
\end{equation}
\end{lemma}

\begin{proof}
\textbf{Use smoothing to find a frequency with modulus close to one.}
Write
\[
 \beta=\beta(P),\qquad \kappa=\EE X^3,\qquad
 f(u)=\EE e^{iuX}.
\]
The violation and \eqref{eq:moment-cutoff} imply
\(\beta<B_*<1.84\), and \(|\kappa|\le\beta<1.84\).
Thus Lemma~\ref{lem:effective-low-frequency} applies.
Let \(X_1,\ldots,X_k\) be independent with law \(P\), and put
\[
 F_k(x)=\PP\left\{\frac{X_1+\cdots+X_k}{\sqrt k}\le x\right\}.
\]
Its characteristic function is \(f(t/\sqrt k)^k\).
Independence, centering, and unit variance give
\[
 \EE\left|\frac{X_1+\cdots+X_k}{\sqrt k}\right|
 \le\left(\frac{\EE(X_1+\cdots+X_k)^2}{k}\right)^{1/2}
 =1
\]
by the Cauchy--Schwarz inequality.
Thus the standardized sum has finite first absolute moment.
The properties established in
\eqref{eq:effective-edgeworth-comparison}--
\eqref{eq:effective-edgeworth-density-bound} verify the
remaining assumptions of Lemma~\ref{lem:signed-smoothing}
for the comparison function \(\mathcal G_{k,\kappa}\),
with Lipschitz constant \(M=1\).

Suppose that
\begin{equation}\label{eq:effective-no-large-frequency}
 |f(u)|\le1-\frac{\log k}{k},
 \qquad 1/2\le|u|\le1000.
\end{equation}
Since \(0<\log k/k<1\), the inequality \(1-s\le e^{-s}\)
gives
\[
 |f(u)|^k\le
 \left(1-\frac{\log k}{k}\right)^k
 \le e^{-\log k}=\frac1k
\]
throughout this interval.
Choose the smoothing cutoff \(L=1000\sqrt k\).
The low-frequency contribution is at most \(6/k\) by
\eqref{eq:effective-fourier-integral}.
On the remaining frequencies, the triangle inequality gives
\[
 \begin{aligned}
 &\int_{\sqrt k/2\le|t|\le1000\sqrt k}
 \frac{|f(t/\sqrt k)^k-b_k(t)|}{|t|}\,dt\\
 &\quad\le
 \int_{\sqrt k/2\le|t|\le1000\sqrt k}
       \frac{|f(t/\sqrt k)|^k}{|t|}\,dt
 +\int_{|t|\ge\sqrt k/2}\frac{|b_k(t)|}{|t|}\,dt\\
 &\quad\le\frac2k\int_{\sqrt k/2}^{1000\sqrt k}\frac{dt}{t}
             +5e^{-k/8}\\
 &\quad=\frac{2\log2000}{k}+5e^{-k/8}
 \le\frac{16}{k}+5e^{-k/8}.
 \end{aligned}
\]
Here we used \eqref{eq:effective-no-large-frequency},
\eqref{eq:effective-gaussian-tail}, and \(\log2000<8\).
Combining the two frequency ranges, and applying
\eqref{eq:signed-smoothing} with the constants in
\eqref{eq:effective-smoothing-constants}, yields
\[
 \begin{aligned}
 \sqrt k\,\|F_k-\mathcal G_{k,\kappa}\|_\infty
 &\le\frac{\sqrt k}{4}
          \left(\frac{22}{k}+5e^{-k/8}\right)
       +\frac{24\sqrt k}{1000\sqrt k}\\
 &=\frac{11}{2\sqrt k}
       +\frac54\sqrt k\,e^{-k/8}+\frac3{125}.
 \end{aligned}
\]
The function \(x\mapsto xe^{-x/8}\) has derivative
\(e^{-x/8}(1-x/8)\), so it is decreasing for \(x\ge8\).
For \(k\ge1024\), it follows that
\[
 \frac54ke^{-k/8}\le1280e^{-128}<1,
\]
where the last inequality also follows from
\(e^{128}\ge128^2/2>1280\).
Consequently,
\begin{equation}\label{eq:effective-nonresonant-smoothing}
 \sqrt k\,\|F_k-\mathcal G_{k,\kappa}\|_\infty
 \le\frac{13}{2\sqrt k}+\frac3{125}.
\end{equation}

By \eqref{eq:jitter-envelope-gaussian-suprema},
\(\sup_x|(1-x^2)\phi(x)|=\phi_0\).
Thus \eqref{eq:effective-edgeworth-comparison} gives
\[
 \sqrt k\,\|\mathcal G_{k,\kappa}-\Phi\|_\infty
 =\frac{|\kappa|\phi_0}{6}.
\]
Also, \(|\kappa|\le\beta\), and
\[
 1=\EE|X|^2\le(\EE|X|^3)^{2/3}=\beta^{2/3}
\]
by H\"older's inequality, so \(\beta\ge1\).
Using the triangle inequality and
\eqref{eq:effective-nonresonant-smoothing}, we obtain
\[
 \begin{aligned}
 R_k(P)
 &=\frac{\sqrt k}{\beta}\|F_k-\Phi\|_\infty\\
 &\le\frac{|\kappa|\phi_0}{6\beta}
       +\frac1\beta
          \left(\frac{13}{2\sqrt k}+\frac3{125}\right)\\
 &\le\frac1{15}+\frac{13}{2\sqrt k}+\frac3{125}\\
 &\le\frac1{15}+\frac{13}{64}+\frac3{125}
   =\frac{7051}{24000}<0.3<\ce.
 \end{aligned}
\]
Here we used \(\phi_0<2/5\) and \(\sqrt k\ge32\).
This contradicts \(R_k(P)>\ce\).
Therefore \eqref{eq:effective-no-large-frequency} fails.
Since \(f(-u)=\overline{f(u)}\), we may choose a positive
frequency \(u\in[1/2,1000]\) such that
\begin{equation}\label{eq:effective-large-frequency}
 |f(u)|>1-\frac{\log k}{k}>0.
\end{equation}

\textbf{Convert concentration of phases into distance from a lattice.}
Choose \(\theta\in\R\) with
\(f(u)=|f(u)|e^{i\theta}\), and set
\[
 a=\frac{\theta}{u},\qquad h=\frac{2\pi}{u}.
\]
The range of \(u\) gives
\[
 \frac{\pi}{500}\le h\le4\pi.
\]
Moreover, \eqref{eq:effective-large-frequency} implies
\begin{equation}\label{eq:effective-phase-concentration}
 \begin{aligned}
 \EE|e^{iuX}-e^{i\theta}|^2
 &=2-2\Re\{e^{-i\theta}f(u)\}\\
 &=2(1-|f(u)|)
 <\frac{2\log k}{k}.
 \end{aligned}
\end{equation}
For \(v\in\R\), let
\(\delta=\dist(v,2\pi\mathbb Z)\in[0,\pi]\).
Then \(|e^{iv}-1|=2\sin(\delta/2)\).
Concavity of sine on \([0,\pi/2]\) gives
\(\sin r\ge2r/\pi\) there, and hence
\[
 |e^{iv}-1|\ge\frac{2\delta}{\pi}.
\]
Equivalently,
\[
 \dist(v,2\pi\mathbb Z)
 \le\frac{\pi}{2}|e^{iv}-1|.
\]
Since \(u>0\), \(ua=\theta\), and \(uh=2\pi\), we have
\[
 \dist(X,a+h\mathbb Z)
 =\frac1u\dist(uX-\theta,2\pi\mathbb Z)
 \le\frac{\pi}{2u}|e^{iuX}-e^{i\theta}|.
\]
Squaring and taking expectations, then using
\eqref{eq:effective-phase-concentration} and \(u\ge1/2\),
gives
\[
 \begin{aligned}
 \EE\dist(X,a+h\mathbb Z)^2
 &\le\frac{\pi^2}{4u^2}
       \EE|e^{iuX}-e^{i\theta}|^2\\
 &<\frac{\pi^2}{2u^2}\frac{\log k}{k}
 \le2\pi^2\frac{\log k}{k}.
 \end{aligned}
\]
This proves \eqref{eq:effective-near-lattice}.
\end{proof}

\subsection{A quantitative two-cluster smoothing expansion}

We now obtain an explicit finite-sample expansion for a sum
of two-cluster variables after adding one uniform variable.
Lemma~\ref{lem:effective-binomial} will use this estimate
with \(U=0\), and Proposition~\ref{prop:effective-clusters}
will use it for \(Y=B+U\).
The proof follows the frequency decomposition used in
Lemma~\ref{lem:jitter}, with numerical bounds on each
part of the Fourier integral.

In this and the next two subsections, let
\begin{equation}\label{eq:effective-cluster-setup}
 \begin{gathered}
 I=[2/5,9/20],\qquad
 B\sim\operatorname{Bernoulli}(p),\quad p\in I,\\
 Y=B+U,\qquad
 \EE(U\mid B)=0,\qquad |U|\le\varepsilon,
 \quad \varepsilon\ge0.
 \end{gathered}
\end{equation}
Put
\[
 v=p(1-p),\qquad s=\EE U^2,\qquad
 \sigma^2=v+s,\qquad
 \gamma=\frac{\EE(Y-p)^3}{\sigma^3}.
\]
Conditional centering gives
\[
 \EE U=\EE[\EE(U\mid B)]=0,\qquad
 \EE[(B-p)U]=\EE[(B-p)\EE(U\mid B)]=0.
\]
Consequently,
\[
 \EE Y=p,\qquad
 \Var Y=\EE(B-p+U)^2=v+s=\sigma^2.
\]
Let \(Y_1,Y_2,\ldots\) be independent copies of \(Y\),
and let \(H\), independent of these variables, be uniform
on \([-1/2,1/2]\).
Define
\[
 J_n(x)=\PP\left\{\sum_{i=1}^nY_i+H\le x\right\}.
\]
Thus the uniform variable is added once to the entire sum.

\begin{lemma}\label{lem:effective-cluster-jitter}
Suppose \(T\ge10\), \(\varepsilon\le(100T)^{-1}\), and
\(n\ge10^6\) is an integer.
With \(z=(x-np)/(\sigma\sqrt n)\), one has
\begin{equation}\label{eq:effective-cluster-jitter}
 \begin{aligned}
 &\sqrt n\sup_{x\in\R}\left|
 J_n(x)-\Phi(z)-\frac{\gamma}{6\sqrt n}(1-z^2)\phi(z)
 \right|\\
 &\quad\le
 \mathcal E(n,\varepsilon,T):=
 \frac{4(1+\log T)}{\sqrt n}
 +2\varepsilon^2T^2(1+\log T)+\frac{50}{T}.
 \end{aligned}
\end{equation}
\end{lemma}

\begin{proof}
The assumptions first give
\begin{equation}\label{eq:effective-cluster-small-noise}
 \varepsilon\le10^{-3},\qquad
 s\le\varepsilon^2\le\frac1{10000T^2},\qquad
 \frac6{25}\le v\le\frac14.
\end{equation}
For the last bounds, \(p(1-p)\) is increasing on \(I\),
so its minimum there is \((2/5)(3/5)=6/25\);
the upper bound \(p(1-p)\le1/4\) holds for every \(p\).

\textbf{Compare the squared characteristic functions and their derivatives.}
Write
\[
 f_Y(u)=\EE e^{iuY},\qquad q_Y(u)=|f_Y(u)|^2,
 \qquad
 f_B(u)=\EE e^{iuB},\qquad q_B(u)=|f_B(u)|^2.
\]
Let \((B',U')\) be an independent copy of \((B,U)\), and set
\[
 D_B=B-B',\qquad D_U=U-U'.
\]
Independence of the two pairs and conditional centering give
\[
 \EE(D_U\mid B,B')=0,\qquad
 \EE D_U^2=2\EE U^2-2(\EE U)^2=2s.
\]
Also, \(|D_B|\le1\) and \(|D_U|\le2\varepsilon\).
Independence of the copies yields
\[
 q_Y(u)=\EE\cos(u(D_B+D_U)),\qquad
 q_B(u)=\EE\cos(uD_B).
\]
The variables inside these expectations are bounded, so
we may differentiate twice under the expectations.

For a fixed \(u\), define functions of \(\xi\) by
\[
 \psi_0(\xi)=\cos(u\xi),\qquad
 \psi_1(\xi)=-\xi\sin(u\xi),\qquad
 \psi_2(\xi)=-\xi^2\cos(u\xi).
\]
Thus, for \(r=0,1,2\),
\[
 q_Y^{(r)}(u)=\EE\psi_r(D_B+D_U),\qquad
 q_B^{(r)}(u)=\EE\psi_r(D_B),
\]
where the superscript \(0\) means the function itself.
The linear Taylor term has expectation zero, since
\[
 \EE[\psi_r'(D_B)D_U]
 =\EE[\psi_r'(D_B)\EE(D_U\mid B,B')]=0.
\]
Every point between \(D_B\) and \(D_B+D_U\) has absolute
value at most \(c:=1+2\varepsilon\).
Taylor's formula with second-order remainder therefore gives
\[
 |q_Y^{(r)}(u)-q_B^{(r)}(u)|
 \le\frac12\EE D_U^2
           \sup_{|\xi|\le c}|\psi_r''(\xi)|
 =s\sup_{|\xi|\le c}|\psi_r''(\xi)|.
\]
The derivatives here are with respect to \(\xi\), and
\[
 \begin{aligned}
 \psi_0''(\xi)&=-u^2\cos(u\xi),\\
 \psi_1''(\xi)&=-2u\cos(u\xi)+u^2\xi\sin(u\xi),\\
 \psi_2''(\xi)&=-2\cos(u\xi)+4u\xi\sin(u\xi)
                           +u^2\xi^2\cos(u\xi).
 \end{aligned}
\]
It follows that
\begin{equation}\label{eq:effective-cluster-derivatives}
 \begin{aligned}
 |q_Y(u)-q_B(u)|
 &\le su^2,\\
 |q_Y'(u)-q_B'(u)|
 &\le s\{2|u|+cu^2\},\\
 |q_Y''(u)-q_B''(u)|
 &\le s\{2+4c|u|+c^2u^2\}
 \le10s(1+u^2).
 \end{aligned}
\end{equation}
For the last inequality, \(c\le1.002<2\) and
\(2|u|\le1+u^2\) imply
\[
 2+4c|u|+c^2u^2
 \le2+8|u|+4u^2
 \le6+8u^2\le10(1+u^2).
\]
This comparison uses conditional centering and does not
require \(B\) and \(U\) to be independent.

\textbf{Locate and bound the peaks near the Bernoulli lattice frequencies.}
Since \(f_B(u)=1-p+pe^{iu}\),
\[
 q_B(u)=1-4v\sin^2(u/2),\qquad
 q_B'(u)=-2v\sin u,\qquad
 q_B''(u)=-2v\cos u.
\]
For each integer \(j\), set
\[
 A_j=[2\pi j-1/4,\,2\pi j+1/4].
\]
Suppose \(A_j\cap[-T,T]\ne\varnothing\).
Then \(|2\pi j|\le T+1/4\), so every \(u\in A_j\)
satisfies
\[
 |u|\le T+1/2\le1.05T.
\]
We first work on this entire interval \(A_j\), including
any part outside the cutoff \([-T,T]\).
Equations \eqref{eq:effective-cluster-small-noise}
and \eqref{eq:effective-cluster-derivatives} give
\[
 \begin{aligned}
 |q_Y''(u)-q_B''(u)|
 &\le10s\{1+(T+1/2)^2\}\\
 &\le0.001\{T^{-2}+(1+1/(2T))^2\}\\
 &\le0.001(0.01+1.05^2)=0.0011125.
 \end{aligned}
\]
For \(u\in A_j\), periodicity and
\(\cos r\ge1-r^2/2\) give
\[
 \cos u=\cos(u-2\pi j)
 \ge1-\frac{(u-2\pi j)^2}{2}\ge\frac{31}{32}.
\]
Using \(2v\ge12/25\), we conclude that
\begin{equation}\label{eq:effective-cluster-concavity}
 q_Y''(u)
 \le-\frac{12}{25}\frac{31}{32}+0.0011125
 =-0.4638875<-\frac25,\qquad u\in A_j.
\end{equation}

For the first derivative, \(T\ge10\) and \(c\le1.002\)
give
\[
 \begin{aligned}
 |q_Y'(u)-q_B'(u)|
 &\le s\{2(1.05T)+1.002(1.05T)^2\}\\
 &\le sT^2\{0.21+1.002(1.05)^2\}\\
 &\le2sT^2\le0.0002<0.001.
 \end{aligned}
\]
Here \(0.21+1.002(1.05)^2=1.314705<2\).
At the left and right endpoints of \(A_j\), the
Bernoulli derivatives are respectively
\(2v\sin(1/4)\) and \(-2v\sin(1/4)\).
Their absolute values satisfy
\[
 2v\sin(1/4)
 \ge\frac{12}{25}
       \left(\frac14-\frac{(1/4)^3}{6}\right)
 =0.11875>0.11.
\]
Thus \(q_Y'\) is positive at the left endpoint and
negative at the right endpoint.
By \eqref{eq:effective-cluster-concavity}, it is strictly
decreasing on \(A_j\).
There is therefore a unique point \(u_j\) in the interior
of \(A_j\) with \(q_Y'(u_j)=0\), and this point is the
unique maximizer of \(q_Y\) on \(A_j\).

Since \(q_B'(2\pi j)=0\), the preceding first-derivative
bound gives \(|q_Y'(2\pi j)|\le2sT^2\).
The mean value theorem and
\eqref{eq:effective-cluster-concavity} imply
\begin{equation}\label{eq:effective-cluster-peak-location}
 |u_j-2\pi j|
 \le\frac52|q_Y'(2\pi j)|
 \le5sT^2\le8sT^2.
\end{equation}
The same concavity bound, applied in the second-order
Taylor formula at \(u_j\), gives
\[
 q_Y(u)\le q_Y(u_j)-\frac15(u-u_j)^2
          \le1-\frac15(u-u_j)^2.
\]
Since \(|u-u_j|\le1/2\), the upper bound is positive.
Using \(1-r\le e^{-r}\), we obtain
\begin{equation}\label{eq:effective-cluster-peak-bound}
 |f_Y(u)|^n=q_Y(u)^{n/2}
 \le e^{-n(u-u_j)^2/10},\qquad u\in A_j.
\end{equation}

\textbf{Use the uniform multiplier to bound the nonzero peaks.}
The characteristic function of \(H\) is
\[
 H_0(u):=\EE e^{iuH}=\sinc(u/2).
\]
It satisfies \(H_0(2\pi j)=0\) for every nonzero integer
\(j\).
Also,
\[
 |H_0'(u)|=|\EE(iH e^{iuH})|
 \le\EE|H|=\int_{-1/2}^{1/2}|x|\,dx=\frac14.
\]
Consequently, for \(j\ne0\) and \(u\in A_j\),
\begin{equation}\label{eq:effective-cluster-multiplier}
 |H_0(u)|
 \le\frac14|u-2\pi j|
 \le\frac14\{|u-u_j|+8sT^2\},
\end{equation}
where the last inequality uses
\eqref{eq:effective-cluster-peak-location}.
For \(j\ge1\), every \(u\in A_j\) satisfies
\[
 |u|\ge2\pi j-\frac14
       \ge(2\pi-1/4)j>5j.
\]
Combining this with
\eqref{eq:effective-cluster-peak-bound} and
\eqref{eq:effective-cluster-multiplier}, and then extending
the integral to the whole real line, gives
\[
 \begin{aligned}
 &\int_{A_j\cap[-T,T]}
       \frac{|f_Y(u)|^n|H_0(u)|}{|u|}\,du\\
 &\quad\le\frac1{20j}
       \int_{\R}(|w|+8sT^2)e^{-nw^2/10}\,dw\\
 &\quad=\frac1{20j}
       \left(\frac{10}{n}
              +8sT^2\sqrt{\frac{10\pi}{n}}\right)\\
 &\quad\le\frac1{2jn}+\frac{12sT^2}{5j\sqrt n}.
 \end{aligned}
\]
Here \(w=u-u_j\), and we used
\[
 \int_{\R}|w|e^{-nw^2/10}\,dw=\frac{10}{n},
 \qquad
 \int_{\R}e^{-nw^2/10}\,dw=\sqrt{\frac{10\pi}{n}},
 \qquad \sqrt{10\pi}<6.
\]
The same integral bound holds for the negative interval
\(A_{-j}\), since \(|f_Y(-u)|=|f_Y(u)|\) and
\(|H_0(-u)|=|H_0(u)|\).

Put
\[
 J=\left\lfloor\frac{T+1/4}{2\pi}\right\rfloor,\qquad
 \mathcal R_T=
 \bigcup_{1\le|j|\le J}(A_j\cap[-T,T]).
\]
These are exactly the nonzero peak intervals that meet
the cutoff.
Since \(T\ge10\), we have \(1\le J<T\), and
\[
 \sum_{j=1}^J\frac1j
 \le1+\int_1^J\frac{dx}{x}
 =1+\log J\le1+\log T.
\]
Summing over positive and negative indices, and using
\(24/5<5\), yields
\begin{equation}\label{eq:effective-cluster-resonances}
 \begin{aligned}
 \int_{\mathcal R_T}
       \frac{|f_Y(u)|^n|H_0(u)|}{|u|}\,du
 &\le
 \left(\frac1n+\frac{24sT^2}{5\sqrt n}\right)
       \sum_{j=1}^J\frac1j\\
 &\le
 \left(\frac1n+\frac{5sT^2}{\sqrt n}\right)(1+\log T).
 \end{aligned}
\end{equation}

\textbf{Bound the remaining frequencies away from zero.}
Define
\[
 \mathcal N_T=
 \{u:1/2\le|u|\le T\}\setminus\mathcal R_T.
\]
For \(u\in\mathcal N_T\),
\(\dist(u,2\pi\mathbb Z)\ge1/4\).
Indeed, a point within \(1/4\) of a nonzero \(2\pi j\)
would belong to \(\mathcal R_T\), while \(|u|\ge1/2\)
excludes the interval around zero.
Writing \(\delta=\dist(u,2\pi\mathbb Z)\), we have
\(1/4\le\delta\le\pi\).
Periodicity and monotonicity of sine on \([0,\pi/2]\)
therefore give
\[
 \sin^2(u/2)=\sin^2(\delta/2)\ge\sin^2(1/8)>
 \left(\frac3{25}\right)^2.
\]
For the strict numerical inequality, use
\(\sin(1/8)\ge1/8-(1/8)^3/6>3/25\).
Equations \eqref{eq:effective-cluster-derivatives}
and \eqref{eq:effective-cluster-small-noise} now give
\[
 \begin{aligned}
 q_Y(u)
 &\le q_B(u)+su^2\\
 &\le1-\frac{24}{25}\left(\frac3{25}\right)^2+10^{-4}\\
 &=0.986276<1-\frac1{100}.
 \end{aligned}
\]
Thus \(|f_Y(u)|^n=q_Y(u)^{n/2}\le e^{-n/200}\).
Using \(|H_0(u)|\le1\), we obtain
\begin{equation}\label{eq:effective-cluster-nonresonances}
 \int_{\mathcal N_T}
       \frac{|f_Y(u)|^n|H_0(u)|}{|u|}\,du
 \le2e^{-n/200}\int_{1/2}^T\frac{du}{u}
 =2\log(2T)e^{-n/200}.
\end{equation}

\textbf{Standardize the summands and estimate the low frequencies.}
Set \(X=(Y-p)/\sigma\), which has mean zero and variance
one by the identities preceding the lemma.
Equation \eqref{eq:effective-cluster-small-noise} gives
\begin{equation}\label{eq:effective-cluster-variance-bounds}
 \frac6{25}\le\sigma^2=v+s\le\frac14+10^{-6}<1,
 \qquad 0.48<\sigma<1.
\end{equation}
Since \(p\in[0.4,0.45]\) and \(|U|\le10^{-3}\),
\[
 |Y-p|\le\max(p,1-p)+\varepsilon\le0.601,
 \qquad
 |X|\le\frac{0.601}{\sqrt{6/25}}<2.
\]
We also need a bound on \(\beta(Y)=\EE|X|^3\).
Put \(A=B-p\).
Because \(|A|\ge0.4>\varepsilon\), \(A+U\) and \(A\)
have the same sign, and hence
\[
 |A+U|^3
 =|A|^3+3A|A|U+3|A|U^2+\operatorname{sgn}(A)U^3.
\]
The term \(3A|A|U\) has expectation zero by conditional
centering.
Moreover,
\[
 \EE|A|^3=(1-p)p^3+p(1-p)^3
          =v\{p^2+(1-p)^2\},
\]
and \(|A|\le1\), \(\EE|U|^3\le\varepsilon s\).
It follows that
\[
 \begin{aligned}
 \rho:=\EE|Y-p|^3
 &\le v\{p^2+(1-p)^2\}+(3+\varepsilon)s\\
 &\le\frac14(0.52)+4\cdot10^{-6}=0.130004.
 \end{aligned}
\]
Here \(p^2+(1-p)^2=1-2v\le0.52\),
\(s\le10^{-6}\), and \(3+\varepsilon<4\).
Consequently,
\begin{equation}\label{eq:effective-cluster-moment-bounds}
 |\gamma|\le\beta(Y)=\frac{\rho}{\sigma^3}
 \le\frac{0.130004}{(6/25)^{3/2}}<1.11.
\end{equation}
The law of \(X\) therefore satisfies the support and
moment assumptions of Lemma~\ref{lem:effective-low-frequency}.

Let \(f_X(w)=\EE e^{iwX}\), and define the standardized
smoothed distribution function
\[
 \mathcal J_n(w)=J_n(np+\sigma\sqrt n\,w).
\]
By independence of the summands and \(H\), the
corresponding probability law has characteristic function
\begin{equation}\label{eq:effective-cluster-transform}
 a_n(t)=f_X(t/\sqrt n)^n
              H_0\left(\frac{t}{\sigma\sqrt n}\right).
\end{equation}
Use the comparison function \(\mathcal G_{n,\gamma}\)
from \eqref{eq:effective-edgeworth-comparison}.
The associated signed measure has Fourier transform
\[
 b_n(t)=e^{-t^2/2}
        \left(1+\frac{\gamma(it)^3}{6\sqrt n}\right).
\]
The original frequency \(u\) and standardized frequency
\(t\) are related by
\begin{equation}\label{eq:effective-cluster-frequency-change}
 t=\sigma\sqrt n\,u,\qquad
 f_X(t/\sqrt n)=e^{-ipu}f_Y(u),\qquad
 \frac{dt}{|t|}=\frac{du}{|u|}.
\end{equation}
In particular, \(|u|\le1/2\) corresponds to
\[
 |t|\le a:=\frac{\sigma\sqrt n}{2}<\frac{\sqrt n}{2}.
\]

From \eqref{eq:effective-cluster-transform},
\[
 \begin{aligned}
 |a_n(t)-b_n(t)|
 &\le|f_X(t/\sqrt n)^n-b_n(t)|\\
 &\quad+
 \left|H_0\left(\frac{t}{\sigma\sqrt n}\right)-1\right|
       |f_X(t/\sqrt n)|^n.
 \end{aligned}
\]
The first term, divided by \(|t|\) and integrated over
\(|t|\le a\), contributes at most \(6/n\) by
\eqref{eq:effective-fourier-integral}.
For the second term, \(\EE H=0\) and
\(\EE H^2=1/12\), so the Taylor remainder bound
\eqref{eq:jitter-proof-taylor-remainder} with \(m=2\) gives
\[
 |H_0(u)-1|
 =|\EE(e^{iuH}-1-iuH)|
 \le\frac{u^2}{2}\EE H^2=\frac{u^2}{24}.
\]
Also, \eqref{eq:effective-characteristic-decay},
applied to \(X\), gives
\[
 |f_X(t/\sqrt n)|^n\le e^{-0.35t^2},
 \qquad |t|\le a.
\]
Thus the second integral is at most
\[
 \begin{aligned}
 &\frac1{24\sigma^2n}
       \int_{\R}|t|e^{-0.35t^2}\,dt\\
 &\qquad=\frac1{24\sigma^2n}\frac{20}{7}
 \le\frac1{24(6/25)n}\frac{20}{7}
 =\frac{125}{252n}<\frac1{2n}.
 \end{aligned}
\]
Combining the two contributions yields
\begin{equation}\label{eq:effective-cluster-low-integral}
 \int_{|t|\le a}\frac{|a_n(t)-b_n(t)|}{|t|}\,dt
 \le\frac6n+\frac1{2n}\le\frac7n.
\end{equation}
The integrand at zero can be defined by continuity.

\textbf{Bound the comparison transform beyond the low-frequency interval.}
We need the Gaussian tail starting at
\(a=\sigma\sqrt n/2\).
Repeating the integral calculation for
\eqref{eq:effective-gaussian-tail}, with this value of
\(a\), gives
\[
 \begin{aligned}
 \int_{|t|\ge a}\frac{|b_n(t)|}{|t|}\,dt
 &\le2\int_a^\infty\frac{e^{-t^2/2}}t\,dt
     +\frac{|\gamma|}{3\sqrt n}
         \int_a^\infty t^2e^{-t^2/2}\,dt\\
 &\le\left\{\frac2{a^2}
      +\frac{|\gamma|}{3\sqrt n}
          \left(a+\frac1a\right)\right\}e^{-a^2/2}\\
 &=\left\{\frac8{\sigma^2n}
      +|\gamma|\left(\frac{\sigma}{6}
                    +\frac2{3\sigma n}\right)\right\}
           e^{-\sigma^2n/8}.
 \end{aligned}
\]
By \eqref{eq:effective-cluster-variance-bounds} and
\eqref{eq:effective-cluster-moment-bounds}, and since
\(n\ge2\), the coefficient in braces is at most
\[
 \frac{100}{3n}
 +1.11\left(\frac16+\frac{25}{18n}\right)
 \le\frac{50}{3}+1.11\frac{31}{36}<20.
\]
Also, \(\sigma^2n/8\ge3n/100\).
We have therefore proved
\begin{equation}\label{eq:effective-cluster-gaussian-tail}
 \int_{|t|\ge a}\frac{|b_n(t)|}{|t|}\,dt
 \le20e^{-3n/100}.
\end{equation}

\textbf{Combine the frequency ranges and apply signed smoothing.}
Take \(L=\sigma T\sqrt n\).
Under the change of variables
\eqref{eq:effective-cluster-frequency-change},
the range \(a\le|t|\le L\) becomes
\(1/2\le|u|\le T\).
This range is the union of \(\mathcal R_T\) and
\(\mathcal N_T\).
Moreover,
\[
 |a_n(t)|=|f_Y(u)|^n|H_0(u)|.
\]
The triangle inequality, followed by
\eqref{eq:effective-cluster-resonances},
\eqref{eq:effective-cluster-nonresonances}, and
\eqref{eq:effective-cluster-gaussian-tail}, therefore gives
\[
 \begin{aligned}
 &\int_{a\le|t|\le L}
       \frac{|a_n(t)-b_n(t)|}{|t|}\,dt\\
 &\quad\le
 \left(\frac1n+\frac{5sT^2}{\sqrt n}\right)(1+\log T)
 +2\log(2T)e^{-n/200}+20e^{-3n/100}.
 \end{aligned}
\]
Adding \eqref{eq:effective-cluster-low-integral}, we obtain
\begin{equation}\label{eq:effective-cluster-total-integral}
 \begin{aligned}
 \int_{-L}^L\frac{|a_n(t)-b_n(t)|}{|t|}\,dt
 &\le\frac7n
 +\left(\frac1n+\frac{5sT^2}{\sqrt n}\right)(1+\log T)\\
 &\quad+2\log(2T)e^{-n/200}+20e^{-3n/100}.
 \end{aligned}
\end{equation}

The probability law with distribution function
\(\mathcal J_n\) has finite first absolute moment because
the summands and \(H\) are bounded.
Since \(|\gamma|<1.11<1.84\), the signed comparison
measure associated with \(\mathcal G_{n,\gamma}\) has
total mass one and finite first absolute moment, as
verified after \eqref{eq:effective-edgeworth-comparison}.
Its primitive is Lipschitz with constant at most one
by \eqref{eq:effective-edgeworth-density-bound}.
Thus Lemma~\ref{lem:signed-smoothing} applies.
Using \eqref{eq:effective-smoothing-constants} and
\(24/\sigma<50\), we obtain from
\eqref{eq:effective-cluster-total-integral}
\begin{equation}\label{eq:effective-cluster-smoothing-budget}
 \begin{aligned}
 \sqrt n\,\|\mathcal J_n-\mathcal G_{n,\gamma}\|_\infty
 &\le\frac7{4\sqrt n}
 +\left(\frac1{4\sqrt n}+\frac54sT^2\right)(1+\log T)\\
 &\quad+\frac{\sqrt n}{2}\log(2T)e^{-n/200}
       +5\sqrt n e^{-3n/100}+\frac{50}{T}.
 \end{aligned}
\end{equation}
Here the last term bounds the scaled smoothing remainder:
\[
 \sqrt n\,\frac{24}{L}
 =\frac{24}{\sigma T}<\frac{50}{T}.
\]

It remains to absorb the two exponential terms in
\eqref{eq:effective-cluster-smoothing-budget}.
Put \(\Lambda_T=1+\log T\), and denote their sum by
\[
 E_{n,T}=\frac{\sqrt n}{2}\log(2T)e^{-n/200}
              +5\sqrt n e^{-3n/100}.
\]
Since \(\log(2T)\le\Lambda_T\), \(\Lambda_T\ge1\), and
\(e^{-3n/100}\le e^{-n/200}\),
\[
 \frac{\sqrt n}{\Lambda_T}E_{n,T}
 \le\frac n2e^{-n/200}+5ne^{-n/200}
 =\frac{11n}{2}e^{-n/200}.
\]
Using \(e^x\ge x^2/2\) with \(x=n/200\) gives
\[
 \frac{11n}{2}e^{-n/200}
 \le\frac{440000}{n}\le0.44<1,
 \qquad n\ge10^6.
\]
Therefore \(E_{n,T}\le\Lambda_T/\sqrt n\), uniformly
over all \(T\ge10\).
Finally,
\[
 \frac7{4\sqrt n}
 +\frac{\Lambda_T}{4\sqrt n}
 +\frac{\Lambda_T}{\sqrt n}
 \le\frac{3\Lambda_T}{\sqrt n}
 \le\frac{4\Lambda_T}{\sqrt n},
 \qquad
 \frac54sT^2\Lambda_T
 \le2\varepsilon^2T^2\Lambda_T.
\]
Substituting these estimates into
\eqref{eq:effective-cluster-smoothing-budget} proves
the claimed bound for
\(\sqrt n\,\|\mathcal J_n-\mathcal G_{n,\gamma}\|_\infty\).
The change of variables \(x=np+\sigma\sqrt n\,w\) is
a bijection of \(\R\), so this supremum is exactly
the one in \eqref{eq:effective-cluster-jitter}.
\end{proof}

\subsection{Binomial estimates and small accumulated variance}

We now derive explicit binomial estimates from
Lemma~\ref{lem:effective-cluster-jitter} and use them to
compare \(Y=B+U\) with \(B\) when
\(n\ge10^{100}\), \(|U|\le\varepsilon\le10^{-12}\), and
\(\lambda=n\EE U^2\le10^{-12}\).
The resulting inequality \(R_n(B+U)\le R_n(B)<\ce\)
will handle the small-variance case in
Proposition~\ref{prop:effective-clusters}.
Near a binomial jump, the comparison uses
Lemma~\ref{lem:one-sided-loss}; thresholds far from the
mean are controlled by the nonuniform bound recorded next.

We use the iid nonuniform bound announced in
\cite[pp.~124--125]{Shevtsova2013} and proved in
\cite[Corollary~3.4.5, p.~90, and Table~3.3, p.~55]{Shevtsova2017}.
Its normalization and the constant \(17.36\) are also recorded
in equation~(1) and Table~1 of the arXiv version of
\cite{Shevtsova2020}.
Let \(Y_1,\ldots,Y_n\) be independent with common law \(Q\), and put
\[
 m_Q=\EE Y_1,\qquad v_Q=\Var Y_1>0,\qquad
 \rho_Q=\EE|Y_1-m_Q|^3<\infty.
\]
For \(Z_i=(Y_i-m_Q)/\sqrt{nv_Q}\), we have
\[
 \sum_{i=1}^n\EE Z_i^2=1,\qquad
 L_{3,n}:=\sum_{i=1}^n\EE|Z_i|^3
 =\frac{\rho_Q}{v_Q^{3/2}\sqrt n}
 =\frac{\beta(Q)}{\sqrt n}.
\]
Thus, with \(z=(x-nm_Q)/\sqrt{nv_Q}\), the bound becomes
\begin{equation}\label{eq:effective-nonuniform}
 \frac{\sqrt n}{\beta(Q)}
 \left|Q^{*n}((-\infty,x])-\Phi(z)\right|
 \le\frac{17.36}{1+|z|^3}.
\end{equation}
The source states the estimate for
\(\PP\{\sum_iY_i<x\}\).
Taking thresholds \(x_r>x\) with \(x_r\downarrow x\)
gives \(\PP\{\sum_iY_i<x_r\}\to
\PP\{\sum_iY_i\le x\}\).
The normal term and the right-hand side are continuous
in the threshold, so the estimate passes to the limit
and yields \eqref{eq:effective-nonuniform}.
Taking thresholds increasing to \(x\) in this
right-continuous version also gives the same bound for
\(Q^{*n}((-\infty,x))\).
For the conditional sums in
Lemma~\ref{lem:effective-local-mass}, we also use
\eqref{eq:independent-berry-esseen} with
\(C_{\mathrm{ind}}=1\).

The following lemma supplies the binomial point-mass
and jump estimates needed in the local comparisons.

\begin{lemma}\label{lem:effective-binomial}
Use the notation of Lemma~\ref{lem:binomial-estimates}, with \(p\in I\).
For every integer \(n\ge10^{100}\),
\begin{align}
 \sup_{k\in\mathbb Z}b_{n,k}(p)
 &\le\frac2{\sqrt n},\label{eq:effective-binomial-mass}\\
 0.99\le\frac{\sqrt{nv}\,b_{n,k}(p)}{\phi(z_k)}
 &\le1.01,\qquad
 b_{n,k}(p)\ge\frac{10^{-6}}{\sqrt n}
 \quad (|z_k|\le5),\label{eq:effective-binomial-central}\\
 b_{n,k}(p)&\ge\frac{e^{-100}}{\sqrt n}
 \quad (|z_k|\le11),\label{eq:effective-binomial-wide}\\
 |R^\pm_{n,k}(B_p)-A_p^\pm(z_k)|
 &\le2\cdot10^{-10},
 \qquad k\in\mathbb Z.\label{eq:effective-binomial-branches}
\end{align}
Moreover, \(R_n(B_p)>0.39\).
\end{lemma}

\begin{proof}
\textbf{Bound the point probabilities and record the parameter bounds.}
Recall
\[
 q=1-p,\qquad v=pq,\qquad
 \tau=p^2+q^2=1-2v,\qquad d=1-2p.
\]
Since \(p\in I=[2/5,9/20]\), these parameters satisfy
\begin{equation}\label{eq:effective-binomial-parameter-bounds}
 \frac6{25}\le v\le\frac14,\qquad
 \frac12\le\tau\le0.52,\qquad
 \frac1{10}\le d\le\frac15.
\end{equation}
The moment identities \eqref{eq:binomial-standardized-moments}
give
\[
 \beta(B_p)=\frac{\tau}{\sqrt v},\qquad
 \gamma=\frac d{\sqrt v},\qquad
 \frac1{\beta(B_p)}=\frac{\sqrt v}{\tau}\le1.
\]
The Fourier estimate established from
\eqref{eq:binomial-fourier-inversion} gives
\[
 \begin{aligned}
 b_{n,k}(p)
 &\le\frac1{2\pi}\int_{\R}e^{-2nvu^2/\pi^2}\,du
 =\frac{\sqrt\pi}{2\sqrt{2nv}}\\
 &\le\sqrt{\frac{25\pi}{48}}\,\frac1{\sqrt n}
 <\frac2{\sqrt n}.
 \end{aligned}
\]
For the final inequality, \(\pi<4\) implies
\(25\pi/48<25/12<4\).
This proves \eqref{eq:effective-binomial-mass}.

\textbf{Expand the distribution function at both sides of each jump.}
Let \(K\sim\operatorname{Bin}(n,p)\), let \(H\) be independent
and uniform on \([-1/2,1/2]\), and write
\[
 J_n(x)=\PP(K+H\le x),\qquad
 \psi(z)=(1-z^2)\phi(z),\qquad
 a=\frac1{2\sqrt{nv}}.
\]
Apply Lemma~\ref{lem:effective-cluster-jitter} with \(U=0\).
For every \(T\ge10\), its small-noise hypothesis holds with
\(\varepsilon=0\), and its sample-size hypothesis follows from
\(n\ge10^{100}>10^6\).
Thus there is a function \(r_{n,p,T}\) such that
\begin{equation}\label{eq:effective-binomial-jitter-error}
 \begin{aligned}
 J_n(np+\sqrt{nv}\,z)
 &=\Phi(z)+\frac{\gamma}{6\sqrt n}\psi(z)+r_{n,p,T}(z),\\
 \sqrt n\,\|r_{n,p,T}\|_\infty
 &\le E_T,\qquad
 E_T:=\frac{4(1+\log T)}{\sqrt n}+\frac{50}{T}.
 \end{aligned}
\end{equation}
The endpoint identities used in
\eqref{eq:binomial-jump-identities} give
\[
 F_{n,p}(k)=J_n(k+1/2),\qquad
 F_{n,p}(k-1)=J_n(k-1/2).
\]
Indeed, \(K\) is integer-valued and \(H\in[-1/2,1/2]\);
the only possible discrepancy at either endpoint requires
\(H=-1/2\), which has probability zero.
The corresponding standardized arguments in
\eqref{eq:effective-binomial-jitter-error} are \(z_k+a\)
and \(z_k-a\).

The calculation in Lemma~\ref{lem:effective-selection} gives
\[
 \|\phi'\|_\infty=\sup_z|z|\phi(z)<\frac14,\qquad
 \|\psi'\|_\infty=\|\phi'''\|_\infty<\frac54,
\]
where \(\psi'(z)=(z^3-3z)\phi(z)=-\phi'''(z)\).
For \(\theta\in\{-1,1\}\), Taylor's theorem and the mean
value theorem therefore yield
\[
 \begin{aligned}
 |\Phi(z+\theta a)-\Phi(z)-\theta a\phi(z)|
 &\le\frac{\|\phi'\|_\infty a^2}{2},\\
 |\psi(z+\theta a)-\psi(z)|
 &\le\|\psi'\|_\infty a.
 \end{aligned}
\]
After multiplication by \(\sqrt n\), the total error from
these two expansions is at most
\begin{equation}\label{eq:effective-binomial-taylor-bound}
 \begin{aligned}
 \sqrt n\,\frac{\|\phi'\|_\infty a^2}{2}
 +\frac{|\gamma|}{6}\|\psi'\|_\infty a
 &=\frac1{\sqrt n}
   \left\{\frac{\|\phi'\|_\infty}{8v}
          +\frac{|d|\|\psi'\|_\infty}{12v}\right\}\\
 &\le\frac1{\sqrt n}\frac{1/32+1/48}{6/25}
 =\frac{125}{576\sqrt n}<\frac1{4\sqrt n}.
 \end{aligned}
\end{equation}
Here we used \(\gamma=d/\sqrt v\), \(a=1/(2\sqrt{nv})\),
and \eqref{eq:effective-binomial-parameter-bounds}.
Combining \eqref{eq:effective-binomial-jitter-error}
and \eqref{eq:effective-binomial-taylor-bound}, we obtain
\begin{equation}\label{eq:effective-binomial-endpoint-expansions}
 \begin{aligned}
 \sqrt n\{F_{n,p}(k)-\Phi(z_k)\}
 &=\frac{\phi(z_k)}{2\sqrt v}
   +\frac{\gamma}{6}\psi(z_k)+e^+_{n,k},\\
 \sqrt n\{F_{n,p}(k-1)-\Phi(z_k)\}
 &=-\frac{\phi(z_k)}{2\sqrt v}
   +\frac{\gamma}{6}\psi(z_k)+e^-_{n,k},\\
 |e^\pm_{n,k}|&\le E_T+\frac1{4\sqrt n}.
 \end{aligned}
\end{equation}
These bounds hold for every integer \(k\), uniformly over
\(p\in I\).

\textbf{Deduce the central point-probability and jump-branch estimates.}
First choose \(T=10^{12}\).
Since \(\log10<3\) and \(n^{-1/2}\le10^{-50}\),
\[
 E_T
 =\frac{4(1+12\log10)}{\sqrt n}+5\cdot10^{-11}
 <148\cdot10^{-50}+5\cdot10^{-11}<10^{-10}.
\]
Thus \(|e^\pm_{n,k}|<2\cdot10^{-10}\).
Subtracting the second identity in
\eqref{eq:effective-binomial-endpoint-expansions} from
the first, and using \(F_{n,p}(k)-F_{n,p}(k-1)=b_{n,k}(p)\),
gives
\begin{equation}\label{eq:effective-binomial-local-error}
 \left|\sqrt n\,b_{n,k}(p)-\frac{\phi(z_k)}{\sqrt v}\right|
 \le4\cdot10^{-10}.
\end{equation}
For \(|z_k|\le5\), we have
\(\phi(z_k)\ge\phi(5)>10^{-6}\).
Also, \(\sqrt v\le1/2\), so
\[
 \left|\frac{\sqrt{nv}\,b_{n,k}(p)}{\phi(z_k)}-1\right|
 \le\frac{4\cdot10^{-10}\sqrt v}{\phi(z_k)}
 <2\cdot10^{-4}<0.01.
\]
For the absolute lower bound in the same region,
\eqref{eq:effective-binomial-local-error} and
\(1/\sqrt v\ge2\) give
\[
 \sqrt n\,b_{n,k}(p)
 \ge\frac{\phi(z_k)}{\sqrt v}-4\cdot10^{-10}
 >2\cdot10^{-6}-4\cdot10^{-10}>10^{-6}.
\]
This proves \eqref{eq:effective-binomial-central}.

To obtain the normalized branches, divide the first identity
in \eqref{eq:effective-binomial-endpoint-expansions} by
\(\beta(B_p)\), and take the negative of the second identity
before dividing by \(\beta(B_p)\).
Since \(\gamma=d/\sqrt v\) and \(\beta(B_p)=\tau/\sqrt v\),
the leading terms become
\[
 \frac1{\beta(B_p)}
 \left\{\frac{\phi(z_k)}{2\sqrt v}
       \mathbin{\pm}\frac{\gamma}{6}\psi(z_k)\right\}
 =\frac{\phi(z_k)}{6\tau}(3\pm d\mp dz_k^2)
 =A_p^\pm(z_k).
\]
Each error has absolute value at most
\(2\cdot10^{-10}/\beta(B_p)\le2\cdot10^{-10}\).
This proves \eqref{eq:effective-binomial-branches},
with the signs in the definitions of
Lemma~\ref{lem:binomial-estimates}.

\textbf{Use a smaller smoothing error on the wider central interval.}
For \eqref{eq:effective-binomial-wide}, choose \(T=e^{200}\)
in \eqref{eq:effective-binomial-jitter-error}.
Then
\[
 E_T=\frac{804}{\sqrt n}+50e^{-200}
 \le804\cdot10^{-50}+50e^{-200}<10^{-47}.
\]
Subtracting the two endpoint expansions again now yields
\[
 \left|\sqrt n\,b_{n,k}(p)-\frac{\phi(z_k)}{\sqrt v}\right|
 \le2E_T+\frac1{2\sqrt n}
 <2\cdot10^{-47}+\frac12\,10^{-50}
 <3\cdot10^{-47}.
\]
If \(|z_k|\le11\), then
\(\phi(z_k)\ge\phi(11)>e^{-62}\), while \(1/\sqrt v\ge2\).
Therefore
\[
 \sqrt n\,b_{n,k}(p)>2e^{-62}-3\cdot10^{-47}>e^{-100}.
\]
The last inequality follows from
\(3\cdot10^{-47}<e^{-100}\) and \(e^{-62}>e^{-100}\).
This proves \eqref{eq:effective-binomial-wide}.

\textbf{Obtain a uniform lower bound for the binomial error.}
Choose an integer \(k\) nearest to \(np\).
Then \(|k-np|\le1/2\), and
\[
 z_k^2\le\frac1{4nv}\le\frac{25}{24n}.
\]
The central value of the positive branch is
\[
 g(p):=A_p^+(0)
 =\frac{\phi_0(3+d)}{6\tau}
 =\frac{\phi_0(2-p)}{3\tau},
\]
as in \eqref{eq:binomial-central-branch}.
Since \(d>0\), \(\phi(z)=\phi_0e^{-z^2/2}\), and
\(0\le1-e^{-r}\le r\) for \(r\ge0\), we have
\[
 \begin{aligned}
 0\le g(p)-A_p^+(z_k)
 &=\frac{\phi_0}{6\tau}
   \left\{(3+d)(1-e^{-z_k^2/2})
               +dz_k^2e^{-z_k^2/2}\right\}\\
 &\le\frac{\phi_0}{6\tau}
       \left(\frac{3+d}{2}+d\right)z_k^2\\
 &\le\frac{2/5}{3}
       \left(\frac{3+1/5}{2}+\frac15\right)
       \frac{25}{24n}
 =\frac1{4n}.
 \end{aligned}
\]
Here we used \(\phi_0<2/5\), \(\tau\ge1/2\), and
\(d\le1/5\).
The definition of \(R_n(B_p)\) as a supremum of absolute
discrepancies implies \(R_n(B_p)\ge R^+_{n,k}(B_p)\).
Consequently, \eqref{eq:effective-binomial-branches} gives
\[
 \begin{aligned}
 R_n(B_p)
 &\ge A_p^+(z_k)-2\cdot10^{-10}\\
 &\ge g(p)-\frac1{4n}-2\cdot10^{-10}
 >g(p)-3\cdot10^{-10},
 \end{aligned}
\]
where \(1/(4n)<10^{-10}\) follows from \(n\ge10^{100}\).
Finally, \(\phi_0>0.3989\), \(2-p\ge1.55\), and
\(\tau\le0.52\) imply
\[
 g(p)>
 \frac{0.3989(1.55)}{3(0.52)}
 =0.3963429487\ldots>0.396.
\]
Since \(0.396-3\cdot10^{-10}>0.39\), the claimed lower
bound follows.
\end{proof}

For the comparison below, we recall the constants in
Lemma~\ref{lem:one-sided-loss}.
Its proof establishes
\eqref{eq:loss-upper} and \eqref{eq:loss-lower} with
\begin{equation}\label{eq:effective-one-sided}
 c_0=\frac38,\qquad C=176000,\qquad
 \kappa=\frac1{160000}.
\end{equation}
Indeed, \eqref{eq:loss-bounded} gives \(c_0=3/8\).
With the clipping level \(a_0=1/20\),
\eqref{eq:loss-clipping-errors} and
\eqref{eq:loss-clipped-first-moment} give
\[
 C=\frac38(16000)+\frac54(8000)+160000=176000,
 \qquad \kappa=a_0^4=\frac1{160000}.
\]
We will apply both strict-event inequalities to the
noise sum conditional on the number of upper-cluster
summands. Their loss term will dominate the fourth-moment
and normalization errors under the bounds on
\(n,\varepsilon,\lambda\) in
Lemma~\ref{lem:effective-small-variance}.

\begin{lemma}\label{lem:effective-small-variance}
If \(p\in I\), \(\EE(U\mid B)=0\),
\(|U|\le\varepsilon\le10^{-12}\), \(n\ge10^{100}\) is an
integer, and \(\lambda=n\EE U^2\le10^{-12}\), then
\(R_n(B+U)\le R_n(B)<\ce\).
For \(\lambda>0\), more precisely,
\begin{equation}\label{eq:effective-small-gap}
 R_n(B+U)\le R_n(B)-
 \min\left\{\frac{5\cdot10^{-9}\lambda}
                  {\sqrt{5\lambda+\varepsilon^2}},\
             \frac1{10}\right\}.
\end{equation}
\end{lemma}

\begin{proof}
If \(\lambda=0\), then \(\EE U^2=0\), so \(U=0\)
almost surely. The conclusion follows from
\eqref{eq:schulz-two-point}.
For the remainder of the proof, assume \(\lambda>0\).
Write
\[
 q=1-p,\qquad v=pq,\qquad \tau=p^2+q^2,\qquad
 s_b=\EE(U^2\mid B=b),\qquad s=qs_0+ps_1=\EE U^2.
\]
In particular,
\begin{equation}\label{eq:effective-small-parameters}
 p,q\ge\frac25,\qquad
 \frac6{25}\le v\le\frac14,\qquad
 \frac12\le\tau\le\frac{13}{25},\qquad
 s=\frac{\lambda}{n}\le\varepsilon^2\le10^{-24}.
\end{equation}
Let \((B_i,U_i)\), \(1\le i\le n\), be independent copies
of \((B,U)\), put \(Y_i=B_i+U_i\), and set
\[
 K=\sum_{i=1}^nB_i,\qquad
 b_k=\PP(K=k),\qquad
 F_B(x)=\PP(K\le x),\qquad
 F(x)=\PP\left\{\sum_{i=1}^nY_i\le x\right\}.
\]
As in \eqref{eq:effective-cluster-setup}, conditional
centering gives \(\EE Y_i=p\) and \(\Var Y_i=v+s>0\).

\textbf{Control the conditional noise and its first absolute moment.}
Conditional on the entire vector \((B_1,\ldots,B_n)\),
the variables \(U_i\) are independent with laws
\(U\mid B=B_i\).
Every vector containing \(k\) ones gives the same
convolution law for their sum.
Thus, for \(0\le k\le n\), conditional on \(K=k\),
the noise sum has the law
of a variable \(W_k\) formed by adding \(n-k\) independent
copies of \(U\mid B=0\) and \(k\) independent copies of
\(U\mid B=1\).
These variables are centered and bounded by \(\varepsilon\).
Since \(p,q\ge2/5\),
\[
 s_0+s_1\le\frac52s,\qquad
 t_k:=\EE W_k^2=(n-k)s_0+ks_1
 \le n(s_0+s_1)\le\frac52\lambda.
\]
Write \(W_k=\sum_{i=1}^n\xi_i\) for these independent
conditional copies.
Expanding the fourth power, independence and
centering eliminate all terms containing an index exactly
once. The remaining terms give, as in
\eqref{eq:cluster-fourth-moment},
\[
 \begin{aligned}
 \EE W_k^4
 &=\sum_{i=1}^n\EE\xi_i^4
   +6\sum_{i<\ell}\EE\xi_i^2\,\EE\xi_\ell^2\\
 &\le3t_k^2+\varepsilon^2t_k\\
 &\le\frac{75}{4}\lambda^2
          +\frac52\varepsilon^2\lambda
 \le20\lambda(\lambda+\varepsilon^2)=:D_4,
 \qquad 0\le k\le n.
 \end{aligned}
\]
Here \(\EE\xi_i^4\le\varepsilon^2\EE\xi_i^2\), and
\(6\sum_{i<\ell}\EE\xi_i^2\,\EE\xi_\ell^2
\le3(\sum_i\EE\xi_i^2)^2=3t_k^2\).
In particular,
\begin{equation}\label{eq:effective-small-fourth-moment}
 \EE W_k^4\le D_4
 \le20\cdot10^{-12}(10^{-12}+10^{-24})
 <3\cdot10^{-23}<\frac1{160000}.
\end{equation}

Put \(z_k=(k-np)/\sqrt{nv}\).
If \(|z_k|\le5\), then
\[
 |k-np|\le5\sqrt{nv}\le\frac52\sqrt n,
 \qquad np,nq\ge\frac25n>\frac52\sqrt n.
\]
Consequently \(0\le k\le n\).
For these indices,
\[
 \begin{aligned}
 |t_k-\lambda|
 &=|k-np|\,|s_1-s_0|\\
 &\le5\sqrt{nv}\,(s_0+s_1)
 \le\frac{25}{4}\frac{\lambda}{\sqrt n}
 \le\frac{\lambda}{2}.
 \end{aligned}
\]
Thus \(\lambda/2\le t_k\le3\lambda/2\), and \(t_k>0\).
H\"older's inequality gives
\[
 t_k=\EE\bigl[|W_k|^{2/3}|W_k|^{4/3}\bigr]
 \le(\EE|W_k|)^{2/3}(\EE W_k^4)^{1/3}.
\]
Rearranging, using the fourth-moment bound in terms of
\(t_k\), and then the two bounds on \(t_k\), yields
\begin{equation}\label{eq:effective-absolute-noise}
 \begin{aligned}
 \EE|W_k|
 &\ge\frac{t_k^{3/2}}{\sqrt{\EE W_k^4}}
 \ge\frac{t_k}{\sqrt{3t_k+\varepsilon^2}}\\
 &\ge\frac{\lambda}{2\sqrt{5\lambda+\varepsilon^2}},
 \qquad |z_k|\le5.
 \end{aligned}
\end{equation}
The last inequality uses \(t_k\ge\lambda/2\) and
\(3t_k\le(9/2)\lambda\le5\lambda\).

\textbf{Separate the nearest binomial count and bound the leakage.}
Write an arbitrary real threshold as \(x=k+u\), where
\(k\in\mathbb Z\) and \(|u|\le1/2\).
For \(0\le k\le n\), the conditional representation above
gives the identities
\eqref{eq:cluster-right-block}--\eqref{eq:cluster-left-block}:
\[
 \begin{aligned}
 F(k+u)&=F_B(k-1)+b_k\PP(W_k\le u)+e^+_{k,u},\\
 F((k+u)-)&=F_B(k-1)+b_k\PP(W_k<u)+e^-_{k,u}.
 \end{aligned}
\]
For a count \(j\ne k\), any difference from its binomial
contribution requires a noise of magnitude at least
\(|j-k|-1/2\).
Markov's inequality and \eqref{eq:effective-binomial-mass}
therefore give
\[
 |e^\pm_{k,u}|
 \le\sum_{\substack{0\le j\le n\\j\ne k}}
       \frac{b_j\EE W_j^4}{(|j-k|-1/2)^4}
 \le\frac{2D_4}{\sqrt n}
       \sum_{\substack{j\in\mathbb Z\\j\ne k}}
          (|j-k|-1/2)^{-4}.
\]
The decreasing-integrand bound for this series is
\[
 \begin{aligned}
 \sum_{\substack{j\in\mathbb Z\\j\ne k}}
       (|j-k|-1/2)^{-4}
 &=2\sum_{m=1}^{\infty}(m-1/2)^{-4}\\
 &\le2\left(16+\int_1^\infty(x-1/2)^{-4}\,dx\right)
 =\frac{112}{3}<40.
 \end{aligned}
\]
It follows that
\begin{equation}\label{eq:effective-leakage}
 |e^\pm_{k,u}|\le\frac{80D_4}{\sqrt n}.
\end{equation}
Using the event \(\{|W_j|\ge |j-k|-1/2\}\) in the
tail bound treats both the strict and non-strict
probabilities, including atoms at the threshold.

\textbf{Control the normalization and the Gaussian derivative near a jump.}
Let \(\rho=\EE|Y-p|^3\).
Since \(\varepsilon<\min(p,q)\), \(Y-p\) is negative
when \(B=0\) and positive when \(B=1\).
Expanding the two conditional cubes as in
\eqref{eq:cluster-third-moment} gives
\[
 \rho=v\tau
   +q\{3ps_0-\EE(U^3\mid B=0)\}
   +p\{3qs_1+\EE(U^3\mid B=1)\}.
\]
The bounds
\(|\EE(U^3\mid B=b)|\le\varepsilon s_b\) show that both
terms in braces are nonnegative, and their weighted sum
is at most \((3+\varepsilon)s\le4s\).
Thus
\[
 0\le\rho-v\tau\le4s,\qquad
 \rho\ge v\tau\ge\frac3{25}.
\]
Define
\[
 A=\frac{(v+s)^{3/2}}{\rho}=\frac1{\beta(Y)},\qquad
 A_B=\frac{\sqrt v}{\tau}=\frac1{\beta(B)}.
\]
Since \(v+s\le1/4+10^{-24}<4/9\), the derivative of
\(r^{3/2}\) is at most one between \(v\) and \(v+s\).
Hence \((v+s)^{3/2}-v^{3/2}\le s\).
Using
\[
 A-A_B
 =\frac{(v+s)^{3/2}-v^{3/2}}{\rho}
   +\frac{v^{3/2}(v\tau-\rho)}{\rho\,v\tau},
\]
we obtain
\[
 |A-A_B|
 \le\frac{s}{3/25}
       +\frac{(1/8)4s}{(3/25)^2}
 =\frac{775}{18}s<50s.
\]
Also,
\(\sqrt{6/25}/(13/25)\le A_B\le(1/2)/(1/2)=1\),
so \(0.9<A_B\le1\).
Together with \(s\le10^{-24}\), this proves
\begin{equation}\label{eq:effective-prefactor}
 \frac12<A<2,\qquad |A-A_B|\le50s.
\end{equation}

For every integer \(k\), define
\[
 G_k(u)=\Phi\left(\frac{k-np+u}{\sqrt{n(v+s)}}\right),
 \qquad
 w_{k,u}=\sqrt{\frac v{v+s}}z_k
               +\frac{u}{\sqrt{n(v+s)}}.
\]
If \(|z_k|\le5\), the preceding central-index argument
gives \(0\le k\le n\), and \(b_k>0\).
For these indices, define
\[
 f_k(u)=\frac{G_k(u)-G_k(0)}{b_k}.
\]
Then \(f_k(0)=0\), and differentiation gives
\[
 f_k'(u)
 =\sqrt{\frac v{v+s}}\,
   \frac{\phi(w_{k,u})}{\phi(z_k)}\,
   \frac{\phi(z_k)}{\sqrt{nv}\,b_k},
\]
as in \eqref{eq:cluster-derivative-ratio}.
For \(|u|\le3/5\), use
\(0\le1-(1+r)^{-1/2}\le r/2\), \(r\ge0\), to obtain
\[
 \begin{aligned}
 |w_{k,u}-z_k|
 &\le\frac{5s}{2v}+\frac{3}{5\sqrt{nv}}\\
 &\le\frac{125}{12}s
       +\frac{3}{5\sqrt{6/25}}\frac1{\sqrt n}
 \le\frac3{\sqrt n}.
 \end{aligned}
\]
The last bound uses \(s=\lambda/n\le10^{-12}/n\).
Consequently,
\[
 \left|\log\frac{\phi(w_{k,u})}{\phi(z_k)}\right|
 =\frac12|w_{k,u}^2-z_k^2|
 \le\frac{3}{\sqrt n}\left(5+\frac{3}{2\sqrt n}\right)
 \le\frac{20}{\sqrt n}.
\]
Furthermore,
\[
 1-3s\le\sqrt{\frac v{v+s}}\le1,\qquad
 \frac1{1.01}\le
       \frac{\phi(z_k)}{\sqrt{nv}\,b_k}\le\frac1{0.99},
\]
where the second pair follows from
\eqref{eq:effective-binomial-central}.
Since \(s\le10^{-24}\) and \(n\ge10^{100}\), we have
\[
 1-3s>0.999,\qquad
 e^{-20/\sqrt n}>0.999,\qquad e^{20/\sqrt n}<1.001.
\]
Multiplication of the three factors therefore gives
\begin{equation}\label{eq:effective-small-derivative}
 \frac34<0.98<
 \frac{0.999^2}{1.01}
 \le f_k'(u)
 \le\frac{1.001}{0.99}<1.02<\frac54,
 \qquad |z_k|\le5,\quad |u|\le3/5.
\end{equation}

\textbf{Obtain a strict loss for both central signed discrepancies.}
Fix \(|z_k|\le5\) and \(|u|\le1/2\).
Using \(F_B(k)=F_B(k-1)+b_k\), rearrange the two
block identities
\eqref{eq:cluster-right-block}--\eqref{eq:cluster-left-block}
to get
\[
 \begin{aligned}
 F(k+u)-G_k(u)
 &=F_B(k)-G_k(0)
   -b_k\{\PP(W_k>u)+f_k(u)\}+e^+_{k,u},\\
 G_k(u)-F((k+u)-)
 &=G_k(0)-F_B(k-1)
   -b_k\{\PP(W_k<u)-f_k(u)\}-e^-_{k,u}.
 \end{aligned}
\]
Equations \eqref{eq:effective-small-fourth-moment}
and \eqref{eq:effective-small-derivative} verify the
assumptions of Lemma~\ref{lem:one-sided-loss}, with the
constants in \eqref{eq:effective-one-sided}.
Applying \eqref{eq:loss-upper} and \eqref{eq:loss-lower},
and then \eqref{eq:effective-leakage}, gives
\[
 \begin{aligned}
 F(k+u)-G_k(u)
 &\le F_B(k)-G_k(0)-\frac38b_k\EE|W_k|
       +176000b_kD_4+\frac{80D_4}{\sqrt n},\\
 G_k(u)-F((k+u)-)
 &\le G_k(0)-F_B(k-1)-\frac38b_k\EE|W_k|
       +176000b_kD_4+\frac{80D_4}{\sqrt n}.
 \end{aligned}
\]
For all integers \(k\) and \(|u|\le1/2\), define
\[
 \begin{aligned}
 D^+_{k,u}&=A\sqrt n\,[F(k+u)-G_k(u)],\\
 D^-_{k,u}&=A\sqrt n\,[G_k(u)-F((k+u)-)],\\
 \widetilde R^+_{n,k}&=A\sqrt n\,[F_B(k)-G_k(0)],\\
 \widetilde R^-_{n,k}&=A\sqrt n\,[G_k(0)-F_B(k-1)].
 \end{aligned}
\]
The central mass bound in
\eqref{eq:effective-binomial-central}, together with
\eqref{eq:effective-prefactor} and
\eqref{eq:effective-absolute-noise}, gives
\[
 \frac38A\sqrt n\,b_k\EE|W_k|
 \ge
 \left(\frac38\cdot\frac12\cdot10^{-6}\cdot\frac12\right)
       \frac{\lambda}{\sqrt{5\lambda+\varepsilon^2}}
 >
 \frac{10^{-8}\lambda}{\sqrt{5\lambda+\varepsilon^2}}.
\]
For the errors, \eqref{eq:effective-binomial-mass}
and \(A<2\) imply
\[
 \begin{aligned}
 A\sqrt n\left(176000b_kD_4+\frac{80D_4}{\sqrt n}\right)
 &\le(176000\cdot4+160)D_4\\
 &\le14{,}083{,}200\,\lambda(\lambda+\varepsilon^2)\\
 &<2\cdot10^7\lambda(\lambda+\varepsilon^2).
 \end{aligned}
\]
We conclude that
\begin{equation}\label{eq:effective-small-central-loss}
 D^\pm_{k,u}
 \le\widetilde R^\pm_{n,k}
      -\frac{10^{-8}\lambda}{\sqrt{5\lambda+\varepsilon^2}}
      +2\cdot10^7\lambda(\lambda+\varepsilon^2),
 \qquad |z_k|\le5,\quad |u|\le1/2.
\end{equation}

\textbf{Restore the Bernoulli variance and normalization.}
For \(r\ge0\), put \(w=z/\sqrt{1+r}\).
Differentiation gives
\[
 \left|\frac{\partial}{\partial r}
       \Phi\left(\frac z{\sqrt{1+r}}\right)\right|
 =\frac{|w|\phi(w)}{2(1+r)}
 \le\frac18,
\]
using \(\sup_w|w|\phi(w)<1/4\), established in
Lemma~\ref{lem:effective-selection}.
Integration from \(r=0\) to \(r=s/v\) gives
\[
 \sup_{z\in\R}
 \left|\Phi\left(\frac z{\sqrt{1+s/v}}\right)-\Phi(z)\right|
 \le\frac{s}{8v}\le\frac{25}{48}s<s.
\]
In particular, \(|G_k(0)-\Phi(z_k)|\le s\).
By \eqref{eq:schulz-two-point},
\(R_n(B)<\ce<1/2\), and \(A_B>0.9>1/2\).
Thus
\[
 \sqrt n\,\Delta_n(B)=\frac{R_n(B)}{A_B}<1.
\]
Both
\(|F_B(k)-\Phi(z_k)|\) and
\(|\Phi(z_k)-F_B(k-1)|\) are at most \(\Delta_n(B)\).
For the second assertion, take thresholds tending to
\(k\) from below and use continuity of the Gaussian CDF.
Comparing with the branches \(R^\pm_{n,k}(B_p)\) from
Lemma~\ref{lem:binomial-estimates}, the exact differences are
\[
 \begin{aligned}
 \widetilde R^+_{n,k}-R^+_{n,k}(B_p)
 &=\sqrt n(A-A_B)[F_B(k)-\Phi(z_k)]
   +A\sqrt n[\Phi(z_k)-G_k(0)],\\
 \widetilde R^-_{n,k}-R^-_{n,k}(B_p)
 &=\sqrt n(A-A_B)[\Phi(z_k)-F_B(k-1)]
   +A\sqrt n[G_k(0)-\Phi(z_k)].
 \end{aligned}
\]
The preceding bounds and the triangle inequality now give
\begin{equation}\label{eq:effective-small-branch-restoration}
 \begin{aligned}
 |\widetilde R^\pm_{n,k}-R^\pm_{n,k}(B_p)|
 &\le\sqrt n\,|A-A_B|\Delta_n(B)
       +A\sqrt n\,|G_k(0)-\Phi(z_k)|\\
 &\le50s+2\sqrt n\,s\\
 &\le\frac{52\lambda}{\sqrt n}
 \le\frac{100\lambda}{\sqrt n}.
 \end{aligned}
\end{equation}

The positive errors in
\eqref{eq:effective-small-central-loss} and
\eqref{eq:effective-small-branch-restoration}
can now be absorbed by half the loss.
Indeed,
\[
 \begin{aligned}
 2\cdot10^7(\lambda+\varepsilon^2)+\frac{100}{\sqrt n}
 &\le2\cdot10^7(10^{-12}+10^{-24})+100\cdot10^{-50}\\
 &<2.1\cdot10^{-5}<\frac1{600},
 \end{aligned}
\]
whereas
\[
 \frac{5\cdot10^{-9}}{\sqrt{5\lambda+\varepsilon^2}}
 \ge\frac{5\cdot10^{-9}}{\sqrt{5\cdot10^{-12}+10^{-24}}}
 >\frac1{600}.
\]
Multiplying this comparison by \(\lambda>0\), and using
\(R^\pm_{n,k}(B_p)\le R_n(B)\), proves
\begin{equation}\label{eq:effective-small-central-gap}
 D^\pm_{k,u}
 \le R_n(B)
       -\frac{5\cdot10^{-9}\lambda}
                    {\sqrt{5\lambda+\varepsilon^2}},
 \qquad |z_k|\le5,\quad |u|\le1/2.
\end{equation}
Every error in this comparison has a factor \(\lambda\);
the argument therefore applies to every \(\lambda>0\)
allowed by the hypotheses.

\textbf{Control the remaining thresholds and take the supremum.}
Suppose \(|z_k|>5\).
For \(|u|\le1/2\), the standardized threshold for \(Y\)
is \(w_{k,u}\), and
\[
 \begin{aligned}
 |w_{k,u}|
 &\ge\sqrt{\frac v{v+s}}\,|z_k|
           -\frac1{2\sqrt{n(v+s)}}\\
 &>5(1-3s)-\frac1{2\sqrt{nv}}\\
 &\ge5-15\cdot10^{-24}
              -\frac{25}{24}10^{-50}>4.9.
 \end{aligned}
\]
Here we used \(\sqrt v\ge\sqrt{6/25}>0.48\).
If \(k<0\) or \(k>n\), then
\[
 |k-np|\ge\min(np,nq)\ge\frac25n,\qquad
 |z_k|\ge\frac45\sqrt n>5,
\]
so these indices are included in the same case.
The nonuniform bound \eqref{eq:effective-nonuniform},
including its left-limit version, now gives
\[
 D^\pm_{k,u}\le\frac{17.36}{1+4.9^3}<0.15.
\]
Lemma~\ref{lem:effective-binomial} gives \(R_n(B)>0.39\),
and hence
\[
 D^\pm_{k,u}<R_n(B)-\frac1{10},
 \qquad |z_k|>5,\quad |u|\le1/2.
\]
Together with \eqref{eq:effective-small-central-gap},
this bounds both signed discrepancies at every real
threshold \(x=k+u\).
For the negative discrepancy of the right-continuous CDF,
we use \(F(x-)\le F(x)\), which implies
\[
 A\sqrt n\left[
 \Phi\left(\frac{x-np}{\sqrt{n(v+s)}}\right)-F(x)
 \right]\le D^-_{k,u}.
\]
The positive discrepancy is \(D^+_{k,u}\).
Taking the supremum over \(x\in\R\) therefore proves
\eqref{eq:effective-small-gap}.
Finally, \(R_n(B)<\ce\) by \eqref{eq:schulz-two-point},
which completes the proof.
\end{proof}

\subsection{An explicit local two-cluster theorem}

We now complete the explicit comparison for the two-cluster
laws in \eqref{eq:effective-cluster-setup}.
Lemma~\ref{lem:effective-small-variance} treats the case
\(\lambda=n\EE U^2\le10^{-12}\).
For the remaining case, we first give a numerical version
of Lemma~\ref{lem:two-cluster-local-mass}.
Its lower bounds for short interval probabilities will
strengthen the comparison between the original sum and
the sum with uniform jitter.
Combining this strict comparison with
Lemma~\ref{lem:effective-cluster-jitter} will give the
parameters \(\eta_*\) and \(N_*\) specified in
\eqref{eq:effective-threshold-parameters}.

\begin{lemma}\label{lem:effective-local-mass}
Suppose \(p\in I\), \(\EE(U\mid B)=0\),
\(|U|\le e^{-\mathcal A}\), and \(n\ge10^{100}\) is an
integer satisfying \(\lambda=n\EE U^2\ge10^{-12}\).
Let \(Y_1,\ldots,Y_n\) be independent copies of
\(Y=B+U\), as in \eqref{eq:effective-cluster-setup}.
For \((a,b)=(1/4,1/2)\) and
\((a,b)=(-1/2,-1/4)\), one has
\begin{equation}\label{eq:effective-local-mass}
 \inf_{|x-np|\le3\sqrt n}
 \PP\left\{\sum_{i=1}^nY_i\in(x+a,x+b)\right\}
 \ge\frac{e^{-4\cdot10^{12}}}{\sqrt n}.
\end{equation}
\end{lemma}

\begin{proof}
\textbf{Select binomial counts near the interval and compare their variances.}
Realize the \(Y_i\) as \(B_i+U_i\), where
\((B_i,U_i)\), \(1\le i\le n\), are independent copies
of \((B,U)\).
Put
\[
 S_n=\sum_{i=1}^nY_i,\qquad K=\sum_{i=1}^nB_i,\qquad
 s_b=\EE(U^2\mid B=b),\qquad
 s=(1-p)s_0+ps_1=\EE U^2=\frac{\lambda}{n}.
\]
Since \(p,1-p\ge2/5\), we have
\begin{equation}\label{eq:effective-local-conditional-variances}
 |s_1-s_0|\le s_0+s_1\le\frac52s.
\end{equation}
For \(0\le k\le n\), conditional on \(K=k\), the sum
\(S_n\) has the law of \(k+W_k\), where \(W_k\) is the
sum of \(n-k\) independent copies of \(U\mid B=0\)
and \(k\) independent copies of \(U\mid B=1\).
This is the conditional representation used in
Lemma~\ref{lem:effective-small-variance}: condition first
on all labels \(B_i\), observe that the noise variables
are then independent, and note that every label vector
with \(k\) ones gives the same convolution law.
The conditional summands are centered, so
\begin{equation}\label{eq:effective-local-noise-variance}
 \begin{aligned}
 \EE W_k&=0,\\
 t_k:=\EE W_k^2
 &=(n-k)s_0+ks_1
 =\lambda+(k-np)(s_1-s_0).
 \end{aligned}
\end{equation}
The last identity follows from
\(\lambda=n[(1-p)s_0+ps_1]\).

Fix \(x\) with \(|x-np|\le3\sqrt n\), and set
\begin{equation}\label{eq:effective-local-selected-counts}
 r=\max(1,\sqrt\lambda),\qquad
 \mathcal K_x=\{k\in\mathbb Z:|k-x|\le r\}.
\end{equation}
The bound on \(U\) gives
\(\lambda=ns\le ne^{-2\mathcal A}\le n\).
Since \(n\ge1\), both \(1\) and \(\sqrt\lambda\) are at
most \(\sqrt n\), and therefore \(r\le\sqrt n\).
For every \(k\in\mathcal K_x\), it follows that
\begin{equation}\label{eq:effective-local-count-location}
 |k-np|\le|k-x|+|x-np|
 \le4\sqrt n\le5\sqrt n.
\end{equation}
Moreover,
\[
 np,n(1-p)\ge\frac25n>5\sqrt n,
\]
where the strict inequality follows from \(n\ge10^{100}\).
Thus all integers in \(\mathcal K_x\) belong to \([0,n]\),
and the conditional representation above applies to them.
Using \eqref{eq:effective-local-conditional-variances}
and \eqref{eq:effective-local-count-location} in
\eqref{eq:effective-local-noise-variance}, and recalling
\(\lambda>0\), gives
\[
 \left|\frac{t_k}{\lambda}-1\right|
 \le\frac{5\sqrt n(5s/2)}{\lambda}
 =\frac{25}{2\sqrt n}\le\frac12.
\]
Consequently,
\begin{equation}\label{eq:effective-local-variance-comparison}
 \frac{\lambda}{2}\le t_k\le\frac{3\lambda}{2},
 \qquad k\in\mathcal K_x.
\end{equation}
In particular, \(t_k>0\).

Put \(v=p(1-p)\) and \(z_k=(k-np)/\sqrt{nv}\).
Since \(v\ge6/25\), equation
\eqref{eq:effective-local-count-location} implies
\[
 |z_k|\le\frac5{\sqrt v}\le\frac{25}{\sqrt6}<11.
\]
The final inequality follows by squaring, since
\(25^2=625<726=11^2\cdot6\).
We may therefore apply \eqref{eq:effective-binomial-wide}
to obtain
\begin{equation}\label{eq:effective-local-binomial-lower}
 \PP(K=k)\ge\frac{e^{-100}}{\sqrt n},
 \qquad k\in\mathcal K_x.
\end{equation}

\textbf{Approximate each conditional noise sum by a normal law.}
For \(k\in\mathcal K_x\), write
\(W_k=\sum_{i=1}^n\xi_{k,i}\) for the independent
centered conditional summands described above.
The bound \(|\xi_{k,i}|\le e^{-\mathcal A}\) gives
\[
 \sum_{i=1}^n\EE|\xi_{k,i}|^3
 \le e^{-\mathcal A}\sum_{i=1}^n\EE\xi_{k,i}^2
 =e^{-\mathcal A}t_k.
\]
All the hypotheses of
\eqref{eq:independent-berry-esseen} hold: these summands
are independent and centered, have finite third absolute
moments, and have total variance \(t_k>0\).
Define
\[
 H_k(z)=\PP\{W_k/\sqrt{t_k}\le z\},\qquad
 \delta_k=\frac{e^{-\mathcal A}}{\sqrt{t_k}}.
\]
Using \(C_{\mathrm{ind}}=1\) in that inequality gives
\begin{equation}\label{eq:effective-local-conditional-be}
 |H_k(z)-\Phi(z)|\le\delta_k,\qquad z\in\R.
\end{equation}
For any fixed \(z\), let \(z_j\uparrow z\) in
\eqref{eq:effective-local-conditional-be}.
Since \(H_k(z_j)\to H_k(z-)\) and \(\Phi\) is continuous,
we also obtain
\begin{equation}\label{eq:effective-local-left-be}
 |H_k(z-)-\Phi(z)|\le\delta_k,\qquad z\in\R.
\end{equation}
We will use this bound at the upper endpoint of the
open interval.

\textbf{Obtain a uniform lower bound for the conditional interval probability.}
For either of the two pairs \((a,b)\) in the statement,
\[
 b-a=\frac14,\qquad \max(|a|,|b|)=\frac12.
\]
For \(k\in\mathcal K_x\), define the standardized endpoints
\[
 \ell_k=\frac{x+a-k}{\sqrt{t_k}},\qquad
 u_k=\frac{x+b-k}{\sqrt{t_k}}.
\]
For \(c\in\{a,b\}\),
\eqref{eq:effective-local-selected-counts} and
\eqref{eq:effective-local-variance-comparison} give
\[
 \begin{aligned}
 \frac{|x+c-k|}{\sqrt{t_k}}
 &\le\frac{r+1/2}{\sqrt{\lambda/2}}\\
 &=\sqrt2\left\{\max(1,\lambda^{-1/2})
                        +\frac12\lambda^{-1/2}\right\}\\
 &\le\sqrt2\left(10^6+\frac12\,10^6\right)
 =\frac{3\cdot10^6}{\sqrt2}=:A_0.
 \end{aligned}
\]
Here \(\lambda\ge10^{-12}\) implies
\(\lambda^{-1/2}\le10^6\), and hence also
\(\max(1,\lambda^{-1/2})\le10^6\).
Thus \(\ell_k,u_k\in[-A_0,A_0]\).
Since \(\phi(z)\ge\phi(A_0)\) throughout this interval,
\begin{equation}\label{eq:effective-local-gaussian-mass}
 \begin{aligned}
 \Phi(u_k)-\Phi(\ell_k)
 &=\int_{\ell_k}^{u_k}\phi(z)\,dz\\
 &\ge\frac{b-a}{\sqrt{t_k}}\phi(A_0)
 \ge\frac{c_1}{\sqrt\lambda},\\
 c_1&:=\frac14\sqrt{\frac23}\,\phi(A_0)
 =\frac{e^{-(9/4)10^{12}}}{4\sqrt{3\pi}}.
 \end{aligned}
\end{equation}
The last lower bound uses \(t_k\le3\lambda/2\), and
the formula for \(c_1\) uses \(A_0^2/2=(9/4)10^{12}\).
To compare this constant with the noise bound, note that
\(4\sqrt{3\pi}<16<e^4\). Consequently,
\[
 c_1>e^{-(9/4)10^{12}-4}>e^{-3\cdot10^{12}}.
\]
Also \(4\sqrt2<e^3\) and
\(\mathcal A=10^{14}\) by
\eqref{eq:effective-threshold-parameters}, so
\begin{equation}\label{eq:effective-local-noise-margin}
 \frac{c_1}{4\sqrt2}
 >e^{-3\cdot10^{12}-3}
 >e^{-10^{14}}=e^{-\mathcal A}.
\end{equation}

The conditional interval probability has the exact form
\[
 \PP\{W_k\in(x+a-k,x+b-k)\}
 =H_k(u_k-)-H_k(\ell_k).
\]
This uses the left limit at the upper endpoint and the
right-continuous value at the lower endpoint, so both
endpoints are excluded.
Apply \eqref{eq:effective-local-left-be} at \(u_k\)
and \eqref{eq:effective-local-conditional-be} at \(\ell_k\).
Together with \eqref{eq:effective-local-gaussian-mass}
and \(t_k\ge\lambda/2\), these inequalities give
\[
 \begin{aligned}
 \PP\{W_k\in(x+a-k,x+b-k)\}
 &\ge\Phi(u_k)-\Phi(\ell_k)-2\delta_k\\
 &\ge\frac{c_1}{\sqrt\lambda}
          -\frac{2e^{-\mathcal A}}{\sqrt{t_k}}\\
 &\ge\frac{c_1-2\sqrt2 e^{-\mathcal A}}{\sqrt\lambda}.
 \end{aligned}
\]
Equation \eqref{eq:effective-local-noise-margin} now gives
\begin{equation}\label{eq:effective-local-conditional-lower}
 \PP\{W_k\in(x+a-k,x+b-k)\}
 \ge\frac{c_1}{2\sqrt\lambda},
 \qquad k\in\mathcal K_x.
\end{equation}

\textbf{Sum over the selected counts.}
The interval \([x-r,x+r]\) contains
\(\lfloor x+r\rfloor-\lceil x-r\rceil+1\) integers.
Using \(\lfloor y\rfloor\ge y-1\),
\(\lceil y\rceil\le y+1\), and \(r\ge1\), we obtain
\[
 \#\mathcal K_x\ge2r-1\ge r\ge\sqrt\lambda.
\]
Condition on \(K\), and keep only the terms with
\(k\in\mathcal K_x\).
Equations \eqref{eq:effective-local-binomial-lower}
and \eqref{eq:effective-local-conditional-lower} yield
\[
 \begin{aligned}
 \PP\{S_n\in(x+a,x+b)\}
 &=\sum_{k=0}^n\PP(K=k)
       \PP\{W_k\in(x+a-k,x+b-k)\}\\
 &\ge\sum_{k\in\mathcal K_x}
       \frac{e^{-100}}{\sqrt n}\frac{c_1}{2\sqrt\lambda}\\
 &\ge\frac{e^{-100}c_1}{2\sqrt n}.
 \end{aligned}
\]
Since \(c_1>e^{-3\cdot10^{12}}\) and \(1/2>e^{-1}\),
\[
 \frac{e^{-100}c_1}{2}
 >e^{-3\cdot10^{12}-101}
 >e^{-4\cdot10^{12}}.
\]
Every bound above is uniform over \(|x-np|\le3\sqrt n\)
and valid for both stated pairs \((a,b)\).
Taking the infimum over these \(x\) proves
\eqref{eq:effective-local-mass}.
\end{proof}

We now use \eqref{eq:effective-local-mass} to obtain a
strict loss in each side of the uniform-jitter comparison.
The resulting bound, together with
Lemma~\ref{lem:effective-small-variance}, proves
\(R_n(B+U)<\ce\) for all \(n\ge N_*\).
We then translate this result into the support and
cluster-probability conditions of
Proposition~\ref{prop:clusters}.

\begin{proposition}\label{prop:effective-clusters}
If \(p\in I\), \(\EE(U\mid B)=0\), and
\(|U|\le e^{-\mathcal A}\), then
\(R_n(B+U)<\ce\) for every integer \(n\ge N_*\).
Consequently, the parameters in Proposition~\ref{prop:clusters}
may be taken to be \(\eta=\eta_*\) and \(N=N_*\).
\end{proposition}

\begin{proof}
\textbf{Separate the small-variance case and fix the error bounds.}
Write
\[
 Y=B+U,\qquad \varepsilon=e^{-\mathcal A},\qquad
 s=\EE U^2,\qquad \lambda=ns,\qquad
 v=p(1-p),\qquad \sigma^2=v+s.
\]
We have \(0\le s\le\varepsilon^2=e^{-2\mathcal A}\).
By \eqref{eq:effective-threshold-parameters}, we have
\(n\ge N_*=\lceil e^{\mathcal A}\rceil\) and
\(\mathcal A=10^{14}>100\log10\). Hence
\[
 n\ge e^{\mathcal A}>10^{100},\qquad
 n^{-1/2}\le e^{-\mathcal A/2},\qquad
 \varepsilon<10^{-12}.
\]
Thus, if \(\lambda\le10^{-12}\),
Lemma~\ref{lem:effective-small-variance} gives
\(R_n(Y)\le R_n(B)<\ce\).
For the rest of the first assertion, assume
\(\lambda>10^{-12}\), and set
\[
 D=e^{-5\cdot10^{12}},\qquad T=e^{10^{13}}.
\]
The assumptions of Lemma~\ref{lem:effective-cluster-jitter}
hold: \(T\ge10\), \(n\ge10^6\), and
\[
 100T\varepsilon=100e^{-9\cdot10^{13}}<1.
\]
For the last inequality, \(e^{9\cdot10^{13}}>100\)
follows from \(e^r\ge1+r\).

Let \(E_n=\mathcal E(n,e^{-\mathcal A},T)\), with
\(\mathcal E\) as in \eqref{eq:effective-cluster-jitter}.
Its three terms satisfy
\[
 \begin{aligned}
 \frac{4(1+\log T)}{\sqrt n}
 &\le4(1+10^{13})e^{-5\cdot10^{13}},\\
 2\varepsilon^2T^2(1+\log T)
 &=2(1+10^{13})e^{-18\cdot10^{13}},\\
 \frac{50}{T}&=50e^{-10^{13}}.
 \end{aligned}
\]
Each displayed prefactor is less than \(e^{100}\).
After division by \(D\), the three exponential factors
are respectively
\[
 e^{-45\cdot10^{12}},\qquad
 e^{-175\cdot10^{12}},\qquad
 e^{-5\cdot10^{12}}.
\]
Hence each term divided by \(D\) is at most
\(e^{100-5\cdot10^{12}}<1/12\).
The latter inequality follows from
\(e^{5\cdot10^{12}-100}\ge1+5\cdot10^{12}-100>12\).
Likewise,
\[
 \frac{2/\sqrt n+100s}{D}
 \le2e^{-45\cdot10^{12}}+100e^{-195\cdot10^{12}}
 \le2e^{100-5\cdot10^{12}}<\frac14.
\]
We have therefore established the two error bounds
\begin{equation}\label{eq:effective-clusters-error-budgets}
 E_n\le\frac D4,\qquad
 \frac2{\sqrt n}+100s\le\frac D4.
\end{equation}

\textbf{Use local interval probabilities to improve both jitter bounds.}
Let \(Y_1,\ldots,Y_n\) be independent copies of \(Y\), and
write
\[
 S_n=\sum_{i=1}^nY_i,\qquad F_n(x)=\PP(S_n\le x),
 \qquad z=\frac{x-np}{\sigma\sqrt n}.
\]
As in Lemma~\ref{lem:effective-cluster-jitter}, let \(H\) be
independent and uniform on \([-1/2,1/2]\), and put
\(J_n(x)=\PP(S_n+H\le x)\).
The exact identities \eqref{eq:cluster-jitter-cdf-identities}
give, for every real \(x\),
\begin{equation}\label{eq:effective-clusters-jitter-identities}
 \begin{aligned}
 J_n(x+1/2)-F_n(x)
 &=\EE[(x+1-S_n)\ind_{\{x<S_n<x+1\}}],\\
 F_n(x)-J_n(x-1/2)
 &=\EE[(S_n-x+1)\ind_{\{x-1<S_n<x\}}]
   +\PP(S_n=x).
 \end{aligned}
\end{equation}
The first weight is at least \(1/2\) on
\((x+1/4,x+1/2)\); the second is at least \(1/2\)
on \((x-1/2,x-1/4)\).
The atom term in the second identity is nonnegative.

The variance bound \eqref{eq:effective-cluster-variance-bounds}
implies \(\sigma^2\le1/4+10^{-6}\), and hence
\(5\sigma<3\).
If \(|z|\le5\), then
\[
 |x-np|=\sigma\sqrt n\,|z|
 \le5\sigma\sqrt n<3\sqrt n.
\]
Lemma~\ref{lem:effective-local-mass} therefore applies
to both quarter intervals: its conditions follow from
\(p\in I\), conditional centering, \(|U|\le e^{-\mathcal A}\),
\(n\ge10^{100}\), and \(\lambda>10^{-12}\).
Consequently, \eqref{eq:effective-clusters-jitter-identities}
gives
\[
 \begin{aligned}
 J_n(x+1/2)-F_n(x)
 &\ge\frac{e^{-4\cdot10^{12}}}{2\sqrt n}
 \ge\frac D{\sqrt n},\\
 F_n(x)-J_n(x-1/2)
 &\ge\frac{e^{-4\cdot10^{12}}}{2\sqrt n}
 \ge\frac D{\sqrt n}.
 \end{aligned}
\]
Here \(e^{-4\cdot10^{12}}/2\ge D\) is equivalent to
\(e^{10^{12}}\ge2\).
Equivalently,
\begin{equation}\label{eq:effective-clusters-strict-sandwich}
 J_n(x-1/2)+\frac D{\sqrt n}
 \le F_n(x)
 \le J_n(x+1/2)-\frac D{\sqrt n},
 \qquad |z|\le5.
\end{equation}

\textbf{Expand the shifted jitter bounds and obtain the central envelope.}
Put
\[
 \gamma=\frac{\EE(Y-p)^3}{\sigma^3},\qquad
 \psi(z)=(1-z^2)\phi(z),\qquad
 a=\frac1{2\sigma\sqrt n}.
\]
Lemma~\ref{lem:effective-cluster-jitter} gives a function
\(r_n\) satisfying
\[
 J_n(x)=\Phi(z)+\frac{\gamma}{6\sqrt n}\psi(z)+r_n(x),
 \qquad \sqrt n\,\|r_n\|_\infty\le E_n.
\]
The standardized argument at \(x+\theta/2\) is
\(z+\theta a\), for \(\theta\in\{-1,1\}\).
Taylor's theorem for \(\Phi\) and the mean value theorem
for \(\psi\) give
\[
 \begin{aligned}
 |\Phi(z+\theta a)-\Phi(z)-\theta a\phi(z)|
 &\le\frac{\|\phi'\|_\infty a^2}{2},\\
 |\psi(z+\theta a)-\psi(z)|
 &\le\|\psi'\|_\infty a.
 \end{aligned}
\]
The derivative bounds from
Lemma~\ref{lem:effective-selection} are
\(\|\phi'\|_\infty<1/4\) and
\(\|\psi'\|_\infty=\|\phi'''\|_\infty<5/4\).
Equations \eqref{eq:effective-cluster-variance-bounds}
and \eqref{eq:effective-cluster-moment-bounds} also give
\[
 \sigma^2\ge6/25,\qquad \sigma>0.48,\qquad
 |\gamma|\le\beta(Y)<1.11.
\]
Thus the Taylor error, after multiplication by \(\sqrt n\),
is bounded by
\[
 \begin{aligned}
 \frac1{\sqrt n}
 \left\{\frac{\|\phi'\|_\infty}{8\sigma^2}
        +\frac{|\gamma|\|\psi'\|_\infty}{12\sigma}\right\}
 &\le\frac1{\sqrt n}
   \left\{\frac{1}{32(6/25)}
          +\frac{1.11(5/4)}{12(0.48)}\right\}\\
 &<\frac2{\sqrt n}.
 \end{aligned}
\]
It follows that
\begin{equation}\label{eq:effective-clusters-shifted-expansions}
 \begin{aligned}
 \sqrt n\{J_n(x+\theta/2)-\Phi(z)\}
 &=\frac{\theta}{2\sigma}\phi(z)
   +\frac{\gamma}{6}\psi(z)+r_{n,\theta}(x),\\
 |r_{n,\theta}(x)|&\le E_n+\frac2{\sqrt n},
 \qquad x\in\R,\quad \theta\in\{-1,1\}.
 \end{aligned}
\end{equation}
Use \(\theta=1\) in the upper bound of
\eqref{eq:effective-clusters-strict-sandwich}, and
\(\theta=-1\) in its lower bound.
For \(|z|\le5\), this gives both signed estimates
\[
 \begin{aligned}
 \sqrt n\{F_n(x)-\Phi(z)\}
 &\le\frac{\phi(z)}{2\sigma}
      +\frac{\gamma}{6}\psi(z)-D+E_n+\frac2{\sqrt n},\\
 \sqrt n\{\Phi(z)-F_n(x)\}
 &\le\frac{\phi(z)}{2\sigma}
      -\frac{\gamma}{6}\psi(z)-D+E_n+\frac2{\sqrt n}.
 \end{aligned}
\]
Take the maximum of these two bounds and use
\(\|\phi\|_\infty=\|\psi\|_\infty=\phi_0\), as proved in
\eqref{eq:jitter-envelope-gaussian-suprema}.
We obtain
\begin{equation}\label{eq:effective-central-envelope}
 \sqrt n\,|F_n(x)-\Phi(z)|
 \le\phi_0\left(\frac1{2\sigma}+\frac{|\gamma|}{6}\right)
       -D+E_n+\frac2{\sqrt n},
 \qquad |z|\le5.
\end{equation}

\textbf{Compare the moment normalization with that of the Bernoulli law.}
Write
\[
 \tau=p^2+(1-p)^2,\qquad d=1-2p,\qquad
 \beta_B=\frac{\tau}{\sqrt v},\qquad
 \gamma_B=\frac d{\sqrt v}.
\]
Set \(A=B-p\).
Conditional centering gives
\[
 \EE(A^2U)=0,\qquad \EE(A|A|U)=0.
\]
Moreover, \(|A|\ge2/5>\varepsilon\), so \(A+U\) has
the same sign as \(A\).
Expanding the two cubic moments therefore gives
\[
 \begin{aligned}
 \EE(Y-p)^3
 &=vd+3\EE(AU^2)+\EE U^3,\\
 \EE|Y-p|^3
 &=v\tau+3\EE(|A|U^2)
            +\EE[\operatorname{sgn}(A)U^3].
 \end{aligned}
\]
Since \(|A|\le1\) and \(\EE|U|^3\le\varepsilon s\),
both raw moment differences are bounded in absolute value by
\((3+\varepsilon)s\le4s\).
The mean value theorem for \(t\mapsto t^{-3/2}\) gives
\[
 |(v+s)^{-3/2}-v^{-3/2}|
 \le\frac32v^{-5/2}s.
\]
Also, \(v\tau\le v\) and \(|vd|\le v\).
Dividing the moment identities by \((v+s)^{3/2}\)
and using the last two bounds, we obtain
\[
 \begin{aligned}
 |\beta(Y)-\beta_B|
 &\le\frac{4s}{v^{3/2}}
       +v\frac32v^{-5/2}s
 =\frac{11s}{2v^{3/2}},\\
 |\gamma-\gamma_B|
 &\le\frac{4s}{v^{3/2}}
       +v\frac32v^{-5/2}s
 =\frac{11s}{2v^{3/2}}.
 \end{aligned}
\]
Because \(v\ge6/25\), we have
\(v^{3/2}\ge(6/25)^{3/2}>1/10\).
Thus
\begin{equation}\label{eq:effective-clusters-moment-comparison}
 |\beta(Y)-\beta_B|\le100s,\qquad
 |\gamma-\gamma_B|\le100s.
\end{equation}

Since \(p\le9/20<1/2\), \(\gamma_B>0\).
Using \(\sigma\ge\sqrt v\),
\(|\gamma|\le\gamma_B+100s\), and \(\phi_0<2/5\),
we conclude that
\begin{equation}\label{eq:effective-clusters-bernoulli-envelope}
 \phi_0\left(\frac1{2\sigma}+\frac{|\gamma|}{6}\right)
 \le L_B+7s,\qquad
 L_B:=\frac{\phi_0(3+d)}{6\sqrt v}.
\end{equation}
Indeed, the additional term is at most
\((100\phi_0/6)s<(20/3)s<7s\).
To compare \(L_B\) with \(\ce\beta_B\), put
\(c_*=3+\sqrt{10}\).
The definition \eqref{eq:esseen-parameters} of \(\pe\) gives
\[
 2c_*\pe=c_*-1,\qquad
 2c_*\pe^2=c_*-4.
\]
Consequently,
\[
 \begin{aligned}
 c_*\tau-(3+d)
 &=2c_*p^2-2(c_*-1)p+c_*-4\\
 &=2c_*(p-\pe)^2\ge0.
 \end{aligned}
\]
Since \(\ce=c_*\phi_0/6\) by \eqref{eq:ce}, multiplication
by \(\phi_0/(6\sqrt v)\) gives \(L_B\le\ce\beta_B\).
Combining this with
\eqref{eq:effective-clusters-moment-comparison},
\eqref{eq:effective-clusters-bernoulli-envelope}, and
\(\ce<1/2\), we obtain
\[
 \begin{aligned}
 \phi_0\left(\frac1{2\sigma}+\frac{|\gamma|}{6}\right)
 &\le\ce\beta_B+7s\\
 &\le\ce\beta(Y)+(100\ce+7)s
 \le\ce\beta(Y)+100s.
 \end{aligned}
\]

\textbf{Combine the central and distant threshold bounds.}
Substitute the preceding estimate into
\eqref{eq:effective-central-envelope} and use
\eqref{eq:effective-clusters-error-budgets}.
For \(|z|\le5\), we obtain
\[
 \begin{aligned}
 \sqrt n\,|F_n(x)-\Phi(z)|
 &\le\ce\beta(Y)-D+E_n+\frac2{\sqrt n}+100s\\
 &\le\ce\beta(Y)-\frac D2.
 \end{aligned}
\]
Since \(0<\beta(Y)<1.11<2\), division by \(\beta(Y)\)
gives
\begin{equation}\label{eq:effective-clusters-central-gap}
 \frac{\sqrt n}{\beta(Y)}|F_n(x)-\Phi(z)|
 \le\ce-\frac{D}{2\beta(Y)}
 \le\ce-\frac D4,\qquad |z|\le5.
\end{equation}
For \(|z|>5\), the nonuniform estimate
\eqref{eq:effective-nonuniform} instead yields
\[
 \frac{\sqrt n}{\beta(Y)}|F_n(x)-\Phi(z)|
 \le\frac{17.36}{1+|z|^3}
 <\frac{17.36}{126}<0.14.
\]
Here the hypotheses of that estimate hold because \(Y\)
is bounded and \(\Var Y=\sigma^2\ge6/25>0\).
Also, \(\ce>2/5\) and \(0<D<1\) imply
\[
 \ce-\frac D4>\frac25-\frac14=0.15>0.14.
\]
Thus the normalized error is at most \(\ce-D/4\) at
every threshold \(x\).
Taking the supremum over \(x\in\R\) gives
\(R_n(Y)\le\ce-D/4<\ce\) in the case
\(\lambda>10^{-12}\).
Together with the small-variance case at the start,
this proves the first assertion.

\textbf{Return to the two intervals in Proposition~\ref{prop:clusters}.}
To establish the second assertion, let \(X\sim Q\) satisfy
\[
 \supp Q\subset[-\eta_*,\eta_*]\cup[1-\eta_*,1+\eta_*],
 \qquad
 \bigl|Q([1-\eta_*,1+\eta_*])-\pe\bigr|<\eta_*.
\]
Define
\[
 B=\ind_{\{X\in[1-\eta_*,1+\eta_*]\}},\qquad
 p=\PP(B=1).
\]
The intervals are disjoint because \(\eta_*<1/4\), and
the support assumption places \(X\) in their union almost
surely.
Equation \eqref{eq:effective-threshold-parameters} gives
\(\eta_*<0.01\).
Since \(3.16^2<10<3.18^2\),
\eqref{eq:esseen-parameters} gives \(0.41<\pe<0.42\).
Consequently,
\[
 |p-\pe|<\eta_*
 \quad\Longrightarrow\quad
 0.4<p<0.43,\qquad p\in I.
\]
In particular, both clusters have positive probability.
Since \(X\) is bounded, the conditional means below are
finite and well defined.
Define
\[
 m_b=\EE(X\mid B=b)\quad(b=0,1),\qquad d_0=m_1-m_0.
\]
The support assumption gives
\[
 m_0\in[-\eta_*,\eta_*],\qquad
 m_1\in[1-\eta_*,1+\eta_*],\qquad
 d_0\ge1-2\eta_*>0.
\]
Set
\[
 \widetilde Y=\frac{X-m_0}{d_0},\qquad
 U=\frac{X-m_B}{d_0}.
\]
Since \(m_B=m_0+d_0B\), we have
\(\widetilde Y=B+U\).
Also,
\[
 \EE(U\mid B)=0,\qquad
 |U|\le\frac{2\eta_*}{d_0}
 \le\frac{2\eta_*}{1-2\eta_*}
 \le4\eta_*\le e^{-\mathcal A}.
\]
The penultimate inequality uses \(\eta_*<1/4\).
The last uses \(\eta_*=e^{-2\mathcal A}\) and
\(e^{\mathcal A}>4\).
Both \(X\) and \(\widetilde Y\) are bounded, and
\[
 \Var(\widetilde Y)=p(1-p)+\EE U^2>0,\qquad
 \Var X=d_0^2\Var(\widetilde Y)>0.
\]
Thus the first assertion applies to \(\widetilde Y\),
and affine invariance \eqref{eq:raw-normalization} gives
\[
 R_n(Q)=R_n(\widetilde Y)<\ce,\qquad n\ge N_*.
\]
Finally, if the upper-cluster probability is within
\(\eta_*\) of \(1-\pe\), replace \(X\) by \(1-X\).
This affine transformation exchanges the two intervals,
and the new upper-cluster probability is \(1-p\), which
satisfies \(|(1-p)-\pe|<\eta_*\).
The preceding argument and affine invariance therefore
give the same bound for \(Q\).
Since \(\eta_*\in(0,1/4)\) and \(N_*<\infty\), these
parameters satisfy the conclusion of
Proposition~\ref{prop:clusters}.
\end{proof}

\subsection{Stability of the lattice moment inequality}

We now quantify the equality case of Esseen's moment
inequality \eqref{eq:esseen-moment} for bounded lattice
laws. If the deficit is small, we will show that the law
is close in \(W_1\) to \(\Pe\) or its reflection, and that
its lattice spacing is close to \(h_{\mathrm E}\).
Lemma~\ref{lem:effective-identification} will apply this
estimate to the standardization of the finite lattice law
constructed in Lemma~\ref{lem:effective-global-jitter}.
The proof first controls the mass outside the two lattice
points nearest zero, and then identifies the probabilities
and spacing of the remaining two-point law.

\begin{lemma}\label{lem:effective-lattice-stability}
Let \(Q\in\Pthree\) be supported in \([-15,15]\), with
\(\beta(Q)\le2\), and suppose that its support is contained
in a translate of \(h\mathbb Z\), where \(h>0\).
Put \(c_*=3+\sqrt{10}\).
If
\[
 0\le D:=c_*\beta(Q)-|\kappa(Q)|-3h
 \le\delta\le10^{-6},
\]
then, for some \(\epsilon\in\{-1,1\}\),
\begin{equation}\label{eq:effective-lattice-stability}
 W_1(Q,\epsilon\Pe)\le1000\sqrt\delta,\qquad
 |h-h_{\mathrm E}|\le1000\sqrt\delta.
\end{equation}
Here \(\epsilon\Pe\) denotes the law of \(\epsilon Z\)
for \(Z\sim\Pe\).
We do not assume that \(h\) is the maximal lattice span.
\end{lemma}

\begin{proof}
\textbf{Choose the orientation and express the deficit through pairs of support points.}
Reflection preserves membership in \(\Pthree\), the support
bound \([-15,15]\), the third absolute moment, and the
given lattice spacing \(h\); it changes the sign of
\(\kappa(Q)\), so the deficit \(D\) is unchanged.
Reflect \(Q\) if necessary so that \(\kappa(Q)\ge0\),
and continue to denote the resulting law by \(Q\).
We will undo this reflection at the end.
Let \(X\sim Q\), write \(c=3+\sqrt{10}\), and fix a
containing lattice
\[
 \Lambda=t+h\mathbb Z,\qquad \supp Q\subset\Lambda.
\]
Put \(X_-=(-X)_+\).
Since \(X=X_+-X_-\) and \(\EE X=0\),
\[
 m:=\EE X_+=\EE X_-.
\]
These expectations are finite because \(X\) is bounded.
If \(m=0\), both nonnegative variables \(X_+\) and
\(X_-\) vanish almost surely, contradicting
\(\EE X^2=1\).
Thus \(m>0\).
For every Borel set \(A\subset(0,\infty)\), define
\[
 \nu_-(A)=\frac{\EE[X_-\ind_{\{X_-\in A\}}]}{m},
 \qquad
 \nu_+(A)=\frac{\EE[X_+\ind_{\{X_+\in A\}}]}{m},
 \qquad \pi=\nu_-\otimes\nu_+.
\]
Both measures are nonnegative and have total mass one;
in particular, \(\pi\) is a probability measure on
\((0,\infty)^2\).
For \(j=1,2\), their definitions and the product structure
give
\[
 m\int (a')^j\,d\pi(a',b')=\EE X_-^{j+1},
 \qquad
 m\int (b')^j\,d\pi(a',b')=\EE X_+^{j+1}.
\]
Adding or subtracting these identities yields
\begin{equation}\label{eq:effective-lattice-pair-moments}
 \begin{aligned}
 m\int(a'+b')\,d\pi&=\EE X^2=1,\\
 m\int((a')^2+(b')^2)\,d\pi&=\EE|X|^3=\beta(Q),\\
 m\int((b')^2-(a')^2)\,d\pi&=\EE X^3=\kappa(Q).
 \end{aligned}
\end{equation}
Since \(c-2>0\), \(c-4>0\), and
\((c-2)(c-4)=9\), expansion of the square gives
\[
 \begin{aligned}
 &(\sqrt{c-2}\,a'-\sqrt{c-4}\,b')^2
       +3(a'+b')(a'+b'-h)\\
 &\qquad=(c-2)(a')^2+(c-4)(b')^2-6a'b'\\
 &\qquad\quad+3(a')^2+6a'b'+3(b')^2-3h(a'+b')\\
 &\qquad=(c+1)(a')^2+(c-1)(b')^2-3h(a'+b').
 \end{aligned}
\]
Integrating and applying
\eqref{eq:effective-lattice-pair-moments}, with
\(|\kappa(Q)|=\kappa(Q)\), proves the exact identity
\begin{equation}\label{eq:effective-lattice-deficit}
 \begin{aligned}
 D=m\int\bigl[
 &(\sqrt{c-2}\,a'-\sqrt{c-4}\,b')^2\\
 &+3(a'+b')(a'+b'-h)\bigr]\,d\pi(a',b').
 \end{aligned}
\end{equation}
For \(\pi\)-almost every pair, \(-a'\) and \(b'\)
belong to \(\Lambda\).
Their positive difference \(a'+b'\) is therefore a
positive integer multiple of \(h\), so \(a'+b'\ge h\).
Both terms in the integrand of
\eqref{eq:effective-lattice-deficit} are consequently
nonnegative.
This conclusion uses only that \(\Lambda\) contains
the support, without assuming that \(h\) is maximal.

\textbf{Locate two adjacent lattice points with positive mass.}
H\"older's inequality and the unit second moment give
\[
 1=\EE|X|^2\le(\EE|X|^3)^{2/3}=\beta(Q)^{2/3},
 \qquad \beta(Q)\ge1.
\]
Also \(0\le\kappa(Q)\le\beta(Q)\).
Since \(3h=c\beta(Q)-\kappa(Q)-D\), the bounds
\(D\le\delta\), \(D\ge0\), and \(\beta(Q)\le2\) imply
\begin{equation}\label{eq:effective-lattice-spacing-bounds}
 \frac{c-1-\delta}{3}\le h\le\frac{2c}{3}.
\end{equation}
Indeed,
\[
 c\beta(Q)-\kappa(Q)-D
 \ge(c-1)\beta(Q)-\delta\ge c-1-\delta,
 \qquad
 c\beta(Q)-\kappa(Q)-D\le c\beta(Q)\le2c.
\]
Using \(3.16<\sqrt{10}<4\) and \(\delta\le10^{-6}\)
in \eqref{eq:effective-lattice-spacing-bounds} gives
\[
 h>\frac{5.16-10^{-6}}{3}>1.7,\qquad
 h<\frac{14}{3}<5.
\]

Suppose first that \(0\in\Lambda\).
Then \(\Lambda=h\mathbb Z\); every negative point of
\(\Lambda\) has absolute value at least \(h\), and every
positive point has value at least \(h\).
Thus \(a'+b'\ge2h\) for \(\pi\)-almost every pair.
Dropping the square in
\eqref{eq:effective-lattice-deficit} and using
\(a'+b'-h\ge h\), we obtain
\[
 D\ge3hm\int(a'+b')\,d\pi=3h>5.
\]
The equality uses
\eqref{eq:effective-lattice-pair-moments}.
This contradicts \(D\le\delta\le10^{-6}\), and proves
that \(0\notin\Lambda\).

Let \(-a\) and \(b\) be the two consecutive points of
\(\Lambda\) immediately below and above zero.
They satisfy
\[
 a>0,\qquad b>0,\qquad a+b=h,
\]
and all negative and positive lattice points are,
respectively,
\[
 -(a+jh),\qquad b+\ell h,\qquad j,\ell\in\{0,1,2,\ldots\}.
\]
If \(Q(\{-a\})=0\), then \(a'\ge a+h\) for
\(\nu_-\)-almost every \(a'\), while \(b'\ge b\) for
\(\nu_+\)-almost every \(b'\).
Consequently \(a'+b'\ge2h\) for \(\pi\)-almost every
pair, and the preceding calculation again gives
\(D\ge3h>5\), a contradiction.
If \(Q(\{b\})=0\), the same conclusion follows from
\(a'\ge a\) and \(b'\ge b+h\).
Therefore
\begin{equation}\label{eq:effective-lattice-principal-atoms}
 Q(\{-a\})>0,\qquad Q(\{b\})>0,\qquad
 0<a,b<h<5.
\end{equation}
This also shows that \(h\) is the maximal lattice span
under the present hypotheses.
Indeed, for any other containing lattice with spacing
\(H>0\), the difference of the two support points
\(b-(-a)=h\) must be a positive integer multiple of
\(H\), and therefore \(H\le h\).
The containing lattice \(\Lambda\) itself has spacing \(h\).

\textbf{Bound the second moment carried by all other lattice points.}
Define
\[
 A_{\mathrm{out}}=\Lambda\setminus\{-a,b\},\qquad
 r=Q(A_{\mathrm{out}}),\qquad
 s_1=\EE[|X|\ind_{\{X\in A_{\mathrm{out}}\}}],\qquad
 s_2=\EE[X^2\ind_{\{X\in A_{\mathrm{out}}\}}].
\]
Write \(L=a'+b'\) for the sum of a pair.
Since \(L\) is a positive integer multiple of \(h\),
\[
 L(L-h)\ge hL\ind_{\{L\ge2h\}}.
\]
If either \(-a'\) or \(b'\) belongs to
\(A_{\mathrm{out}}\), the lattice representation above
shows that \(L\ge2h\).
It follows for \(\pi\)-almost every pair that
\[
 L\ind_{\{L\ge2h\}}
 \ge a'\ind_{\{-a'\in A_{\mathrm{out}}\}}
      +b'\ind_{\{b'\in A_{\mathrm{out}}\}}.
\]
The size-biased definitions and the fact that each
factor of \(\pi\) has total mass one give
\[
 \begin{aligned}
 m\int a'\ind_{\{-a'\in A_{\mathrm{out}}\}}\,d\pi
 &=\EE[X_-^2\ind_{\{X\in A_{\mathrm{out}}\}}],\\
 m\int b'\ind_{\{b'\in A_{\mathrm{out}}\}}\,d\pi
 &=\EE[X_+^2\ind_{\{X\in A_{\mathrm{out}}\}}].
 \end{aligned}
\]
Their sum is \(s_2\).
Dropping the square in
\eqref{eq:effective-lattice-deficit} and then applying
these inequalities yields
\[
 D\ge3hm\int L\ind_{\{L\ge2h\}}\,d\pi\ge3h s_2.
\]
Since \(h>1.7\), we conclude that
\begin{equation}\label{eq:effective-extra-lattice-mass}
 s_2\le\frac{D}{3h}\le\delta.
\end{equation}
Every point of \(A_{\mathrm{out}}\) has absolute value
at least \(\min(a+h,b+h)>h>1\).
Thus the pointwise inequalities
\[
 \ind_{\{X\in A_{\mathrm{out}}\}}
 \le |X|\ind_{\{X\in A_{\mathrm{out}}\}}
 \le X^2\ind_{\{X\in A_{\mathrm{out}}\}}
\]
and \eqref{eq:effective-extra-lattice-mass} give
\begin{equation}\label{eq:effective-lattice-outside-moments}
 r\le s_1\le s_2\le\delta.
\end{equation}

\textbf{Condition on the two principal atoms and control the normalization.}
By \eqref{eq:effective-lattice-outside-moments},
\[
 0\le r\le s_1\le s_2\le\delta\le10^{-6}.
\]
In particular, \(1-r>0\).
Let \(W\) have the conditional law of \(X\) on
\(\{-a,b\}\), and put
\(\mu=\EE W\) and \(v_0=\Var W\).
Since \(\EE X=0\) and \(\EE X^2=1\), conditioning on
the principal atoms gives
\begin{equation}\label{eq:effective-lattice-conditional-moments}
 \mu=-\frac{\EE[X\ind_{\{X\in A_{\mathrm{out}}\}}]}{1-r},
 \qquad
 \EE W^2=\frac{1-s_2}{1-r}.
\end{equation}
Consequently,
\[
 |\mu|
 \le\frac{s_1}{1-r}
 \le\frac{\delta}{1-\delta}\le2\delta.
\]
Here \(\delta\le10^{-6}<1/2\), so \(1/(1-\delta)\le2\).
Also, \(s_2\ge r\), and therefore
\[
 \left|\EE W^2-1\right|
 =\frac{s_2-r}{1-r}
 \le\frac{\delta}{1-\delta}\le2\delta.
\]
Using \(v_0=\EE W^2-\mu^2\), we obtain
\begin{equation}\label{eq:effective-lattice-conditional-variance}
 |\mu|\le2\delta,\qquad
 |v_0-1|
 \le2\delta+4\delta^2\le5\delta,\qquad
 v_0\ge1-5\delta>0.
\end{equation}
The second inequality uses \(4\delta^2\le3\delta\),
which holds for \(0\le\delta\le10^{-6}\).

\textbf{Transfer the moment deficit to the centered conditional law.}
Define
\[
 F_h(x)=c|x|^3-x^3-3hx^2.
\]
For \(|x|\ge h\), the inequality \(x^3\le|x|^3\) gives
\[
 \begin{aligned}
 F_h(x)
 &\ge(c-1)|x|^3-3hx^2\\
 &=x^2\{(c-1)|x|-3h\}
 \ge(c-4)hx^2\ge0.
 \end{aligned}
\]
Here \(c=3+\sqrt{10}>4\).
Every point in \(A_{\mathrm{out}}\) has absolute value
at least \(h\), so \(F_h(X)\ind_{\{X\in A_{\mathrm{out}}\}}\ge0\).
Since the law has already been reflected so that
\(\kappa(Q)\ge0\), its deficit is
\[
 \EE F_h(X)=c\beta(Q)-\kappa(Q)-3h=D.
\]
Separating the contribution of the principal atoms,
we have
\[
 D=(1-r)\EE F_h(W)
       +\EE[F_h(X)\ind_{\{X\in A_{\mathrm{out}}\}}].
\]
The second term is nonnegative. It follows that
\begin{equation}\label{eq:effective-lattice-conditional-deficit}
 \EE F_h(W)\le\frac{D}{1-r}
 \le\frac{\delta}{1-\delta}\le2\delta.
\end{equation}

Because \(a,b>0\) and \(a+b=h<5\), we have \(|W|<5\).
By \eqref{eq:effective-lattice-conditional-variance},
every point between \(W\) and \(W-\mu\) has absolute
value at most \(5+2\delta\).
Each of the functions \(g(x)=|x|^3\) and \(g(x)=x^3\)
is differentiable on \(\R\), with
\(|g'(x)|=3x^2\).
The mean value theorem therefore gives
\[
 |g(W-\mu)-g(W)|
 \le3(5+2\delta)^2|\mu|
 \le6\delta(5+2\delta)^2.
\]
Since \(5+2\delta\le5.1\) and
\(6(5.1)^2=156.06<200\), taking expectations yields
\begin{equation}\label{eq:effective-lattice-centering-cubics}
 \begin{aligned}
 \left|\EE|W-\mu|^3-\EE|W|^3\right|&\le200\delta,\\
 \left|\EE(W-\mu)^3-\EE W^3\right|&\le200\delta.
 \end{aligned}
\end{equation}
For the quadratic term, \(v_0=\EE W^2-\mu^2\) gives
\[
 -3hv_0=-3h\EE W^2+3h\mu^2,\qquad
 0\le3h\mu^2\le60\delta^2\le\delta.
\]
Here we used \(h<5\), \(|\mu|\le2\delta\), and
\(60\delta\le60\cdot10^{-6}<1\).
Combining this identity with
\eqref{eq:effective-lattice-conditional-deficit}
and \eqref{eq:effective-lattice-centering-cubics}, we obtain
\begin{equation}\label{eq:effective-lattice-centered-deficit}
 \begin{aligned}
 C_W&:=c\EE|W-\mu|^3-\EE(W-\mu)^3-3hv_0\\
 &=\EE F_h(W)
      +c\{\EE|W-\mu|^3-\EE|W|^3\}\\
 &\quad-\{\EE(W-\mu)^3-\EE W^3\}+3h\mu^2\\
 &\le2\delta+200c\delta+200\delta+60\delta^2\\
 &\le\{3+200(c+1)\}\delta\le2000\delta.
 \end{aligned}
\end{equation}
For the last inequality, \(c<7\) implies
\(3+200(c+1)<1603<2000\).

\textbf{Identify the mass of the upper principal atom.}
Both principal atoms have positive mass, so
\[
 p:=\PP(W=b)\in(0,1),\qquad q:=1-p\in(0,1).
\]
The identity \(h=a+b\) gives
\[
 \mu=-aq+bp=-a+hp.
\]
Thus \(W-\mu\) takes the values \(-hp\) and \(hq\)
with probabilities \(q\) and \(p\), respectively.
Computing its moments directly gives
\[
 \begin{aligned}
 v_0&=q h^2p^2+p h^2q^2=h^2pq,\\
 \EE|W-\mu|^3
 &=q h^3p^3+p h^3q^3=h^3pq(p^2+q^2),\\
 \EE(W-\mu)^3
 &=-q h^3p^3+p h^3q^3=h^3pq(q^2-p^2)
 =h^3pq(1-2p).
 \end{aligned}
\]
Since \(v_0^{3/2}=h^3(pq)^{3/2}>0\),
division of \eqref{eq:effective-lattice-centered-deficit}
by \(v_0^{3/2}\) yields
\[
 \frac{C_W}{v_0^{3/2}}
 =\frac{c(p^2+q^2)-(1-2p)-3}{\sqrt{pq}}
 \le\frac{2000\delta}{v_0^{3/2}}.
\]
Moreover, \eqref{eq:effective-lattice-conditional-variance}
gives
\[
 v_0^{3/2}\ge(1-5\delta)^{3/2}
 \ge(9/10)^{3/2}>2/3.
\]
The last inequality follows, for example, from
\((9/10)^{3/2}>(9/10)^2=81/100>2/3\).
Thus \(2000\delta/v_0^{3/2}\le3000\delta\).

By \eqref{eq:esseen-parameters},
\(\pe=(4-\sqrt{10})/2\).
Together with \(c=3+\sqrt{10}\), this gives
\[
 2c\pe=c-1,\qquad 2c\pe^2=c-4.
\]
Since \(p^2+q^2=2p^2-2p+1\), the numerator above
therefore has the exact representation
\[
 \begin{aligned}
 c(p^2+q^2)-(1-2p)-3
 &=2cp^2-2(c-1)p+c-4\\
 &=2c(p-\pe)^2.
 \end{aligned}
\]
Consequently,
\begin{equation}\label{eq:effective-lattice-bernoulli-deficit}
 \frac{2c(p-\pe)^2}{\sqrt{pq}}\le3000\delta.
\end{equation}
Because \(pq\le1/4\), multiplication by \(\sqrt{pq}\)
in \eqref{eq:effective-lattice-bernoulli-deficit} gives
\[
 2c(p-\pe)^2\le3000\delta\sqrt{pq}\le1500\delta.
\]
Taking square roots yields
\begin{equation}\label{eq:effective-lattice-probability-stability}
 |p-\pe|\le\sqrt{\frac{750}{c}}\sqrt\delta
 \le12\sqrt\delta.
\end{equation}
For the final constant, \(c>6\) implies
\(750/c<125<144\).
To locate \(p\), the inequalities
\[
 3.162^2=9.998244<10<10.010896=3.164^2
\]
give \(3.162<\sqrt{10}<3.164\), and hence
\(0.418<\pe<0.419\).
Also, \(\delta\le10^{-6}\) implies
\(12\sqrt\delta\le0.012\).
It follows from
\eqref{eq:effective-lattice-probability-stability} that
\[
 0.406<p<0.431,\qquad p\in I=[2/5,9/20].
\]

\textbf{Construct explicit couplings to bound \(W_1(Q,\Pe)\).}
Let \(W'\), independent of \(X\), have the same law as \(W\),
and define
\[
 \widetilde W=
 \begin{cases}
 X,&X\in\{-a,b\},\\
 W',&X\in A_{\mathrm{out}}.
 \end{cases}
\]
The support of \(Q\) is contained in
\(\{-a,b\}\cup A_{\mathrm{out}}\), so this specifies
\(\widetilde W\) almost surely.
For every Borel set \(E\), independence and the definition
of the conditional law of \(W\) give
\[
 \begin{aligned}
 \PP(\widetilde W\in E)
 &=\PP(X\in E,\ X\in\{-a,b\})
       +\PP(X\in A_{\mathrm{out}})\PP(W'\in E)\\
 &=(1-r)\PP(W\in E)+r\PP(W\in E)
 =\PP(W\in E).
 \end{aligned}
\]
Thus \((X,\widetilde W)\) is a coupling of \(Q\) and
\(\mathcal L(W)\).
It moves \(X\) only on \(A_{\mathrm{out}}\).
Since \(|X|\le15\) and \(|W'|<h<5\) almost surely,
each such move has distance less than \(20\).
Using \(r\le\delta\) from
\eqref{eq:effective-lattice-outside-moments}, we obtain
\begin{equation}\label{eq:effective-lattice-conditional-coupling}
 W_1(Q,\mathcal L(W))
 \le\EE|X-\widetilde W|
 \le20r\le20\delta.
\end{equation}

Now put
\[
 Y=\frac{W-\mu}{\sqrt{v_0}}.
\]
The lower bound \(v_0\ge1-5\delta>0\) in
\eqref{eq:effective-lattice-conditional-variance}
ensures that this variable is well defined.
It satisfies \(\EE Y=0\) and \(\EE Y^2=1\), so the
Cauchy--Schwarz inequality gives \(\EE|Y|\le1\).
The bound \(|v_0-1|\le5\delta\) from
\eqref{eq:effective-lattice-conditional-variance} also gives
\begin{equation}\label{eq:effective-lattice-scale-error}
 |\sqrt{v_0}-1|
 =\frac{|v_0-1|}{1+\sqrt{v_0}}
 \le\frac{5\delta}{1+\sqrt{1-5\delta}}
 \le5\delta.
\end{equation}
Use \(W=\mu+\sqrt{v_0}Y\) to couple \(W\) and \(Y\).
The triangle inequality and \(|\mu|\le2\delta\) from
\eqref{eq:effective-lattice-conditional-variance} give
\begin{equation}\label{eq:effective-lattice-standardization-coupling}
 \begin{aligned}
 W_1(\mathcal L(W),\mathcal L(Y))
 &\le\EE|W-Y|\\
 &\le|\mu|+|\sqrt{v_0}-1|\,\EE|Y|\\
 &\le2\delta+5\delta\le10\delta.
 \end{aligned}
\end{equation}

To compare \(Y\) with \(\Pe\), first identify its atoms.
Since \(\mu=pb-qa\), \(a+b=h\), and \(v_0=h^2pq\),
we have
\[
 \frac{-a-\mu}{\sqrt{v_0}}
 =-\frac{ph}{h\sqrt{pq}}=-\sqrt{\frac pq},\qquad
 \frac{b-\mu}{\sqrt{v_0}}
 =\frac{qh}{h\sqrt{pq}}=\sqrt{\frac qp}.
\]
For \(\theta\in I\), define
\[
 l(\theta)=-\sqrt{\frac{\theta}{1-\theta}},\qquad
 u(\theta)=\sqrt{\frac{1-\theta}{\theta}}.
\]
Then \(Y\) has law
\((1-p)\delta_{l(p)}+p\delta_{u(p)}\).
The definitions \eqref{eq:esseen-parameters} and
\eqref{eq:esseen-law} similarly give
\[
 \Pe=(1-\pe)\delta_{l(\pe)}+\pe\delta_{u(\pe)}.
\]
Differentiation yields
\[
 l'(\theta)
 =-\frac1{2\sqrt\theta(1-\theta)^{3/2}},\qquad
 u'(\theta)
 =-\frac1{2\theta^{3/2}\sqrt{1-\theta}}.
\]
Both \(\theta\) and \(1-\theta\) belong to
\([2/5,3/5]\) when \(\theta\in I\).
Therefore
\[
 |l'(\theta)|,\ |u'(\theta)|
 \le\frac1{2(2/5)^2}=\frac{25}{8}<5,\qquad
 |l(\theta)|,\ |u(\theta)|
 \le\sqrt{\frac32}<\frac54.
\]
In particular, the distance between any two atoms with
parameters in \(I\) is less than \(5/2\).

Let \(V\) be uniform on \((0,1)\), and for
\(\theta\in\{p,\pe\}\) define
\[
 B_\theta=\ind_{\{V\le\theta\}},\qquad
 Z_\theta=l(\theta)(1-B_\theta)+u(\theta)B_\theta.
\]
The laws of \(Z_p\) and \(Z_{\pe}\) are
\(\mathcal L(Y)\) and \(\Pe\), respectively.
Their indicators disagree with probability
\(|p-\pe|\).
On the event where the indicators agree, the mean value
theorem and the derivative bounds give
\(|Z_p-Z_{\pe}|\le5|p-\pe|\), because the interval
between \(p\) and \(\pe\) is contained in \(I\).
On the event where they disagree, the atom bounds give
\(|Z_p-Z_{\pe}|<5/2\).
Hence
\begin{equation}\label{eq:effective-lattice-bernoulli-coupling}
 \begin{aligned}
 W_1(\mathcal L(Y),\Pe)
 &\le\EE|Z_p-Z_{\pe}|\\
 &\le5|p-\pe|\PP(B_p=B_{\pe})
       +\frac52\PP(B_p\ne B_{\pe})\\
 &\le\frac{15}{2}|p-\pe|\le8|p-\pe|.
 \end{aligned}
\end{equation}
Combining
\eqref{eq:effective-lattice-conditional-coupling},
\eqref{eq:effective-lattice-standardization-coupling},
and \eqref{eq:effective-lattice-bernoulli-coupling}
by the triangle inequality, and using
\eqref{eq:effective-lattice-probability-stability}, gives
\[
 W_1(Q,\Pe)
 \le20\delta+10\delta+8(12\sqrt\delta)
 =30\delta+96\sqrt\delta
 \le1000\sqrt\delta.
\]
The last inequality uses \(0\le\delta\le1\), which implies
\(\delta\le\sqrt\delta\), and \(30+96<1000\).

\textbf{Control the lattice span and undo the initial reflection.}
Define
\[
 s(\theta)=[\theta(1-\theta)]^{-1/2},
 \qquad \theta\in I.
\]
Since \(\theta\in[2/5,9/20]\),
\[
 \theta(1-\theta)
 =\frac14-\left(\theta-\frac12\right)^2
 \ge\frac14-\frac1{100}=\frac6{25}.
\]
Consequently,
\[
 s(\theta)\le(6/25)^{-1/2}<3,\qquad
 s'(\theta)
 =-\frac{1-2\theta}{2[\theta(1-\theta)]^{3/2}}.
\]
Since \(0<1-2\theta\le1/5\) and
\([\theta(1-\theta)]^{3/2}\ge(6/25)^{3/2}>1/10\),
we have \(|s'(\theta)|<1\) throughout \(I\).
The mean value theorem therefore gives
\(|s(p)-s(\pe)|\le|p-\pe|\).
The identities \(v_0=h^2p(1-p)\) and
\(h_{\mathrm E}=1/\se=s(\pe)\) imply
\[
 h=\sqrt{v_0}s(p),\qquad h_{\mathrm E}=s(\pe).
\]
Using \eqref{eq:effective-lattice-scale-error} and
\eqref{eq:effective-lattice-probability-stability}, we obtain
\begin{equation}\label{eq:effective-lattice-span-comparison}
 \begin{aligned}
 |h-h_{\mathrm E}|
 &\le s(p)|\sqrt{v_0}-1|+|s(p)-s(\pe)|\\
 &\le3(5\delta)+|p-\pe|\\
 &\le15\delta+12\sqrt\delta
 \le1000\sqrt\delta.
 \end{aligned}
\end{equation}
Again, the final inequality follows from
\(\delta\le\sqrt\delta\) and \(15+12<1000\).

Finally, let \(Q_{\mathrm{orig}}\) denote the law before
the initial reflection.
Choose \(\epsilon=1\) if no reflection was made and
\(\epsilon=-1\) otherwise.
The law denoted by \(Q\) above is the law of
\(\epsilon X_{\mathrm{orig}}\) for
\(X_{\mathrm{orig}}\sim Q_{\mathrm{orig}}\).
Multiplication of both coordinates of a coupling by
\(\epsilon\) preserves the distance between them and
is its own inverse.
It therefore gives
\[
 W_1(Q_{\mathrm{orig}},\epsilon\Pe)=W_1(Q,\Pe).
\]
Reflection also preserves the lattice spacing \(h\).
Thus the two bounds just proved give
\eqref{eq:effective-lattice-stability} for the original
law, with this choice of \(\epsilon\).
\end{proof}

\subsection{Finite lattice approximation and uniform resonance control}

We now use Lemma~\ref{lem:effective-near-lattice} to
construct a finite lattice law \(Q\) close to a violating
law \(P\).
The maximal span \(h\) of \(Q\) will determine the width
of the uniform variable added to sums with the original
law \(P\).
We obtain an explicit smoothing expansion, uniformly
over the integers \(\ell\in[n/2,n]\).
Lemma~\ref{lem:effective-identification} will combine
this expansion and the moment bounds for \(Q\) with
Lemma~\ref{lem:effective-lattice-stability} to locate
the maximizing law.

For this and the next subsection, fix
\begin{equation}\label{eq:effective-global-parameters}
 T=e^{20\mathcal A},\qquad \tau=e^{-100\mathcal A},\qquad
 g=10^{-20}\tau^2,\qquad \xi=10^{10}\tau T^2.
\end{equation}
Here \(\mathcal A\) and \(N_{\mathrm{conf}}\) are defined
in \eqref{eq:effective-threshold-parameters}.

\begin{lemma}\label{lem:effective-global-jitter}
Suppose that \(P\in\Pthree\) is supported in \([-6,6]\),
\(n\ge N_{\mathrm{conf}}\), and \(R_n(P)>\ce\).
There is a finite lattice law \(Q\), with mean
\(\mu_Q\), variance \(v_Q>0\), central third absolute
moment \(\rho_Q\), central third moment \(\kappa_Q\),
and maximal span \(h\), such that
\begin{gather}
 \pi/500\le h\le5,\qquad W_1(P,Q)\le10^5\tau,\notag\\
 |\mu_Q|\le10^5\tau,\qquad
 |v_Q-1|\le10^9\tau,\label{eq:effective-rounded-moments}\\
 |\rho_Q-\beta(P)|+|\kappa_Q-\kappa(P)|
 \le2\cdot10^9\tau.\notag
\end{gather}
The standardization of \(Q\) is supported in \([-15,15]\)
and has third absolute moment less than two.
Let \(X_1,X_2,\ldots\) be independent with law \(P\),
and let \(H_h\), independent of these variables, be
uniform on \([-h/2,h/2]\).
Then, for every integer \(\ell\in[n/2,n]\),
\begin{equation}\label{eq:effective-global-jitter}
\begin{split}
 \sqrt\ell\sup_{x\in\R}\Bigg|
 &\PP\left\{\sum_{i=1}^{\ell}X_i+H_h\le x\right\}
 -\Phi(x/\sqrt\ell)\\
 &-\frac{\kappa(P)}{6\sqrt\ell}
       (1-x^2/\ell)\phi(x/\sqrt\ell)
 \Bigg|\le e^{-19\mathcal A}.
\end{split}
\end{equation}
\end{lemma}

\begin{proof}
\textbf{Round the summand to a nearby lattice.}
Let \(X\sim P\).
The violation \(R_n(P)>\ce\) and
\eqref{eq:moment-cutoff} give
\(\beta(P)<B_*<1.84\).
Also, \(n\ge N_{\mathrm{conf}}\ge e^{1000\mathcal A}>1024\),
so Lemma~\ref{lem:effective-near-lattice} applies.
It supplies a lattice
\(\Lambda_0=a_0+h_0\mathbb Z\), where
\[
 \frac{\pi}{500}\le h_0\le4\pi,\qquad
 \EE\dist(X,\Lambda_0)^2\le2\pi^2\frac{\log n}{n}.
\]
Let \(Z\) be a nearest point of \(\Lambda_0\) to \(X\),
choosing the smaller point in case of a tie.
This defines a measurable function of \(X\).
Since successive lattice points are \(h_0\) apart,
\[
 |X-Z|=\dist(X,\Lambda_0)\le h_0/2,\qquad
 |Z|\le6+h_0/2\le6+2\pi<13.
\]
The Cauchy--Schwarz inequality therefore gives
\[
 e_n:=\EE|X-Z|
 \le\left(\EE|X-Z|^2\right)^{1/2}
 \le\sqrt2\,\pi\sqrt{\frac{\log n}{n}}
 \le5\sqrt{\frac{\log n}{n}}.
\]
The function \(x\mapsto(\log x)/x\) is decreasing for
\(x\ge e\), because its derivative is
\((1-\log x)/x^2\).
Using \(n\ge e^{1000\mathcal A}\), we obtain
\begin{equation}\label{eq:effective-rounding-error}
 e_n\le5\sqrt{1000\mathcal A}\,e^{-500\mathcal A}
 \le e^{-100\mathcal A}=\tau.
\end{equation}
For the last comparison,
\(\mathcal A=10^{14}\) gives
\[
 5\sqrt{1000\mathcal A}=5\sqrt{10^{17}}
 <10^{10}<1+400\mathcal A\le e^{400\mathcal A}.
\]
The value of \(\tau\) is the one fixed in
\eqref{eq:effective-global-parameters}.

All values of \(Z\) lie in
\([-6-h_0/2,6+h_0/2]\), an interval of length \(12+h_0\).
If its occupied points have integer indices \(k\) in
\(\Lambda_0=a_0+h_0\mathbb Z\), the difference between
the largest and smallest such indices is at most
\[
 \frac{12+h_0}{h_0}
 \le\frac{6000}{\pi}+1.
\]
Their number is at most one more than this difference.
Thus both their number and their index range are less
than \(2000\), since \(\pi>3.1\) gives
\(6000/\pi+2<2000\).

\textbf{Remove atoms with small probability and construct a coupling.}
Retain precisely those atoms \(z\) of \(Z\) for which
\(\PP(Z=z)\ge\tau\).
Let \(E\) be the event that \(Z\) is retained, and put
\[
 r=\PP(E^c)\le2000\tau.
\]
The inequality follows because fewer than \(2000\)
atoms are discarded and each has probability less
than \(\tau\).
We will use the following numerical bounds:
\begin{equation}\label{eq:effective-retention-parameters}
 r\le2000\tau<\frac1{100},\qquad
 w:=10^5\tau\le10^{-6}.
\end{equation}
Indeed, \(100\mathcal A>11\log10\) implies
\(\tau=e^{-100\mathcal A}<10^{-11}\).
It follows in particular that \(r<1\), so the
conditional law
\[
 Q=\mathcal L(Z\mid E)
\]
is well defined and has finite support in \((-13,13)\).

On an extension of the probability space, let \(Z'\)
be independent of \((X,Z)\) with law \(Q\), and define
\[
 Y=
 \begin{cases}
 Z,&E\text{ occurs},\\
 Z',&E^c\text{ occurs}.
 \end{cases}
\]
For any Borel set \(A\), independence gives
\[
 \begin{aligned}
 \PP(Y\in A)
 &=\PP(Z\in A,E)+\PP(E^c)\PP(Z'\in A)\\
 &=(1-r)Q(A)+rQ(A)=Q(A).
 \end{aligned}
\]
Thus \((X,Y)\) is a coupling of \(P\) and \(Q\).
On \(E\), \(Y=Z\), whereas on \(E^c\) both \(Z\)
and \(Y=Z'\) have absolute value less than \(13\).
Consequently,
\begin{equation}\label{eq:effective-retained-coupling}
 \begin{aligned}
 W_1(P,Q)
 &\le\EE|X-Y|
 \le\EE|X-Z|+\EE|Z-Y|\\
 &\le e_n+26r
 \le\tau+26(2000\tau)
 =52001\tau\le w.
 \end{aligned}
\end{equation}
This explicit coupling will also be used to compare
characteristic functions and their derivatives below.

\textbf{Control the moments and the standardization of the retained law.}
Write
\[
 \mu_Q=\EE Y,\qquad v_Q=\Var Y,\qquad
 \rho_Q=\EE|Y-\mu_Q|^3,\qquad
 \kappa_Q=\EE(Y-\mu_Q)^3.
\]
Since \(\EE X=0\), \(\EE X^2=1\), and
\(|X|,|Y|\le13\), the coupling gives
\[
 \begin{aligned}
 |\mu_Q|&=|\EE(Y-X)|\le w,\\
 |\EE Y^2-1|
 &\le\EE\{|Y-X|\,|Y+X|\}\le26w.
 \end{aligned}
\]
Therefore
\begin{equation}\label{eq:effective-retained-moment-bounds}
 |\mu_Q|\le w,\qquad
 |v_Q-1|
 \le26w+w^2\le27w,
 \qquad 0.99<v_Q<1.01.
\end{equation}
Here \(w\le10^{-6}<1\) gives \(w^2\le w\), and
\(27w\le27\cdot10^{-6}<0.01\).
In particular, \(v_Q>0\), so \(Q\) has at least two
distinct atoms.

Centering \(Y\) increases the mean coupling distance
by at most \(|\mu_Q|\):
\[
 \EE|(Y-\mu_Q)-X|
 \le\EE|Y-X|+|\mu_Q|\le2w.
\]
Also,
\[
 |Y-\mu_Q|<13+w<14,\qquad |X|\le6<14.
\]
On \([-14,14]\), the derivatives of \(x\mapsto|x|^3\)
and \(x\mapsto x^3\) have absolute value at most
\(3(14)^2=588\).
The mean value theorem, followed by the coupling bound,
therefore yields
\begin{equation}\label{eq:effective-retained-cubic-errors}
 |\rho_Q-\beta(P)|\le1176w,\qquad
 |\kappa_Q-\kappa(P)|\le1176w.
\end{equation}
Since \(27w=2.7\cdot10^6\tau<10^9\tau\) and
\(1176w=1.176\cdot10^8\tau<10^9\tau\),
\eqref{eq:effective-retained-coupling},
\eqref{eq:effective-retained-moment-bounds}, and
\eqref{eq:effective-retained-cubic-errors} prove the
moment and Wasserstein bounds in
\eqref{eq:effective-rounded-moments}.

For later use, \(\beta(P)<1.84\) and \(w\le10^{-6}\)
also give
\[
 \rho_Q\le\beta(P)+1176w
 <1.84+0.001176<1.85.
\]
Let \(\widetilde Q\) be the law of
\((Y-\mu_Q)/\sqrt{v_Q}\).
Its support and third absolute moment satisfy
\begin{equation}\label{eq:effective-retained-standardization}
 \supp\widetilde Q\subset[-15,15],\qquad
 \beta(\widetilde Q)=\frac{\rho_Q}{v_Q^{3/2}}<2.
\end{equation}
Indeed, \(v_Q>0.99\) and \(|Y-\mu_Q|<14\) imply
\[
 \frac{|Y-\mu_Q|}{\sqrt{v_Q}}
 <\frac{14}{\sqrt{0.99}}<15.
\]
For the moment bound,
\[
 \frac{\rho_Q}{v_Q^{3/2}}
 <\frac{1.85}{(0.99)^{3/2}}
 <\frac{1.85}{(0.99)^2}<2,
\]
where \((0.99)^2=0.9801>1.85/2\).
The law \(\widetilde Q\) has mean zero and variance one,
so it belongs to \(\Pthree\).

\textbf{Identify the maximal lattice span of the retained support.}
List the retained atoms in increasing order as
\(z_0<z_1<\cdots<z_m\), where \(m\ge1\).
Since they belong to \(\Lambda_0\), there are integers
\[
 b_i=\frac{z_i-z_0}{h_0},\qquad
 b_0=0,\qquad 1\le b_i<2000\quad(1\le i\le m).
\]
Let
\[
 G=\left\{\sum_{i=1}^m a_i b_i:a_i\in\mathbb Z\right\}
 \subset\mathbb Z.
\]
The set \(G\) contains a positive integer, so it has a
least positive element \(g_0\).
Division of any \(b_i\) by \(g_0\) leaves a remainder
in \(G\cap[0,g_0)\), which must be zero by minimality.
Thus \(g_0\) divides every \(b_i\).
Conversely, any common divisor of the \(b_i\) divides
every element of \(G\), and in particular divides \(g_0\).
Hence \(g_0=\gcd(b_1,\ldots,b_m)\).
Set
\[
 h=g_0h_0,\qquad d_i=b_i/g_0,\qquad
 z_i=z_0+hd_i.
\]
Then
\[
 d_0=0,\qquad 0\le d_i\le2000,\qquad
 \gcd(d_1,\ldots,d_m)=1.
\]
The retained support is contained in \(z_0+h\mathbb Z\).
Furthermore, \(g_0\in G\) means that \(h\) is an integer
linear combination of the differences \(z_i-z_0\).
If the same support were contained in a lattice with
spacing \(H>0\), every such difference would belong to
\(H\mathbb Z\); their integer combination \(h\) would
then be a positive integer multiple of \(H\).
Thus \(H\le h\), proving that \(h\) is the maximal
lattice span of \(Q\).
In particular, \(h\ge h_0\ge\pi/500\).

The centering and scaling map divides lattice spacings
by \(\sqrt{v_Q}\), and its inverse takes any containing
lattice for \(\widetilde Q\) to one for \(Q\).
Thus the maximal span of \(\widetilde Q\) is
\(h/\sqrt{v_Q}\).
Apply Esseen's moment inequality \eqref{eq:esseen-moment}
to \(\widetilde Q\), whose third moments are
\(\rho_Q/v_Q^{3/2}\) and \(\kappa_Q/v_Q^{3/2}\).
Multiplication by \(v_Q^{3/2}\) gives
\[
 |\kappa_Q|+3hv_Q\le c_*\rho_Q.
\]
Using \(c_*<7\), \(\rho_Q<1.85\), and \(v_Q>0.99\),
we conclude that
\begin{equation}\label{eq:effective-retained-span-bounds}
 \frac{\pi}{500}\le h
 \le\frac{c_*\rho_Q-|\kappa_Q|}{3v_Q}
 <\frac{7(1.85)}{3(0.99)}<5.
\end{equation}
The last numerical comparison follows from
\(7(1.85)=12.95<14.85=5\cdot3(0.99)\).

\textbf{Find a bounded integer representation of the lattice step.}
We next prove that the integers \(d_i\) admit coefficients
\(c_i\in\mathbb Z\) satisfying
\begin{equation}\label{eq:effective-bezout}
 \sum_{i=0}^m c_i d_i=1,\qquad
 \sum_{i=0}^m|c_i|\le4000.
\end{equation}
Choose any positive \(d_j=d_*\).
If \(d_*=1\), take \(c_j=1\) and all other coefficients
equal to zero.
Suppose instead that \(d_*\ge2\).
Consider the graph on the \(d_*\) residue classes
\(\mathbb Z/d_*\mathbb Z\), joining each class to the
classes obtained by adding or subtracting one of the
\(d_i\).
The condition \(\gcd(d_1,\ldots,d_m)=1\) implies that
these residues generate the entire group.
More explicitly, \(g_0\in G\), divided by \(g_0\),
expresses \(1\) as an integer linear combination of
the \(d_i\); following the corresponding signed steps
reaches the residue \(1\), and repetition reaches
every residue.
Thus the graph is connected.

Choose a shortest path from \(0\) to \(1\).
Such a path has no repeated vertices, so it uses at
most \(d_*-1\) edges.
Let \(S\) be the signed integer sum of its steps.
Since every \(d_i\le2000\), we have
\[
 S\equiv1\pmod{d_*},\qquad
 |S|\le2000(d_*-1).
\]
Write \(S-1=M d_*\) for an integer \(M\).
Then
\[
 |M|
 \le\frac{|S|+1}{d_*}
 \le\frac{2000(d_*-1)+1}{d_*}
 \le2000.
\]
Subtracting \(M d_j\) from the signed path sum gives
an integer combination of the \(d_i\) equal to \(1\).
After combining coefficients of repeated indices, the
sum of their absolute values cannot increase, and is
therefore at most
\[
 (d_*-1)+|M|\le1999+2000=3999<4000.
\]
This proves \eqref{eq:effective-bezout}.

\textbf{Deduce a spectral gap away from the lattice frequencies.}
For \(u\in\R\), write
\[
 f_P(u)=\EE e^{iuX},\qquad f_Q(u)=\EE e^{iuY},
 \qquad \omega=\frac{2\pi}{h},\qquad q_i=Q(\{z_i\}).
\]
Every retained unconditional mass is at least \(\tau\),
and conditioning divides it by \(1-r\le1\). Hence
\[
 q_i=\frac{\PP(Z=z_i)}{1-r}\ge\tau,
 \qquad 0\le i\le m.
\]
Put
\[
 \delta_i(u)=\dist(uhd_i,2\pi\mathbb Z).
\]
For each \(i\), choose an integer \(k_i\) attaining
this distance.
Equation \eqref{eq:effective-bezout} gives
\[
 \begin{aligned}
 \dist(uh,2\pi\mathbb Z)
 &\le\left|uh-2\pi\sum_{i=0}^m c_i k_i\right|\\
 &=\left|\sum_{i=0}^m c_i(uhd_i-2\pi k_i)\right|\\
 &\le\sum_{i=0}^m|c_i|\delta_i(u)
 \le4000\max_{1\le i\le m}\delta_i(u).
 \end{aligned}
\]
The maximum may be taken over \(i\ge1\) because
\(d_0=0\) and \(\delta_0(u)=0\).

Since \(\sum_iq_i=1\), direct multiplication of the
finite characteristic-function sum gives
\[
 1-|f_Q(u)|^2
 =\sum_{i,j=0}^m q_iq_j
       [1-\cos(u(z_i-z_j))].
\]
Every term is nonnegative.
Choose \(i\in\{1,\ldots,m\}\) attaining the largest
\(\delta_i(u)\), and retain the two distinct terms
\((i,0)\) and \((0,i)\).
For any real \(t\), with
\(\delta_t=\dist(t,2\pi\mathbb Z)\in[0,\pi]\),
concavity of sine on \([0,\pi/2]\) gives
\[
 1-\cos t=2\sin^2(\delta_t/2)
 \ge\frac{2\delta_t^2}{\pi^2}.
\]
Consequently,
\[
 \begin{aligned}
 1-|f_Q(u)|^2
 &\ge2q_iq_0[1-\cos(uhd_i)]\\
 &\ge\frac{4\tau^2}{\pi^2}\delta_i(u)^2\\
 &\ge\frac{4\tau^2}{\pi^2\,4000^2}
          \dist(uh,2\pi\mathbb Z)^2.
 \end{aligned}
\]
The argument is also valid if all circular distances
are zero, because both the resulting lower bound and
the selected two terms then vanish.
Since \(h>0\),
\(\dist(uh,2\pi\mathbb Z)=h\,\dist(u,\omega\mathbb Z)\).
Using \(h\ge\pi/500\) and
\[
 \frac{4(\pi/500)^2}{\pi^2\,4000^2}=10^{-12},
\]
we have proved
\begin{equation}\label{eq:effective-spectral-gap}
 1-|f_Q(u)|^2
 \ge10^{-12}\tau^2\dist(u,\omega\mathbb Z)^2,
 \qquad u\in\R.
\end{equation}

\textbf{Transfer the spectral gap through conditioning and rounding.}
Put \(\rho=10^{-3}\).
If \(\dist(u,\omega\mathbb Z)\ge\rho\), then
\eqref{eq:effective-spectral-gap} and the definition of \(g\)
give
\[
 1-|f_Q(u)|^2
 \ge10^{-12}\tau^2\rho^2
 =10^{-18}\tau^2=100g.
\]
Since \(|f_Q(u)|\le1\), it follows that
\[
 1-|f_Q(u)|
 =\frac{1-|f_Q(u)|^2}{1+|f_Q(u)|}
 \ge50g.
\]
Thus \(|f_Q(u)|\le1-50g\) at these frequencies.

The retention event \(E\) has probability \(1-r\), and
the conditional law of \(Z\) given \(E\) is \(Q\).
Consequently,
\[
 \begin{aligned}
 f_P(u)
 &=(1-r)f_Q(u)
   +\EE[(e^{iuX}-e^{iuZ})\ind_E]
   +\EE[e^{iuX}\ind_{E^c}].
 \end{aligned}
\]
The bound \(|e^{it}-e^{is}|\le|t-s|\) gives
\begin{equation}\label{eq:effective-mixture-gap}
 \begin{aligned}
 |f_P(u)|
 &\le(1-r)|f_Q(u)|
      +|u|\EE[|X-Z|\ind_E]+r\\
 &\le(1-r)|f_Q(u)|+r+|u|e_n.
 \end{aligned}
\end{equation}
We next check that the rounding error is small enough
to preserve this spectral gap throughout the cutoff.
The function \(x\mapsto(\log x)/x\) is decreasing for
\(x\ge e\).
Since \(n\ge N_{\mathrm{conf}}\ge e^{1000\mathcal A}\),
the bound \eqref{eq:effective-rounding-error} on \(e_n\) implies
\[
 Te_n
 \le5\sqrt{1000\mathcal A}\,e^{-480\mathcal A}<g.
\]
For the last inequality, division by
\(g=10^{-20}e^{-200\mathcal A}\) leaves
\[
 5\cdot10^{20}\sqrt{1000\mathcal A}\,e^{-280\mathcal A}<1.
\]
Indeed, \(\mathcal A=10^{14}\) gives
\(5\cdot10^{20}\sqrt{1000\mathcal A}<10^{30}\), whereas
\[
 e^{280\mathcal A}
 \ge\frac{(280\mathcal A)^2}{2}>10^{30}.
\]
Also, \(r\le2000\tau<1/100\), because
\(e^{100\mathcal A}\ge1+100\mathcal A>200000\).
Therefore \eqref{eq:effective-mixture-gap} yields
\[
 |f_P(u)|
 \le(1-r)(1-50g)+r+g
 \le1-\frac{97}{2}g\le1-g
\]
whenever \(|u|\le T\) and
\(\dist(u,\omega\mathbb Z)\ge\rho\).
We record this as
\begin{equation}\label{eq:effective-global-nonresonant-gap}
 |f_P(u)|\le1-g,\qquad
 |u|\le T,\quad \dist(u,\omega\mathbb Z)\ge\rho.
\end{equation}
Keeping the retention factor \(1-r\) in
\eqref{eq:effective-mixture-gap} preserves the gap after
both conditioning and rounding.

\textbf{Compare the derivatives near the lattice frequencies.}
Write
\[
 q_P(u)=|f_P(u)|^2,\qquad q_Q(u)=|f_Q(u)|^2.
\]
Use the coupling \((X,Y)\) of \(P\) and \(Q\) from
\eqref{eq:effective-retained-coupling}, satisfying
\(\EE|X-Y|\le10^5\tau\), and take an
independent copy \((X',Y')\).
Then
\[
 D_P=X-X',\qquad D_Q=Y-Y'
\]
are both supported in \([-26,26]\), and
\[
 \EE|D_P-D_Q|
 \le\EE|X-Y|+\EE|X'-Y'|
 \le2\cdot10^5\tau.
\]
Independence gives
\[
 q_P(u)=\EE\cos(uD_P),\qquad
 q_Q(u)=\EE\cos(uD_Q).
\]
All the variables are bounded, so differentiation under
the expectations is justified.
The first and second derivatives are expectations of
\(-x\sin(ux)\) and \(-x^2\cos(ux)\), respectively.
As functions of \(x\in[-26,26]\), their derivatives
have absolute values at most
\[
 1+26|u|,\qquad 52+676|u|.
\]
For \(|u|\le T+1\), the coupling therefore gives
\[
 \begin{aligned}
 |q_P'(u)-q_Q'(u)|
 &\le2\cdot10^5\tau\{1+26(T+1)\},\\
 |q_P''(u)-q_Q''(u)|
 &\le2\cdot10^5\tau\{52+676(T+1)\}.
 \end{aligned}
\]
Since \(T\ge1\), both right-hand sides are at most
\(3\cdot10^8\tau T\), and hence at most
\(\xi=10^{10}\tau T^2\).
We have proved
\begin{equation}\label{eq:effective-global-derivatives}
 |q_P'(u)-q_Q'(u)|\le\xi,\qquad
 |q_P''(u)-q_Q''(u)|\le\xi,\qquad |u|\le T+1.
\end{equation}

The variance estimate in
\eqref{eq:effective-retained-moment-bounds}
and the parameter choices give
\begin{equation}\label{eq:effective-global-small-parameters}
 0.99\le v_Q\le1.01,\qquad \xi<10^{-4}.
\end{equation}
For the bound on \(\xi\), use \(e^x\ge1+x\) to obtain
\[
 \xi=10^{10}e^{-60\mathcal A}
 \le\frac{10^{10}}{1+60\mathcal A}<10^{-4},
\]
where \(\mathcal A=10^{14}\).

For every integer \(j\), put
\[
 u_j^0=j\omega,\qquad
 I_j=[u_j^0-\rho,u_j^0+\rho].
\]
We consider only intervals that intersect \([-T,T]\).
For such an interval, \(|u_j^0|\le T+\rho\), so every
point of the full interval satisfies
\[
 |u|\le T+2\rho<T+1.
\]
Thus \eqref{eq:effective-global-derivatives} applies on
the entire interval, including any part beyond the cutoff.
Because \(D_Q\) is a multiple of \(h\) and
\(\omega=2\pi/h\), we have
\[
 \sin(u_j^0D_Q)=0,\qquad \cos(u_j^0D_Q)=1
 \quad\hbox{almost surely}.
\]
It follows that
\[
 q_Q'(u_j^0)=0,\qquad
 q_Q''(u_j^0)=-\EE D_Q^2=-2v_Q.
\]
The third derivative satisfies
\[
 \|q_Q'''\|_\infty
 \le\EE|D_Q|^3
 \le26\,\EE D_Q^2=52v_Q<60.
\]
For \(u\in I_j\), the mean value theorem and
\eqref{eq:effective-global-small-parameters} consequently give
\[
 q_Q''(u)\le-2v_Q+60\rho
 \le-1.98+0.06=-1.92\le-\frac32.
\]
Combining this with \eqref{eq:effective-global-derivatives},
we obtain
\begin{equation}\label{eq:effective-global-peak-concavity}
 q_P''(u)\le-\frac32+\xi\le-1,\qquad u\in I_j.
\end{equation}
Integration of \(q_Q''\le-3/2\) from the center to the
two endpoints gives
\[
 q_Q'(u_j^0-\rho)\ge\frac32\rho,\qquad
 q_Q'(u_j^0+\rho)\le-\frac32\rho.
\]
Equation \eqref{eq:effective-global-derivatives} then implies
\[
 q_P'(u_j^0-\rho)\ge\frac32\rho-\xi>0,\qquad
 q_P'(u_j^0+\rho)\le-\frac32\rho+\xi<0,
\]
since \((3/2)\rho=0.0015>\xi\).
By continuity, \(q_P'\) has a zero \(u_j\) in the
interior of \(I_j\).
Equation \eqref{eq:effective-global-peak-concavity} makes
\(q_P'\) strictly decreasing, so this zero is unique
and is the unique maximizer of \(q_P\) on \(I_j\).

At the center, \(q_Q'(u_j^0)=0\), so
\(|q_P'(u_j^0)|\le\xi\).
Applying the mean value theorem between \(u_j^0\) and
\(u_j\), and using \(-q_P''\ge1\), yields
\begin{equation}\label{eq:effective-global-peak-location}
 |u_j-u_j^0|
 \le|q_P'(u_j)-q_P'(u_j^0)|
 =|q_P'(u_j^0)|\le\xi.
\end{equation}
Taylor's theorem at \(u_j\) gives
\[
 q_P(u)
 \le q_P(u_j)-\frac12(u-u_j)^2
 \le1-\frac12(u-u_j)^2,\qquad u\in I_j.
\]
Here \(|u-u_j|\le2\rho\), so the last upper bound
is positive.
For every integer \(\ell\ge1\), using
\(1-t\le e^{-t}\) and
\(|f_P(u)|^\ell=q_P(u)^{\ell/2}\), we conclude that
\begin{equation}\label{eq:effective-global-peak-bound}
 |f_P(u)|^\ell
 \le e^{-\ell(u-u_j)^2/4},\qquad
 u\in I_j.
\end{equation}

\textbf{Use the jitter multiplier to integrate the nonzero peaks.}
For the independent uniform variable \(H_h\), define
\[
 H_0(u)=\EE e^{iuH_h}=\sinc(hu/2).
\]
For every nonzero integer \(j\), \(H_0(u_j^0)=0\).
Differentiation under the expectation gives
\[
 |H_0'(u)|
 =|\EE(iH_h e^{iuH_h})|
 \le\EE|H_h|=\frac h4\le\frac54<\frac32.
\]
The mean value bound for this complex-valued function
follows by integrating its derivative.
Together with \eqref{eq:effective-global-peak-location},
it gives
\[
 |H_0(u)|
 \le\frac32|u-u_j^0|
 \le\frac32\{|u-u_j|+\xi\},
 \qquad u\in I_j,\quad j\ne0.
\]
Since \(h\le5\), we have \(\omega=2\pi/h\ge2\pi/5>6/5\).
For \(j\ge1\) and \(u\in I_j\),
\[
 |u|\ge j\omega-\rho
 \ge j(\omega-\rho)>j.
\]
Combining these bounds with
\eqref{eq:effective-global-peak-bound}, and extending the
integral to the whole real line, yields
\[
 \begin{aligned}
 &\int_{I_j\cap[-T,T]}
       \frac{|f_P(u)|^\ell|H_0(u)|}{|u|}\,du\\
 &\quad\le\frac3{2j}
       \int_{\R}(|w|+\xi)e^{-\ell w^2/4}\,dw\\
 &\quad=\frac3{2j}
       \left(\frac4\ell+\frac{2\sqrt\pi\,\xi}{\sqrt\ell}\right)\\
 &\quad\le\frac6j
       \left(\frac1\ell+\frac{\xi}{\sqrt\ell}\right).
 \end{aligned}
\]
Here \(w=u-u_j\), and the last step uses \(\sqrt\pi<2\).
The same bound holds for \(I_{-j}\), since
\(|f_P(-u)|=|f_P(u)|\), \(|H_0(-u)|=|H_0(u)|\), and
the cutoff is symmetric.

Define
\[
 J=\left\lfloor\frac{T+\rho}{\omega}\right\rfloor,\qquad
 \mathcal R_T=\bigcup_{1\le|j|\le J}(I_j\cap[-T,T]).
\]
These are exactly the nonzero peak intervals that meet
the cutoff, including intervals only partially contained
in \([-T,T]\).
The bounds \(h\ge\pi/500\) and \(h\le5\) give
\(\omega\le1000\) and \(\omega>6/5\).
Since \(T=e^{20\mathcal A}\ge1+20\mathcal A>1000\),
they imply \(1\le J<T\).
Indeed,
\[
 J\le\frac{T+\rho}{\omega}
 <\frac56(T+\rho)<T,
\]
where the last inequality uses \(T>5\rho\).
Therefore
\[
 \sum_{j=1}^J\frac1j
 \le1+\int_1^J\frac{dx}{x}
 =1+\log J\le1+\log T.
\]
Summing over both signs of \(j\), we obtain
\begin{equation}\label{eq:effective-global-resonance-integral}
 \int_{\mathcal R_T}
       \frac{|f_P(u)|^\ell|H_0(u)|}{|u|}\,du
 \le12(1+\log T)
       \left(\frac1\ell+\frac{\xi}{\sqrt\ell}\right).
\end{equation}

For the remaining Fourier integral, put
\[
 \mathcal N_T=
 \{u:1/2\le|u|\le T\}\setminus\mathcal R_T.
\]
Every nonzero \(I_j\) lies outside \((-1,1)\), so
\(\mathcal R_T\subset\{u:1/2\le|u|\le T\}\).
Moreover, \(I_0=[-\rho,\rho]\) lies within
\((-1/2,1/2)\).
Thus each \(u\in\mathcal N_T\) satisfies
\(\dist(u,\omega\mathbb Z)\ge\rho\), and
\eqref{eq:effective-global-nonresonant-gap} applies there.
The sets \(\mathcal R_T\) and \(\mathcal N_T\) partition
the frequencies \(1/2\le|u|\le T\).

\textbf{Define the standardized Fourier comparison and bound its low frequencies.}
Fix an integer \(\ell\in[n/2,n]\), and put
\(\kappa=\kappa(P)\).
Since \(R_n(P)>\ce\), equation \eqref{eq:moment-cutoff}
gives
\[
 |\kappa|\le\beta(P)<B_*<1.84.
\]
Also \(\ell\ge n/2\ge2\).
The assumptions of Lemma~\ref{lem:effective-low-frequency}
therefore hold for \(P\) and sample size \(\ell\).
Define
\[
 \mathcal J_\ell(z)
 =\PP\left\{\frac{\sum_{i=1}^{\ell}X_i+H_h}{\sqrt\ell}
                  \le z\right\},\qquad
 \mathcal G_{\ell,\kappa}(z)
 =\Phi(z)+\frac{\kappa}{6\sqrt\ell}(1-z^2)\phi(z).
\]
The characteristic function of the law with CDF
\(\mathcal J_\ell\) and the Fourier transform of the
signed comparison measure from
\eqref{eq:effective-edgeworth-comparison} are, respectively,
\begin{equation}\label{eq:effective-global-fourier-comparison}
 \begin{aligned}
 a_\ell(t)&=f_P(t/\sqrt\ell)^\ell H_0(t/\sqrt\ell),\\
 b_\ell(t)&=e^{-t^2/2}
             \left(1+\frac{\kappa(it)^3}{6\sqrt\ell}\right).
 \end{aligned}
\end{equation}
The first identity uses independence of the summands
and of \(H_h\), together with
\(H_0(u)=\EE e^{iuH_h}\).
The original frequency \(u\) and standardized frequency
\(t\) are related by
\begin{equation}\label{eq:effective-global-frequency-change}
 t=\sqrt\ell\,u,\qquad
 \frac{dt}{|t|}=\frac{du}{|u|}.
\end{equation}

Since \(H_h\) is centered and
\(\EE H_h^2=h^2/12\), the Taylor remainder bound
\eqref{eq:jitter-proof-taylor-remainder} with \(m=2\)
gives
\[
 |H_0(u)-1|
 =|\EE(e^{iuH_h}-1-iuH_h)|
 \le\frac{u^2}{2}\EE H_h^2=\frac{h^2u^2}{24}.
\]
For \(|t|\le\sqrt\ell/2\), the triangle inequality yields
\[
 \begin{aligned}
 |a_\ell(t)-b_\ell(t)|
 &\le|f_P(t/\sqrt\ell)^\ell-b_\ell(t)|\\
 &\quad+|f_P(t/\sqrt\ell)|^\ell
              |H_0(t/\sqrt\ell)-1|.
 \end{aligned}
\]
The first term contributes at most \(6/\ell\) after
division by \(|t|\) and integration, by
\eqref{eq:effective-fourier-integral}.
Equation \eqref{eq:effective-characteristic-decay}
also gives
\[
 |f_P(t/\sqrt\ell)|^\ell\le e^{-0.35t^2},
 \qquad |t|\le\sqrt\ell/2.
\]
Thus the contribution from the second term is bounded by
\[
 \begin{aligned}
 \frac{h^2}{24\ell}
       \int_{\R}|t|e^{-0.35t^2}\,dt
 &=\frac{h^2}{24(0.35)\ell}\\
 &\le\frac{25}{24(0.35)\ell}<\frac3\ell,
 \end{aligned}
\]
where we used \(h\le5\) and
\(\int_{\R}|t|e^{-\alpha t^2}\,dt=1/\alpha\)
for \(\alpha>0\).
Combining the two contributions and enlarging the bound
gives
\begin{equation}\label{eq:effective-global-low-integral}
 \int_{|t|\le\sqrt\ell/2}
       \frac{|a_\ell(t)-b_\ell(t)|}{|t|}\,dt
 \le\frac9\ell\le\frac{11}{\ell}.
\end{equation}
The integrand at zero is defined by continuity; the
Taylor bounds above and
\eqref{eq:effective-fourier-error} show that its limit
there is zero.

\textbf{Combine the remaining frequencies and apply signed smoothing.}
Choose the standardized cutoff
\[
 L=T\sqrt\ell.
\]
Under \eqref{eq:effective-global-frequency-change},
the remaining interval
\(\sqrt\ell/2\le|t|\le L\) becomes
\(1/2\le|u|\le T\).
The nonzero resonance intervals comprising
\(\mathcal R_T\) lie in this range: the preceding bound
\(|u|\ge|j|\) on \(I_j\), \(j\ne0\), gives
\(|u|\ge1\) there.
Consequently, the remaining frequencies are the union
of \(\mathcal R_T\) and \(\mathcal N_T\), and no part
of the interval around zero is left out.
On \(\mathcal N_T\),
\eqref{eq:effective-global-nonresonant-gap} gives
\(|f_P(u)|\le1-g\).
Since \(0<g<1\), \(1-g\le e^{-g}\), and
\(|H_0(u)|\le1\), we obtain
\begin{equation}\label{eq:effective-global-nonresonance-integral}
 \begin{aligned}
 \int_{\mathcal N_T}
        \frac{|f_P(u)|^\ell|H_0(u)|}{|u|}\,du
 &\le e^{-\ell g}
          \int_{1/2\le|u|\le T}\frac{du}{|u|}\\
 &=2\log(2T)e^{-\ell g}.
 \end{aligned}
\end{equation}
For the comparison transform, the condition
\(|\kappa|<1.84\) allows us to use
\eqref{eq:effective-gaussian-tail}:
\[
 \int_{|t|\ge\sqrt\ell/2}
       \frac{|b_\ell(t)|}{|t|}\,dt
 \le5e^{-\ell/8}.
\]
The triangle inequality, the frequency change, and
\eqref{eq:effective-global-resonance-integral}
therefore give
\[
 \begin{aligned}
 &\int_{\sqrt\ell/2\le|t|\le L}
       \frac{|a_\ell(t)-b_\ell(t)|}{|t|}\,dt\\
 &\quad\le
 12(1+\log T)\left(\frac1\ell+\frac{\xi}{\sqrt\ell}\right)
 +2\log(2T)e^{-\ell g}+5e^{-\ell/8}.
 \end{aligned}
\]
Adding \eqref{eq:effective-global-low-integral}, we conclude
that
\begin{equation}\label{eq:effective-global-total-integral}
 \begin{aligned}
 \int_{-L}^{L}\frac{|a_\ell(t)-b_\ell(t)|}{|t|}\,dt
 &\le\frac{11}{\ell}
 +12(1+\log T)\left(\frac1\ell+\frac{\xi}{\sqrt\ell}\right)\\
 &\quad+2\log(2T)e^{-\ell g}+5e^{-\ell/8}.
 \end{aligned}
\end{equation}

The law with CDF \(\mathcal J_\ell\) has finite first
absolute moment because \(P\) and \(H_h\) have bounded
support.
As verified after \eqref{eq:effective-edgeworth-comparison},
the signed measure with primitive \(\mathcal G_{\ell,\kappa}\)
has total mass one and finite first absolute moment.
Equation \eqref{eq:effective-edgeworth-density-bound}
gives a Lipschitz constant at most one for this primitive.
We may therefore apply Lemma~\ref{lem:signed-smoothing}
with the constants in
\eqref{eq:effective-smoothing-constants}.
Using \eqref{eq:effective-global-total-integral} and
putting \(H_T=1+\log T\), we obtain
\begin{equation}\label{eq:effective-global-smoothing-budget}
 \begin{aligned}
 \sqrt\ell\,\|\mathcal J_\ell-\mathcal G_{\ell,\kappa}\|_\infty
 &\le\frac{11}{4\sqrt\ell}
       +\frac{3H_T}{\sqrt\ell}+3\xi H_T\\
 &\quad+\frac12\sqrt\ell\log(2T)e^{-\ell g}
       +\frac54\sqrt\ell e^{-\ell/8}+\frac{24}{T}.
 \end{aligned}
\end{equation}
The last term is the scaled smoothing remainder,
since \(24\sqrt\ell/L=24/T\).

\textbf{Check the numerical bounds uniformly over \(\ell\in[n/2,n]\).}
By \eqref{eq:effective-threshold-parameters} and
\eqref{eq:effective-global-parameters},
\[
 n\ge N_{\mathrm{conf}}\ge e^{1000\mathcal A},
 \qquad
 g=10^{-20}e^{-200\mathcal A}.
\]
Hence
\[
 \frac{\ell g}{\sqrt n}
 \ge\frac{\sqrt n\,g}{2}
 \ge\frac12\,10^{-20}e^{300\mathcal A}>1.
\]
For the final inequality, the elementary bound
\(e^x\ge x^2/2\), \(x\ge0\), gives
\[
 \frac12\,10^{-20}e^{300\mathcal A}
 \ge22500\cdot10^{-20}\mathcal A^2
 =2.25\cdot10^{12}>1,
\]
because \(\mathcal A=10^{14}\).
Thus \(\ell g\ge\sqrt n\).
Also \(\sqrt n\ge16\), so
\(\ell/8\ge n/16\ge\sqrt n\).
Since \(T\ge1\), we have \(H_T\ge1\) and
\(\log(2T)\le H_T\).
It follows that the two exponential terms in
\eqref{eq:effective-global-smoothing-budget} satisfy
\[
 \begin{aligned}
 &\frac12\sqrt\ell\log(2T)e^{-\ell g}
       +\frac54\sqrt\ell e^{-\ell/8}\\
 &\quad\le\left(\frac12H_T+\frac54\right)
                    \sqrt\ell e^{-\sqrt n}\\
 &\quad\le2H_T\sqrt\ell e^{-\sqrt n}
 \le\frac{2H_T}{\sqrt\ell}.
 \end{aligned}
\]
To justify the last step, set \(s=\sqrt n\ge16\).
The exponential series gives
\(e^s\ge s^3/6\ge s^2=n\), and therefore
\(\ell e^{-\sqrt n}\le n e^{-\sqrt n}\le1\).
Since \(H_T\ge1\), the remaining terms proportional
to \(\ell^{-1/2}\) are bounded by
\[
 \frac{11}{4\sqrt\ell}
 +\frac{3H_T}{\sqrt\ell}+\frac{2H_T}{\sqrt\ell}
 \le\frac{31H_T}{4\sqrt\ell}
 \le\frac{10H_T}{\sqrt\ell}.
\]
Consequently,
\begin{equation}\label{eq:effective-global-reduced-budget}
 \sqrt\ell\,\|\mathcal J_\ell-\mathcal G_{\ell,\kappa}\|_\infty
 \le\frac{10H_T}{\sqrt\ell}+3\xi H_T+\frac{24}{T}.
\end{equation}

Finally,
\[
 H_T=1+20\mathcal A,\qquad
 \xi=10^{10}e^{-60\mathcal A},\qquad
 \ell^{-1/2}\le\sqrt2\,e^{-500\mathcal A}.
\]
Thus the three terms on the right of
\eqref{eq:effective-global-reduced-budget} are at most
\[
 10\sqrt2 H_T e^{-500\mathcal A},\qquad
 3\cdot10^{10}H_T e^{-60\mathcal A},\qquad
 24e^{-20\mathcal A},
\]
respectively.
Since \(\mathcal A=10^{14}\), we have
\(H_T\le21\mathcal A\) and
\[
 30\sqrt2 H_T
 <9\cdot10^{10}H_T
 \le1.89\cdot10^{12}\mathcal A
 <2\cdot10^{12}\mathcal A
 <\frac{\mathcal A^2}{2}\le e^{\mathcal A}.
\]
Also \(72<\mathcal A^2/2\le e^{\mathcal A}\).
The three terms are therefore bounded, respectively, by
\[
 \frac13e^{-499\mathcal A},\qquad
 \frac13e^{-59\mathcal A},\qquad
 \frac13e^{-19\mathcal A}.
\]
Their sum is less than \(e^{-19\mathcal A}\).
Combining this with
\eqref{eq:effective-global-reduced-budget} proves the
required bound for every allowed \(\ell\).
Finally, \(x=\sqrt\ell\,z\) is a bijection of \(\R\),
so the supremum in
\eqref{eq:effective-global-smoothing-budget} is exactly
the supremum in \eqref{eq:effective-global-jitter}.
This proves the lemma.
\end{proof}

\subsection{Effective identification of a maximizing law}

We now combine the smoothing expansion from
Lemma~\ref{lem:effective-global-jitter} with the stability
estimate in Lemma~\ref{lem:effective-lattice-stability}.
At a positive maximizing discrepancy, the expansion forces
the auxiliary lattice law to have a small deficit in
Esseen's moment inequality.
The stability estimate then places the maximizing law
near \(\Pe\) or its reflection.
We use the positive sign of the discrepancy to exclude
the reflection, and then show that the standardized
maximizing threshold is close to zero.
These quantitative bounds will be used in
Proposition~\ref{prop:effective-confinement} to control the
entire support.

For a law \(Q\) with finite second moment, write
\[
 M(Q)=\int x|x|\,Q(dx).
\]

\begin{lemma}\label{lem:effective-identification}
Let \(P\in\Pthree\) and \(t\in\R\) attain the positive
maximum \(R=C_n>\ce\), where
\[
 n\ge N_{\mathrm{conf}},\qquad \supp P\subset[-6,6].
\]
Put \(z=t/\sqrt n\) and \(\delta_0=e^{-6\mathcal A}\).
Then the maximal span \(h\) of the auxiliary law in
Lemma~\ref{lem:effective-global-jitter} satisfies
\begin{equation}\label{eq:effective-identification}
 W_1(P,\Pe)\le e^{-7\mathcal A},\qquad
 |h-h_{\mathrm E}|\le e^{-7\mathcal A},
\end{equation}
and
\begin{equation}\label{eq:effective-moment-identification}
 \begin{split}
 |\beta(P)-\be|&\le\delta_0,\qquad
 |\kappa(P)-\ke|\le\delta_0,\\
 |M(P)-M(\Pe)|&\le\delta_0,\qquad z^2\le e^{-5\mathcal A}.
 \end{split}
\end{equation}
\end{lemma}

\begin{proof}
\textbf{Use the positive maximizing discrepancy to bound the moment deficit.}
Write \(\beta=\beta(P)\) and \(\kappa=\kappa(P)\).
The maximizing assumption and \eqref{eq:cn} give
\[
 R=\frac{\sqrt n}{\beta}\{F_{n,P}(t)-\Phi(z)\},
 \qquad R\le R_n(P)\le C_n=R.
\]
Thus \(R_n(P)=R>\ce\), and
\eqref{eq:moment-cutoff} gives
\[
 0<\beta<B_*<1.84,\qquad |\kappa|\le\beta<1.84.
\]
Let \(Q,h\) be supplied by
Lemma~\ref{lem:effective-global-jitter}.
Its assumptions hold because
\(n\ge N_{\mathrm{conf}}\), \(\supp P\subset[-6,6]\),
and \(R_n(P)>\ce\).
In particular, \(0<h\le5\).
For independent \(X_1,\ldots,X_n\) with law \(P\), put
\[
 S_n=\sum_{i=1}^nX_i,\qquad
 J_n(x)=\PP(S_n+H_h\le x),
\]
where \(H_h\) is independent and uniform on \([-h/2,h/2]\).
Because \(H_h\le h/2\), the event \(S_n\le t\) is contained
in \(S_n+H_h\le t+h/2\).
Thus
\[
 F_{n,P}(t)\le J_n(t+h/2).
\]
The positive maximizing identity and \(z=t/\sqrt n\)
give
\begin{equation}\label{eq:effective-identification-positive-maximum}
 \ce\beta<R\beta
 =\sqrt n\{F_{n,P}(t)-\Phi(z)\}
 \le\sqrt n\{J_n(t+h/2)-\Phi(z)\}.
\end{equation}

Put \(\psi(x)=(1-x^2)\phi(x)\).
Taking \(\ell=n\) in \eqref{eq:effective-global-jitter}
gives
\[
 J_n(x)
 =\Phi(x/\sqrt n)
     +\frac{\kappa}{6\sqrt n}\psi(x/\sqrt n)+r_n(x),
 \qquad
 \sqrt n\,\|r_n\|_\infty\le e^{-19\mathcal A}.
\]
Set \(a=h/(2\sqrt n)\), so that the standardized argument
at \(t+h/2\) is \(z+a\).
Taylor's theorem and the mean value theorem give
\[
 \begin{aligned}
 |\Phi(z+a)-\Phi(z)-a\phi(z)|
 &\le\frac12\|\phi'\|_\infty a^2,\\
 |\psi(z+a)-\psi(z)|
 &\le\|\psi'\|_\infty a.
 \end{aligned}
\]
The bounds proved in Lemma~\ref{lem:effective-selection}
give \(\|\phi'\|_\infty<1/4\) and
\(\|\psi'\|_\infty=\|\phi'''\|_\infty<5/4\).
Using \(h\le5\) and \(|\kappa|<1.84\), the total Taylor
error after multiplication by \(\sqrt n\) is at most
\[
 \begin{aligned}
 \frac1{\sqrt n}
 \left\{\frac{\|\phi'\|_\infty h^2}{8}
          +\frac{|\kappa|\|\psi'\|_\infty h}{12}\right\}
 &\le\frac1{\sqrt n}
       \left(\frac{25}{32}+\frac{23}{24}\right)\\
 &=\frac{167}{96\sqrt n}<\frac{10}{\sqrt n}.
 \end{aligned}
\]
Consequently,
\begin{equation}\label{eq:effective-identification-shifted-jitter}
 \begin{aligned}
 \sqrt n\{J_n(t+h/2)-\Phi(z)\}
 &=\frac h2\phi(z)+\frac{\kappa}{6}\psi(z)+\eta_n,\\
 |\eta_n|&\le e^{-19\mathcal A}+\frac{10}{\sqrt n}.
 \end{aligned}
\end{equation}
Since \(n\ge N_{\mathrm{conf}}\ge e^{1000\mathcal A}\),
\[
 e^{-19\mathcal A}+\frac{10}{\sqrt n}
 \le e^{-19\mathcal A}+10e^{-500\mathcal A}
 \le2e^{-19\mathcal A}\le e^{-18\mathcal A}.
\]
Here \(10e^{-500\mathcal A}\le e^{-19\mathcal A}\)
uses \(e^{481\mathcal A}\ge10\), and the final inequality
uses \(e^{\mathcal A}\ge2\).
Combining \eqref{eq:effective-identification-positive-maximum}
and \eqref{eq:effective-identification-shifted-jitter}
yields
\begin{equation}\label{eq:effective-maximizer-envelope}
 \ce\beta<
 \phi(z)\left\{\frac h2+\frac{\kappa}{6}(1-z^2)\right\}
 +e^{-18\mathcal A}.
\end{equation}

The Gaussian bounds
\(\|\phi\|_\infty=\|\psi\|_\infty=\phi_0\) from
\eqref{eq:jitter-envelope-gaussian-suprema} imply
\[
 \ce\beta
 <\phi_0\left(\frac h2+\frac{|\kappa|}{6}\right)
       +e^{-18\mathcal A}.
\]
Using \(\ce=c_*\phi_0/6\) from \eqref{eq:ce} and
multiplying by \(6/\phi_0\), we obtain
\begin{equation}\label{eq:effective-identification-summand-deficit}
 c_*\beta-|\kappa|-3h
 \le16e^{-18\mathcal A}.
\end{equation}
Indeed, \(\phi_0>0.39\) gives
\(6/\phi_0<6/0.39<16\).

\textbf{Transfer the deficit to the standardized lattice approximation.}
Write the mean, variance, and central third moments of
\(Q\) as \(\mu_Q,v_Q,\rho_Q,\kappa_Q\), as in
Lemma~\ref{lem:effective-global-jitter}.
Let
\[
 \sigma_Q=\sqrt{v_Q},\qquad
 \widetilde Q=\mathcal L((Y-\mu_Q)/\sigma_Q),\quad Y\sim Q.
\]
By \eqref{eq:effective-retained-moment-bounds}
and \eqref{eq:effective-retained-standardization},
\[
 0.99<v_Q<1.01,\qquad
 \widetilde Q\in\Pthree,\qquad
 \supp\widetilde Q\subset[-15,15],\qquad
 \beta(\widetilde Q)=\frac{\rho_Q}{v_Q^{3/2}}<2.
\]
The maximal lattice span of \(\widetilde Q\) is
\(h/\sigma_Q\), and its third moment is
\(\kappa(\widetilde Q)=\kappa_Q/v_Q^{3/2}\).
Thus Esseen's moment inequality \eqref{eq:esseen-moment}
for \(\widetilde Q\), multiplied by \(v_Q^{3/2}\), gives
\[
 D_Q:=c_*\rho_Q-|\kappa_Q|-3hv_Q\ge0.
\]
On the other hand,
\[
 \begin{aligned}
 D_Q-(c_*\beta-|\kappa|-3h)
 &=c_*(\rho_Q-\beta)
      -(|\kappa_Q|-|\kappa|)-3h(v_Q-1)\\
 &\le c_*|\rho_Q-\beta|
        +|\kappa_Q-\kappa|+3h|v_Q-1|.
 \end{aligned}
\]
Since \(c_*<7\), \(h\le5\), and
\eqref{eq:effective-rounded-moments} gives
\[
 |\rho_Q-\beta|+|\kappa_Q-\kappa|\le2\cdot10^9\tau,
 \qquad |v_Q-1|\le10^9\tau,
\]
the preceding difference is bounded above by
\[
 7(2\cdot10^9\tau)+15(10^9\tau)
 =29\cdot10^9\tau\le10^{12}\tau.
\]
Together with
\eqref{eq:effective-identification-summand-deficit},
this proves
\begin{equation}\label{eq:effective-identification-lattice-deficit}
 0\le D_Q\le16e^{-18\mathcal A}+10^{12}\tau.
\end{equation}

The deficit for the standardized law is exactly
\[
 \begin{aligned}
 \widetilde D
 &:=c_*\beta(\widetilde Q)
      -|\kappa(\widetilde Q)|-\frac{3h}{\sigma_Q}\\
 &=\frac{c_*\rho_Q-|\kappa_Q|-3hv_Q}{v_Q^{3/2}}
 =\frac{D_Q}{v_Q^{3/2}}.
 \end{aligned}
\]
Because \(v_Q>0.99\),
\(v_Q^{3/2}>(0.99)^{3/2}>(0.99)^2>1/2\).
Using \(\tau=e^{-100\mathcal A}\) and
\eqref{eq:effective-identification-lattice-deficit},
we therefore obtain
\[
 0\le\widetilde D
 \le32e^{-18\mathcal A}
      +2\cdot10^{12}e^{-100\mathcal A}
 \le e^{-16\mathcal A}.
\]
For the last inequality, division by \(e^{-16\mathcal A}\)
gives
\[
 32e^{-2\mathcal A}
       +2\cdot10^{12}e^{-84\mathcal A}\le1.
\]
Each term is at most \(1/2\):
\(e^{2\mathcal A}\ge1+2\mathcal A>64\), and
\(e^{84\mathcal A}\ge1+84\mathcal A>4\cdot10^{12}\).

Set \(\delta=e^{-16\mathcal A}\).
Then \(0<\delta\le10^{-6}\), since
\(16\mathcal A>6\log10\).
The law \(\widetilde Q\) has mean zero, variance one,
support in \([-15,15]\), and third absolute moment
less than two, as verified above.
Its positive lattice span \(h/\sigma_Q\) has deficit
\(\widetilde D\in[0,\delta]\).
All the hypotheses of
Lemma~\ref{lem:effective-lattice-stability} therefore hold.
That lemma gives a sign \(\epsilon\in\{-1,1\}\) such that
\begin{equation}\label{eq:effective-identification-standardized-stability}
 W_1(\widetilde Q,\epsilon\Pe)\le1000e^{-8\mathcal A},
 \qquad
 \left|\frac h{\sigma_Q}-h_{\mathrm E}\right|
 \le1000e^{-8\mathcal A}.
\end{equation}

\textbf{Transfer the stability bounds back to \(P\) and to the original span.}
Let \(Y\sim Q\), and put
\[
 \widetilde Y=\frac{Y-\mu_Q}{\sigma_Q},
 \qquad \mathcal L(\widetilde Y)=\widetilde Q.
\]
Then \(\EE\widetilde Y=0\) and
\(\EE\widetilde Y^2=1\).
The Cauchy--Schwarz inequality therefore gives
\(\EE|\widetilde Y|\le1\).
Using the moment bounds \eqref{eq:effective-rounded-moments},
we have
\[
 |\sigma_Q-1|
 =\frac{|v_Q-1|}{1+\sigma_Q}\le10^9\tau.
\]
The deterministic coupling
\(Y=\mu_Q+\sigma_Q\widetilde Y\) now yields
\begin{equation}\label{eq:effective-identification-standardization-cost}
 \begin{aligned}
 W_1(Q,\widetilde Q)
 &\le\EE|Y-\widetilde Y|\\
 &\le|\mu_Q|+|\sigma_Q-1|\,\EE|\widetilde Y|\\
 &\le(10^5+10^9)\tau.
 \end{aligned}
\end{equation}
Combining \eqref{eq:effective-rounded-moments},
\eqref{eq:effective-identification-standardization-cost},
and \eqref{eq:effective-identification-standardized-stability}
by the triangle inequality gives
\begin{equation}\label{eq:effective-identification-law-transfer}
 \begin{aligned}
 W_1(P,\epsilon\Pe)
 &\le W_1(P,Q)+W_1(Q,\widetilde Q)
                  +W_1(\widetilde Q,\epsilon\Pe)\\
 &\le1000e^{-8\mathcal A}+(2\cdot10^5+10^9)\tau\\
 &\le1000e^{-8\mathcal A}+2\cdot10^9\tau.
 \end{aligned}
\end{equation}

For the lattice spacing, the bound from
Lemma~\ref{lem:effective-lattice-stability} concerns
\(h/\sigma_Q\), the spacing of the standardized law.
We now convert it to a bound on \(h\).
The variance estimate
\eqref{eq:effective-retained-moment-bounds} gives
\(\sigma_Q^2=v_Q<1.01\), so \(\sigma_Q<2\).
Since \(3.1^2<10<3.2^2\),
\eqref{eq:esseen-parameters} gives
\(2/5<\pe<9/20\). Consequently,
\[
 \se^2=\frac14-\left(\pe-\frac12\right)^2
 \ge\frac14-\frac1{100}=\frac6{25},\qquad
 h_{\mathrm E}=\frac1{\se}\le\sqrt{\frac{25}{6}}<3.
\]
It follows from
\eqref{eq:effective-identification-standardized-stability}
that
\begin{equation}\label{eq:effective-identification-span-transfer}
 \begin{aligned}
 |h-h_{\mathrm E}|
 &=\left|\sigma_Q
       \left(\frac h{\sigma_Q}-h_{\mathrm E}\right)
        +(\sigma_Q-1)h_{\mathrm E}\right|\\
 &\le\sigma_Q\left|\frac h{\sigma_Q}-h_{\mathrm E}\right|
       +h_{\mathrm E}|\sigma_Q-1|\\
 &\le2000e^{-8\mathcal A}+3\cdot10^9\tau.
 \end{aligned}
\end{equation}
Together,
\eqref{eq:effective-identification-law-transfer} and
\eqref{eq:effective-identification-span-transfer} give
\begin{equation}\label{eq:effective-unoriented-identification}
 \max\{W_1(P,\epsilon\Pe),\,|h-h_{\mathrm E}|\}
 \le4000e^{-8\mathcal A}+10^{10}\tau
 <e^{-7\mathcal A}.
\end{equation}
To verify the last comparison, use
\(\tau=e^{-100\mathcal A}\) and divide by
\(e^{-7\mathcal A}\).
The resulting sum satisfies
\[
 4000e^{-\mathcal A}+10^{10}e^{-93\mathcal A}
 <\frac12+\frac12=1,
\]
because
\(e^{\mathcal A}\ge1+\mathcal A>8000\) and
\(e^{93\mathcal A}\ge1+93\mathcal A>2\cdot10^{10}\).
The sign \(\epsilon\) in
\eqref{eq:effective-unoriented-identification}
has not yet been determined.

\textbf{Convert the Wasserstein bound into bounds on the three moments.}
Both \(P\) and \(\epsilon\Pe\) are supported in
\([-6,6]\).
For \(P\), this is an assumption of the lemma.
For \(\epsilon\Pe\), the two atom magnitudes are
\(\aE,\bE>0\), and \(\aE+\bE=h_{\mathrm E}<3\), by
\eqref{eq:esseen-law} and the preceding span bound.
Reflection does not change their magnitudes.

On \([-6,6]\), the functions
\[
 g_1(x)=|x|^3,\qquad g_2(x)=x^3,\qquad g_3(x)=x|x|
\]
are differentiable and satisfy
\[
 |g_1'(x)|=|g_2'(x)|=3x^2\le108,\qquad
 |g_3'(x)|=2|x|\le12\le108.
\]
Thus each is Lipschitz with constant at most \(108\)
on this interval.
For any coupling \((V,Z)\) with laws \(P\) and
\(\epsilon\Pe\), respectively, we consequently have
\[
 |\EE g_i(V)-\EE g_i(Z)|
 \le\EE|g_i(V)-g_i(Z)|
 \le108\,\EE|V-Z|,\qquad i=1,2,3.
\]
The left-hand side depends only on the marginal laws.
Taking the infimum over all such couplings and using
\eqref{eq:effective-unoriented-identification} gives
\[
 |\EE g_i(V)-\EE g_i(Z)|
 \le108\,W_1(P,\epsilon\Pe)
 \le108e^{-7\mathcal A}
 \le e^{-6\mathcal A}=\delta_0.
\]
The final inequality follows from
\(e^{\mathcal A}\ge1+\mathcal A>108\).
The moments of the reflected Esseen law satisfy
\[
 \beta(\epsilon\Pe)=\be,\qquad
 \kappa(\epsilon\Pe)=\epsilon\ke,\qquad
 M(\epsilon\Pe)=\epsilon M(\Pe).
\]
These identities follow because \(g_1\) is even and
\(g_2,g_3\) are odd, together with
\eqref{eq:esseen-identities}.
We therefore obtain
\begin{equation}\label{eq:effective-unoriented-moments}
 |\beta(P)-\be|\le\delta_0,\qquad
 |\kappa(P)-\epsilon\ke|\le\delta_0,\qquad
 |M(P)-\epsilon M(\Pe)|\le\delta_0.
\end{equation}

\textbf{Use the positive maximizing discrepancy to determine the orientation.}
For \(\epsilon\in\{-1,1\}\), define
\[
 E_\epsilon(w)
 =\phi(w)\left\{\frac{h_{\mathrm E}}2
           +\frac{\epsilon\ke}{6}(1-w^2)\right\},
 \qquad w\in\R.
\]
The parameter identities
\eqref{eq:esseen-parameters}--\eqref{eq:esseen-identities}
give
\[
 \begin{aligned}
 (3+\sqrt{10})(\pe^2+\qe^2)
 &=(3+\sqrt{10})(10-3\sqrt{10})\\
 &=\sqrt{10}=3+\qe-\pe.
 \end{aligned}
\]
Dividing by \(\se\) and using the formulas for
\(\be,\ke,h_{\mathrm E}\), we obtain
\((3+\sqrt{10})\be=3h_{\mathrm E}+\ke\).
Together with the definition of \(\ce\) in
\eqref{eq:ce}, this yields
\begin{equation}\label{eq:effective-esseen-envelope-height}
 E_+(0)=\phi_0\left(\frac{h_{\mathrm E}}2+\frac{\ke}{6}\right)
       =\ce\be.
\end{equation}

We next compare the original envelope with \(E_\epsilon\),
where \(\epsilon\) is supplied by
\eqref{eq:effective-unoriented-identification}.
By \eqref{eq:jitter-envelope-gaussian-suprema},
\(\sup_w\phi(w)=\phi_0\) and
\(\sup_w|(1-w^2)\phi(w)|=\phi_0\).
Equations \eqref{eq:effective-unoriented-identification}
and \eqref{eq:effective-unoriented-moments} therefore imply
\[
 \begin{aligned}
 &\left|
 \phi(w)\left\{\frac h2+\frac{\kappa(P)}6(1-w^2)\right\}
                -E_\epsilon(w)\right|\\
 &\qquad\le\frac{\phi_0}{2}|h-h_{\mathrm E}|
           +\frac{\phi_0}{6}|\kappa(P)-\epsilon\ke|\\
 &\qquad\le\frac{\phi_0}{2}e^{-7\mathcal A}
           +\frac{\phi_0}{6}\delta_0,
 \qquad w\in\R.
 \end{aligned}
\]
We must also compare the third absolute moment on the
left side of \eqref{eq:effective-maximizer-envelope}:
\[
 \ce\be\le\ce\beta(P)+\ce|\beta(P)-\be|
          \le\ce\beta(P)+\ce\delta_0.
\]
Combining these inequalities with
\eqref{eq:effective-maximizer-envelope} gives
\[
 \ce\be
 < E_\epsilon(z)
   +\frac{\phi_0}{2}e^{-7\mathcal A}
   +\left(\frac{\phi_0}{6}+\ce\right)\delta_0
   +e^{-18\mathcal A}.
\]
Since \(e^{-7\mathcal A}\le\delta_0=e^{-6\mathcal A}\),
\(\phi_0<2/5\), and \(\ce<1/2\), the coefficient of
\(\delta_0\) in these perturbation terms is at most
\[
 \frac{\phi_0}{2}+\frac{\phi_0}{6}+\ce
 <\frac15+\frac1{15}+\frac12=\frac{23}{30}<2.
\]
Thus
\begin{equation}\label{eq:effective-oriented-envelope-comparison}
 \ce\be<E_\epsilon(z)+2\delta_0+e^{-18\mathcal A}.
\end{equation}
The estimates leading to this comparison are uniform
over the standardized threshold \(z\).

Suppose that \(\epsilon=-1\).
Then
\[
 E_-(w)=\phi(w)
       \left\{\frac{h_{\mathrm E}}2-\frac{\ke}{6}
                          +\frac{\ke w^2}{6}\right\}.
\]
The identities
\[
 \frac{\ke}{h_{\mathrm E}}=\qe-\pe=\sqrt{10}-3\in(0,1)
\]
show that \(h_{\mathrm E}>\ke>0\).
Since \(\phi'(w)=-w\phi(w)\), differentiation gives
\[
 E_-'(w)
 =w\phi(w)\left\{\frac{\ke-h_{\mathrm E}}2
                           -\frac{\ke w^2}{6}\right\}<0,
 \qquad w>0.
\]
The function \(E_-\) is even, so its maximum on \(\R\)
is attained at zero.
By \eqref{eq:effective-esseen-envelope-height},
\[
 \sup_{w\in\R}E_-(w)
 =E_-(0)
 =\ce\be-\frac{\phi_0\ke}{3}.
\]
To bound this deficit numerically, use
\(\se=\sqrt{\pe\qe}\le1/2\) and \(\sqrt{10}>3.15\):
\[
 \ke=\frac{\sqrt{10}-3}{\se}
 \ge2(\sqrt{10}-3)>0.3.
\]
Since \(\phi_0>0.39\), it follows that
\begin{equation}\label{eq:effective-reflected-envelope-gap}
 \frac{\phi_0\ke}{3}>
 \frac{0.39\cdot0.3}{3}=0.039.
\end{equation}
On the other hand, \(\mathcal A=10^{14}\) and
\(e^{6\mathcal A}\ge6\mathcal A\) give
\[
 2\delta_0+e^{-18\mathcal A}
 \le3\delta_0
 \le\frac1{2\mathcal A}<0.039.
\]
Substituting the maximum of \(E_-\) and these two
bounds into
\eqref{eq:effective-oriented-envelope-comparison}
would imply \(\ce\be<\ce\be\), a contradiction.
Therefore \(\epsilon=1\).
Equations \eqref{eq:effective-unoriented-identification}
and \eqref{eq:effective-unoriented-moments} now give
\eqref{eq:effective-identification} and the first three
bounds in \eqref{eq:effective-moment-identification}.

\textbf{Locate the standardized maximizing threshold.}
Because \((P,t)\) attains the positive maximum,
\[
 R=\frac{\sqrt n}{\beta(P)}
       \{F_{n,P}(t)-\Phi(z)\}>\ce>\frac25,
 \qquad z=\frac{t}{\sqrt n}.
\]
Applying \eqref{eq:effective-nonuniform} at this
threshold yields
\[
 R\le\frac{17.36}{1+|z|^3}.
\]
If \(|z|\ge4\), the right side is at most
\(17.36/65<2/5\), contradicting the preceding strict
lower bound. Hence \(|z|<4\).

The positive envelope is
\[
 E_+(w)=\phi(w)
       \left\{\frac{h_{\mathrm E}}2+\frac{\ke}{6}
                          -\frac{\ke w^2}{6}\right\}.
\]
For \(|w|\le4\), put \(u=w^2/2\in[0,8]\).
The inequality \(e^u\ge1+u\) implies
\[
 1-e^{-w^2/2}
 =1-e^{-u}\ge\frac{u}{1+u}\ge\frac{u}{10}
 =\frac{w^2}{20}.
\]
Also \(h_{\mathrm E}=1/\se\ge2\), and \(\ke>0\).
Using \(\phi(w)=\phi_0e^{-w^2/2}\), we obtain
\begin{equation}\label{eq:effective-positive-envelope-gap}
 \begin{aligned}
 E_+(0)-E_+(w)
 &=\phi_0\left[
   \left(\frac{h_{\mathrm E}}2+\frac{\ke}{6}\right)
                         (1-e^{-w^2/2})
     +\frac{\ke w^2}{6}e^{-w^2/2}\right]\\
 &\ge\phi_0\frac{h_{\mathrm E}}2
                         (1-e^{-w^2/2})\\
 &\ge\frac{0.39}{20}w^2
 \ge\frac{w^2}{100},
 \qquad |w|\le4.
 \end{aligned}
\end{equation}
The first lower bound drops nonnegative terms, and
the second uses \(\phi_0>0.39\) and
\(h_{\mathrm E}/2\ge1\).
We may apply \eqref{eq:effective-positive-envelope-gap}
at \(w=z\).
Equations \eqref{eq:effective-esseen-envelope-height}
and \eqref{eq:effective-oriented-envelope-comparison},
now with \(\epsilon=1\), give
\[
 \frac{z^2}{100}
 \le E_+(0)-E_+(z)
 <2\delta_0+e^{-18\mathcal A}
 \le3\delta_0.
\]
It follows that \(z^2\le300\delta_0\).
Finally \(300<1+\mathcal A\le e^{\mathcal A}\), so
\[
 z^2\le300e^{-6\mathcal A}\le e^{-5\mathcal A}.
\]
This proves the remaining assertion in
\eqref{eq:effective-moment-identification}.
\end{proof}

\subsection{An explicit bound on the entire support}

Lemma~\ref{lem:effective-identification} places the maximizing
law close to \(\Pe\) in \(W_1\), but this still permits additional
support points in regions of very small probability.
We now use the exact contact equations to exclude such points.
We first obtain a coarse interval for the support, and then
compare upper and lower bounds for the change in
\(F_{n-1,P}\) between two support contacts.
This confines the entire support near the two Esseen atoms.
After the affine change of coordinates sending these atoms
to \(0\) and \(1\), we also check the cluster probabilities.
The resulting bounds allow us to apply
Proposition~\ref{prop:effective-clusters} and complete the
proof of the explicit threshold.

\begin{proposition}\label{prop:effective-confinement}
Let \(P\in\Pthree\) and \(t\in\R\) attain the positive
maximum \(R=C_n>\ce\), where
\[
 n\ge N_{\mathrm{conf}},\qquad
 \supp P\subset[-6,6],\qquad d_n\le20/\sqrt n.
\]
Let \(\widetilde P\) be the law of \(\pe+\se X\),
where \(X\sim P\), and put \(\zeta=\eta_*/2\).
Then
\begin{equation}\label{eq:effective-full-support}
 \supp\widetilde P\subset[-\zeta,\zeta]\cup[1-\zeta,1+\zeta],
 \qquad
 |\widetilde P([1-\zeta,1+\zeta])-\pe|<\eta_*.
\end{equation}
\end{proposition}

\begin{proof}
Write
\[
 \beta=\beta(P),\qquad \kappa=\kappa(P),\qquad
 M=M(P),\qquad z=\frac{t}{\sqrt n},\qquad
 \delta_0=e^{-6\mathcal A}.
\]
Let \(h\) be the auxiliary lattice span supplied by
Lemma~\ref{lem:effective-global-jitter}.
The hypotheses of Lemma~\ref{lem:effective-identification}
are satisfied, so we may use both
\eqref{eq:effective-identification} and
\eqref{eq:effective-moment-identification} below.
Together with \(n\ge N_{\mathrm{conf}}\) and
\eqref{eq:effective-threshold-parameters}, these bounds give
\[
 n^{-1/2}\le e^{-500\mathcal A},\qquad
 z^2\le e^{-5\mathcal A}.
\]
We first prove the coarser support bound
\(\pe+\se\supp P\subset[-1/2,3/2]\).
The contact argument for this bound requires an explicit
version of the Gaussian cancellation estimate, which we
establish next.

\textbf{Make the Gaussian remainder explicit on the bounded support.}
The function \(H_n\) is defined in \eqref{eq:gaussian-H}.
Since \(n\ge N_{\mathrm{conf}}\ge2\), the bound
\eqref{eq:gaussian-compact-error-bound} in
Lemma~\ref{lem:gaussian-expansion} applies.
It gives, with \(M_3=\|\phi'''\|_\infty<5/4\)
as established in Lemma~\ref{lem:effective-selection},
\begin{equation}\label{eq:effective-gaussian-three-remainders}
 \begin{aligned}
 &\left|H_n(z,y)-\phi''(z)
             \left(\frac y2-\frac{y^3}{6}\right)\right|\\
 &\quad\le\frac{M_3}{\sqrt n}
       \left(\frac{|y|^4}{24}+\frac12+\frac{y^2}{4}\right)\\
 &\quad\le\frac{67.5+0.625+11.25}{\sqrt n}
       =\frac{79.375}{\sqrt n},\qquad |y|\le6.
 \end{aligned}
\end{equation}
Indeed, the three numerical terms are
\((5/4)6^4/24=67.5\), \((5/4)/2=0.625\), and
\((5/4)6^2/4=11.25\).
They bound, respectively, the translation, scale, and
\(g\)-expansion remainders in
\eqref{eq:gaussian-translation-expansion},
\eqref{eq:gaussian-scale-expansion}, and
\eqref{eq:gaussian-g-expansion}.

We next replace \(\phi''(z)\) by \(\phi''(0)\).
Differentiation gives
\[
 \phi^{(4)}(r)=(r^4-6r^2+3)\phi(r),\qquad
 \phi'''(0)=0.
\]
The derivatives of \(ve^{-v/2}\) and \(v^2e^{-v/2}\)
are \(e^{-v/2}(1-v/2)\) and
\(e^{-v/2}(2v-v^2/2)\), respectively.
Thus the functions \(ve^{-v/2}\) and \(v^2e^{-v/2}\)
attain their maxima on \([0,\infty)\) at \(v=2\) and
\(v=4\), with values \(2e^{-1}\) and \(16e^{-2}\).
Using \(v=r^2\), we obtain
\begin{equation}\label{eq:effective-gaussian-fourth-derivative}
 \|\phi^{(4)}\|_\infty
 \le\phi_0(16e^{-2}+12e^{-1}+3)<4.
\end{equation}
For the numerical inequality, \(\phi_0<2/5\) and
\(e>8/3\) give
\[
 \phi_0(16e^{-2}+12e^{-1}+3)
 <\frac25\left(\frac94+\frac92+3\right)
 =3.9<4.
\]
Taylor's theorem at zero therefore yields
\[
 |\phi''(z)-\phi''(0)|
 \le\frac12\|\phi^{(4)}\|_\infty z^2\le2z^2.
\]
For \(|y|\le6\),
\[
 \left|\frac y2-\frac{y^3}{6}\right|
 \le\frac{|y|}{2}+\frac{|y|^3}{6}\le39,
\]
so this replacement changes the main term by at most
\(78z^2\).
Since \(\phi''(0)=-\phi_0\), combining these estimates gives
\begin{equation}\label{eq:effective-H-local}
 \begin{aligned}
 \left|H_n(z,y)-\frac{\phi_0}{6}(y^3-3y)\right|
 &\le\frac{79.375}{\sqrt n}+78z^2\\
 &\le100\left(z^2+\frac1{\sqrt n}\right),
 \qquad |y|\le6.
 \end{aligned}
\end{equation}

\textbf{Keep the constant cancellation in the influence bound.}
Let \(I(y)\) be the contamination derivative from
Lemma~\ref{lem:influence} for the maximizing pair \((P,t)\).
The moment cutoff \eqref{eq:moment-cutoff} and the bounds
\eqref{eq:effective-cn} give
\begin{equation}\label{eq:effective-confinement-basic-bounds}
 \begin{gathered}
 0<\beta<2,\qquad
 |M|\le\int x^2\,P(dx)=1,\\
 \ce<R\le0.4690<\frac12,\qquad
 0<R-\ce\le\frac{4.75}{\sqrt n}.
 \end{gathered}
\end{equation}
In particular, \(R\beta<1\).
The comparison \eqref{eq:bounded-influence-comparison}
applies to the present maximizing pair.
Expanding its last two terms yields
\begin{equation}\label{eq:effective-influence-polynomial-bound}
 \begin{aligned}
 \beta I(y)\le{}&
 H_n(z,y)-R|y|^3+\frac32R\beta y^2+3RMy+K_n,\\
 K_n:={}&
 \beta\left\{\frac{n^{3/2}}{\sqrt{n-1}}C_{n-1}-nR\right\}
                  -\frac12R\beta.
 \end{aligned}
\end{equation}
Here the term \(-R\beta/2\) comes from
\(-3R\beta/2+R\beta\); retaining it is necessary for
the following cancellation.

Put
\[
 a_n=\sqrt{\frac n{n-1}},\qquad
 \theta_n=n(a_n-1)-\frac12.
\]
Since \(n\ge2\), we have
\(1<a_n\le\sqrt2<3/2\).
Rationalizing \(a_n-1\) gives
\[
 n(a_n-1)=\frac{n}{(n-1)(a_n+1)}
          =\frac{a_n^2}{a_n+1}.
\]
Consequently,
\[
 \theta_n
 =\frac{(a_n-1)(2a_n+1)}{2(a_n+1)}\ge0,
 \qquad
 \theta_n
 \le\frac{n}{2(n-1)}-\frac12
 =\frac1{2(n-1)}\le\frac2n.
\]
Using \(R=C_n\) and the definition
\eqref{eq:effective-increment} of \(d_n\), we can write
\[
 \begin{aligned}
 K_n
 &=\beta a_n\,n(C_{n-1}-R)
       +\beta R\left\{n(a_n-1)-\frac12\right\}\\
 &\le\beta a_n d_n+\beta R\theta_n.
 \end{aligned}
\]
Because \(\beta a_n<3\), \(R\beta<1\), and
\(d_n,\theta_n\ge0\), it follows that
\begin{equation}\label{eq:effective-recursive-constant}
 K_n\le3d_n+\frac2n.
\end{equation}

\textbf{Compare the contact inequality with its limiting polynomial.}
Set \(M_{\mathrm E}=M(\Pe)\), and define
\begin{equation}\label{eq:effective-contact-polynomial-definition}
 L(y)=\frac{\phi_0}{6}(y^3-3y)-\ce|y|^3
        +\frac32\ce\be y^2+3\ce M_{\mathrm E}y.
\end{equation}
We compare the right side of
\eqref{eq:effective-influence-polynomial-bound} with
\(L(y)\), uniformly for \(|y|\le6\).
The Gaussian term is controlled by
\eqref{eq:effective-H-local}, and the constant term by
\eqref{eq:effective-recursive-constant}.
Replacing \(R\) by \(\ce\) in the remaining terms has
absolute error at most
\[
 \begin{aligned}
 &(R-\ce)
   \left\{|y|^3+\frac32\beta y^2+3|M|\,|y|\right\}\\
 &\qquad\le\frac{4.75}{\sqrt n}(216+108+18)
   =\frac{1624.5}{\sqrt n}
   \le\frac{1625}{\sqrt n},
 \end{aligned}
\]
by \eqref{eq:effective-confinement-basic-bounds}.
Next, \eqref{eq:effective-moment-identification} gives
\(|\beta-\be|\le\delta_0\) and
\(|M-M_{\mathrm E}|\le\delta_0\).
Thus replacing \(\beta,M\) by \(\be,M_{\mathrm E}\)
in the terms now carrying the coefficient \(\ce\) has
absolute error at most
\[
 \ce\left\{\frac32|\beta-\be|y^2
                  +3|M-M_{\mathrm E}|\,|y|\right\}
 \le\ce(54+18)\delta_0\le36\delta_0,
\]
where \(\ce<1/2\) was used in the last step.
The four error bounds are therefore
\[
\begin{array}{ll}
\text{Error source}&\text{Upper bound}\\[2pt]
\hline
\text{Gaussian cancellation remainder}
 &100(z^2+n^{-1/2})\\
\text{Recursive constant, including }-R\beta/2
 &3d_n+2/n\\
\text{Replacement of }R\text{ by }\ce
 &1625n^{-1/2}\\
\text{Replacement of }\beta,M\text{ by }\be,M_{\mathrm E}
 &36\delta_0.
\end{array}
\]
Since \(d_n\le20n^{-1/2}\) by assumption and
\(n^{-1}\le n^{-1/2}\), their sum is at most
\[
 (100+60+2+1625)n^{-1/2}+100z^2+36\delta_0
 \le2000n^{-1/2}+1000(\delta_0+z^2).
\]
Also, \(n\ge N_{\mathrm{conf}}\ge e^{1000\mathcal A}\),
\(\delta_0=e^{-6\mathcal A}\), and
\(z^2\le e^{-5\mathcal A}\) by
\eqref{eq:effective-moment-identification}.
Hence
\[
 \begin{aligned}
 2000n^{-1/2}+1000(\delta_0+z^2)
 &\le2000e^{-500\mathcal A}
       +1000e^{-6\mathcal A}+1000e^{-5\mathcal A}\\
 &=e^{-4\mathcal A}
   \left(2000e^{-496\mathcal A}
        +1000e^{-2\mathcal A}+1000e^{-\mathcal A}\right)
 <e^{-4\mathcal A}.
 \end{aligned}
\]
For the last inequality, each of the three terms in
parentheses is less than \(1/3\):
\[
 e^{496\mathcal A}\ge1+496\mathcal A>6000,\quad
 e^{2\mathcal A}\ge1+2\mathcal A>3000,\quad
 e^{\mathcal A}\ge1+\mathcal A>3000.
\]
At every \(y\in\supp P\), the contact condition
\eqref{eq:contact} gives \(I(y)=0\).
The preceding estimates therefore yield
\begin{equation}\label{eq:effective-contact-polynomial}
 0\le L(y)+\frac{2000}{\sqrt n}
                  +1000(\delta_0+z^2)
 \le L(y)+e^{-4\mathcal A},
 \qquad y\in\supp P.
\end{equation}

\textbf{Use the polynomial to obtain a coarse interval for the support.}
Recall \(c_*=3+\sqrt{10}=6\ce/\phi_0\).
The linear cancellation
\eqref{eq:support-linear-cancellation} gives
\(c_*M_{\mathrm E}=1\).
Consequently,
\begin{equation}\label{eq:effective-contact-polynomial-factorization}
 \begin{aligned}
 L(y)
 &=\frac{\phi_0}{6}
     \left\{y^3-c_*|y|^3+\frac32c_*\be y^2\right\}\\
 &=-\frac{\phi_0}{6}y^2
    \left\{(c_*-\sgn y)|y|-\frac32c_*\be\right\}.
 \end{aligned}
\end{equation}
By \eqref{eq:support-third-moment-identity} and
\eqref{eq:support-parameter-identities},
\[
 c_*\be=\frac{\sqrt{10}}{\se},\qquad
 \frac3{c_*+1}=\pe,\qquad
 \frac3{c_*-1}=\qe.
\]
Thus the magnitudes of the negative and positive
nonzero roots are, respectively,
\begin{equation}\label{eq:effective-contact-polynomial-roots}
 \begin{aligned}
 b_-&:=\frac{3c_*\be}{2(c_*+1)}
       =\frac{\sqrt{10}}2\,\frac{\pe}{\se}
       =\frac{\sqrt{10}}2\aE,\\
 b_+&:=\frac{3c_*\be}{2(c_*-1)}
       =\frac{\sqrt{10}}2\,\frac{\qe}{\se}
       =\frac{\sqrt{10}}2\bE.
 \end{aligned}
\end{equation}
In particular, the factored polynomial can be written as
\[
 L(y)=
 \begin{cases}
  -\dfrac{\phi_0}{6}(c_*+1)y^2(|y|-b_-),&y<0,\\[4pt]
  -\dfrac{\phi_0}{6}(c_*-1)y^2(y-b_+),&y>0.
 \end{cases}
\]

We claim that \(L(y)<-0.1\) whenever
\(r=\pe+\se y\notin[-1/2,3/2]\).
To check the numerical bounds explicitly, note that
\(3.16^2<10<3.163^2\).
The definitions \eqref{eq:esseen-parameters} therefore give
\[
 0.418\le\pe\le0.420,\qquad
 0.580\le\qe\le0.582,\qquad
 0.49\le\se\le0.5.
\]
For the last pair of bounds, use
\[
 \se^2=\pe\qe\ge0.418\cdot0.580>0.49^2,
 \qquad
 \se^2=\frac14-\left(\pe-\frac12\right)^2\le\frac14.
\]
Equations \eqref{eq:effective-contact-polynomial-roots}
then imply
\[
 b_-\le\frac{3.163\cdot0.420}{2\cdot0.49}<1.363,
 \qquad
 b_+\le\frac{3.163\cdot0.582}{2\cdot0.49}<1.89.
\]
If \(r<-1/2\), then
\[
 y<0,\qquad
 |y|>\frac{1/2+\pe}{\se}
       \ge\frac{0.918}{0.5}=1.836.
\]
Using \(\phi_0>0.39\) and \(c_*+1>7\), we get
\[
 -L(y)>
 \frac{0.39}{6}\,7(1.836)^2(1.836-1.363)>0.7>0.1.
\]
If \(r>3/2\), then
\[
 y>\frac{3/2-\pe}{\se}
       \ge\frac{1.080}{0.5}=2.16,
\]
and \(c_*-1>5\) gives
\[
 -L(y)>
 \frac{0.39}{6}\,5(2.16)^2(2.16-1.89)>0.4>0.1.
\]
This proves the claim.
Since \(e^{-4\mathcal A}<0.1\), any support point
with \(r\notin[-1/2,3/2]\) would contradict
\eqref{eq:effective-contact-polynomial}.
We conclude that
\begin{equation}\label{eq:effective-coarse-support}
 \pe+\se\supp P\subset[-1/2,3/2].
\end{equation}

\textbf{Use a common jitter width for the contact comparison.}
Put \(m=n-1\) and \(G=F_{m,P}\).
Since \(n\ge N_{\mathrm{conf}}>10^{100}\), we have
\(m\in[n/2,n]\) and \(m\ge2\).
Let \(V\) be uniform on \([-1/2,1/2]\), independently
of \(m\) independent summands with law \(P\), and let \(J_h\) and
\(J\) be the CDFs of their sum plus \(hV\) and
\(h_{\mathrm E}V\), respectively.
Write
\[
 \psi(w)=(1-w^2)\phi(w),\qquad
 \mathcal G_m(v)=\Phi(v/\sqrt m)
           +\frac{\kappa}{6\sqrt m}\psi(v/\sqrt m).
\]
Lemma~\ref{lem:effective-global-jitter}, applied with
\(\ell=m\), gives
\[
 \sqrt m\,\|J_h-\mathcal G_m\|_\infty
 \le e^{-19\mathcal A}.
\]
By \eqref{eq:effective-edgeworth-density-bound},
\(\|\mathcal G_m'\|_\infty\le1/\sqrt m\);
the hypotheses of that bound hold because
\(|\kappa|\le\beta<B_*<1.84\) and \(m\ge2\).
Set \(a_0=|h-h_{\mathrm E}|/2\).
The coupling by \(V\) satisfies
\(|hV-h_{\mathrm E}V|\le a_0\), and consequently
\[
 J_h(v-a_0)\le J(v)\le J_h(v+a_0),\qquad v\in\R.
\]
Apply the approximation for \(J_h\) to both shifted
arguments, and use the derivative bound for
\(\mathcal G_m\).
We obtain
\begin{equation}\label{eq:effective-contact-jitter-width}
 \sqrt m\,\|J-\mathcal G_m\|_\infty
 \le e^{-19\mathcal A}+a_0
 \le e^{-19\mathcal A}+\frac12e^{-7\mathcal A}
 \le2\delta_0.
\end{equation}
Here \eqref{eq:effective-identification} bounds \(a_0\),
and the last inequality uses \(\delta_0=e^{-6\mathcal A}\).

\textbf{Express the discrepancy at each support contact exactly.}
Define
\[
 \mathcal P(y)=|y|^3-\beta-3My-\frac32\beta(y^2-1).
\]
For every \(y\in\supp P\), the contact equation
\eqref{eq:contact-equation} reads
\[
 n\{G(t-y)-F_{n,P}(t)\}
 =-\sqrt n\,\phi(z)y-\frac{z\phi(z)}2(y^2-1)
        +\frac{R}{\sqrt n}\mathcal P(y).
\]
Also, the positive maximizing identity is
\(F_{n,P}(t)-\Phi(z)=R\beta/\sqrt n\).
Subtracting \(\Phi((t-y)/\sqrt m)\), and using
\eqref{eq:gaussian-H}, whose definition here gives
\[
 H_n(z,y)
 =n^{3/2}\left\{\Phi\left(\frac{t-y}{\sqrt m}\right)-\Phi(z)\right\}
      +ny\phi(z)+\frac{\sqrt n}{2}z\phi(z)(y^2-1),
\]
proves the exact identity
\begin{equation}\label{eq:effective-contact-exact}
 \begin{aligned}
 &\sqrt m\left\{G(t-y)-\Phi\left(\frac{t-y}{\sqrt m}\right)\right\}\\
 &\qquad=\sqrt{\frac mn}
       \left\{R\beta+\frac{R\mathcal P(y)-H_n(z,y)}n\right\},
 \qquad y\in\supp P.
 \end{aligned}
\end{equation}
For \(|y|\le6\), the bounds in
\eqref{eq:effective-confinement-basic-bounds},
in particular \(\beta<2\) and \(|M|\le1\), imply
\[
 |\mathcal P(y)|
 \le216+2+18+\frac32(2)(35)=341,
\]
where \(|y^2-1|\le35\).
Equations \eqref{eq:effective-H-local} and
\eqref{eq:effective-moment-identification} also give
\[
 \begin{aligned}
 |H_n(z,y)|
 &\le\frac{\phi_0}{6}|y^3-3y|
           +100\left(z^2+\frac1{\sqrt n}\right)\\
 &<\frac{2/5}{6}(216+18)
           +100(e^{-5\mathcal A}+e^{-500\mathcal A})
 <15.6+1<17.
 \end{aligned}
\]
The penultimate inequality uses
\(200e^{-5\mathcal A}<1\), which follows from
\(e^{5\mathcal A}\ge1+5\mathcal A>200\).
Since \(R<1/2\) and \(R\beta<1\),
\[
 |R\mathcal P(y)-H_n(z,y)|
 \le\frac{341}{2}+17=187.5.
\]
Furthermore,
\[
 0\le1-\sqrt{\frac mn}
 =\frac1{n(1+\sqrt{1-1/n})}\le\frac1n.
\]
Applying these estimates to
\eqref{eq:effective-contact-exact} gives
\[
 \sqrt m\left\{G(t-y)-\Phi\left(\frac{t-y}{\sqrt m}\right)\right\}
 \ge R\beta-\frac1n-\frac{187.5}{n}
 \ge\ce\beta-\frac{200}{n}.
\]
Thus
\begin{equation}\label{eq:effective-contact-saturation}
 \sqrt m\left\{G(t-y)-\Phi\left(\frac{t-y}{\sqrt m}\right)\right\}
 \ge\ce\beta-\frac{200}{n},
 \qquad y\in\supp P.
\end{equation}

\textbf{Bound the change caused by jitter at each contact.}
Fix \(y\in\supp P\), put
\[
 v=t-y,\qquad w=\frac v{\sqrt m},\qquad
 a=\frac{h_{\mathrm E}}2.
\]
Since \(\pe\in[2/5,9/20]\), we have
\(\se^2=\pe\qe\ge6/25>1/9\), and hence \(a<3/2\).
Taylor's theorem and the mean value theorem yield
\[
 \begin{aligned}
 |\Phi(w+a/\sqrt m)-\Phi(w)-a\phi(w)/\sqrt m|
 &\le\frac{\|\phi'\|_\infty a^2}{2m},\\
 |\psi(w+a/\sqrt m)-\psi(w)|
 &\le\frac{\|\psi'\|_\infty a}{\sqrt m}.
 \end{aligned}
\]
The derivative bounds from
Lemma~\ref{lem:effective-selection} give
\(\|\phi'\|_\infty<1/4\) and
\(\|\psi'\|_\infty=\|\phi'''\|_\infty<5/4\),
since \(\psi=-\phi''\).
Using these bounds and \(|\kappa|\le\beta<2\),
the total scaled Taylor remainder is at most
\[
 \frac1{\sqrt m}\left(\frac{a^2}{8}+\frac{5a}{12}\right)
 <\frac1{\sqrt m}\left(\frac9{32}+\frac58\right)
 =\frac{29}{32\sqrt m}<\frac2{\sqrt n}.
\]
The last inequality uses \(m\ge n/2\).
It follows from \eqref{eq:effective-contact-jitter-width}
that
\[
 \begin{aligned}
 \sqrt m\{J(v+a)-\Phi(w)\}
 &\le a\phi(w)+\frac{\kappa}{6}\psi(w)
       +e^{-19\mathcal A}+\frac12e^{-7\mathcal A}
       +\frac2{\sqrt n}\\
 &\le\phi_0\left(\frac{h_{\mathrm E}}2+\frac{|\kappa|}{6}\right)
       +e^{-19\mathcal A}+\frac12e^{-7\mathcal A}
       +\frac2{\sqrt n}.
 \end{aligned}
\]
Here \(\|\phi\|_\infty=\|\psi\|_\infty=\phi_0\),
as in \eqref{eq:jitter-envelope-gaussian-suprema}.
Since \(\ke>0\), the triangle inequality gives
\(|\kappa|-\ke\le|\kappa-\ke|\).
By \eqref{eq:effective-esseen-envelope-height}
and \eqref{eq:effective-moment-identification},
\[
 \begin{aligned}
 \phi_0\left(\frac{h_{\mathrm E}}2+\frac{|\kappa|}{6}\right)
 -\ce\beta
 &\le\frac{\phi_0}{6}|\kappa-\ke|
           +\ce|\beta-\be|\\
 &\le\left(\ce+\frac{\phi_0}{6}\right)\delta_0
 <\frac{\delta_0}{2}.
 \end{aligned}
\]
Indeed,
\[
 \ce+\frac{\phi_0}{6}
 =\frac{(4+\sqrt{10})\phi_0}{6}
 <\frac{(4+3.2)(2/5)}6=0.48<\frac12.
\]
Subtracting \eqref{eq:effective-contact-saturation}
therefore gives
\[
 \begin{aligned}
 \sqrt m\{J(t-y+h_{\mathrm E}/2)-G(t-y)\}
 &\le e^{-19\mathcal A}+\frac12e^{-7\mathcal A}
       +\frac12\delta_0+\frac2{\sqrt n}+\frac{200}{n}\\
 &\le4\delta_0+\frac{1000}{\sqrt n}\\
 &\le e^{-5\mathcal A}
          \left(4e^{-\mathcal A}+1000e^{-495\mathcal A}\right)
 <e^{-5\mathcal A}.
 \end{aligned}
\]
For the second inequality,
\(e^{-19\mathcal A},e^{-7\mathcal A}\le\delta_0\)
and \(200/n\le200/\sqrt n\).
For the third, use \(n^{-1/2}\le e^{-500\mathcal A}\).
Finally,
\(1000e^{-495\mathcal A}\le e^{-\mathcal A}\) follows
from \(e^{494\mathcal A}\ge1+494\mathcal A>1000\),
and \(5e^{-\mathcal A}<1\).
For the lower bound, the jitter is at most
\(h_{\mathrm E}/2\), so the event that the unjittered sum
is at most \(t-y\) is contained in the event defining
\(J(t-y+h_{\mathrm E}/2)\).
We have proved
\begin{equation}\label{eq:effective-contact-flatness}
 0\le\sqrt m\{J(t-y+h_{\mathrm E}/2)-G(t-y)\}
 \le4\delta_0+\frac{1000}{\sqrt n}\le e^{-5\mathcal A},
 \qquad y\in\supp P.
\end{equation}

\textbf{Compare the increment between two support contacts.}
Let \(x<y\) belong to \(\supp P\), and suppose
\(d=y-x\in(0,h_{\mathrm E})\).
Conditioning on the uniform jitter and changing variables
gives the exact formula
\[
 J(t-y+h_{\mathrm E}/2)
 =\frac1{h_{\mathrm E}}
        \int_0^{h_{\mathrm E}}G(t-y+u)\,du.
\]
The integrand minus \(G(t-y)\) is nonnegative for
\(u\ge0\), and for \(u\in[d,h_{\mathrm E}]\) it is at
least \(G(t-y+d)-G(t-y)\).
Since \(t-y+d=t-x\), we obtain
\[
 J(t-y+h_{\mathrm E}/2)-G(t-y)
 \ge\frac{h_{\mathrm E}-d}{h_{\mathrm E}}
          \{G(t-x)-G(t-y)\}.
\]
The denominator \(h_{\mathrm E}-d\) is positive.
Using \eqref{eq:effective-contact-flatness} therefore gives
\begin{equation}\label{eq:effective-increment-upper}
 \sqrt m\{G(t-x)-G(t-y)\}
 \le\frac{h_{\mathrm E}}{h_{\mathrm E}-d}e^{-5\mathcal A}.
\end{equation}

On the other hand, subtracting the exact contact equations
\eqref{eq:contact-equation} at \(x\) and \(y\) gives
\[
 \begin{aligned}
 \sqrt m\{G(t-x)-G(t-y)\}
 &=\sqrt{\frac mn}\phi(z)d
   +\frac{\sqrt m}{2n}z\phi(z)(y^2-x^2)\\
 &\quad+\frac{R\sqrt m}{n^{3/2}}
                    \{\mathcal P(x)-\mathcal P(y)\}.
 \end{aligned}
\]
Since \(n\ge100\) and
\(z^2\le e^{-5\mathcal A}<1/100\),
\[
 \sqrt{\frac mn}\ge1-\frac1n\ge0.99,\qquad
 e^{-z^2/2}\ge1-\frac{z^2}{2}>0.99.
\]
Together with \(\phi_0>0.39\), these bounds give
\[
 \sqrt{\frac mn}\phi(z)
 =\sqrt{\frac mn}\phi_0e^{-z^2/2}
 >0.99^2(0.39)>\frac13.
\]
Also \(|x|,|y|\le6\), so \(|y^2-x^2|\le36\).
Using \(\sup_z|z\phi(z)|<1/4\) and
\(\sqrt m\le\sqrt n\), the second term has absolute
value at most
\[
 \frac{\sqrt m}{2n}\frac14(36)\le\frac{4.5}{\sqrt n}.
\]
The bound \(|\mathcal P|\le341\) and \(R<1/2\)
bound the last term by \(341/n\).
Since \(n\ge10^{100}\),
\[
 \frac{4.5}{\sqrt n}+\frac{341}{n}
 =\frac{4.5+341/\sqrt n}{\sqrt n}
 \le\frac{30}{\sqrt n}.
\]
We conclude that
\begin{equation}\label{eq:effective-increment-lower}
 \sqrt m\{G(t-x)-G(t-y)\}
 \ge\frac d3-\frac{30}{\sqrt n}.
\end{equation}
These bounds use the contact equation at every point of
\(\supp P\); no lower bound on the probability of either
support point is required.

\textbf{Use Wasserstein distance to locate a support point near each Bernoulli atom.}
Finally, the affine map
\[
 T(x)=\pe+\se x
\]
sends \(\Pe\) to \(B_{\pe}\): indeed,
\(T(-\aE)=0\) and \(T(\bE)=1\), with probabilities
\(\qe\) and \(\pe\), respectively.
For every coupling \((X,Y)\) of \(P\) and \(\Pe\),
\[
 |T(X)-T(Y)|=\se|X-Y|.
\]
Conversely, applying \(T^{-1}\) to a coupling of
\(\widetilde P\) and \(B_{\pe}\) gives a coupling of
\(P\) and \(\Pe\).
Taking infima in both directions proves the affine scaling
identity. Thus \eqref{eq:effective-identification} gives
\begin{equation}\label{eq:effective-confinement-affine-W1}
 W_1(\widetilde P,B_{\pe})
 =\se W_1(P,\Pe)
 \le\se e^{-7\mathcal A}
 \le e^{-7\mathcal A}.
\end{equation}
Here \(0<\se\le1/2\) follows from
\(\se^2=\pe(1-\pe)\).

We claim that there are \(r_0,r_1\in\supp\widetilde P\)
such that
\begin{equation}\label{eq:effective-confinement-anchors}
 |r_0|<\frac\zeta{10},\qquad
 |r_1-1|<\frac\zeta{10}.
\end{equation}
To prove this, fix \(b\in\{0,1\}\) and suppose that
\(\supp\widetilde P\) contains no point of
\((b-\zeta/10,b+\zeta/10)\).
Then \(|V-b|\ge\zeta/10\) almost surely for
\(V\sim\widetilde P\).
For every coupling \((V,B)\) with \(B\sim B_{\pe}\),
\[
 \EE|V-B|
 \ge\EE\bigl[|V-b|\ind_{\{B=b\}}\bigr]
 \ge\frac\zeta{10}\PP(B=b)
 \ge\frac{0.4\zeta}{10}.
\]
The last inequality uses \(\min(\pe,\qe)>0.4\), which
follows from \eqref{eq:esseen-parameters}.
Since \(\zeta=\eta_*/2=e^{-2\mathcal A}/2\),
\[
 \frac{0.4\zeta}{10}
 =\frac1{50}e^{-2\mathcal A}
 >e^{-7\mathcal A},
\]
where the strict inequality is equivalent to
\(e^{5\mathcal A}>50\).
Taking the infimum over couplings contradicts
\eqref{eq:effective-confinement-affine-W1} and proves
\eqref{eq:effective-confinement-anchors}.

\textbf{Exclude every support point outside the two small clusters.}
Suppose that \(r\in\supp\widetilde P\) satisfies
\(\dist(r,\{0,1\})\ge\zeta\).
By \eqref{eq:effective-coarse-support},
\(r\in[-1/2,3/2]\).
Choose \(b\in\{0,1\}\) nearest to \(r\), choosing either
one in case of a tie. Then
\[
 \zeta\le |r-b|\le\frac12.
\]
Choose the point \(r_b\) from
\eqref{eq:effective-confinement-anchors}.
The reverse triangle inequality and the triangle inequality give
\begin{equation}\label{eq:effective-confinement-anchor-distance}
 0.9\zeta
 \le |r-r_b|
 \le\frac12+\frac\zeta{10}
 <0.6.
\end{equation}
The last inequality uses \(0<\zeta<1\).
In particular, the points \(r\) and \(r_b\) are distinct.
Since \(T\) is a homeomorphism,
\(\supp\widetilde P=T(\supp P)\).
Let \(x<y\) be their inverse images under \(T\), listed
in increasing order, and put \(d=y-x\).
Because \(h_{\mathrm E}=1/\se\),
\eqref{eq:effective-confinement-anchor-distance} implies
\begin{equation}\label{eq:effective-confinement-standardized-distance}
 0.9\zeta h_{\mathrm E}
 \le d\le0.6h_{\mathrm E}<h_{\mathrm E}.
\end{equation}
The two increment estimates therefore apply to this pair.
From \eqref{eq:effective-increment-upper},
\[
 \sqrt m\{G(t-x)-G(t-y)\}
 \le\frac{h_{\mathrm E}}{h_{\mathrm E}-d}e^{-5\mathcal A}
 \le\frac52e^{-5\mathcal A}
 \le3e^{-5\mathcal A}.
\]
On the other hand, \(h_{\mathrm E}\ge2\), so
\eqref{eq:effective-increment-lower} and
\eqref{eq:effective-confinement-standardized-distance} give
\[
 \sqrt m\{G(t-x)-G(t-y)\}
 \ge\frac d3-\frac{30}{\sqrt n}
 \ge0.6\zeta-\frac{30}{\sqrt n}.
\]
Since \(n\ge N_{\mathrm{conf}}\ge e^{1000\mathcal A}\),
\begin{equation}\label{eq:effective-confinement-gap-contradiction}
 \begin{aligned}
 0.6\zeta-\frac{30}{\sqrt n}-3e^{-5\mathcal A}
 &\ge0.3e^{-2\mathcal A}
       -30e^{-500\mathcal A}-3e^{-5\mathcal A}\\
 &=e^{-2\mathcal A}
   \left(0.3-30e^{-498\mathcal A}-3e^{-3\mathcal A}\right)
 >0.
 \end{aligned}
\end{equation}
For the last inequality, \(e^{\mathcal A}>10\) implies
\(30e^{-498\mathcal A}<0.1\) and
\(3e^{-3\mathcal A}<0.1\).
The upper and lower bounds for the same increment contradict
\eqref{eq:effective-confinement-gap-contradiction}.
There is therefore no such point \(r\).
This proves, in particular, the support assertion in
\eqref{eq:effective-full-support}.

\textbf{Collapse the clusters and control their probabilities.}
The intervals \([-\zeta,\zeta]\) and
\([1-\zeta,1+\zeta]\) are disjoint because \(2\zeta<1\).
Let \(V\sim\widetilde P\), set
\[
 B=\ind_{\{V\in[1-\zeta,1+\zeta]\}},\qquad
 p=\PP(B=1)
   =\widetilde P([1-\zeta,1+\zeta]),
\]
and note that \(V\) belongs to the union of the two
clusters almost surely.
Thus \(|V-B|\le\zeta\) almost surely, and this coupling gives
\[
 W_1(\widetilde P,B_p)\le\EE|V-B|\le\zeta.
\]
For Bernoulli laws, \(W_1(B_p,B_{\pe})=|p-\pe|\).
Indeed, every coupling has cost at least the absolute difference
of the means, whereas coupling their indicators by the same
uniform random variable attains cost \(|p-\pe|\).
The triangle inequality and
\eqref{eq:effective-confinement-affine-W1} now yield
\begin{equation}\label{eq:effective-confinement-cluster-mass}
 \begin{aligned}
 |p-\pe|
 &=W_1(B_p,B_{\pe})\\
 &\le W_1(B_p,\widetilde P)+W_1(\widetilde P,B_{\pe})\\
 &\le\zeta+e^{-7\mathcal A}
 <2\zeta=\eta_*.
 \end{aligned}
\end{equation}
The strict inequality follows from
\(e^{-7\mathcal A}<e^{-2\mathcal A}/2=\zeta\).
This proves the second assertion in
\eqref{eq:effective-full-support} and completes the proof.
\end{proof}

\begin{proof}[Proof of the explicit bound in Theorem~\ref{thm:main}]
\textbf{Select a violating index to which support confinement applies.}
Suppose that \(C_N>\ce\) for an integer
\(N\ge2N_{\mathrm{conf}}\).
Lemma~\ref{lem:effective-selection} supplies an integer
\(n\in[\lceil N/2\rceil,N]\) with \(C_n>\ce\) and
\(d_n\le10/\sqrt{N-2}\).
Lemma~\ref{lem:attainment} then supplies a positive maximizing
pair \((P,t)\) at this index, and the support conclusion of
Lemma~\ref{lem:effective-selection} applies to this pair.
Thus
\begin{equation}\label{eq:effective-final-selected-index}
 \begin{gathered}
 R_n(P)=C_n>\ce,\qquad
 n\ge\lceil N/2\rceil\ge N_{\mathrm{conf}},\qquad
 \supp P\subset[-6,6],\\
 d_n\le\frac{10}{\sqrt{N-2}}
       \le\frac{10}{\sqrt{n-2}}
       \le\frac{20}{\sqrt n}.
 \end{gathered}
\end{equation}
For the middle inequality in the last line, use \(n\le N\).
For the last one, \(n\ge N_{\mathrm{conf}}>3\) gives
\(n\le4(n-2)\); taking square roots and reciprocals yields
\(1/\sqrt{n-2}\le2/\sqrt n\).
All denominators are positive.

\textbf{Apply the strict local comparison to the confined law.}
Let \(\widetilde P\) be the law of \(\pe+\se X\) for
\(X\sim P\).
The hypotheses of Proposition~\ref{prop:effective-confinement}
follow from \eqref{eq:effective-final-selected-index}.
Its conclusion \eqref{eq:effective-full-support}, with
\(\zeta=\eta_*/2\), implies
\[
 \supp\widetilde P
 \subset[-\eta_*,\eta_*]\cup[1-\eta_*,1+\eta_*].
\]
The upper-cluster probabilities are unchanged when the radius
is enlarged from \(\zeta\) to \(\eta_*\): the entire support
already lies in the two intervals of radius \(\zeta\), the
upper one is contained in \([1-\eta_*,1+\eta_*]\), and the
lower one is disjoint from that interval because
\(\zeta+\eta_*<1\).
Consequently,
\[
 \left|\widetilde P([1-\eta_*,1+\eta_*])-\pe\right|<\eta_*.
\]
These are the two-cluster hypotheses in
Proposition~\ref{prop:effective-clusters}.
Also, by \eqref{eq:effective-threshold-parameters} and
monotonicity of the ceiling function,
\[
 n\ge N_{\mathrm{conf}}
   =\lceil e^{1000\mathcal A}\rceil
   \ge\lceil e^{\mathcal A}\rceil=N_*.
\]

To record explicitly the strict inequality furnished by that
proposition, let \(V\sim\widetilde P\), and define
\[
 B=\ind_{\{V\in[1-\eta_*,1+\eta_*]\}},\qquad
 p=\PP(B=1),\qquad
 m_b=\EE(V\mid B=b)\quad(b=0,1).
\]
Since \(|p-\pe|<\eta_*<0.01\) and
\(0.41<\pe<0.42\), we have \(0.4<p<0.43\), hence
\(p\in I=[2/5,9/20]\).
Both conditional means are therefore defined.
The cluster supports imply
\[
 m_0\in[-\eta_*,\eta_*],\qquad
 m_1\in[1-\eta_*,1+\eta_*],\qquad
 d_0:=m_1-m_0\ge1-2\eta_*>0.
\]
Set \(Y=(V-m_0)/d_0\) and \(U=(V-m_B)/d_0\).
Then
\begin{equation}\label{eq:effective-final-cluster-normalization}
 Y=B+U,\qquad
 \EE(U\mid B)=0,\qquad
 |U|\le\frac{2\eta_*}{1-2\eta_*}
       \le4\eta_*\le e^{-\mathcal A}.
\end{equation}
Here \(\eta_*<1/4\) gives the penultimate inequality, and
\(\eta_*=e^{-2\mathcal A}\), \(e^{\mathcal A}>4\) give
the last.
The first assertion of Proposition~\ref{prop:effective-clusters}
now applies to \(Y\) and gives \(R_n(Y)<\ce\).
The laws of \(X\), \(V\), and \(Y\) are bounded and
nondegenerate: their variances are \(1\), \(\se^2\), and
\(\se^2/d_0^2>0\), respectively.
Affine invariance \eqref{eq:raw-normalization} therefore gives
\begin{equation}\label{eq:effective-final-strict-contradiction}
 C_n=R_n(P)=R_n(\widetilde P)=R_n(Y)<\ce,
\end{equation}
contrary to \eqref{eq:effective-final-selected-index}.

It follows that \(C_N\le\ce\) for every integer
\(N\ge2N_{\mathrm{conf}}\).
By \eqref{eq:cn} and \eqref{eq:raw-normalization}, this is the
asserted Berry--Esseen bound for every standardized summand
law with finite third absolute moment at all such sample sizes.
Finally, \eqref{eq:effective-threshold-parameters} gives
\[
 2N_{\mathrm{conf}}
 =2\left\lceil\exp(1000\mathcal A)\right\rceil
 =2\left\lceil\exp(10^{17})\right\rceil.
\]
Thus the explicit threshold \eqref{eq:explicit-threshold} is valid.
\end{proof}

\newpage
% ===== references =====

\end{document}